\RequirePackage{fix-cm}

\documentclass[leqno, fleqn, centertags, 12pt]{article}

\usepackage{latexsym}
\usepackage{amsmath}
\usepackage{amssymb}
\usepackage{verbatim}

\usepackage[T1]{fontenc}

\usepackage[upright, widespace]{fourier}

\usepackage{bm}

\usepackage{fullpage}

\usepackage{tikz-cd}

\tikzcdset{row sep/special/.initial=12ex}
\tikzcdset{row sep/spe/.initial=10.5ex}
\tikzcdset{row sep/speh/.initial=13ex}
\tikzcdset{row sep/speh/.initial=13ex}
\tikzcdset{row sep/sma/.initial=1.2em}
\tikzcdset{row sep/smaa/.initial=1.7em}

\tikzcdset{column sep/llarge/.initial=3.6em}
\tikzcdset{column sep/spec/.initial=5em}
\tikzcdset{column sep/spech/.initial=6em}
\tikzcdset{column sep/speci/.initial=6em}
\tikzcdset{column sep/specc/.initial=8em}
\tikzcdset{column sep/speccc/.initial=10em}
\tikzcdset{column sep/sboom/.initial=5.6em}

\tikzcdset{row sep/red/.initial=3.79em}

\tikzcdset{row sep/tri/.initial=5em}
\tikzcdset{column sep/tri/.initial=2.4em}

\tikzcdset{row sep/boom/.initial=5em}
\tikzcdset{column sep/boom/.initial=5.1em}

\tikzcdset{row sep/sboom/.initial=5.6em}
\tikzcdset{column sep/sboom/.initial=5.712em}

\tikzcdset{row sep/qboom/.initial=6.4em}
\tikzcdset{column sep/qboom/.initial=6.528em}
\tikzcdset{row sep/hqboom/.initial=1.65em}
\tikzcdset{row sep/shqboom/.initial=1.0em}

\tikzcdset{row sep/rboom/.initial=7.5em}

\tikzcdset{row sep/boomm/.initial=5.6em}
\tikzcdset{column sep/boomm/.initial=4.25em}

\tikzcdset{row sep/bboom/.initial=8em}
\tikzcdset{column sep/bboom/.initial=8.16em}

\tikzcdset{row sep/booms/.initial=4.25em}
\tikzcdset{column sep/booms/.initial=4.335em}

\tikzcdset{row sep/boomss/.initial=4.75em}
\tikzcdset{column sep/boomss/.initial=4.845em}

\tikzcdset{row sep/free/.initial=12ex}
\tikzcdset{column sep/freee/.initial=4em}

\usepackage{comment}

\newskip\nineskipamount \nineskipamount=9pt plus 0pt minus 0pt
\newskip\zeroskipamount \zeroskipamount=0pt plus 0pt minus 0pt
\usepackage[noorphans,vskip=\nineskipamount]{quoting}

\makeatletter
\renewcommand{\@makefntext}[1]{\vspace*{0.5ex}\parindent=0em
\hspace*{-0.4em}
\hbox to 0.4em{\hss\@makefnmark}\hspace*{0.4em}{#1}
}
\makeatother

\newcounter{mysectionnumber}
\newcommand{\mysection}[2]{\setcounter{footnote}{0}
\setcounter{myparnum}{0}
\refstepcounter{mysectionnumber}
\vspace{27pt}{\Large {\themysectionnumber.} {#1}}\label{#2}\vspace*{15pt}}

\newcommand{\myuppar}[1]{\vspace{\medskipamount}\textbf{#1}\hspace*{0.5em}}

\newcommand{\myit}[1]{\textbf{\textit{#1}}\hspace{0.0em}}

\newcounter{myparnum}[mysectionnumber]
\renewcommand{\themyparnum}{\arabic{mysectionnumber}.\arabic{myparnum}}
\newcommand{\mypar}[2]{\refstepcounter{myparnum}{\vspace{\medskipamount}\textbf{{\themyparnum. #1}\label{#2}}\hspace{0.5em}}}

\newcounter{mylemmanum}[myparnum]
\newcounter{myappendnumber}
\newcounter{myaparnum}[myappendnumber]
\newcounter{myapparnum}[mysectionnumber]
\newcommand{\proof}{\vspace{\medskipamount}{\textbf{{\emph{Proof}.}}\hspace*{1em}}}

\newcommand{\eproof}{ $\blacksquare$}

\newcommand{\dis}{\displaystyle}

\def\sss{\hspace{0.05em}\ }
\def\dss{\hspace{0.1em}\ }
\def\trs{\hspace{0.15em}\ }
\def\qss{\hspace{0.2em}\ }
\def\pss{\hspace{0.3em}\ }
\def\oss{\hspace{0.4em}\ }

\def\halfff{\hspace*{0.025em}}
\def\fff{\hspace*{0.05em}}
\def\dff{\hspace*{0.1em}}
\def\trf{\hspace*{0.15em}}
\def\qff{\hspace*{0.2em}}
\def\pff{\hspace*{0.3em}}
\def\off{\hspace*{0.4em}}

\newcommand{\hnsp}{\hspace*{-0.05em}}
\newcommand{\nsp}{\hspace*{-0.1em}}
\newcommand{\nnsp}{\hspace*{-0.15em}}

\newcommand{\dnsp}{\hspace*{-0.2em}}

\renewcommand{\leq}{\leqslant}
\renewcommand{\geq}{\geqslant}

\newcommand{\zzz}{\mathbf{Z}}

\newcommand{\ccc}{\mathbf{C}}
\newcommand{\rrr}{\mathbf{R}}
\newcommand{\nnn}{\mathbf{N}}

\newcommand{\image}{\operatorname{Im}\trf}
\newcommand{\kernel}{\operatorname{Ker}\trf}

\newcommand{\pr}{\operatorname{p{\fff}r}}
\newcommand{\inte}{\operatorname{int}\trf}

\newcommand{\core}{\operatorname{\textbf{\textup{c{\fff}o}}}\trf}

\newcommand{\ssk}{\operatorname{S{\fff}k}}
\newcommand{\id}{\operatorname{id}}

\newcommand{\gr}{\operatorname{G{\fff}r}}

\newcommand{\ob}{\operatorname{Ob}}
\newcommand{\mor}{\operatorname{Mor}}

\newcommand{\fin}{\operatorname{fin}}

\newcommand{\ev}{\operatorname{e{\fff}v}}
\newcommand{\vect}{\mbox{Vect}}
\newcommand{\ffin}{^{\dff \mathrm{fin}}}
\newcommand{\comp}{^{\dff \mathrm{comp}}}

\newcommand{\num}[1]{|\qff #1 \qff|}

\newcommand{\bbnum}[1]{\llbracket\qff #1 \qff\rrbracket}

\newcommand{\norm}[1]{\|\qff #1 \qff\|}
\newcommand{\sco}[1]{\langle\trf #1 \trf\rangle}

\newcommand{\bsl}{\backslash}
\newcommand{\sa}{\mathrm{sa}}
\newcommand{\odd}{\mathrm{odd}}
\newcommand{\inv}{\mathrm{inv}}

\newcommand{\ttoo}{\hspace*{0.2em}\longrightarrow\hspace*{0.2em}}

\newcommand{\degenerate}{de\-gen\-er\-ate\ }

\begin{document}

\setlength{\baselineskip}{12pt plus 0pt minus 0pt}
\setlength{\parskip}{12pt plus 0pt minus 0pt}
\setlength{\abovedisplayskip}{12pt plus 0pt minus 0pt}
\setlength{\belowdisplayskip}{12pt plus 0pt minus 0pt}

\newskip\smallskipamount \smallskipamount=3pt plus 0pt minus 0pt
\newskip\medskipamount   \medskipamount  =6pt plus 0pt minus 0pt
\newskip\bigskipamount   \bigskipamount =12pt plus 0pt minus 0pt

\author{Nikolai\qss V.\qss Ivanov}
\title{Topological\qss categories\qss related\qss to\qss
Fredholm\qss operators\fff:\oss
I.\oss Classifying\qss spaces}
\date{}

\footnotetext{\hspace*{-0.65em}\copyright\oss 
Nikolai\qss V.\qss Ivanov,\oss 2021.\oss}

\maketitle

\renewcommand{\baselinestretch}{1}
\selectfont

\myit{\hspace*{0em}\large Contents}  \vspace*{1ex} \vspace*{\bigskipamount}\\ 
\hbox to 0.8\textwidth{\myit{\phantom{1}1.}\hspace*{0.5em} Introduction \hfil 2}\hspace*{0.5em} 
\vspace*{1.5ex}\\
\hbox to 0.8\textwidth{\myit{\phantom{1}2.}\hspace*{0.5em} Simplicial\sss spaces \hfil 14}\hspace*{0.5em} \vspace*{0.25ex}\\
\hbox to 0.8\textwidth{\myit{\phantom{1}3.}\hspace*{0.5em} Topological\sss categories \hfil 20}\hspace*{0.5em} \vspace*{0.25ex}\\
\hbox to 0.8\textwidth{\myit{\phantom{1}4.}\hspace*{0.5em} Partially\sss ordered spaces \hfil 22}\hspace*{0.5em} \vspace*{0.25ex}\\
\hbox to 0.8\textwidth{\myit{\phantom{1}5.}\hspace*{0.5em} Topological\sss simplicial\sss complexes \hfil 23}  \hspace*{0.5em} \vspace*{0.25ex}\\
\hbox to 0.8\textwidth{\myit{\phantom{1}6.}\hspace*{0.5em} Segal's\dss introductory\sss example \hfil 27}  \hspace*{0.5em} \vspace*{0.25ex}\\
\hbox to 0.8\textwidth{\myit{\phantom{1}7.}\hspace*{0.5em} Coverings and simplicial\sss spaces \hfil 30}  \hspace*{0.5em} \vspace*{0.25ex}\\
\hbox to 0.8\textwidth{\myit{\phantom{1}8.}\hspace*{0.5em} Self-adjoint\dss Fredholm\sss operators \hfil 31}  \hspace*{0.5em} \vspace*{1.5ex}\\
\hbox to 0.8\textwidth{\myit{\phantom{1}9.}\hspace*{0.5em} Classifying spaces for self-adjoint\dss Fredholm\sss operators \hfil 34}  \hspace*{0.5em} \vspace*{0.25ex}\\ 
\hbox to 0.8\textwidth{\myit{10.}\hspace*{0.5em} Polarizations and\sss splittings \hfil 47}  \hspace*{0.25em} \vspace*{0.25ex}\\
\hbox to 0.8\textwidth{\myit{11.}\hspace*{0.5em} Restricted\dss Grassmannians \hfil 55} \hspace*{0.5em} \vspace*{0.25ex}\\
\hbox to 0.8\textwidth{\myit{12.}\hspace*{0.5em} Categories related\sss to restricted\dss Grassmannians \hfil 58} \hspace*{0.5em} \vspace*{0.25ex}\\
\hbox to 0.8\textwidth{\myit{13.}\hspace*{0.5em} Grassmannian\sss bundle and\dss Grassmannian\dss quasi-fibration \hfil 67} \hspace*{0.5em} \vspace*{1.5ex}\\
\hbox to 0.8\textwidth{\myit{14.}\hspace*{0.5em} Classifying spaces for Fredholm\dss operators \hfil 76} \hspace*{0.5em} \vspace*{0.25ex}\\ 
\hbox to 0.8\textwidth{\myit{15.}\hspace*{0.5em} Finite-unitary\dss bundle
and\sss finite-unitary\sss quasi-fibration \hfil 88}  \hspace*{0.5em} \vspace*{1.5ex}\\
\hbox to 0.8\textwidth{\myit{16.}\hspace*{0.5em} Atiyah--Singer\dss map \hfil 93}  \hspace*{0.5em} 
\vspace*{1.5ex}\\
\hbox to 0.8\textwidth{\myit{17.}\hspace*{0.5em} Hilbert\dss bundles \hfil 102}  \hspace*{0.5em} \vspace*{1.5ex}\\ 
\hbox to 0.8\textwidth{\myit{References}  \hfil 105}\hspace*{0.5em}  \hspace*{0.5em}  \vspace*{1.5ex}

\renewcommand{\baselinestretch}{1}
\selectfont

{\small
Throughout\sss this paper we work over\sss complex numbers\sss $\ccc$\nnsp.\oss
There\dss is\dss no doubt\sss that\sss similar\sss results hold over\sss
the real\sss numbers\sss $\rrr$\nnsp,\oss
but,\oss as usual,\oss the main\sss ideas are more\sss transparent\sss over $\ccc$\nnsp.\oss
}

\newpage
\mysection{Introduction}{introduction}

\myuppar{Fredholm\dss operators and classifying spaces.}
By well\dss known\sss theorems of\trs Atiyah\qss \cite{a} and of\trs
Atiyah\dss and\dss Singer\qss \cite{as}\qss 
the space $\mathcal{F}$ of\dss bounded\dss Fredholm\dss operators and\sss a component\sss $\hat{\mathcal{F}}$
of\dss the space of\dss bounded self-adjoint\trs Fredholm\sss operators
are\sss the classifying spaces for\sss the\sss topological\sss complex
$K^{\trf 0}$\nsp\dnsp-theory and $K^{\fff 1}$\nsp\dnsp-theory\sss respectively.\oss
Segal\qss \cite{s4}\qss explained\sss how\sss these spaces 
are related\dss to Quillen's\dss constructions\qss \cite{q}\qss in\sss higher algebraic $K$\dnsp-theory.\oss
Namely,\pss $\mathcal{F}$ and\sss $\hat{\mathcal{F}}$ are homotopy equivalent\sss to\sss the classifying spaces
of\dss two categories,\oss which\dss Segal\dss denoted\sss by\sss $\hat{C}$\sss
and $QC$\sss respectively\qss 
(Segal\dss attributed\sss the definitions of\dss these categories\sss to\dss Quillen).\oss
The idea behind\sss these relations\dss is\dss very attractive and\dss
described\sss by\dss Segal\dss as follows.\vspace{-9pt}\vspace{1.25pt}

\begin{quoting}
The idea\dss is\dss that\sss a\dss Fredholm\sss operator\dss is\dss
determined\sss by\sss its kernel\sss and cokernel\dss --\dss a pair
of\dss finite dimensional\sss vector spaces\dss --\dss in\sss the sense\sss
that\sss the operators with a presrcibed\sss kernel\sss and cokernel\sss form
a contractible space.\oss
When a\dss Fredholm\sss operator\dss is\dss deformed continuously\sss its
kernel\sss and cokernel\sss can\sss jump,\oss
but\sss only\sss by\sss adding\sss isomorphic pieces\sss to each:\oss
i.e.\qss the\sss jumps correspond\sss to morphisms in\sss $\hat{C}$\dnsp.\oss
In\sss the same sense a self-adjoint\sss operator\dss is\dss determined\sss by\sss
its kernel,\oss and when\sss the kernel\sss jumps\sss the piece added\sss to\sss it\dss
is\dss the sum of\dss a part\sss on which\sss the operator was positive and a part\sss on
which\sss it\sss was negative,\oss so\sss that\sss the\sss jump corresponds\sss to a
morphism in\sss $QC$\nnsp.\oss
\end{quoting}

\vspace{-9pt}\vspace{1.25pt}
Segal's\dss paper\qss \cite{s4}\qss is\dss an\sss inspiring overview of\dss
a wide range of\dss ideas related\sss to $K$\dnsp-theory,\oss but\sss
hardly contains any detailed\sss proofs.\oss 
Segal\dss devoted one page\sss to an outline of\dss proofs\sss that\sss 
$\mathcal{F}$ and\sss $\hat{\mathcal{F}}$ are homotopy equivalent\sss to\sss the classifying spaces
of\dss $\hat{C}$\sss and $QC$\sss respectively,\oss
and\sss indicated\sss that\sss a delicate argument\dss is\dss similar\sss to 
an\sss argument\sss in\qss \cite{s2},\oss Proposition\qss 2.7.\oss  
There\dss is\dss no doubt\sss that\sss the proofs are well\dss known\sss in some quarters,\oss
but\sss author's attempt\sss to reconstruct\sss these proofs was only\sss partially successful.\oss
The resulting\sss proofs are inspired\sss by\sss the ideas of\trs Segal\pss 
\cite{s1}\trs --\trs \cite{s4}\qss and\trs Quillen\trs \cite{q},\oss 
but\sss follow a different\sss route.\oss
In particular,\oss key\sss roles are played\sss by\sss the notion of\qss \emph{enhanced operators}\qss 
and\sss categories of\qss \emph{subspace models}. 

After\sss proving\sss these\sss homotopy equivalence\sss theorems 
we continue\sss in\sss the spirit\sss of\qss \cite{s4}\qss
and\sss transplant\sss some further\sss ideas of\trs Quillen\qss \cite{q},\qss \cite{gr}\qss
to\sss the realm of\trs Fredholm\dss and\sss self-adjoint\trs Fredholm\dss operators.\oss
Along\sss the way\sss we relate\sss the spaces $\mathcal{F}$ and $\hat{\mathcal{F}}$
and\sss the classifying spaces of\dss related categories\sss 
with more concrete classifying spaces,\oss 
such as a particular kind of\trs Grassmannians and\sss unitary\sss groups.\oss
The main results may\sss be considered as partial\sss analogues of\dss 
the equivalence of\qss Quillen's\trs two definitions of\dss higher algebraic $K$\dnsp-theory,\oss
the one based on plus-construction and\sss the one based on\sss the $Q$\dnsp-construction.\oss
This\dss is\dss one of\dss the main\sss theorems of\trs Quillen\qss \cite{q},\pss \cite{gr}.\oss
The proofs are simpler\sss than\dss Quillen's\dss ones.\oss
In\sss particular,\pss no analogue of\trs the $S^{\fff -1}\dff S$\sss construction is\dss needed.\oss
Morally,\oss the reason\dss is\dss that\dss the\dss Hilbert\dss space\dss
is\dss already stable\sss and no stabilization\sss procedure\dss is\dss need\-ed,\oss
in contrast\sss with\sss the algebraic $K$\dnsp-theory.\oss
At\dss the same\sss time\sss the proofs reveal\sss a parallel\dss
between\sss the methods of\trs Quillen\dss and\sss 
the methods of\trs Atiyah\dss and\dss Singer\qss \cite{as}.\oss

\myuppar{Topological\sss categories.}
Following\dss Segal\qss \cite{s1},\oss
we define a\qss \emph{topological\dss category}\pss as a small\sss category\sss
in\sss which\sss the sets\sss of\dss objects and\sss morphisms
have\sss topologies for which\sss the usual\sss structure maps are continuous.\oss
We denote\sss the space of\dss objects of\dss a\sss topological\sss
category $\mathcal{C}$ by\sss $\ob\dff \mathcal{C}$\dnsp,\oss
and\sss the space of\dss morphisms\sss by\sss $\mor\dff \mathcal{C}$\dnsp.\oss
A\sss topological\sss category $\mathcal{C}$ leads\sss to a\sss 
topological\sss space $\num{\mathcal{C}}$\sss called\sss
its\qss \emph{classifying\dss space}\pss or\qss \emph{geometric realization}.\oss
Sections\qss \ref{simplicial-spaces} -- \ref{topological-simplicial-complexes}\qss
include a review of\dss related definitions and\sss basic properties,\oss
while\dss Section\qss \ref{segal-example}\qss is\dss devoted\sss to an\sss
example due\sss to\dss Segal\qss \cite{s4},\oss
which\dss is\dss instructive and\sss used\sss later.\oss 
Section\qss \ref{coverings}\qss is\dss devoted\sss to a fundamental\dss
theorem of\qss Segal\qss \cite{s1}\qss about\sss topological\sss categories
related\sss to coverings.\oss

A simple\sss example has a\sss topological\sss space\sss $X$\sss
as\sss the space of\dss objects and\sss the space of\dss morphisms,\oss
and only\sss identity\sss morphisms.\oss
We will\sss often\sss identify\sss $X$\sss with\sss this category.\oss
Its geometric realization\dss is\sss $\num{X}\off =\off X$\nnsp.\oss
A partial\sss order\sss $\leq$\sss on a space\sss $X$\sss
defines a category\sss having\sss $X$\sss as\sss the space of\dss objects,\oss
and\sss the space of\dss pairs\sss $(\trf x\fff,\qff y\trf)$\sss such\sss that\sss 
$x\qff \leq\qff y$\sss as\sss the space of\dss morphisms.\oss 
There\dss is\dss a unique morphism\sss $x\qff \ttoo\qff y$\trs
if\trs $x\qff \leq\qff y$\sss and\sss no other morphisms.\oss

\myuppar{Simplicial\sss spaces and\sss topological\sss simplicial\sss complexes.}
The classifying space $\num{\mathcal{C}}$\sss is,\oss by\sss the definition,\oss
the geometric realization of\dss a simplicial\sss space $\mathit{N}\trf \mathcal{C}$\nnsp,\oss
the nerve of\sss $\mathcal{C}$\dnsp.\oss
Working\sss with geometric realizations of\dss simplicial\sss spaces\dss is\dss complicated\sss by
a well\sss known\sss technical\sss difficulty caused\sss 
by\sss collapsing of\dss degenerate simplices.\oss
In order\sss to deal\sss with\sss this difficulty one needs\sss either\sss 
to modify\sss the definition of\dss the geometric realization,\oss
or\sss to work only with sufficiently\sss nice\sss simplicial\sss spaces.\oss
It\dss turns out\sss that\sss all\sss simplicial\sss spaces\qss
(the nerves of\dss topological\sss categories)\qss 
arising\sss in\sss this paper
are as nice in\sss this respect\sss as one may wish.\oss

In\sss this paper we deal\sss mostly\sss with\sss 
topological\sss spaces $S$ 
partially\sss ordered\sss by a relation\sss $\leq$\sss
and considered as a\sss topological\sss category.\oss
We denote\sss its classifying space by\sss $\num{S}$\nnsp.\oss
At\sss the same\sss time\sss $S$\sss together\sss with\sss $\leq$\sss
define a structure of\dss a\qss \emph{topological\sss simplicial\sss complex}\pss on\sss $S$\nnsp.\oss
Its simplices are\sss finite\sss linearly ordered subsets of\sss $S$\nnsp.\oss
This structure\sss leads\sss to another\sss 
geometric realization\sss $\bbnum{S}$\nnsp.\oss
It\dss is\dss defined\sss in\sss the same way as\sss the classical\sss geometric
realization of\dss simplicial\sss complexes,\oss 
but\sss takes\sss the\sss topology of\sss $S$\sss into account.\oss
See\dss Section\qss \ref{topological-simplicial-complexes}.\oss
The definition of\sss $\bbnum{S}$\sss involves no
degenerate simplices and\sss avoids\sss the difficulties caused\sss by\sss them.\oss
In\sss general\sss $\bbnum{S}$\sss is\dss different\sss from\sss $\num{S}$\nnsp,\oss
but\sss $\num{S}$\sss is\dss canonically\sss homeomorphic\sss to\sss $\bbnum{S}$\sss 
if\dss the partially ordered space\sss $S$\sss has\qss \emph{free equalities}.\oss
The\sss latter property\sss means\sss that\sss the\sss three subspaces of\sss
$S\dff \times\dff S$\sss defined\sss by\sss the conditions\sss
$x\off =\off y$\nnsp,\qss $x\qff <\qff y$\nnsp,\oss 
and\sss $x\qff >\qff y$\sss on\sss the pairs\sss 
$(\trf x\fff,\qff y\trf)\qff \in\qff S\dff \times\dff S$\sss are closed.\oss
See\dss Section\qss \ref{pos-section}\qss and\dss 
Corollary\qss \ref{free-full-realizations}.\oss

It\sss turns out\sss that\sss all\sss partially\sss ordered spaces
used\sss in\sss this paper have\sss free equalities,\oss
but\sss some of\dss the\sss topological\sss categories used
do not\sss arise from\sss partial\sss orders.\oss
Still,\oss they are close\sss to partially ordered\sss 
spaces with\sss free equalities.\oss
Their\sss spaces of\dss objects are partially ordered,\oss
and morphisms\sss $x\qff \ttoo\qff y$\sss exist\dss if\trs and\dss only\trs if\dss
$x\qff \leq\qff y$\nnsp.\oss
Moreover,\oss they\sss have\qss \emph{free units}\pss in\sss the sense\sss
that\sss the subspace of\dss identity\sss morphisms\dss is\dss equal\sss to\sss
the union of\dss several\sss components of\dss the space of\dss all\dss morphisms.\oss
The classifying spaces of\dss such categories can be constructed\sss
in\sss terms of\dss non-degenerate simplices only.\oss
See\dss Lemmas\qss \ref{free-d-categories}\qss and\qss \ref{ndc-spaces}\qss
for precise statements.\oss
This allows us\sss to avoid\sss difficulties caused\sss by degenerate simplices.\oss

\myuppar{Topological\sss categories related\sss to self-adjoint\trs Fredholm\dss operators.}
Let\sss $H$\sss be a separable infinite dimensional\dss Hilbert\dss space over\sss $\ccc$\nnsp.\oss
By\sss $\hat{\mathcal{F}}$\sss we will\sss denote\sss the space of\dss
self-adjoint\trs Fredholm\dss operators\sss $H\qff \ttoo\qff H$\sss
which are neither essentially\sss positive,\oss
nor essentially\sss negative.\oss
The\sss last\sss condition\dss is\dss well\sss known and\sss motivated\sss by\qss \cite{as}.\oss
For\sss the purposes of\dss this introduction one can assume\sss
the all\sss operators are bounded and\sss $\hat{\mathcal{F}}$\sss
is\dss equipped\sss with\sss the usual\dss topology defined\sss by\sss
the norm of\dss operators.\oss
But\sss our\sss theory applies also\sss to unbounded\sss operators
and depends only on\sss few basic properties of\dss 
the\sss topology\sss of\sss $\hat{\mathcal{F}}$\dnsp.\oss
See\dss Section\qss \ref{spaces-operators}\qss for\sss the details.\oss
We will\sss consider\sss $\hat{\mathcal{F}}$\sss also as a\sss
topological\sss category,\oss in\sss the manner\sss indicated above.\oss

The most\sss basic\sss property of\dss
self-adjoint\trs Fredholm\dss operators\dss
is\dss the following.\oss
If\sss $A\qff \in\qff \hat{\mathcal{F}}$\nnsp,\oss
then\sss there exists an\sss $\varepsilon\qff >\qff 0$\sss
such\sss that\sss $-\qff \varepsilon\fff,\qff \varepsilon$\sss do not\sss belong\sss
to\sss the spectrum\sss $\sigma\trf(\trf A\trf)$\sss of\sss $A$\sss
and\sss the image of\dss the spectral\sss projection\sss 
$P_{\dff [\trf -\qff \varepsilon\fff,\dff \varepsilon \trf]}\dff(\trf A\trf)$\sss
is\dss finitely\sss dimensional.\oss
We will\sss treat\sss the image\sss
$\image\dff P_{\dff [\trf -\qff \varepsilon\fff,\dff \varepsilon \trf]}\dff(\trf A\trf)$\sss
as an approximation\sss to\sss the kernel\sss $\kernel\dff A$\nnsp.\oss

We define an\qss \emph{enhanced\sss operator}\trs
(which\dss is\dss silently assumed\sss to be self-adjoint\sss and\dss Fredholm)\qss
as a pair\sss $(\trf A\fff,\qff \varepsilon\trf)$\sss such\sss that\sss
$A\qff \in\qff \hat{\mathcal{F}}$\sss and\sss $\varepsilon$\sss
is\dss related\sss to\sss $A$\sss as above.\oss
Let\sss 
$\hat{\mathcal{E}}
\qff \subset\qff 
\hat{\mathcal{F}}\dff \times\dff \rrr_{\qff >\dff 0}$\sss 
be\sss the set\sss of\dss enhanced operators.\oss
We equip\sss $\hat{\mathcal{E}}$\sss with\sss topology\sss
induced\sss from\sss the product\sss of\dss the\sss topology of\sss
$\hat{\mathcal{F}}$\sss and\sss the\qss \emph{discrete}\pss 
topology on\sss $\rrr_{\qff >\dff 0}$\nsp.\oss
There\dss is\dss natural\dss partial\sss order on\sss $\hat{\mathcal{E}}$\dnsp.\oss
Namely,\pss
$(\trf A\fff,\qff \varepsilon\trf)
\qff \leq\qff 
(\trf A'\fff,\qff \delta\trf)$\sss
if\sss $A\off =\off A'$\sss and\sss 
$\varepsilon\qff \leq\qff \delta$\nnsp.\oss
This partial\sss order\sss turns\sss $\hat{\mathcal{E}}$\sss 
into a\sss topological\sss category.\oss 
The rule
$(\trf A\fff,\qff \varepsilon\trf)
\qff \longmapsto\qff
A$\sss
defines a continuous\qss \emph{forgetting\qss functor}\trs
$\hat{\mathcal{E}}
\qff \ttoo\qff
\hat{\mathcal{F}}$\dnsp.\oss

We would\dss like\sss to\sss turn also\sss the rule\sss
$(\trf A\fff,\qff \varepsilon\trf)
\qff \longmapsto\qff
\image\dff P_{\dff [\trf -\qff \varepsilon\fff,\qff \varepsilon \trf]}\dff(\trf A\trf)$\sss
into a continuous functor.\oss
To\sss this end,\oss we need\sss to define\sss first\sss the\sss target\sss
category of\dss this functor.\oss
Let\sss $\hat{\mathcal{S}}$\sss be\sss the category\sss having\sss
finitely\sss dimensional\sss subspaces\sss $V\qff \subset\qff H$\sss
as\sss its objects,\oss
with morphisms\sss $V\qff \ttoo\qff V\fff'$\sss
being\sss pairs\sss $U_{\dff -}\dff,\qff U_{\dff +}$\sss
of\dss finitely\sss dimensional\sss subspaces of\sss $H$\sss
such\sss that\sss
\[
\quad
V\fff'
\off =\off 
U_{\dff -}\dff \oplus\dff V\dff \oplus\dff U_{\dff +}
\pff,
\]

\vspace{-12pt}
where\sss $\oplus$\sss denotes\sss the sum of\pss \emph{orthogonal}\pss subspaces.\oss
The subspaces\sss $U_{\dff -}$ and\sss $U_{\dff +}$\sss are called\sss the\qss
\emph{negative}\qss and\qss \emph{positive parts}\oss of\dss the morphism\sss in question.\oss
The composition of\dss morphisms\dss is\dss defined\sss by\sss taking\sss the sum of\dss
negative parts and\sss the sum of\dss positive parts\sss to get,\oss
respectively,\oss the negative and\sss
the positive parts of\dss the composition.\oss 
The standard\sss topology on\sss the set\sss of\dss finitely dimensional\sss subspaces\sss
turns\sss $\hat{\mathcal{S}}$\sss into a\sss topological\sss category.\oss

The above rule\sss takes objects of\sss $\hat{\mathcal{E}}$\sss
to objects of\sss $\hat{\mathcal{S}}$\dnsp.\oss
It\sss naturally extends\sss to morphisms.\oss
Indeed,\oss suppose\sss that\sss 
$(\trf A\fff,\qff \varepsilon\trf)
\qff \ttoo\qff 
(\trf A\fff,\qff \delta\trf)$\sss
is\dss a morphism of\dss the category\sss $\hat{\mathcal{E}}$\dnsp.\oss
Then\sss $\varepsilon\qff \leq\qff \delta$\sss and\vspace{1.5pt}\vspace{0.75pt}
\[
\quad
\image\dff P_{\dff [\trf -\qff \delta\fff,\qff \delta \trf]}\dff(\trf A\trf)
\off\qff =\off\qff
\image\dff P_{\dff [\trf -\qff \delta\fff,\qff -\qff \varepsilon \trf]}\dff(\trf A\trf)
\qff \oplus\qff
\image\dff P_{\dff [\trf -\qff \varepsilon\fff,\qff \varepsilon \trf]}\dff(\trf A\trf)
\qff \oplus\qff
\image\dff P_{\dff [\trf \varepsilon\fff,\qff \delta \trf]}\dff(\trf A\trf)
\pff.
\]

\vspace{-12pt}\vspace{1.5pt}\vspace{0.75pt}
This orthogonal\sss decomposition defines a morphism\vspace{1.5pt}\vspace{0.75pt}
\[
\quad
\image\dff P_{\dff [\trf -\qff \varepsilon\fff,\qff \varepsilon \trf]}\dff(\trf A\trf)
\off \ttoo\off
\image\dff P_{\dff [\trf -\qff \delta\fff,\qff \delta \trf]}\dff(\trf A\trf)
\pff
\]

\vspace{-12pt}\vspace{1.5pt}\vspace{0.75pt}
of\dss the category\sss $\hat{\mathcal{S}}$\dnsp.\oss
This defines a functor\sss
$\hat{\mathcal{E}}\qff \ttoo\qff \hat{\mathcal{S}}$\dnsp.\oss 

Finally,\oss let\sss us\sss define\dss Quillen--Segal\dss category\sss $QC$\nnsp,\oss
which we will\sss from\sss now on denote simply\sss by\sss $Q$\nnsp.\oss
It\dss is\dss an abstract\sss analogue of\sss $\hat{\mathcal{S}}$\dnsp.\oss
The objects of\sss $Q$\sss are\sss finitely dimensional\dss Hilbert\sss spaces
over\sss $\ccc$\nnsp.\oss
A morphism\sss $V\qff \ttoo\qff V\fff'$\sss 
is\dss a\dss triple\sss 
$(\trf U_{\dff -}\dff,\qff U_{\dff +}\dff,\qff f\qff)$\nnsp,\oss
where\sss
$f\dff \colon\dff
V\qff \ttoo\qff V\fff'$\dss
is\dss an\sss isometric embedding\dss and\dss
$U_{\dff -}\dff,\pff U_{\dff +}$\dss
are subspaces of\trs $V\fff'$\sss
such\sss that
\[
\quad
V\fff'
\off =\off
U_{\dff -}\qff \oplus\qff
f\dff(\trf V\trf)\qff \oplus\qff U_{\dff +}
\qff.
\]

\vspace{-12pt}
The composition of\dss morphisms\dss is\dss defined\sss in\sss the same
way as\sss for\sss $\hat{\mathcal{S}}$\dnsp.\oss 
The\sss topology on\sss the set\sss of\dss objects\dss is\dss discrete,\oss
and\sss the\sss topology on\sss the set\sss of\dss
morphisms\dss is\dss the obvious one.\oss
The categories\sss $\hat{\mathcal{S}}$\sss and\sss $Q$\sss are related only\sss indirectly.\oss
There\dss is\dss an\sss intermediate category $Q/\fff H$\sss
having as objects\sss isometric embeddings\sss $V\qff \ttoo\qff H$\nnsp,\oss
and\sss canonical\dss functors\sss
$\hat{\mathcal{S}}
\off \longleftarrow\off 
Q/\fff H
\qff \ttoo\qff
Q$\dnsp.\oss
See\dss Section\qss \ref{classifying-spaces-saf}.\oss
The main\sss result\sss about\sss all\dss these categories\dss
is\dss the following\dss theorem.\oss

\myuppar{Theorem\qss A.}
\emph{The maps of\dss classifying\sss spaces}\vspace{1.5pt}\vspace{-0.375pt}
\[
\quad
\hat{\mathcal{F}}
\off =\off
\num{\hat{\mathcal{F}}}
\off \longleftarrow\off
\num{\hat{\mathcal{E}}}
\qff \ttoo\qff
\num{\hat{\mathcal{S}}}
\off \longleftarrow\off
\num{Q/\fff H}
\qff \ttoo\qff
\num{Q}
\]

\vspace{-12pt}\vspace{1.5pt}\vspace{-0.375pt}
\emph{induced\sss by\dss the\sss functors defined or mentioned above,\oss
are homotopy equivalences.\oss}

\vspace{6pt}
See\dss Theorems\qss \ref{intermediary}\qss and\qss \ref{operators-categories}.\oss
The proofs involve\sss two intermediate categories\sss
$\mathcal{E}\hat{\mathcal{O}}$\sss and\sss $\hat{\mathcal{O}}$\sss 
and\sss forgetting\sss functors\sss
$\hat{\mathcal{E}}\qff \ttoo\qff
\mathcal{E}\hat{\mathcal{O}}\qff \ttoo\qff
\hat{\mathcal{O}}\qff \ttoo\qff
\hat{\mathcal{S}}$\dnsp,\oss
which also\sss induce homotopy equivalences of\dss
classifying spaces.\oss
See\dss Section\qss \ref{classifying-spaces-saf}\qss for\sss the definitions and\sss proofs.\oss
Except\sss of\sss $\hat{\mathcal{F}}$\sss and\sss $\hat{\mathcal{E}}$\dnsp,\oss
the above categories do not\sss involve\sss genuine\dss Hilbert\dss space operators
and\dss belong\sss to\sss the\sss linear algebra of\dss finitely dimensional\sss
subspaces of\sss $H$\nnsp.\oss
Even\sss to deal\sss with\sss $\hat{\mathcal{F}}$\sss and\sss $\hat{\mathcal{E}}$\sss
we need only\sss the first\sss
ideas of\dss the\sss theory of\dss operators and\sss 
the contractibility\dss theorem of\trs Kuiper\qss \cite{ku}.\oss

While\sss the category\sss $Q$\sss has a definite conceptual\sss advantage
of\dss involving only\sss finitely dimensional\sss vector spaces,\oss
the category\sss $\hat{\mathcal{S}}$\sss turns out\sss to be\sss the most\sss
powerful\dss tool,\oss and after\trs Theorem\qss A\qss we work\sss mostly\sss
with\sss $\hat{\mathcal{S}}$\sss and\dss similar\sss categories.\oss 
Somewhat\sss surprisingly,\oss the classifying\sss space\sss 
$\num{\hat{\mathcal{S}}}$\sss turns out\sss to be\sss homeomorphic\sss to\sss
the space\sss $U\ffin$\sss of\dss isometries\sss $H\qff \ttoo\qff H$\sss
equal\sss to\sss the\sss identity\sss $\id_{\trf H}$\sss on a subspace\sss
of\dss finite codimension.\oss
See\dss Theorem\qss \ref{harris-h},\oss the proof\dss of\dss which\dss
is\dss based on a beautiful\dss idea of\qss Harris\qss \cite{h}.\oss
There\dss is\dss a\sss little subtle point\fff:\oss
$U\ffin$\sss should\sss be equipped\sss 
with a colimit\qss (direct\dss limit\halfff)\qss topology.\oss

\myuppar{Polarizations.}
A\qss \emph{polarization}\qss of\dss a separable infinitely dimensional\dss
Hilbert\sss space\sss $K$\sss is\dss a pair\sss 
of\dss infinitely dimensional\dss subspaces\sss 
$K_{\dff -}\dff,\pff  K_{\dff +}\off \subset\off H$\sss 
such\sss that\sss
$K\off =\off K_{\dff -}\dff \oplus\dff K_{\dff +}$\nsp.\oss
A version\sss of\dss the category\sss $\hat{\mathcal{S}}$\sss involving\sss
polarizations\sss turns out\sss to be very\sss useful.\oss
It\dss is\dss the\sss topological\sss category\sss 
$\mathcal{P}{\nsp}\hat{\mathcal{S}}$\sss 
having as objects\sss triples\sss 
$(\trf V\fff,\pff H_{\dff -}\dff,\pff H_{\dff +}\trf)$\sss
such\sss that\sss $V$\sss is\dss an object\sss of\sss $\hat{\mathcal{S}}$\sss 
and\sss $H_{\dff -}\dff,\pff  H_{\dff +}$\sss is\dss a polarization
of\dss the orthogonal\sss complement\sss $H\dff \ominus\dff V$\sss
of\sss $V$\sss in\sss $H$\nnsp.\oss
The\sss topology on\sss the set\sss of\dss objects\dss is\dss
defined\sss by\sss the usual\sss topology on\sss subspaces\sss $V$\sss
and\sss norm\sss topology on orthogonal\dss projections\sss
$H\qff \ttoo\qff H_{\dff -}$\nsp.\oss
Morphisms of\dss $\mathcal{P}{\nsp}\hat{\mathcal{S}}$\sss 
have\sss the form\vspace{1.5pt}
\[
\quad
(\trf V\fff,\off H_{\dff -}\dff,\off H_{\dff +}\dff)
\qff \ttoo\qff
(\trf U_{\dff -}\qff \oplus\qff
V\qff \oplus\qff U_{\dff +}\dff,\off 
H_{\dff -}\dff \ominus\dff U_{\dff -}\dff,\off 
H_{\dff +}\dff \ominus\dff U_{\dff +}
\trf)
\qff,
\]

\vspace{-12pt}\vspace{1.5pt}
where\sss $U_{\dff -}\dff,\pff U_{\dff +}$\sss
are finitely dimensional\sss subspaces of\dss  $H_{\dff -}\dff,\pff H_{\dff +}$\sss
respectively.\oss
Clearly,\oss the pair\sss $U_{\dff +}\dff,\pff U_{\dff -}$\sss is\dss uniquely\sss
determined and defines a morphism\dss
$V\qff \ttoo\qff U_{\dff -}\qff \oplus\qff
V\qff \oplus\qff U_{\dff +}$\dss
of\dss $\hat{\mathcal{S}}$\dnsp,\oss
and\dss hence\sss there\dss is\dss a canonical\sss functor\sss
$\mathcal{P}{\nsp}\hat{\mathcal{S}}
\qff \ttoo\qff
\hat{\mathcal{S}}$\sss
discarding\sss polarizations.\oss
One can easily show\sss that\sss the\sss induced\sss map
$\num{\mathcal{P}{\nsp}\hat{\mathcal{S}}}
\dff \ttoo\dff
\num{\hat{\mathcal{S}}}$\sss
is\dss a homotopy equivalence.\oss
See\dss Theorem\qss \ref{forgetting-equivalence}.\oss
An\sss important\sss technical\sss advantage of\dss the category\sss
$\mathcal{P}{\nsp}\hat{\mathcal{S}}$\sss
is\dss the fact\sss that\sss it\sss can\sss be defined also 
in\sss terms of\dss a partial\sss order\sss on\sss the space of\dss objects.\oss
See\qss Section\qss \ref{polarizations-splittings}.\oss\vspace{-1.25pt}

\myuppar{Restricted\dss Grassmannians and\dss the\dss Grassmannian\sss bundle.}
Another advantage of\dss the category\sss $\mathcal{P}{\nsp}\hat{\mathcal{S}}$\sss
is\dss its close relation\sss with\sss the\sss
restricted\dss Grassmannians.\oss
There are several\sss versions of\dss this notion.\oss
See\dss the book\qss \cite{ps}\qss by\dss Pressley\dss and\dss Segal,\oss Chapter\qss 7.\oss
We will\sss need\sss the following one.\oss
Suppose\sss that\sss a polarization\sss
$H\off =\off K_{\dff -}\dff \oplus\dff K_{\dff +}$\sss
of\sss $H$\sss is\dss fixed.\oss
The corresponding\qss \emph{restricted\dss Grassmannian}\trs $\gr$\sss
is\dss the space of\dss subspaces\sss $K\qff \subset\qff H$\sss which are\qss
\emph{commensurable}\pss with\sss $K_{\dff -}$\sss in\sss the sense\sss
that\sss the intersection\sss $K\dff \cap\dff K_{\dff -}$\sss
has finite codimension\sss in\sss $K$\sss and\sss in\sss $K_{\dff -}$\nsp.\vspace{-1.25pt}

This definition\sss naturally\sss extends\sss from\sss polarizations of\sss $H$\sss
to objects of\dss the category\sss $\mathcal{P}{\nsp}\hat{\mathcal{S}}$\nsp\dnsp.\oss
Namely,\oss if\dss $P\off =\off (\trf V\fff,\pff H_{\dff -}\dff,\pff H_{\dff +}\trf)$\sss
is\dss an object\sss of\sss $\mathcal{P}{\nsp}\hat{\mathcal{S}}$\nsp\dnsp,\oss
then we define\sss the restricted\dss Grassmannian\sss
$\gr\trf(\trf P\trf)$\sss as\sss the space of\dss subspaces\sss
$K\qff \subset\qff H$\sss which are
commensurable with\sss $H_{\dff -}$\nnsp.\oss
Obviously,\oss if\dss there\dss is\dss a morphism\sss
$P\qff \ttoo\qff P'$\sss of\dss the category\sss
$\mathcal{P}{\nsp}\hat{\mathcal{S}}$\dnsp,\oss
then\sss
$\gr\trf(\trf P\trf)\off =\off \gr\trf(\trf P'\trf)$\nnsp.\oss
Using\sss this fact\sss one can easily\sss construct\sss a\sss 
locally\sss trivial\dss bundle
\[
\quad
\bm{\pi}\dff \colon\dff
\mathbf{G}
\qff \ttoo\qff 
\num{\mathcal{P}{\nsp}\hat{\mathcal{S}}}
\pff
\]

\vspace{-12pt}
with\sss the fiber\sss $\gr\trf(\trf P\trf)$\sss 
over simplices\sss having\sss $P$\sss as a vertex.\oss
See\dss Section\qss \ref{grassmannian-fibrations}.\oss
The methods of\trs Quillen\qss \cite{q}\qss allow\sss 
to construct\sss a categorical\dss model\sss of\dss this bundle,\oss
related\sss to\sss it\sss in approximately\sss the same way as
$\num{\hat{\mathcal{S}}}$\sss is\dss related\sss to\sss $\hat{\mathcal{F}}$\dnsp.\oss
This model\sss turns out\sss to be not\sss a bundle,\oss
but\sss only\sss a\qss \emph{quasi-fibration}.\oss
Unfortunately,\oss this notion\dss is\dss still\dss less known\sss
than\sss it\sss deserves.\oss\vspace{-1.25pt}

\myuppar{Quasi-fibrations.}
A continuous surjective map\sss 
$p\dff \colon\dff E\qff \ttoo\qff B$\sss
is\dss called\qss \emph{quasi-fibration}\pss
if\dss for every\sss point\sss $b\qff \in\qff B$\nnsp,\oss
every\sss point\sss $e\qff \in\qff p^{\dff -\dff 1}\dff(\trf b\trf)$\sss
and every\sss $i\qff \geq\qff 0$\sss the induced\sss map\sss
\[
\quad
p_{\dff *}\dff \colon\dff
\pi_{\dff i}\dff\left(\trf E\fff,\qff p^{\dff -\dff 1}\dff(\trf b\trf)\dff,\qff e\trf\right)
\qff \ttoo\qff 
\pi_{\dff i}\trf(\qff B\dff,\qff b\qff)
\]

\vspace{-12pt}
is\dss an\sss isomorphism.\oss
This notion\dss is\dss due\sss to\dss Dold\dss and\dss Thom\qss \cite{dt}\qss
and\sss plays an\sss important\sss role in\dss Quillen's\dss work\qss \cite{q}.\oss
If\dss $p\dff \colon\dff E\qff \ttoo\qff B$\sss is\dss a quasi-fibration,\oss
then\sss every\sss fiber\sss $p^{\dff -\dff 1}\dff(\trf b\trf)$\sss is\dss
weakly\sss homotopy equivalent\sss to\sss 
the homotopy\sss fiber of\dss $p$\nnsp.\oss\vspace{-1.25pt}

Quasi-fibrations which are not\sss fibrations usually\sss 
arise in\sss the following\sss context.\oss
We will\sss assume\sss that\sss all\sss involved\sss spaces are compactly\sss generated\qss
(this notion\dss is\dss recalled\sss in\dss Section\qss \ref{simplicial-spaces})\qss
and\sss paracompact.\oss
Let\sss $p\dff \colon\dff E\qff \ttoo\qff B$\sss
be a surjective continuous map.\oss 
Suppose\sss that\sss $B$\sss is\dss the union of\dss an\sss increasing\sss 
sequence of\dss closed subspaces\sss
$B_{\trf 0}\qff \subset\qff
B_{\dff 1}\qff \subset\qff
B_{\trf 2}\qff \subset\qff
\ldots$\sss
and\sss that\sss the\sss topology of\sss $B$\sss is\dss the colimit\sss topology
defined\sss by\dss this sequence.\oss
Let\sss $E_{\dff i}\off =\off p^{\dff -\dff 1}\dff(\trf B_{\dff i}\trf)$\nnsp.\oss
Then\sss $p\dff \colon\dff E\qff \ttoo\qff B$\sss is\dss a quasi-fibration\dss
if\dss the following\sss three conditions hold.\oss
First,\oss the map
\[
\quad
E_{\dff i\dff +\dff 1}\qff \smallsetminus\pff E_{\dff i}
\off \ttoo\off
B_{\dff i\dff +\dff 1}\qff \smallsetminus\pff B_{\dff i}
\pff
\]

\vspace{-12pt}
induced\sss by\sss $p$\sss is\dss a\sss locally\sss trivial\dss bundle
for every\sss $i\qff \geq\qff 0$\nnsp.\oss
Second,\oss for every\sss $i\qff \geq\qff 0$\sss
there exists an open\sss neighborhood\sss $U_{\dff i}$\sss
of\trs $B_{\dff i}$\sss in\sss $B_{\dff i\dff +\dff 1}$\sss
and\sss a deformation\sss retraction\sss
$r_{\dff i}\dff \colon\dff 
U_{\dff i}\qff \ttoo\qff B_{\dff i}$\sss
which\sss lifts\sss together\sss with\sss the deformation\sss
to a deformation\sss retraction\sss
$\widetilde{r}_{\dff i}\dff \colon\dff 
p^{\dff -\dff 1}\dff(\trf U_{\dff i}\trf)\qff \ttoo\qff E_{\dff i}$\nsp.\oss
Finally,\oss for every\sss $u\qff \in\qff U_{\dff i}$\sss
the map\sss\vspace{1.5pt}
\[
\quad
p^{\dff -\dff 1}\dff(\trf u\trf)
\qff \ttoo\qff
p^{\dff -\dff 1}\dff\left(\trf r_{\dff i}\trf(\trf u\trf)\trf\right)
\]

\vspace{-12pt}\vspace{1.5pt}
induced\dss by\sss $\widetilde{r}_{\dff i}$\sss is\dss a weak\sss homotopy equivalence.\oss
One can\sss find detailed\sss expositions of\trs Dold--Thom\dss theory\sss in\qss \cite{sr},\oss 
Appendix\qss A,\oss and\sss in\qss \cite{agp},\oss Appendix\qss A.\oss
See also\qss \cite{ha},\oss Section\qss 4.K.\oss
Theorem\qss A.1.19\qss of\qss \cite{agp}\qss differs\sss from\sss the above
claim\sss only\sss by\sss the superfluous assumption\sss 
that\sss the bundles over  
$B_{\dff i\dff +\dff 1}\trf \smallsetminus\qff B_{\dff i}$\sss 
are\sss trivial.\oss

\myuppar{Quasi-fibrations arising\sss from\sss functors.}
Let\sss 
$f\dff \colon\dff \mathcal{C}\qff \ttoo\qff \mathcal{C}'$\sss
be a functor.\oss
Let\sss us define a category\sss $S\trf(\trf f\trf)$\sss
having as objects\sss triples\sss $(\trf X\fff,\qff Y\fff,\qff v\trf)$\sss
with\sss $X$\sss and\sss $Y$\sss being objects of\sss
$\mathcal{C}$\sss and\sss $\mathcal{C}'$\sss respectively,\oss
and\sss $v$\sss being a morphism\sss 
$Y\qff \ttoo\qff f\trf(\trf X\trf)$\sss of\sss $\mathcal{C}'$\dnsp.\oss
Morphisms with\sss the source\sss $(\trf X\fff,\qff Y\fff,\qff v\trf)$\sss
are in\sss one-to-one correspondence with\sss the pairs\sss of\dss
morphisms\dss $u\dff \colon\dff X\qff \ttoo\qff X\fff'$\sss and\sss
$w\dff \colon\dff Y\fff'\qff \ttoo\qff Y$\sss of\sss
$\mathcal{C}$\sss and\sss $\mathcal{C}'$\sss respectively.\oss
The\sss target\sss of\dss the morphism corresponding\sss to\sss 
$(\trf u\fff,\qff w\trf)$\sss
is\dss the\sss triple\dss
$\left(\trf X\fff'\fff,\pff Y\fff'\fff,\pff 
f\trf(\trf u\trf)\dff \circ\dff v\dff \circ\dff w\trf\right)$\nnsp.\oss

There\dss is\dss a covariant\sss functor\sss
$p_{\dff 1}\dff \colon\dff
S\trf(\trf f\trf)\qff \ttoo\qff \mathcal{C}$\sss
taking\sss $(\trf X\fff,\qff Y\fff,\qff v\trf)$\sss to\sss $X$\nnsp,\oss
and a\qss \emph{contravariant}\pss functor\sss
$p_{\trf 2}\dff \colon\dff
S\trf(\trf f\trf)\qff \ttoo\qff \mathcal{C}'$\sss
taking\sss $(\trf X\fff,\qff Y\fff,\qff v\trf)$\sss to\sss $Y$\nnsp.\oss
A contravariant\sss functor can\sss be considered as a covariant\sss functor\sss
to\sss the opposite category.\oss
Since\sss the classifying space of\dss a category\sss and of\dss its opposite are
the same,\oss contravariant\sss functors also induce maps of\dss classifying spaces.\oss
This\sss leads\sss to continuous maps\vspace{1.5pt}\vspace{1.5pt}
\[
\quad
\num{p_{\dff 1}}\dff \colon\dff
\num{S\trf(\trf f\trf)}\qff \ttoo\qff \num{\mathcal{C}}
\quad
\mbox{and}\quad\dff
\num{p_{\trf 2}}\dff \colon\dff
\num{S\trf(\trf f\trf)}\qff \ttoo\qff \num{\mathcal{C}'}
\pff.
\]

\vspace{-12pt}\vspace{1.5pt}\vspace{1.5pt}
For every object\sss $Z$\sss of\sss $\mathcal{C}'$\sss let\sss us denote\sss
by\sss $Z\fff\bsl f$\sss the subcategory\sss of\sss $S\trf(\trf f\trf)$\sss
having as objects\sss triples of\dss the form\sss
$(\trf X\fff,\qff Z\fff,\qff v\trf)$\sss
and as morphisms pairs of\dss the form\sss
$(\trf u\fff,\qff \id_{\trf Z}\qff)$\nnsp.\oss
Each\sss morphism\sss
$w\dff \colon\dff Z\fff'\qff \ttoo\qff Z$\sss
of\sss $\mathcal{C}'$\sss
induces a functor\sss
$w^{\dff *}\dff \colon\dff
Z\fff\bsl f
\qff \ttoo\qff
Z\fff'\fff\bsl f$\sss
and\sss hence a continuous map\vspace{1.5pt}\vspace{1.5pt}
\[
\quad
\num{w^{\dff *}}\dff \colon\dff
\num{Z\fff\bsl f}
\qff \ttoo\qff
\num{Z\fff'\fff\bsl f}
\pff.
\]

\vspace{-12pt}\vspace{1.5pt}\vspace{1.5pt}
Quillen\qss \cite{q}\qss showed\sss that\sss if\dss all\dss
induced\sss maps\sss $\num{w^{\dff *}}$\sss are homotopy\sss equivalences,\oss
then\sss the map\sss
$\num{p_{\trf 2}}\dff \colon\dff
\num{S\trf(\trf f\trf)}\qff \ttoo\qff \num{\mathcal{C}'}$\sss
is\dss a quasi-fibration.\oss
This\dss is\dss a key\sss step in\sss the proofs
of\trs Theorems\qss A\qss and\qss B\qss
of\trs Quillen\qss \cite{q}.\oss
We will\dss need\sss an analogue of\dss this result\sss for some\sss
topological\sss categories.\oss
These categories are nice and concrete,\oss
and direct\sss proofs based on\dss Quillen's\dss ideas are
more simple and\sss transparent\dss than\sss proving a\sss
general\dss theorem and\sss applying\sss it.\oss

The original\sss proofs were based\sss on
a\sss general\sss version of\trs Quillen's\trs
Theorems\qss A\qss and\qss B\qss due\sss to\dss 
Ebert\sss and\dss Randal-Williams\qss \cite{er}.\oss
This version\dss is\dss concerned\sss with a generalization of\dss topological\sss
categories,\oss but\sss only\sss with\sss fat\sss geometric realizations,\oss and
quasi-fibrations are present\sss only\sss implicitly.\oss
Partly\sss by\sss these reasons\sss the direct\sss proofs appear\sss to be more attractive.\oss

\myuppar{Splittings.}
Let\sss $\mathcal{P}$\sss
be\sss the full\sss subcategory of\sss
$\mathcal{P}{\nsp}\hat{\mathcal{S}}$\sss
having as objects\sss the objects of\sss
$\mathcal{P}{\nsp}\hat{\mathcal{S}}$\sss
of\dss the form\sss
$(\dff 0\fff,\qff H_{\dff -}\fff,\qff H_{\dff +}\dff)$\nnsp.\oss
This subcategory\sss has only\sss identity\sss morphisms,\oss
and\dss its space of\dss objects\dss is\dss the
space of\dss polarizations of\sss $H$\nnsp.\oss 
This implies\sss that\sss $\num{\mathcal{P}}$\sss is\dss contractible.\oss
A\qss \emph{splitting}\pss of\dss an object\sss $M$\sss of\sss
$\mathcal{P}{\nsp}\hat{\mathcal{S}}$\sss is\dss a morphism\sss
$s\dff \colon\dff N\qff \ttoo\qff M$\sss such\sss that\sss $N$\sss
is\dss an object\sss of\sss $\mathcal{P}$\dnsp.\oss

There\dss is\dss a\sss topological\sss category\sss
$s\dff \hat{\mathcal{S}}$\sss having as objects 
splittings of\dss objects of\sss $\mathcal{P}{\nsp}\hat{\mathcal{S}}$\dnsp.\oss
Morphisms of\dss the category\sss $s\dff \hat{\mathcal{S}}$\sss having as\sss the source\sss
some splitting\sss
$s\dff \colon\dff N\qff \ttoo\qff M$\sss 
are in one-to-one correspondence\sss with\sss morphisms\sss
$f\dff \colon\dff M\qff \ttoo\qff M\fff'$\dnsp.\oss
The\sss target\sss of\dss the morphism corresponding\sss to\sss 
$f\dff \colon\dff M\qff \ttoo\qff M\fff'$\sss
is\dss the splitting\sss 
$f\dff \circ\dff s\dff \colon\dff N\qff \ttoo\qff M\fff'$\dnsp.\oss
There\dss is\dss an obvious\sss forgetting\sss functor\sss
$\phi\dff \colon\dff
s\dff \hat{\mathcal{S}}
\qff \ttoo\qff
\mathcal{P}{\nsp}\hat{\mathcal{S}}$\sss
taking a splitting\sss 
$s\dff \colon\dff N\qff \ttoo\qff M$\sss 
of\dss $M$\sss to\sss $M$\nnsp.\oss

\myuppar{The main quasi-fibration.}
For\sss the functor\sss
$\phi\dff \colon\dff
s\dff \hat{\mathcal{S}}
\qff \ttoo\qff
\mathcal{P}{\nsp}\hat{\mathcal{S}}$\dss
the objects of\sss $S\trf(\trf \phi\trf)$\sss
can\sss be identified\sss with diagrams in\sss 
$\mathcal{P}{\nsp}\hat{\mathcal{S}}$\sss
of\dss the form\sss
$P\qff \ttoo\qff M\off \longleftarrow\off N$\nnsp,\oss
where\sss $P$\sss is\dss an object\sss of\dss
$\mathcal{P}{\nsp}\hat{\mathcal{S}}$\sss
and\sss $N$\sss is\dss an object\sss of\sss
$\mathcal{P}$\dnsp.\oss
We will\sss denote\sss the functor\sss
$p_{\trf 2}\dff \colon\dff
S\trf(\trf f\trf)\qff \ttoo\qff \mathcal{P}{\nsp}\hat{\mathcal{S}}$\sss
simply\sss by\sss $p$\nnsp.\oss
Clearly,\oss the functor\sss 
$p$\sss takes\sss
$P\qff \ttoo\qff M\off \longleftarrow\off N$\sss to\sss $P$\dnsp.\oss
It\dss induces a continuous map\vspace{1.5pt}\vspace{-0.85pt}
\[
\quad
\num{p}\dff \colon\dff
\num{S\dff(\trf \phi\trf)}
\qff \ttoo\qff
\num{\mathcal{P}{\nsp}\hat{\mathcal{S}}}
\pff.
\]

\vspace{-12pt}\vspace{1.5pt}\vspace{-0.85pt}
This map\dss is\dss our categorical\dss model\sss of\dss 
$\bm{\pi}\dff \colon\dff
\mathbf{G}
\qff \ttoo\qff 
\num{\mathcal{P}{\nsp}\hat{\mathcal{S}}}$\nnsp.\oss
It\sss turns out\sss to be a quasi-fibration.

Let\sss $P$\sss be an object\sss of\sss $\mathcal{P}{\nsp}\hat{\mathcal{S}}$\dnsp.\oss
It\dss turns out\sss that\sss $P\bsl \phi$\sss is\dss a categorial\sss model\sss
of\dss the restricted\dss Grassmannian\sss $\gr\trf(\trf P\trf)$\nnsp.\oss
The objects of\dss  $P\bsl \phi$\sss
can\sss be identified\sss with diagrams in\sss 
$\mathcal{P}{\nsp}\hat{\mathcal{S}}$\sss
of\dss the form\sss
$P\qff \ttoo\qff M\off \longleftarrow\off N$\nnsp,\oss
where\sss $N$\sss is\dss an object\sss of\sss
$\mathcal{P}$\dnsp,\oss
i.e.\qss a polarization\sss 
$H\off =\off N_{\dff -}\dff \oplus\dff N_{\dff +}$\sss 
of\sss $H$\nnsp.\oss
It\dss is\dss easy\sss to see\sss that\sss
$N_{\dff -}\qff \in\qff \gr\trf(\trf P\trf)$\nnsp.\oss
This\sss leads\sss to a map\sss
$\ob\dff P\bsl \phi\qff \ttoo\qff \gr\trf(\trf P\trf)$\nnsp,\oss
which naturally extends\sss to a functor\sss
$h\trf(\trf P\trf)\dff \colon\dff
P\bsl \phi \ttoo\qff \gr\trf(\trf P\trf)$\nnsp.\oss
The map\sss
$\num{h\trf(\trf P\trf)}\dff \colon\dff
\num{P\bsl \phi}\qff \ttoo\qff \gr\trf(\trf P\trf)$\sss
turns out\sss to be a homotopy equivalence,\oss
justifying\sss the idea\sss that\sss $P\bsl \phi$\sss 
is\dss a model\sss of\dss $\gr\trf(\trf P\trf)$\nnsp.\oss

As we already\sss pointed\sss out,\oss if\dss
$w\dff \colon\dff P'\qff \ttoo\qff P$\sss is\dss a morphism\sss
$\mathcal{P}{\nsp}\hat{\mathcal{S}}$\dnsp,\oss then\sss
$\gr\trf(\trf P'\trf)\off =\off \gr\trf(\trf P\trf)$\nnsp.\oss
It\sss follows\sss that\sss
$\num{w^{\dff *}}\dff \colon\dff
\num{P'\bsl \phi}
\qff \ttoo\qff
\num{P\bsl \phi}$\sss
is\dss a homotopy equivalence.\oss
With\sss this result\sss at\sss hand,\oss
we are ready\sss to apply\dss Quillen's\dss methods and\sss
prove\sss the following\dss theorem.\oss

\myuppar{Theorem\qss B.}
\emph{The map\sss
$\num{p}\dff \colon\dff
\num{S\dff(\trf \phi\trf)}
\qff \ttoo\qff
\num{\mathcal{P}{\nsp}\hat{\mathcal{S}}}$\sss
is\dss a quasi-fibration.\oss}

\vspace{6pt}
In\sss the main\sss part\sss of\dss the paper\sss we proceed\sss
in a different\sss order.\oss
Namely,\oss we define first\sss the categories\sss $P\bsl \phi$\sss
as analogues of\dss the restricted\dss Grassmannians,\oss
and denote\sss them\sss by\sss $\mathcal{G}\dff (\trf P\trf)$\nnsp.\oss
Then\sss we use\sss the bundle\sss 
$\bm{\pi}\dff \colon\dff
\mathbf{G}
\qff \ttoo\qff 
\num{\mathcal{P}{\nsp}\hat{\mathcal{S}}}$\sss
and an analogy\sss between\sss these categories
and\sss restricted\dss Grassmannians\sss
to motivate\sss the definition of\dss the category\sss $S\dff(\trf \phi\trf)$\nnsp.\oss
In\dss Section\qss \ref{restricted-grassmannians}\qss we prove\sss that\sss
restricted\dss Grassmannians are homotopy equivalent\sss to a more classical\sss
version of\dss infinite dimensional\dss Grassmannians,\oss
which we denote\sss by\sss $\gr\trf(\dff \infty\dff)$\nnsp.\oss
See\dss Theorem\qss \ref{two-grassmannians}.\oss
In\dss Section\qss \ref{categories-grassmannians}\qss we prove\sss that\sss
the spaces\sss $\num{\mathcal{G}\dff (\trf P\trf)}$\sss
are also homotopy equivalent\sss to\sss $\gr\trf(\dff \infty\dff)$\nnsp.\oss
See\dss Theorem\qss \ref{s-sigma-homotopy-type}.\oss
In\dss Section\qss \ref{grassmannian-fibrations}\qss we prove\sss that\sss
for a morphism\sss 
$w\dff \colon\dff P'\qff \ttoo\qff P$\sss
the map\sss
$\num{w^{\dff *}}\dff \colon\dff
\num{\mathcal{G}\dff (\trf P\trf)}
\qff \ttoo\qff
\num{\mathcal{G}\dff (\trf P'\trf)}$\sss
is\dss a homotopy equivalence and\sss
then prove\dss Theorem\qss B.\oss
See\dss Lemma\qss \ref{fibers-are-he}\qss 
and\trs Theorem\qss \ref{p-is-quasi-fibration}.\oss
Similar\sss ideas can\sss be used\sss to prove\sss
the following\sss theorem.\oss

\myuppar{Theorem\qss C.}
\emph{The classifying\sss space\sss $\num{S\dff(\trf \phi\trf)}$\sss 
is\dss contractible.\oss}

\vspace{6pt}
The proof\dss of\trs Theorem\qss B\qss can\sss be easily\sss modified\sss
to prove\sss that\sss the map\sss
$\num{S\dff(\trf \phi\trf)}
\qff \ttoo\qff
\num{s\dff \hat{\mathcal{S}}}$\sss
induced\sss by\sss the forgetting\sss functor\sss
$p_{\dff 1}\dff \colon\dff
S\dff(\trf \phi\trf)\qff \ttoo\qff s\dff \hat{\mathcal{S}}$\sss
is\dss a quasi-fibration.\oss
It\dss is\dss easy\sss to see\sss that\sss $\num{s\dff \hat{\mathcal{S}}}$\sss
is\dss contractible and\sss that\sss for every splitting\sss
$s\dff \colon\dff N\qff \ttoo\qff M$\nnsp,\oss i.e.\qss every object\sss
$s$\sss of\sss $s\dff \hat{\mathcal{S}}$\dnsp,\oss 
the space\sss $\num{s\bsl p_{\dff 1}}$ is\dss also contractible.\oss
Together\sss these facts imply\sss that\sss $\num{S\dff(\trf \phi\trf)}$\sss
is\dss weakly contractible.\oss
One can also prove directly\sss that\sss $\num{S\dff(\trf \phi\trf)}$\sss
is\dss homotopy equivalent\sss to a\dss CW-complex\sss and deduce\sss that\sss
$\num{S\dff(\trf \phi\trf)}$\sss is\dss contractible.\oss

But\sss the more concrete\sss arguments used\sss to prove\dss
Theorem\qss B\qss lead\sss to a direct\sss proof\dss 
of\dss the fact\sss that\sss
$\num{S\dff(\trf \phi\trf)}
\qff \ttoo\qff
\num{s\dff \hat{\mathcal{S}}}$\sss
is\dss a homotopy equivalence,\oss
and\sss hence\sss to a direct\sss proof\dss of\dss
the contractibility\sss of\dss $\num{S\dff(\trf \phi\trf)}$\nnsp.\oss
In\sss fact,\oss this sort\sss of\dss arguments\dss is\dss used\sss
already\sss in\sss the proof\dss of\trs Theorem\qss A.\oss
See\sss the proofs of\qss Theorems\qss \ref{to-models}\qss and\qss \ref{intermediary}.\oss

\myuppar{Comparing\sss  
$\bm{\pi}\dff \colon\dff
\mathbf{G}
\qff \ttoo\qff 
\num{\mathcal{P}{\nsp}\hat{\mathcal{S}}}$ 
and\sss 
$\num{p}\dff \colon\dff
\num{S\dff(\trf \phi\trf)}
\qff \ttoo\qff
\num{\mathcal{P}{\nsp}\hat{\mathcal{S}}}$\nnsp.}
The above continuous maps\sss
$\num{h\trf(\trf P\trf)}\dff \colon\dff
\num{P\bsl \phi}\qff \ttoo\qff \gr\trf(\trf P\trf)$\sss
can\sss be assembled\sss into\sss a continuous map\sss
$h\dff \colon\dff
\num{S\dff(\trf \phi\trf)}
\qff \ttoo\qff
\mathbf{G}$\sss 
such\sss that\sss the\sss following\dss triangle\dss
is\dss commutative.\oss
\[
\quad
\begin{tikzcd}[column sep=tri, row sep=sboom]
\protect{\num{S\dff(\trf \phi\trf)}}
\arrow[dr, "\dis \protect{\num{p}}"']
\arrow[rr, "\dis h"]
&
&
\mathbf{G}
\arrow[dl, "\dis \bm{\pi}"]
\\
&
\protect{\num{\mathcal{P}{\nsp}\hat{\mathcal{S}}}}
&
\end{tikzcd}
\]

\vspace{-9pt}
\myuppar{Theorem\qss D.}
\emph{The map\sss
$h\dff \colon\dff
\num{S\dff(\trf \phi\trf)}
\qff \ttoo\qff
\mathbf{G}$\sss
is\dss a homotopy\sss equivalence.\oss}

\vspace{6pt}
The proof\dss of\trs Theorem\qss B\qss shows\sss that\sss the homotopy\sss fiber of\dss
$\num{p}\dff \colon\dff
\num{S\dff(\trf \phi\trf)}
\qff \ttoo\qff
\num{\mathcal{P}{\nsp}\hat{\mathcal{S}}}$\dss
is\dss the geometric realizations\sss $\num{P\bsl \phi}$\sss for every object\sss
$P$\sss of\sss $\mathcal{P}{\nsp}\hat{\mathcal{S}}$\nsp\dnsp.\oss
On\sss the other\sss hand,\oss the spaces\sss $\gr\trf(\trf P\trf)$\sss
are\sss the fibers,\oss and\sss hence\sss the homotopy fibers of\dss 
$\num{p}\dff \colon\dff
\num{S\dff(\trf \phi\trf)}
\qff \ttoo\qff
\num{\mathcal{P}{\nsp}\hat{\mathcal{S}}}$\nnsp.\oss
Comparing\sss the homotopy sequences of\sss $\num{p}$\sss
and\sss $\bm{\pi}$\sss shows\sss that\sss
$h\dff \colon\dff
\num{S\dff(\trf \phi\trf)}
\qff \ttoo\qff
\mathbf{G}$\sss 
is\dss a weak\sss homotopy equivalence.\oss
In order\sss to prove\dss Theorem\qss D\qss it\dss remains\sss to
check\sss that\sss $\mathbf{G}$\sss is\dss homotopy equivalent\sss
to a\dss CW-complex.\oss
This\dss is\dss verified\sss in\trs Lemma\qss \ref{g-is-cw}\qss
with\sss the help of\dss some general\sss results of\qss
Morita\qss \cite{mo}\qss and\trs tom\dss Dieck\qss \cite{td1}.\oss
It\sss seems\sss that\sss the information about\sss the structure of\dss
the map\sss $\num{p}$\sss contained\sss 
in\sss the proof\dss of\trs Theorem\qss B\qss
is\dss sufficient\sss to prove\dss Theorem\qss D\qss directly,\oss
without\sss proving\sss first\sss that\sss $\mathbb{G}$\sss 
is\dss homotopy equivalent\sss to a\dss CW-complex.\oss

\myuppar{Theorem\qss E.}
\emph{The\sss total\sss space\sss $\mathbf{G}$\sss of\trs the bundle\sss
$\bm{\pi}\dff \colon\dff
\mathbf{G}
\qff \ttoo\qff 
\num{\mathcal{P}{\nsp}\hat{\mathcal{S}}}$\sss
is\dss contractible.\oss}

\vspace{6pt}
This immediately\sss follows\sss from\trs Theorems\qss C\qss and\qss D.\oss

\myuppar{Bott\dss periodicity.}
Theorem\qss E\qss implies\sss that\sss the fibers\sss
$\gr\trf(\trf P\trf)$\sss of\sss $\bm{\pi}$\sss are homotopy equivalent\sss
to\sss the\sss loop space\sss 
$\Omega\trf \num{\mathcal{P}{\nsp}\hat{\mathcal{S}}}$\nnsp.\oss
By\sss the very\sss definition\sss $\gr\trf(\trf P\trf)$\sss
is\dss homeomorphic\sss to\sss $\gr$\nnsp.\oss
As was mentioned above,\pss $\gr$\sss is\dss homotopy equivalent\sss to\sss
$\gr\trf(\dff \infty\dff)$\nnsp.\oss
But\sss $\gr\trf(\dff \infty\dff)$\sss
is\dss nothing else but\sss the  
classifying space $\zzz\dff \times\dff B\fff U$\nnsp.\oss
Here\sss $U\off =\off U\dff(\dff \infty\dff)$\sss is\dss the colimit\sss 
of\dss unitary groups\sss $U\dff(\dff n\trf)$\nnsp.\oss
The same methods show\sss that\sss $U\dff(\dff \infty\dff)$\sss is\dss
homotopy equivalent\sss to\sss $U\ffin$\dnsp.\oss
But\sss $U\ffin$\sss is\dss homeomorphic\sss to\sss
$\num{\hat{\mathcal{S}}}$\sss and\sss hence\dss 
is\dss homotopy equivalent\sss to\sss 
$\num{\mathcal{P}{\nsp}\hat{\mathcal{S}}}$\nnsp.\oss
It\sss follows\sss that\sss $\zzz\dff \times\dff B\fff U$\sss 
is\dss homotopy equivalent\sss to\sss
$\Omega\dff U$\nnsp.\oss
This\dss is\dss one of\dss the classical\sss forms of\qss Bott\dss periodicity.\oss

This\dss is\dss not\sss the simplest\sss 
proof\dss of\qss Bott\dss periodicity.\oss
Its\sss interest\dss is\dss in\sss relating\dss the functional\sss analytic\sss
methods of\trs Atiyah\sss and\sss Singer\qss \cite{as}\qss
with algebraic and\sss topological\sss ideas of\trs Quillen\qss \cite{q}.\oss
While\sss the definitions of\dss $\hat{\mathcal{S}}$\nsp\dnsp,\qss 
$\mathcal{P}{\nsp}\hat{\mathcal{S}}$\sss and\sss related categories are
motivated\sss by\sss the properties of\dss self-adjoint\trs
Fredholm\dss operators,\oss Fredholm\dss operators 
are not\sss used\sss in\sss this proof.\oss 
But\sss there\dss is\dss a similarity,\oss unexpected at\sss least\sss
to\sss the present\sss author,\oss between\sss the
proof\trs in\trs \cite{as}\qss and\sss this one,\oss
and\sss hence between\sss the methods of\qss
Atiyah--Singer\qss \cite{as}\qss and\dss Quillen\qss \cite{q}.\oss
A\sss key steps of\trs Atiyah--Singer\trs proof\dss
amounts\sss to proving\sss that\sss the map\sss
$\exp \pi\dff \colon\dff \hat{F}_{\fff *}\qff \ttoo\qff -\qff C$\sss
from\qss \cite{as},\oss Proposition\qss 3.3,\oss is\dss a quasi-fibration.\pss
See\sss the proof\dss of\trs Proposition\qss 3.5\qss in\qss \cite{as}.\oss
The key step in our proof\trs is\dss the fact\sss that\sss 
$\num{p}$\sss is\dss quasi-fibration.\oss
There\dss is\dss even a close analogue of\dss 
the\dss Atiyah--Singer\trs map\sss $\exp \pi$\nnsp,\oss
namely\sss the composition of\dss the map\sss
$\num{\hat{\mathcal{E}}}
\qff \ttoo\qff
\num{\hat{\mathcal{S}}}$\sss
from\trs Theorem\qss A\qss with\sss the homeomorphism\sss
$\num{\hat{\mathcal{S}}}
\qff \ttoo\qff
U\ffin$\sss
from\trs Theorem\qss \ref{harris-h}.\oss

The above proof\dss of\qss Bott\dss periodicity can\sss be somewhat\sss shortened\sss by\sss
using\trs Theorems\qss B\qss and\qss C\qss instead of\trs Theorem\qss E.\oss
One can go even\sss further and\sss replace\sss the Hilbert\dss space\sss $H$\sss
by\sss the union\sss $\ccc^{\dff \infty}$\sss of\dss the sequence\sss
$\ccc^{\dff 0}\qff \subset\pff
\ccc^{\dff 1}\qff \subset\pff
\ccc^{\trf 2}\qff \ldots \pff$\nnsp.\oss
The resulting\sss proof\trs is\dss not\sss included\sss in\sss this paper\sss
because similar\sss proofs are well\sss known.\oss 
It\dss is\dss fairly close\sss to\sss the proof\dss
of\qss McDuff\pss \cite{mc}\qss and especially\sss to one of\qss Giffen\qss \cite{g}.\oss
McDuff's\dss proof\trs is\dss an ultimate simplification of\qss
Atiyah--Singer\qss \cite{as}\qss proof,\oss 
but\sss the exposition\sss in\qss \cite{mc}\qss is\dss quite condensed.\oss
A detailed version was provided\dss by\dss Aguilar\dss and\dss Prieto\qss \cite{ap}\qss
and\sss then streamlined\sss by\dss Behrens\qss \cite{b1},\qss \cite{b2}.\oss
Giffen's\dss proof\trs
is\dss based on\sss ideas of\trs Quillen\qss \cite{q},\qss \cite{gr}\qss
and\dss Segal\qss \cite{s3}.

\myuppar{Topological\sss categories related\sss to\dss Fredholm\dss operators.}
All\sss these results have natural\sss analogues for general\qss
(not\sss necessarily\sss self-adjoint\halfff)\qss Fredholm\dss operators.\oss
Let\sss $\mathcal{F}$\sss be\sss the space of\qss
Fredholm\dss operators\sss $H\qff \ttoo\qff H$\nnsp.\oss
Recall\dss that\sss
$\num{A}
\off =\off
\sqrt{\dff A^* A\dff}$\nnsp.\oss
and\sss that\sss $\num{A}$\sss may be not\sss equal\sss to\sss $\num{A^*}$\dnsp.\oss
Recall\sss also\sss that\sss by\sss $W\dff \ominus\dff V$\sss we denote\sss
the orthogonal\sss complement\sss of\sss $V$\sss in\sss $W$\dnsp.\oss

We define an\qss \emph{enhanced\trs Fredholm\dss operator}\pss 
as a pair\sss $(\trf A\dff,\qff \varepsilon\trf)$\sss
consisting of\dss an operator\sss $A\qff \in\qff \mathcal{F}$\sss 
and a number\sss $\varepsilon\qff >\qff 0$\sss
such\sss that\sss 
$\varepsilon\qff \not\in\qff \sigma\dff(\trf \num{A}\trf)$\nnsp,\oss 
the interval\dss $[\trf 0\fff,\qff \varepsilon\trf]$\sss 
is\dss disjoint\sss from\sss the essential\sss spectrum of\sss $\num{A}$\nnsp,\oss
and\sss the same condition holds for\sss $A^*$\sss in\sss the role of\sss $A$\nnsp.\oss
Equivalently and,\oss perhaps,\oss more intuitively,\oss one can\sss require\sss that\sss
$\norm{A\trf(\trf v\trf)}\off \neq\off \varepsilon\dff \norm{v}$\sss
for every\sss non-zero\sss $v\qff \in\qff H$\nnsp,\oss
the space of\dss vectors\sss $v\qff \in\qff H$\sss
such\sss that\sss
$\norm{A\trf(\trf v\trf)}\off \leq\off \varepsilon\dff \norm{v}$\sss
is\dss finitely dimensional,\oss
and\sss the same condition holds for\sss $A^*$\sss in\sss the role of\sss $A$\nnsp.\oss
We equip\sss the set\sss $\mathcal{E}$\sss of\dss enhanced\dss Fredholm\dss operators\dss
with\sss topology induced\sss by\sss the\sss topology of\dss
$\mathcal{F}$\sss and\sss the discrete\sss topology on\sss $\rrr_{\qff >\dff 0}$\nsp.\oss
Similarly\sss to\sss $\hat{\mathcal{E}}$\dnsp,\oss
the space\sss $\mathcal{E}$\sss is\dss partially ordered,\oss
and\sss the partial\sss order\sss turns\sss $\mathcal{E}$\sss
into a\sss topological\sss category.\oss
There\dss is\dss an obvious forgetting functor\sss
$\mathcal{E}\qff \ttoo\qff \mathcal{F}$\dnsp.\oss

The analogue\sss $\mathcal{S}$\sss of\dss the category\sss $\hat{\mathcal{S}}$\sss
is\dss defined as follows.\oss
The objects of\sss $\mathcal{S}$\sss are pairs\sss
$(\trf E_{\dff 1}\fff,\qff E_{\trf 2}\trf)$ of\dss finitely\sss
dimensional\sss subspaces\sss of\sss $H$\nnsp.\oss
Morphisms\sss
$(\trf E_{\dff 1}\fff,\qff E_{\trf 2}\trf)
\qff \ttoo\qff 
(\trf E\fff'_{\dff 1}\fff,\qff E\fff'_{\trf 2}\trf)$\sss
are defined as\sss triples\sss
$(\trf F_{\dff 1}\fff,\pff F_{\dff 2}\dff,\pff f\trf)$\nnsp,\oss
where\sss
$F_{\dff 1}\dff F_{\trf 2}$
are subspaces of\sss $H$\sss
such\sss that\vspace{1.5pt}
\[
\quad
E_{\dff 1}\dff \oplus\dff F_{\dff 1}
\off =\off
E\fff'_{\dff 1}
\quad
\mbox{and}\quad
E_{\trf 2}\dff \oplus\dff F_{\trf 2}
\off =\off
E\fff'_{\trf 2}
\]

\vspace{-12pt}\vspace{1.5pt}
and\sss $f$\sss is\dss an\sss isometry\sss
$F_{\dff 1}\qff \ttoo\qff F_{\dff 2}$\nsp.\oss
The composition\dss is\dss defined\sss by\sss taking\sss the direct\sss sums
of\dss the corresponding subspaces\sss $F_{\dff 1}\fff,\qff F_{\trf 2}$
and of\dss the isometries.\oss
The category $\mathcal{S}$\sss is\dss a\sss topological\sss category\sss 
in\sss an obvious way.\oss
There\dss is\dss a natural\sss functor\sss
$\mathcal{E}\qff \ttoo\qff \mathcal{S}$\nsp\dnsp.\oss
See\dss Section\qss \ref{classifying-spaces-odd-saf}.\oss

Finally,\oss let\sss us\sss define\dss Quillen--Segal\dss category\sss $\hat{C}$\nnsp,\oss
which we will\sss from\sss now on denote by\sss $S$\nnsp.\oss
It\dss is\dss an abstract\sss analogue of\sss $\mathcal{S}$\dnsp.\oss
The objects of\sss $S$\sss are\sss pairs\sss $(\trf V_{\fff 1}\dff,\qff V_{\dff 2}\trf)$\sss
of\dss finitely dimensional\dss Hilbert\sss spaces over\sss $\ccc$\nnsp.\oss
A morphism\sss
$(\trf V_{\dff 1}\dff,\pff V_{\dff 2}\trf)
\qff \ttoo\qff 
(\trf W_{\dff 1}\dff,\pff W_{\dff 2}\trf)$\sss
of\sss $S$\sss
is\dss a\sss triple\sss
$(\trf i_{\trf 1}\dff,\pff i_{\trf 2}\dff,\pff g\trf)$\nnsp,\oss
where\dss
$i_{\trf 1}\dff \colon\dff V_{\fff 1}\qff \ttoo\qff W_{\fff 1}$\dss
and\dss
$i_{\trf 2}\dff \colon\dff V_{\fff 2}\qff \ttoo\qff W_{\fff 2}$\dss
are isometric embeddings and\vspace{1.5pt}
\[
\quad
g\dff \colon\dff
W_{\fff 1}\qff \ominus\qff i_{\trf 1}\dff(\trf V_{\fff 1}\trf)
\qff \ttoo\qff
W_{\fff 2}\qff \ominus\qff i_{\trf 2}\dff(\trf V_{\fff 2}\trf)
\]

\vspace{-12pt}\vspace{1.5pt}
is\dss an\sss isometry.\oss 
The composition\dss is\dss defined\sss in\sss an obvious way
and amounts\sss to\sss taking\sss the direct\sss sum of\dss
the corresponding\sss isometries\sss $g$\nnsp.\oss
There\dss is\dss an\sss intermediate category 
$S/\fff H$\sss
between\sss $\mathcal{S}$\sss and\sss $S$\sss
and\sss canonical\dss functors\sss
$\mathcal{S}
\off \longleftarrow\off 
S/\fff H
\qff \ttoo\qff
S$\dnsp.\oss
See\dss Section\qss \ref{classifying-spaces-odd-saf}.\oss

\myuppar{Theorem\qss F.}
\emph{The maps of\dss classifying\sss spaces}\vspace{1.5pt}
\[
\quad
\mathcal{F}
\off =\off
\num{\mathcal{F}}
\off \longleftarrow\off
\num{\mathcal{E}}
\qff \ttoo\qff
\num{\mathcal{S}}
\off \longleftarrow\off
\num{S/\fff H}
\qff \ttoo\qff
\num{S}
\]

\vspace{-12pt}\vspace{1.5pt}
\emph{induced\sss by\dss the\sss functors defined or mentioned above,\oss
are homotopy equivalences.\oss}

\vspace{6pt}
See\dss Theorems\qss \ref{forgetting-odd},\oss \ref{two-fredholm},\oss
and\qss \ref{intermediary-f}.\oss

\myuppar{Bundles and\sss quasi-fibrations related\sss to\dss Fredholm\dss operators.}
There are analogues of\dss polarizations,\oss splittings,\oss
restricted\dss Grassmannians,\oss etc.\qss for\dss Fredholm\dss operators.\oss
The category\sss $\mathcal{P}{\nsp}\mathcal{S}$\sss has as objects\sss triples\sss
$(\trf E_{\dff 1}\dff,\qff E_{\trf 2}\dff,\qff i\trf)$\sss
such\sss that\sss 
$(\trf E_{\dff 1}\dff,\qff E_{\trf 2}\trf)$\sss
is\dss an object\sss of\sss $\mathcal{S}$\sss and\sss
$i\dff \colon\dff
H\dff \ominus\dff E_{\dff 1}
\qff \ttoo\qff
H\dff \ominus\dff E_{\trf 2}$\sss
is\dss an\sss isometry.\oss
Morphisms of\sss $\mathcal{P}{\nsp}\mathcal{S}$\sss 
have\sss the form\vspace{1.5pt}
\[
\quad
\bigl(\trf E_{\dff 1}\dff,\off E_{\dff 2}\dff,\off i\trf\bigr)
\qff \ttoo\qff 
\bigl(\trf E_{\dff 1}\dff \oplus\dff F_{\dff 1}\dff,\off 
E_{\trf 2}\dff \oplus\dff i\trf(\trf F_{\dff 1}\trf)\dff,\off i\fff'\trf\bigr)
\pff,
\]

\vspace{-12pt}\vspace{1.5pt}
where\sss $F_{\dff 1}$\sss is\dss a finitely dimensional\sss subspace of\sss
$H\dff \ominus\dff E_{\dff 1}$\sss and\sss 
$i'$\sss is\dss induced\sss by\sss $i$\nnsp.\oss
There\dss is\dss a canonical\sss functor\sss 
$\mathcal{P}{\nsp}\mathcal{S}\qff \ttoo\qff \mathcal{S}$\sss
taking\sss
$(\trf E_{\dff 1}\dff,\qff E_{\trf 2}\dff,\qff i\trf)$\sss
to\sss 
$(\trf E_{\dff 1}\dff,\qff E_{\trf 2}\trf)$\nnsp.\oss
One can easily see\sss that\sss
$\num{\mathcal{P}{\nsp}\mathcal{S}}\qff \ttoo\qff \num{\mathcal{S}}$\sss
is\dss a homotopy equivalence.\oss
Like\sss the category\sss $\mathcal{P}{\nsp}\hat{\mathcal{S}}$\nnsp,\oss
the category\sss $\mathcal{P}{\nsp}\mathcal{S}$\sss
can\sss be defined also in\sss terms of\dss a 
partial\sss order\sss on\sss the space of\dss objects.\oss

The analogue of\dss the restricted\dss Grassmannian\sss $\gr$\sss
is\dss the group\sss $U\ffin$\sss which we already encountered,\oss
but\sss with\sss the norm\sss topology\sss instead of\dss the colimit\sss one.\oss
There are also analogues of\dss spaces\sss $\gr\trf(\trf P\trf)$\nnsp.\oss
For an object\sss
$P\off =\off (\trf E_{\dff 1}\dff,\qff E_{\trf 2}\dff,\qff i\trf)$\sss
of\trs $\mathcal{P}{\nsp}\mathcal{S}$\sss
let\sss $U\ffin\dff(\trf P\trf)$\sss
be\sss the space of\dss isometries\sss
$H\qff \ttoo\qff H$\sss
equal\sss to\sss $i$\sss on a subspace of\dss finite codimension\sss 
in\sss $H\dff \ominus\dff E_{\dff 1}$\sss
considered\sss with\sss the norm\sss topology.\oss
Clearly,\pss $U\ffin\dff(\trf P\trf)$\sss is\dss non-empty\dss if\trs and\dss only\trs if\dss
$\dim\trf E_{\dff 1}\off =\off \dim\trf E_{\trf 2}$\nsp.\oss
Let\sss $\mathcal{P}{\nsp}\mathcal{S}_{\qff 0}$\sss 
be\sss the full\sss subcategory\sss of\sss
$\mathcal{P}{\nsp}\mathcal{S}$\sss
having as objects\sss triples\sss
$(\trf E_{\dff 1}\dff,\qff E_{\trf 2}\dff,\qff i\trf)$\sss
such\sss that\sss
$\dim\trf E_{\dff 1}\off =\off \dim\trf E_{\trf 2}$\nsp.\oss
Obviously,\oss
if\dss there\dss is\dss a morphism\sss
$P\dff \ttoo\qff P\fff'$\nsp\dnsp,\oss then\sss
$U\ffin\dff(\trf P\trf)\off =\off U\ffin\dff(\trf P\fff'\trf)$\nnsp.\oss
This fact\sss allows\sss to construct\sss a\sss 
locally\sss trivial\dss bundle\vspace{1.5pt}\vspace{0.9pt}
\[
\quad
\bm{\pi}\dff \colon\dff
\mathbf{U}
\qff \ttoo\qff 
\num{\mathcal{P}{\nsp}\mathcal{S}_{\qff 0}}
\pff
\]

\vspace{-12pt}\vspace{1.5pt}\vspace{0.9pt}
with\sss the fiber\sss $U\ffin\dff(\trf P\trf)$\sss 
over simplices\sss having\sss $P$\sss as a vertex.\oss
See\dss Section\qss \ref{unitary-quasi-fibrations}.\oss

Let\sss $\mathcal{U}$\sss be\sss the full\sss subcategory\sss of\dss
$\mathcal{P}{\nsp}\mathcal{S}$\sss having as objects\sss
triples of\dss the form\sss $(\trf 0\fff,\qff 0\fff,\qff i\trf)$\nnsp.\oss
The category\sss $\mathcal{U}$\sss has only\sss identity\sss morphisms,\oss
and\sss $\num{\mathcal{U}}$\sss is\dss
homeomorphic\sss to\sss the unitary group of\dss a\sss Hilbert\sss space,\oss
and\sss hence\dss is\dss contractible by\trs Kuiper\trs theorem.\oss
A\qss \emph{splitting}\pss of\dss an object\sss $M$\sss of\sss
$\mathcal{P}{\nsp}\mathcal{S}$\sss is\dss a morphism\sss
$s\dff \colon\dff N\qff \ttoo\qff M$\sss such\sss that\sss $N$\sss
is\dss an object\sss of\sss $\mathcal{U}$\dnsp.\oss
Clearly,\oss a splitting of\dss $M$\sss exists only\dss if\dss $M$\sss
is\dss an object\sss of\sss $\mathcal{P}{\nsp}\mathcal{S}_{\qff 0}$\nsp.\oss

The\sss topological\sss category\sss
$s\dff \mathcal{S}_{\qff 0}$\sss having as objects 
splittings of\dss objects of\sss $\mathcal{P}{\nsp}\mathcal{S}_{\qff 0}$\sss
and\sss the forgetting\sss functor\sss
$\phi\dff \colon\dff
s\dff \mathcal{S}_{\qff 0}
\qff \ttoo\qff
\mathcal{P}{\nsp}\mathcal{S}_{\qff 0}$\sss
are defined exactly as in\sss the self-adjoint\sss case,\oss
as also\sss the category\sss $S\dff(\trf \phi\trf)$\nnsp,\oss
the\sss forgetting\sss functor\sss
$p\dff \colon\dff
S\dff(\trf \phi\trf)
\qff \ttoo\qff
\mathcal{P}{\nsp}\mathcal{S}_{\qff 0}$,\oss
and\sss the continuous map\vspace{1.5pt}\vspace{0.9pt}
\[
\quad
\num{p}\dff \colon\dff
\num{S\dff(\trf \phi\trf)}
\qff \ttoo\qff
\num{\mathcal{P}{\nsp}\mathcal{S}_{\qff 0}}
\pff.
\]

\vspace{-12pt}\vspace{1.5pt}\vspace{0.9pt}
This\dss is\dss our categorical\dss model\sss of\dss the bundle\sss
$\bm{\pi}\dff \colon\dff
\mathbf{U}
\qff \ttoo\qff 
\num{\mathcal{P}{\nsp}\mathcal{S}_{\qff 0}}$\nnsp.\oss
The main properties of\dss this model\sss are\sss the same as in\sss the self-adjoint\sss case.\oss
They are summarized\sss in\sss the following\sss
theorem and\sss proved\sss in\dss Section\qss \ref{unitary-quasi-fibrations}.\oss

\myuppar{Theorem\qss G.}
\emph{The map\sss
$\num{p}\dff \colon\dff
\num{S\dff(\trf \phi\trf)}
\qff \ttoo\qff
\num{\mathcal{P}{\nsp}\mathcal{S}_{\qff 0}}$\sss
is\dss a quasi-fibration,\oss
there\dss is\dss a natural\dss homotopy equivalence\dss
$h\dff \colon\dff
\num{S\dff(\trf \phi\trf)}
\qff \ttoo\qff
\mathbf{U}$\nnsp,\oss
and\dss the spaces\dss $\num{S\dff(\trf \phi\trf)}$\sss
and\trs $\mathbf{U}$\sss are contractible.\oss}

\myuppar{The other\sss half\dss of\qss Bott\dss periodicity.}
The above proof\dss of\trs Bott\dss periodicity,\oss
like most\sss of\dss other proofs,\oss
assumed\sss as a known\sss fact\sss that\sss
$\gr\trf(\dff \infty\dff)$\sss
is\dss the classifying space $\zzz\dff \times\dff B\fff U$\nnsp.\oss
In order\sss to prove\sss this one needs\sss 
a\sss locally\sss trivial\sss bundle over\sss a space homotopy equivalent\sss 
to a connected\sss component\sss of\sss $\gr\trf(\dff \infty\dff)$\sss
with contractible\sss total\sss space and\sss fibers homotopy equivalent\sss
to\sss $U\off =\off U\dff(\dff \infty\dff)$\nnsp.\oss
The\sss total\sss space of\dss the bundle\sss 
$\bm{\pi}\dff \colon\dff
\mathbf{U}
\qff \ttoo\qff 
\num{\mathcal{P}{\nsp}\mathcal{S}_{\qff 0}}$\sss
is\dss contractible by\trs Theorem\qss G,\oss
and\sss its\sss fibers are homotopy equivalent\sss to\sss $U\ffin$\sss
and\sss hence\sss to\sss $U\dff(\dff \infty\dff)$\nnsp.\oss
Its base\sss $\num{\mathcal{P}{\nsp}\mathcal{S}_{\qff 0}}$\sss
is\dss indeed\sss homotopy equivalent\sss to a connected component\sss 
of\sss $\gr\trf(\dff \infty\dff)$\nnsp.\oss
Moreover,\pss $\num{\mathcal{P}{\nsp}\mathcal{S}}$\sss
is\dss homotopy equivalent\sss to\sss $\gr\trf(\dff \infty\dff)$\nnsp.\oss

But\dss it\sss seems\sss that\sss there\dss is\dss no proof\dss of\dss this
fact\sss which\dss is\dss as nice as\sss the proof\dss that\sss 
$\num{\mathcal{P}{\nsp}\hat{\mathcal{S}}}$\sss 
is\dss homotopy equivalent\sss to\sss 
$U\dff(\dff \infty\dff)$\nnsp,\oss
because\sss there\dss is\dss no analogue of\dss the
homeomorphism between\sss $\num{\hat{\mathcal{S}}}$\sss and\sss $U\ffin$\dnsp.\oss
A similar difficulty\sss was encountered\sss by\sss Harris\qss \cite{h},\oss Section\qss 4.\oss
Following\sss his example,\oss we omit\sss the details.\oss\vspace{0.45pt}

\myuppar{The\sss homotopy equivalence\sss 
$\hat{\mathcal{F}}\qff \ttoo\qff \Omega\trf \mathcal{F}$\dnsp.}
Theorem\qss G\qss implies\sss that\sss the\sss loop space\sss
$\Omega\trf \mathcal{F}$\sss is\dss canonically\sss homotopy equivalent\sss
to\sss $U\ffin$\dnsp,\oss and\sss hence\sss to\sss $\num{\hat{\mathcal{S}}}$\sss
and\sss $\hat{\mathcal{F}}$\nnsp.\oss
The homotopy equivalence of\dss the spaces\sss $\hat{\mathcal{F}}$\sss
and\sss $\Omega\trf \mathcal{F}$\sss is\dss one of\dss the main\sss results
of\trs Atiyah\dss and\dss Singer\qss \cite{as}.\oss
In\sss fact,\oss Atiyah\dss and\dss Singer\qss \cite{as}\qss
used a simple\sss formula\sss to construct\sss a map\sss
$\hat{\mathcal{F}}\qff \ttoo\qff \Omega\trf \mathcal{F}$\sss
and\sss proved\sss that\sss this map\dss is\dss a homotopy equivalence.\oss
The proof\dss based on\dss Theorem\qss G\qss also\sss leads\sss to a map\sss
$\hat{\mathcal{F}}\qff \ttoo\qff \Omega\trf \mathcal{F}$\dnsp,\oss 
which\dss is\dss well\sss defined\sss up\sss to homotopy.\oss\vspace{0.45pt}

\myuppar{Theorem\qss H.}
\emph{These\sss two homotopy equivalences\sss
$\hat{\mathcal{F}}\qff \ttoo\qff \Omega\trf \mathcal{F}$\sss
are\sss the same up\sss to homotopy.\oss}\vspace{0.45pt}

\vspace{6pt}
See\dss Theorem\qss \ref{alpha-alpha-beta-beta}.\oss
Surprisingly,\oss one needs\sss to go
fairly deep into\dss Atiyah--Singer\dss proof\dss in order\sss to prove\sss
that\sss these maps are homotopic.\oss
A key\sss role\sss in\dss Atiyah--Singer\dss proof\trs is\dss played\sss by\sss
the group\sss $G$\sss of\dss unitary elements of\trs Calkin\dss algebra.\oss
It\dss turns out\sss that\sss the classifying space\sss
$\num{\mathcal{P}{\nsp}\mathcal{S}}$\sss
may serve as a categorical\dss model\sss of\sss $G$\nnsp.\oss
In\sss more details,\oss there\dss is\dss a natural\sss map\sss
$\num{\mathcal{P}{\nsp}\mathcal{S}}\qff \ttoo\qff G$\nnsp,\oss
and\sss this map\dss is\dss a homotopy equivalence.\oss
See\dss Corollary\qss \ref{kappa-equivalence}.\oss\vspace{0.45pt}

\myuppar{Hilbert\dss bundles.}
A big\sss part\sss of\dss this\sss theory easily extends\sss from a fixed\dss
Hilbert\dss space\sss $H$\sss to\sss Hilbert\dss bundles,\oss 
thought\sss as families $H_{\dff x}\dff,\qff x\qff \in\qff X$ 
of\dss Hilbert\sss spaces parameterized\sss by
a\sss topological\sss space\sss $X$\nnsp.\oss
More precisely,\oss for a\dss Hilbert\dss bundle\sss $\mathbb{H}$\sss
one can define analogues of\dss categories not\sss involving\dss
Hilbert\sss space operators or\sss polarizations.\oss
The categories and\sss topological\sss spaces involving\dss Hilbert\dss space operators
can\sss be defined only\sss when\sss $\mathbb{H}$\sss is\dss trivial\sss
and\sss trivialized,\oss i.e.\qss when one deals with operators
in a fixed\dss Hilbert\sss space.\oss\vspace{0.45pt}

For example,\oss the analogue\sss
$\hat{\mathcal{S}}\dff(\trf \mathbb{H}\trf)$\sss
of\dss the category\sss $\hat{\mathcal{S}}$\sss has as objects
finitely dimensional\sss subspaces of\dss the fibers 
$H_{\dff x}\dff,\qff x\qff \in\qff X$\sss of\sss $\mathbb{H}$\nnsp.\oss
Morphisms\sss $V\qff \ttoo\qff V\fff'$ exists only\dss if\sss
$V\fff,\qff V\fff'$\sss are contained\sss in\sss the same fiber\sss $H_{\dff x}$\nsp,\oss
and\sss in\sss this case morphisms are defined exactly as morphisms of\sss
$\hat{\mathcal{S}}$\sss with\sss $H$\sss replaced\sss by\sss $H_{\dff x}$\nsp.\oss
The definition of\sss $Q$\sss does not\sss involve\sss $H$\sss at\sss all,\oss
and\sss by\sss this reason\sss the analogue\sss 
$Q\dff(\trf \mathbb{H}\trf)$\sss of\sss $Q$\sss depends on\sss $X$\sss
but\sss otherwise\dss is\dss independent\sss from\sss $\mathbb{H}$\nnsp.\oss
All\sss this\dss is\dss explained\sss in\dss Section\qss \ref{hilbert-bundles}.\oss
The situation\sss in\sss the non-self-adjoint\sss case\dss 
is\dss completely\sss similar.\oss
There are no surprises,\oss and\sss this case\dss 
is\dss left\sss to\sss the reader.\oss

\myuppar{The analytical\dss index.}
The partial\sss extension of\dss the\sss theory\sss to\sss
Hilbert\dss bundles\dss is\dss motivated\sss by\sss applications\sss
to\sss the index\sss theory,\oss
discussed\sss in\sss the companion\sss paper\qss \cite{i}.\oss
While\sss there\dss is\dss a satisfactory definition of\dss
the analytical\dss index\sss for\sss families of\dss general\dss
Fredholm\dss operators,\oss this\dss is\dss not\sss quite\sss the case
for\sss self-adjoint\dss Fredholm\dss operators.\oss
The standard definition,\oss
going\sss back\sss to\dss Atiyah\dss and\dss Singer\qss \cite{as5},\oss \cite{as}\qss
and\sss to\dss Atiyah,\pss Patodi,\oss and\dss Singer\qss \cite{aps},\oss 
requires\sss the bundle\sss $\mathbb{H}$\sss
to be\sss trivial\sss and even\sss trivialized.\oss
Using categories\sss $\hat{\mathcal{S}}\dff(\trf \mathbb{H}\trf)$\sss
and\sss $Q\dff(\trf \mathbb{H}\trf)$\sss one can\sss give a natural\sss
definition of\dss the analytical\dss index\sss for\sss general\sss bundles
and under\sss minimal\sss continuity assumptions about\sss the family.\oss
Theorem\qss H\qss plays a key\sss role\sss in\sss proving\sss
that\sss this definition of\dss the analytical\dss index\sss agrees\sss
with\sss the\dss Atiyah--Singer\dss one.\oss

One can also define partial\sss analogues of\dss categories
involving\sss polarizations.\oss
This clarifies some subtle aspects of\trs
the notion of\qss \emph{spectral\dss section}\pss
of\qss Melrose\dss and\dss Piazza\qss \cite{mp}.\oss

\newpage
\mysection{Simplicial\qss spaces}{simplicial-spaces}

\myuppar{Compactly\sss generated\sss spaces.}
A\sss topological\sss space $X$\sss is\dss said\sss to be\qss
\emph{compactly\sss generated}\pss if\dss a subset\sss of\sss $X$\sss
intersecting every compact\sss subset\sss in\sss closed\sss subset\dss
is\dss itself\sss closed.\oss
While working with simplicial\sss spaces,\oss
all\sss topological\sss spaces are usually\sss assumed\sss to be
compactly\sss generated,\oss and\sss the definition of\dss products\dss
is\dss modified\sss in such a way\sss that\dss products are also compactly generated.\oss
We will\dss largely\sss follow\sss this\sss tradition,\oss
which goes back\sss to\dss Segal\qss \cite{s1}.\oss
In\sss the context\sss of\dss the present\sss paper\sss this\dss is\dss hardly
a restriction.\oss 
In fact,\oss the class of\dss compactly generated spaces\dss is\dss very\sss broad.\oss
It\sss includes spaces satisfying\sss the first\sss countability axiom,\oss
in\sss particular\sss metric spaces,\oss
as also locally\sss compact\sss spaces.\oss
See\qss \cite{ke},\oss Theorem\qss 7.13.\oss

In general,\oss working\sss with\sss compactly generated spaces requires\sss to 
redefine\sss the\sss topology of\dss the products,\oss
but\dss this\dss is\dss again\sss hardly an\sss issue in\sss the context\sss
of\dss the present\sss paper.\oss
For example,\oss the product\sss of\dss a compactly generated space
and\sss a\sss locally compact\sss space\dss is\dss automatically compactly generated,\oss
as\dss is\dss the product\sss of\dss two metric spaces.\oss
At\sss the same\sss time\sss it\dss is\dss well\sss known\sss 
that\sss the usual\sss topological\sss
constructions preserve\sss the property of\dss being compactly generated,\oss 
if\dss the products are redefined.\oss
See\qss \cite{st}.\oss

\myuppar{Simplicial\sss spaces and $\Delta$\dnsp-spaces.}
Let\sss $\bm{\Delta}$ be\sss the category having sets\sss
$[\halfff n\dff]\off =\off\{\trf 0\fff,\qff 1\fff,\qff \ldots\fff,\qff n\trf\}$\sss
as objects and\sss non-decreasing maps\sss
$[\halfff m\dff]\qff \ttoo\qff [\halfff n\dff]$
as morphisms.\oss
A\qss \emph{simplicial\sss set}\pss is\dss a contravariant\sss functor\sss
from $\bm{\Delta}$\sss to\sss the category of\dss sets,\oss
and a\qss \emph{simplicial\sss space}\pss is\dss a contravariant\sss functor\sss
from $\bm{\Delta}$\sss to\sss the category of\dss topological\sss spaces.\oss
A simplicial\sss set\sss can\sss be considered as a simplicial\sss space\sss
if\dss we equip sets with\sss the discrete\sss topology.\oss
For a simplicial\sss set\sss or\sss space\sss $K$\sss we will\sss denote by\sss $K_{\dff n}$\sss
the value of\dss the functor $K$ on\sss the object\sss $[\halfff n\dff]$\nnsp,\oss
and\sss for every non-decreasing map\sss
$\theta\dff \colon\dff
[\halfff m\dff]\qff \ttoo\qff [\halfff n\dff]$\sss
we will\sss denote by\sss
$\theta^{\dff *}\dff \colon\dff
K_{\dff n}\qff \ttoo\qff K_{\dff m}$\sss
the value of\dss $K$\sss on\sss $\theta$\nnsp.\oss
The points of\sss $K_{\dff n}$ are called\sss the\sss \emph{$n$\dnsp-simplices}\pss
of\sss $K$\nnsp,\oss and\sss the maps $\theta^{\dff *}$ the\qss
\emph{structure maps}\pss of\dss $K$\nnsp.\oss
The\qss \emph{simplicial\dss maps}\qss
$K\qff \ttoo\qff L$\sss are defined as\sss the natural\sss transformations of\dss functors.\oss

There\dss is\dss a version of\dss these notions which appears\sss to be simpler.\oss
Let\sss $\Delta$ be\sss the category having\sss the same objects as\sss $\bm{\Delta}$\nnsp,\oss
but\sss only strictly\sss increasing maps\sss
$[\halfff m\dff]\qff \ttoo\qff [\halfff n\dff]$
as morphisms.\oss
A\dss \emph{$\Delta$\dnsp-set}\pss is\dss a contravariant\sss functor from $\Delta$\sss
to\sss the category of\dss sets,\oss
and a\dss \emph{$\Delta$\dnsp-space}\pss is\dss
a contravariant\sss functor from $\Delta$\sss
to\sss the category of\dss topological\sss spaces.\oss
A simplicial\sss set\sss or space defines a $\Delta$\dnsp-set\sss
or $\Delta$\dnsp-space respectively\sss by\sss restricting\sss
functors\sss to $\Delta$\nnsp.\oss
While $\Delta$\dnsp-sets and $\Delta$\dnsp-spaces
are actually simpler in some respects\sss than\sss
simplicial\sss sets and spaces,\oss
the\sss theories of\dss simplicial\sss things are\sss the right\sss
ones,\oss and we will\sss use $\Delta$\dnsp-things only occasionally.\oss

\myuppar{Geometric realizations.}
Let\sss $\Delta^n$\sss
be\sss the\qss \emph{standard\dss geometric $n$\dnsp-simplex}\vspace{-0.8pt}
\[
\quad
\Delta^n
\off =\off
\left\{\off (\trf t_{\dff 0}\dff,\qff t_{\dff 1}\dff,\qff \ldots\dff,\qff t_{\dff n}\trf)
\qff \in\qff
\rrr^{\dff n\dff +\dff 1}
\off \bigl|\off
t_{\dff 0}\qff +\qff t_{\dff 1}\qff +\qff \ldots\qff +\qff t_{\dff n}
\off =\off 1\qff,\off
t_{\dff 0}\dff,\qff t_{\dff 1}\dff,\qff \ldots\dff,\qff t_{\dff n}
\qff \geq\qff
0
\off\right\}
\pff.
\]

\vspace{-12pt}\vspace{-0.8pt}
Every non-decreasing map\sss
$\theta\dff \colon\dff
[\halfff m\dff]\qff \ttoo\qff [\halfff n\dff]$\sss
defines a\sss linear\sss map\sss
$\theta_{\dff *}\dff \colon\dff
\rrr^{\dff m\dff +\dff 1}\qff \ttoo\qff \rrr^{\dff n\dff +\dff 1}$\sss
taking\sss the $i${\dnsp}th vector of\dss 
the standard\sss basis of\sss $\rrr^{\dff m\dff +\dff 1}$\sss
to\sss the $\theta\trf(\dff i\trf)${\dnsp}th vector of\dss
the standard\sss basis of\sss $\rrr^{\dff n\dff +\dff 1}$\dnsp.\oss
This map\sss induces a map\sss $\Delta^m\qff \ttoo\qff \Delta^n$\dnsp,\oss
also denoted\sss by\sss $\theta_{\dff *}$\nsp.\oss
The\qss \emph{geometric\sss realization}\qss $\num{K}$ of\dss
a simplicial\sss space $K$\sss is\dss defined as\sss the quotient\sss
of\dss the disjoint\sss union\vspace{1.5pt}
\begin{equation}
\label{disjoint-union}
\quad
\coprod\nolimits_{\qff n\qff =\qff 0\fff,\dff 1\fff,\dff \ldots}\pff K_{\dff n}\dff \times\dff \Delta^n
\end{equation}

\vspace{-12pt}\vspace{1.5pt}
by\sss the equivalence relation\dss $\sim$\dss generated\sss by\vspace{1.5pt}
\[
\quad
\bigl(\trf \sigma\fff,\qff \theta_{\dff *}\dff(\trf x\trf)\trf\bigr)
\off \sim\off
\left(\trf \theta^{\dff *}\dff(\trf \sigma\trf)\fff,\qff x\trf\right)
\qff,
\]

\vspace{-12pt}\vspace{1.5pt}
where\sss $\sigma\qff \in\qff K_{\dff n}$\nsp,\pss
$x\qff \in\qff \Delta^m$\dnsp,\oss
and\sss
$\theta\dff \colon\dff
[\halfff m\dff]\qff \ttoo\qff [\halfff n\dff]$\sss
is\dss a non-decreasing map,\oss
i.e.\qss a morphism of\dss the category $\bm{\Delta}$\nnsp.\oss
For\sss $\sigma\qff \in\qff K_{\dff n}$\sss the image of\sss
$\sigma\dff \times\dff \Delta^n$\sss in\sss $\num{K}$\sss is\dss
denoted\sss by\sss $\num{\sigma}$\nnsp.\oss
Clearly,\oss a simplicial\sss map\sss
$K\qff \ttoo\qff L$\sss
defines a map\sss
$\num{K}\qff \ttoo\qff \num{L}$\nnsp,\oss
and\sss in\sss this way we get\sss a functor from 
simplicial\sss spaces\sss to\sss to\sss topological\sss spaces.\oss
This functor\sss respects\sss the\sss products in\sss the sense\sss that\sss
$\num{K\dff \times\dff L}
\off =\off 
\num{K}\dff \times\dff \num{L}$\nnsp,\oss
where\sss $K\dff \times\dff L$\sss is\dss defined dimension-wise,\oss
i.e.\dss 
$(\trf K\dff \times\dff L\trf)_{\dff n}
\off =\off
K_{\dff n}\dff \times\dff L_{\dff n}$\nsp.

The\qss \emph{geometric\sss realization}\qss of\dss
a $\Delta$\dnsp-space is\dss defined\sss in\sss
the same way,\oss except\sss only\sss strictly\sss increasing maps $\theta$
are involved.\oss
We will\sss denote\sss the geometric realization of\dss a $\Delta$\dnsp-space $K$\sss
by\sss $\num{K}_{\dff \Delta}$\nsp.\oss
Since $\Delta$\sss is\dss a subcategory of\sss $\bm{\Delta}$\nnsp,\oss
every simplicial\sss space $K$ defines a $\Delta$\dnsp-space,\oss
which we will\sss denote by\sss $\Delta\dff K$\nnsp.\oss
The geometric realization\sss $\num{\Delta\dff K}_{\dff \Delta}$\sss of\sss
$\Delta\dff K$\sss is\dss often called\sss the\qss
\emph{fat\sss geometric realization}\pss of\sss $K$ and\dss is\dss denoted\sss by\sss
$\norm{K}$\nnsp.\oss
There\dss is\dss a natural\sss quotient\sss
map\sss $\norm{K}\qff \ttoo\qff \num{K}$\nnsp,\oss
which\dss ia\dss a homotopy equivalence under some mild assumptions about\sss $K$\nnsp.\oss

\myuppar{Geometric realizations and\dss level-wise homotopy equivalences.}
Suppose\sss that\dss $K\fff,\qff L$\sss are simplicial\sss spaces and\sss
$f\dff \colon\dff K\qff \ttoo\qff L$\sss
is\dss a simplicial\sss map.\oss
The map\sss $f$\sss is\dss a\qss \emph{level-wise homotopy equivalence}\pss
if\trs the corresponding maps\sss 
$K_{\dff n}\qff \ttoo\qff L_{\dff n}$\sss
are homotopy equivalences.\oss
We would\sss like\sss to be able\sss to deduce from\sss this\sss that\sss the map\sss
$\num{f}\dff \colon\dff \num{K}\qff \ttoo\qff \num{L}$\sss
is\dss also a homotopy equivalence.\oss
Actually,\oss such arguments are going\sss to be one of\dss our\sss main\sss tools.\oss
Similarly,\oss suppose\sss that\sss the maps\sss
$K_{\dff n}\qff \ttoo\qff L_{\dff n}$\sss
are weak\sss homotopy equivalences.\oss
In\sss this case we would\sss like\sss to be able\sss to 
deduce from\sss this\sss that\sss the map\sss
$\num{f}\dff \colon\dff \num{K}\qff \ttoo\qff \num{L}$\sss
is\dss also a weak\sss homotopy equivalence.\oss
Unfortunately,\oss these implications are not\sss true 
without\sss additional\sss assumptions
about $K$ and\sss $L$\nnsp,\oss
and\sss this\dss is\dss a somewhat\sss annoying aspect\sss of\dss
working with simplicial\sss spaces.\oss

There are\sss two standard\sss ways\sss to deal\sss with\sss this difficulty.\oss
First,\oss one can use\sss the fat\sss geometric realization\sss $\norm{\bullet}$\sss
instead of\sss $\num{\bullet}$\nnsp.\oss
In\sss fact,\oss if\dss $f$\sss is\dss a\sss
level-wise homotopy equivalence,\oss then\sss
$\norm{f}\dff \colon\dff \norm{K}\qff \ttoo\qff \norm{L}$\sss
is\dss also a homotopy equivalence.\oss
See\sss \cite{s3},\oss Proposition\qss A.1{\fff}(ii).\oss
The same\dss is\dss true for\sss weak\sss homotopy equivalences.\oss
Unfortunately,\pss $\norm{K}$\sss is\dss a much\sss bigger space\sss than one expects,\oss
even\sss when\sss $K\off =\off [\halfff n\dff]$\nnsp.\oss

The second\sss way\dss is\dss to impose some\qss ``niceness''\qss conditions,\oss 
ensuring\sss that\sss the above implications hold,\oss
on\sss the considered simplicial\sss spaces.\oss
The simplicial\sss spaces considered\sss in\sss this paper are,\oss
in a sense,\oss given,\oss and\sss we cannot\sss replace\sss them\sss by\sss
better ones.\oss
Fortunately,\oss it\dss turns out\sss that\sss they\sss have\sss a very\sss strong\sss
niceness property,\oss to be defined\sss in\sss the next\sss subsection.\oss
Moreover,\oss their\sss geometric realizations admit\sss a\sss simpler
description\sss than\sss in\sss general.\oss

\myuppar{Simplicial\sss spaces with\sss free degeneracies.}
For a simplicial\sss space $K$\sss let\sss\vspace{1.5pt} 
\[
\quad
N_{\dff n}
\off =\off
N_{\dff n}\dff(\trf K\trf)
\off =\off
K_{\dff n}\qff \smallsetminus\off \bigcup\pff \theta^{\dff *}\dff(\trf K_{\dff k}\trf)
\pff,
\]

\vspace{-12pt}\vspace{1.5pt}
where\sss the union\dss is\dss taken over all\sss morphisms\sss
$\theta\dff \colon\dff
[\dff n\trf]\qff \ttoo\qff [\dff k\trf]$\sss
with\sss $n\qff >\qff k$\nnsp,\oss
or,\oss equivalently,\oss
over all\sss surjective morphisms\sss
$\theta\dff \colon\dff
[\dff n\trf]\qff \ttoo\qff [\dff k\trf]$\sss
with\sss $n\qff >\qff k$\nnsp.\oss
The $n$\dnsp-simplices belonging\sss to\sss 
$N_{\dff n}\dff(\trf K\trf)$\sss are called\dss the\qss
\emph{non-degenerate $n$\dnsp-simplices}\pss of\sss $K$\nnsp,\oss
and\sss other $n$\dnsp-simplices are called\dss the\qss \emph{degenerate $n$\dnsp-simplices}.\oss
The simplicial\sss space\sss $K$\sss is\dss said\sss to be\qss
\emph{split},\oss or\sss 
to have\qss \emph{free degeneracies},\oss
if\dss structure maps\sss 
$\theta^{\dff *}\dff \colon\dff
K_{\dff k}\qff \ttoo\qff K_{\dff n}$\sss induce a homeomorphism\vspace{1.5pt}
\begin{equation}
\label{free-d}
\quad
\coprod\nolimits_{\qff \theta}\qff N_{\dff k}
\off \ttoo\off
K_{\dff n}
\pff,
\end{equation}

\vspace{-12pt}\vspace{1.5pt}
where\sss the disjoint\sss union\dss is\dss taken over all\sss surjective maps\sss
$\theta\dff \colon\dff
[\dff n\trf]\qff \ttoo\qff [\dff k\trf]$\sss
and\sss to each such $\theta$ corresponds its own copy of\dss $N_{\dff k}$\nsp.\oss
The decompositions\qss (\ref{free-d})\qss of\dss spaces $K_{\dff n}$\sss 
into disjoint\sss unions\dss
is\dss called\sss the\qss \emph{splitting}\pss of\dss $K$\nnsp.\oss
These notions\sss go back\sss to\dss Artin\dss and\trs Mazur\qss \cite{am}\qss
and\sss were used\sss in\sss the context\sss of\dss simplicial\sss spaces by\trs
Dugger\sss and\dss Isaksen\qss \cite{di}.\oss

\myuppar{Skeletons.}
By a classical\sss lemma of\qss Eilenberg\dss and\trs Zilber\qss \cite{ez}\qss
every simplex $\sigma$\sss of\sss a simplicial\sss space\sss $K$
admits a unique presentation\sss in\sss the form\sss
$\sigma\off =\off \theta^{\dff *}\dff(\trf \tau\trf)$\sss
with surjective $\theta$ and\sss non-degenerate $\tau$\dss
({\fff}the\sss topological\sss structure\dss is\dss irrelevanlt\sss here).\oss
We will\sss call\sss such\sss presentation\sss the\qss 
\emph{Eilenberg--Zilber\dss presentation}\pss of\sss $\sigma$\nnsp.\oss
The\sss \emph{$n${\dnsp}th skeleton}\qss 
$\ssk_{\dff n}\dff K$\sss
of\dss $K$\sss 
is\dss the simplicial\sss subspace of\sss $K$\sss
having as\sss its $m$\dnsp-simplices\sss the $m$\dnsp-simplices $\sigma$ of\sss $K$\sss
such\sss that\sss in\sss the\dss Eilenberg--Zilber\dss presentation\sss
$\sigma\off =\off \theta^{\dff *}\dff(\trf \tau\trf)$\sss
the simplex $\tau$\sss is\dss a $k$\dnsp-simplex with\sss $k\qff \leq\qff n$\nnsp.\oss

Clearly,\pss $\num{K}$\sss is\dss equal\sss to\sss 
the union of\dss subspaces $\num{\nsp\ssk_{\dff n}\dff K}$\nnsp.\oss
Moreover,\oss the\sss topology of\sss $\num{K}$\sss is\dss the same as\sss
the direct\sss limit\sss topology of\dss this union.\oss
Indeed, a map from\sss the direct\sss limit\sss to a\sss topological\sss space $X$\sss
is\dss continuous\sss if\trs and\dss only\trs if\dss its restriction\sss to each
subspace $\num{\nsp\ssk_{\dff n}\dff K}$\sss is\dss continuous,\oss
and\sss hence\sss if\trs and\dss only\trs if\dss the induced\sss maps\sss
$K_{\dff n}\dff \times\dff \Delta^n\qff \ttoo\qff X$\sss are continuous.\oss
The\sss latter condition\dss is\dss equivalent\sss to continuity of\dss
$\num{K}\qff \ttoo\qff X$\nnsp.\oss

\myuppar{Geometric realizations of\dss simplicial\sss spaces with\sss free degeneracies.}
Let\sss $K$\sss be a simplicial\sss space with free degeneracies
and\sss let\sss $N_{\dff n}$\sss be\sss as above.\oss
The diagram\vspace{1.5pt}
\begin{equation}
\label{push-out-free-d}
\quad
\begin{tikzcd}[column sep=boom, row sep=boomm]
N_{\dff n}\dff \times\dff \partial\dff \Delta^n
\arrow[r]
\arrow[d]
&
\protect{\num{\nsp\ssk_{\dff n\dff -\dff 1}\dff K}}
\arrow[d, "\dis i\dff"']
\\
N_{\dff n}\dff \times\dff \Delta^n
\arrow[r, "\dis \pi"]
&
\protect{\num{\nsp\ssk_{\dff n}\dff K}}\dff,
\end{tikzcd}
\end{equation}

\vspace{-10.5pt}\vspace{1.5pt}
where\sss the horizontal\sss arrows are induced\sss by\sss the quotient\sss
map from\sss the disjoint\sss union\qss (\ref{disjoint-union})\qss
to\sss $\num{K}$\sss and\sss the vertical\sss arrows are inclusions,\oss
is\dss obviously\sss commutative.\oss 

Clearly,\oss this diagram\dss is\dss
a\sss push-out\sss square of\dss sets.\oss
In\sss fact,\oss it\dss is\dss also a\sss 
push-out\sss square of\dss topological\sss spaces.\oss
In order\sss to see\sss this,\oss suppose\sss that\sss
$f\dff \colon\dff 
\num{\nsp\ssk_{\dff n}\dff K} 
\qff \ttoo\qff
X$\sss
is\dss a map such\sss that\sss $f\dff \circ\dff i$\dss
and\sss $f\dff \circ\dff \pi$\sss are are continuous.\oss
The continuity of\dss $f\dff \circ\dff i$\sss implies\sss
continuity of\dss the induced\sss maps\sss
$N_{\dff k}\dff \times\dff \Delta^m\qff \ttoo\qff X$\nnsp,\oss
where subspaces\sss $N_{\dff k}\qff \subset\qff K_{\dff m}$\sss
correspond\sss to surjective maps\sss
$\theta\dff \colon\dff
[\dff m\trf]\qff \ttoo\qff [\dff k\trf]$\sss
with\sss $k\qff \leq\qff n\qff -\qff 1$\nnsp.\oss
The continuity of\dss $f\dff \circ\dff \pi$\sss implies\sss
continuity of\dss such\sss induced\sss maps 
with\sss $k\off =\off n$\nnsp.\oss
In view of\pss (\ref{free-d}),\oss the definition of\sss
$\num{\nsp\ssk_{\dff n}\dff K}$ as\sss the quotient\sss space
implies\sss that\sss $f$\sss is\dss continuous.\oss
Hence\sss the above diagram\dss is\dss a\sss 
push-out\sss square of\dss topological\sss spaces.\oss
This conclusion does not\sss require\sss that\sss the spaces $K_{\dff n}$\sss
are compactly\sss generated.\oss

\myuppar{Closed cofibrations.}
Let\sss $A$\sss be a subspace of\sss $X$\nnsp.\oss 
The\sss inclusion\sss $A\qff \ttoo\qff X$\sss
is\dss said\sss to be a\qss \emph{closed\sss cofibration}\oss
if\dss $A$\sss is\dss a closed\sss
and\sss has\sss the homotopy extension\sss property\sss for all\sss spaces.\oss

\mypar{Proposition.}{level-heq}
\emph{Let\qss $f\dff \colon\dff K\qff \ttoo\qff L$\sss
be\dss a map of\dss simplicial\sss spaces with\sss free degeneracies.\oss
Suppose\sss that\trs $K_{\dff n}\dff,\qff L_{\dff n}$
are\sss locally\sss path connected\sss for every\sss $n$\nnsp.\oss
If\trs the maps\sss 
$f_{\dff n}\dff \colon\dff K_{\dff n}\qff \ttoo\qff L_{\dff n}$\sss
are\qss (weak)\qss homotopy\sss equivalences for all\dss $n$\nnsp,\oss
then\sss the map\dss 
$\num{f}\dff \colon\dff
\num{K}\qff \ttoo\qff \num{L}$\sss
is\dss a\qss (weak)\qss homotopy\sss equivalence.\oss}

\proof
The case of\dss weak\sss homotopy equivalences\dss is\dss due\sss to\dss 
D.\dss Dugger\sss and\dss D.\dss Isaksen\qss \cite{di},\oss 
Corollary\qss A.6.\oss
The proofs for both cases are similar.\oss
Clearly,\pss  $f_{\dff n}$\sss takes degenerate simplices\sss
to degenerate simplices even\sss without\sss the\sss free degeneracies assumption.\oss
Since\sss $K_{\dff n}\dff,\qff L_{\dff n}$
are\sss locally\sss path connected,\oss their path components are open and\sss closed.\oss
Since\sss  
$f_{\dff n}\dff \colon\dff 
K_{\dff n}\qff \ttoo\qff L_{\dff n}$\sss 
is\dss a\sss weak\sss homotopy equivalence,\pss
$f_{\dff n}$ 
induces a bijection of\dss the sets of\dss path components.\oss
Since $K$ and\sss $L$ have free degeneracies,\oss
it\sss follows\sss that\sss $f_{\dff n}$\sss
maps\sss $N_{\dff n}\dff(\trf K\trf)$\sss to\sss
$N_{\dff n}\dff(\trf L\trf)$\nnsp.\oss

It\dss is\dss well\sss known\sss that\sss the inclusion\sss
$\partial\dff \Delta^n
\qff \ttoo\qff 
\Delta^n$\sss 
is\dss a closed cofibration.\oss
It\dss follows\sss that\sss
$N_{\dff n}\dff(\trf K\trf)\dff \times\dff \partial\dff \Delta^n
\qff \ttoo\qff 
N_{\dff n}\dff(\trf K\trf)\dff \times\dff \Delta^n$\dss 
is\dss a closed cofibration\qss
(see\qss \cite{br},\oss Corollary\qss 2\qss
of\trs Theorem\qss 7.2.4).\oss
Since\qss (\ref{push-out-free-d})\qss is\dss a push-out\sss square,\oss
it\sss follows\sss that\sss the inclusion\sss
$\num{\nsp\ssk_{\dff n\dff -\dff 1}\dff K}
\qff \ttoo\qff 
\num{\nsp\ssk_{\dff n}\dff K}$\sss
is\dss a closed cofibration.\oss
Of\dss course,\pss $L$\sss has\sss the same properties.\oss

The simplicial\sss map\sss $f$\sss induces a map\sss from\sss the
push-out\sss square\qss (\ref{push-out-free-d})\qss 
to\sss the similar push-out\sss square with\sss $L$
in\sss the role of\dss $K$\nnsp.\oss
Suppose\sss that\sss the maps $f_{\dff n}$ are homotopy equivalences.\oss
In view of\dss the previous paragraph,\oss the gluing\sss theorem\sss
for adjunction spaces applies\qss
(see\qss \cite{br},\oss Theorem\qss 7.5.7).\oss
By\sss this\sss theorem,\oss if\qss 
$\num{\nsp\ssk_{\dff n\dff -\dff 1}\dff K}
\qff \ttoo\qff 
\num{\nsp\ssk_{\dff n\dff -\dff 1}\dff L}$\sss
is\dss a homotopy equivalence,\oss then\sss
$\num{\nsp\ssk_{\dff n}\dff K}
\qff \ttoo\qff 
\num{\nsp\ssk_{\dff n}\dff L}$\sss
is\dss also a homotopy equivalence.\oss
An\sss induction shows\sss that\sss
$\num{\nsp\ssk_{\dff n}\dff K}
\qff \ttoo\qff 
\num{\nsp\ssk_{\dff n}\dff L}$\sss
is\dss a homotopy equivalence for every $n$\nnsp.\oss 
Since\sss the inclusions\sss
$\num{\nsp\ssk_{\dff n\dff -\dff 1}\dff K}
\qff \ttoo\qff 
\num{\nsp\ssk_{\dff n}\dff K}$\sss
are closed cofibration and\sss
$\num{K}$\sss is\dss the direct\sss limit\sss of\dss spaces\sss
$\num{\nsp\ssk_{\dff n}\dff K}$\dss
(and\sss the same\dss is\dss true for\sss $L$\nsp),\oss
it\sss follows\sss that\sss
$\num{K}\qff \ttoo\qff \num{L}$\sss
is\dss a\sss homotopy\sss equivalence.\oss
See\qss \cite{td1},\oss Lemma\qss 6.\oss
The case of\dss weak\sss homotopy equivalences\dss is\dss similar.\oss  \eproof

\myuppar{Good and\sss proper simplicial\sss spaces.}
The material\sss of\dss this subsection\sss will\sss be used\sss
only\sss a couple of\dss times.\oss
We will\sss discuss\sss
what\sss happens without\sss the free degeneracies assumption.\oss
Let\sss
$s\trf(\dff i\trf)\dff \colon\dff
[\halfff n\dff]\qff \ttoo\qff [\halfff n\qff -\qff 1\dff]$
be\sss the unique non-decreasing surjective map\sss assuming\sss the value
$i$\sss twice,\oss
and\sss let\sss $s_{\dff i}\off =\off s\trf(\dff i\trf)^{\dff *}$\dnsp,\oss
where\sss $i\qff \in\qff [\halfff n\qff -\qff 1\dff]$\nnsp.\oss
For a simplicial\sss space $K$\sss let\vspace{3pt}
\begin{equation}
\label{push-out-proper}
\quad
K_{\dff n\fff,\dff i}
\off =\off
s_{\dff i}\dff(\trf K_{\dff n\dff -\dff 1}\trf)
\quad
\mbox{and}\dff\quad
s\fff K_{\dff n}
\off =\off 
\bigcup\nolimits_{\qff i\qff \in\qff [\halfff n\dff -\dff 1\dff]}\qff K_{\dff n\fff,\dff i}
\pff.
\end{equation} 

\vspace{-12pt}\vspace{3pt}
Without\sss the free degeneracies assumption\sss the square\qss (\ref{push-out-free-d})\qss
needs\sss to be replaced\sss by\sss the square\vspace{3pt}
\[
\quad
\begin{tikzcd}[column sep=boom, row sep=boomm]
\bigl(\trf 
s\fff K_{\dff n}\dff \times\dff \Delta^n 
\trf\bigr)
\qff \cup\qff
\bigl(\trf 
K_{\dff n}\dff \times\dff \partial\dff \Delta^n
\trf\bigr)
\arrow[r]
\arrow[d]
&
\protect{\num{\nsp\ssk_{\dff n\dff -\dff 1}\dff K}}
\arrow[d, "\dis i\dff"']
\\
K_{\dff n}\dff \times\dff \Delta^n
\arrow[r, "\dis \pi"]
&
\protect{\num{\nsp\ssk_{\dff n}\dff K}}\dff,
\end{tikzcd}
\]

\vspace{-10.5pt}\vspace{3pt}
where,\oss as in\qss (\ref{push-out-free-d}),\oss
the horizontal\sss arrows are induced\sss by\sss the quotient\sss
map from\qss (\ref{disjoint-union})\qss
to\sss $\num{K}$\sss and\sss the vertical\sss arrows are inclusions.\oss
This\dss is\dss also a push-out\sss square of\dss topological\sss spaces,\oss
but\sss in order\sss to proceed as in\sss the proof\dss of\trs
Proposition\qss \ref{level-heq}\qss one needs\sss to know\sss that\sss
the\sss left\sss vertical\sss arrow\dss is\dss a closed cofibration.\oss
One can easily see\sss that\sss this arrow\dss is\dss a closed cofibration\sss
if\dss the inclusion\sss
$s\fff K_{\dff n}\qff \ttoo\qff K_{\dff n}$\sss
is\dss a closed cofibration.\oss
The simplicial\sss space $K$\sss is\dss said\sss to be\qss \emph{proper}\pss
if\dss this property\sss holds for all\sss $n$\nnsp.\oss
In\sss the same spirit,\pss $K$\sss is\dss said\sss to be\qss \emph{good}\oss
if\dss the inclusions\sss
$K_{\dff n\fff,\dff i}\qff \ttoo\qff K_{\dff n}$\sss are closed cofibrations.\oss
It\dss is\dss easy\sss to see\sss that\sss simplicial\sss spaces with\sss
free degeneracies are both\sss proper and\sss good.\oss\vspace{-0.62pt}

If\dss for proper simplicial\sss spaces one replaces\sss 
the square\qss (\ref{push-out-free-d})\qss by\sss
the square\qss (\ref{push-out-proper}),\oss
then\sss the arguments\sss in\sss the proof\dss of\trs
Proposition\qss \ref{level-heq}\qss still\sss work\sss and\sss
lead\sss to\sss the following\sss result.\oss
Suppose\sss that\dss $f\dff \colon\dff K\qff \ttoo\qff L$\sss
is\dss a map of\dss proper simplicial\sss spaces
such\sss that\sss the maps\sss 
$f_{\dff n}\dff \colon\dff K_{\dff n}\qff \ttoo\qff L_{\dff n}$\sss
are\sss homotopy\sss equivalences for all\dss $n$\nnsp.\oss
Then\sss 
$\num{f}\dff \colon\dff
\num{K}\qff \ttoo\qff \num{L}$\sss
is\dss a\sss homotopy\sss equivalence.\vspace{-0.62pt}

This\sss theorem\dss is\dss due\sss to\dss Segal\qss \cite{s2}\qss and\qss \cite{s3}.\oss
Segal's\dss proof\dss is\dss arranged differently.\oss
He proves\sss that\sss if\dss $K$\sss is\dss proper,\oss
then\sss the quotient\sss map\sss $\norm{K}\qff \ttoo\qff \num{K}$\sss
is\dss a homotopy equivalence.\oss
Then one can use\sss the nearly obvious fact\sss that\sss 
$\norm{f}\dff \colon\dff
\norm{K}\qff \ttoo\qff \norm{L}$\sss
is\dss a\sss homotopy\sss equivalence\sss
if\dss the maps\sss
$f_{\dff n}\dff \colon\dff K_{\dff n}\qff \ttoo\qff L_{\dff n}$\sss
are\sss homotopy\sss equivalences.\oss
See\dss Segal\qss \cite{s3},\oss Proposition\qss A.1,\oss parts\qss (iv)\qss and\qss (ii).\oss
Segal's\dss proof\qss \cite{s3}\qss 
uses\sss properness and\sss refers\sss to\dss Lillig\qss \cite{l}\qss
for\sss the proof\dss that\sss a good\sss simplicial\sss space\dss is\dss proper.\oss
It\sss seems\sss that\sss applying\qss \cite{l}\qss requires either\sss a slightly stronger
assumptions\sss than\sss being\dss good,\pss
which still\sss holds in all\sss applications,\oss
or some additional\sss argument.\oss
Another proof\dss for proper simplicial\sss spaces\dss is\dss
due\trs to\dss tom\sss Dieck\qss \cite{td2}.\oss\vspace{-0.62pt}

\myuppar{Simplicial\sss spaces and\sss $\Delta$\dnsp-spaces.}
An important\sss source of\dss simplicial\sss spaces with\sss free degeneracies
are $\Delta$\dnsp-spaces.\oss
A $\Delta$\dnsp-space\sss $D$
gives rise\sss to simplicial\sss space
$\bm{\Delta}\dff D$ as\sss follows.\oss
The space {\dnsp}$(\trf \bm{\Delta}\dff D\trf)_{\dff n}${\nsp}
of\dss $n$\dnsp-simplices of\trs $\bm{\Delta}\dff D$\dss is\dss the space of\dss pairs\sss
$(\dff \sigma\fff,\pff \rho\trf)$\dss such\dss that\dss
$\sigma\qff \in\qff D_{\dff k}$\dss 
for some\dss $k\qff \leq\qff n$\dss and\dss
$\rho\dff \colon\dff
[\halfff n\dff]\qff \ttoo\qff [\dff k\trf]$\dss
is\dss a surjective non-decreasing\sss map.\oss
In order\dss to define\sss the structure map\sss $\theta^{\dff *}$
for a non-decreasing\sss map\sss
$\theta\dff \colon\dff [\halfff m\dff]\qff \ttoo\qff [\halfff n\dff]$
we represent\dss $\rho\dff \circ\dff \theta$\sss in\sss the form\dss
$\rho\trf \circ\trf \theta\off =\off \tau\trf \circ\trf \varphi$\nnsp,\oss
where\dss $\tau$\dss is\dss strictly\dss increasing and\sss
$\varphi$\dss is\dss surjective and\sss non-decreasing,\pss
and\dss set\dss 
$\theta^{\fff *}\dff(\dff \sigma\fff,\pff \rho\trf)
\off =\off
(\trf \tau^{\dff *}\dff(\dff \sigma\trf)\fff,\pff \varphi\trf)$\nnsp.\oss
One can easily\sss check\dss that\dss 
$(\dff \theta\dff \circ\trf \eta\trf)^{\dff *}
\off =\off
\theta^{\dff *}\fff \circ\qff\fff \eta^{\dff *}$\dss
and\dss hence\sss $\bm{\Delta}\fff D$\dss is\dss 
a simplicial\sss space.\oss
The construction\dss
$D\off \longmapsto\off \bm{\Delta}\dff D$\dss
naturally\sss extends\sss to simplicial\sss maps
and defines a functor from\sss  
$\Delta$\dnsp-spaces\sss 
to\sss simplicial\sss spaces.\oss
By\sss the definition,\oss the canonical\sss map\vspace{1.5pt}
\[
\quad
\coprod\nolimits_{\qff \theta}\qff D_{\dff k}
\off \ttoo\off
(\trf \bm{\Delta}\dff D\trf)_{\dff n}
\]

\vspace{-12pt}\vspace{1.5pt}
is\dss a homeomorphism.\oss
Hence\sss 
$D_{\dff n}$\sss is\dss the space of\dss non-\degenerate $n$\dnsp-simplices of\sss
$\bm{\Delta}\dff D$\sss and\dss $\bm{\Delta}\dff D$\sss has free degeneracies,\oss
i.e.\qss
$N_{\dff n}\dff(\trf \bm{\Delta}\dff D\trf)
\off =\off
D_{\dff n}$\nsp.\oss

The\sss \emph{$n${\dnsp}th skeleton}\qss 
$\ssk_{\dff n}\dff D$\sss
of\dss a $\Delta$\dnsp-space\sss $D$\sss 
is\dss the simplicial\sss $\Delta$\dnsp-subspace of\dss $D$\sss
having as\sss its simplices\sss the $k$\dnsp-simplices of\dss $D$\sss
with\sss $k\qff \leq\qff n$\nnsp.\oss
Clearly,\pss
$\ssk_{\dff n}\dff \bm{\Delta}\dff D
\off =\off
\bm{\Delta}\dff \ssk_{\dff n}\dff D$\nnsp.\oss
Similarly\sss to\sss the case of\dss simplicial\sss spaces,\oss the geometric realization\sss
$\num{D}_{\dff \Delta}$\sss is\dss equal\sss to\sss the direct\sss limit\sss
of\dss subspaces\sss $\num{\nsp\ssk_{\dff n}\dff D}_{\dff \Delta}$\nsp.\oss
Also,\oss similarly\sss to\sss the case of\dss simplicial\sss spaces,\oss
one can describe\sss the subspaces\sss
$\num{\nsp\ssk_{\dff n}\dff D}_{\dff \Delta}$\sss
in\sss terms of\dss push-out\sss diagrams.\oss
Namely,\oss the diagram\vspace{1.5pt}
\begin{equation}
\label{push-out-delta}
\quad
\begin{tikzcd}[column sep=boom, row sep=boomm]
D_{\dff n}\dff \times\dff \partial\dff \Delta^n
\arrow[r]
\arrow[d]
&
\protect{\num{\nsp\ssk_{\dff n\dff -\dff 1}\dff D}_{\dff \Delta}}
\arrow[d, "\dis i\dff"']
\\
D_{\dff n}\dff \times\dff \Delta^n
\arrow[r, "\dis \pi"]
&
\protect{\num{\nsp\ssk_{\dff n}\dff D}_{\dff \Delta}}\dff
\end{tikzcd}
\end{equation}

\vspace{-10.5pt}\vspace{1.5pt}
is\dss a\sss push-out\sss square of\dss topological\sss spaces.\oss
The proof\dss is\dss completely similar\sss to\sss the proof\dss for\sss
the diagram\qss (\ref{push-out-free-d}).\oss
By comparing\sss for each $n$ this diagram\sss  
with\sss the diagram\qss (\ref{push-out-free-d})\qss
for\dss $K\off =\off \bm{\Delta}\dff D$\sss
one can see\sss that\sss the canonical\sss map\sss
$\num{D}_{\Delta}
\qff \ttoo\qff
\num{\bm{\Delta}\dff D}$\sss 
is\dss a\sss homeomorphism.\oss

\myuppar{Simplicial\sss spaces with non-degenerate core.}
The following notions were introduced\sss
in\sss the case of\dss simplicial\sss sets by\dss
Rourke\sss and\trs Sanderson\qss \cite{rosa}.\oss
Suppose\sss that\sss $K$\sss is\dss a simplicial\sss space and\sss let\sss
$N_{\dff n}\off =\off N_{\dff n}\dff(\trf K\trf)$\nnsp.\oss
The\qss \emph{core}\qss of\sss $K$\sss is\sss the $\Delta$\dnsp-subspace $\core K$
of\sss $K$\sss having as its space\sss $(\dff \core K\trf)_{\dff n}$ of\sss
$n$\dnsp-simplices\sss the union\vspace{1.5pt}
\[
\quad
(\dff \core K\trf)_{\dff n}
\off =\off
\bigcup\nolimits_{\qff \theta}\pff \theta^{\dff *}\dff(\trf N_{\dff k}\trf)
\pff,
\]

\vspace{-12pt}\vspace{1.5pt}
where $\theta$ runs over strictly\sss increasing morphisms\sss
$[\halfff n\dff]\qff \ttoo\qff [\dff k\trf]$\nnsp.\oss
In other\sss terms,\oss the simplices of\sss $\core K$
are\sss the non-degenerate simplices of\sss $K$ and\sss their\qss \emph{faces}.\oss
There\dss is\dss
a canonical\sss simplicial\sss map\sss
$\bm{\Delta}\dff \core K\qff \ttoo\qff K$\sss
taking\sss the simplex\sss $(\dff \sigma\fff,\pff \rho\trf)$\sss
of\sss $\bm{\Delta}\dff \core K$\sss
to\sss the simplex\sss $\rho^{\dff *}\dff(\dff \sigma\dff)$ of\sss $K$\nnsp.\oss
The existence part\sss of\qss Eilenberg--Zilber\dss lemma\sss implies\sss that\sss
this map\dss is\dss always surjective\qss
({\fff}i.e.\qss the corresponding maps of\dss the spaces of\sss
$n$\dnsp-simplices are surjective).\oss

The simplicial\sss space $K$\sss is\dss said\sss to be a
simplicial\sss space\qss \emph{with non-degenerate core}\qss
if\dss $\core K$ contains no degenerate simplices of\sss $K$\nnsp,\oss
or,\oss equivalently,\oss if\sss
$\theta^{\dff *}\dff(\trf N_{\dff k}\trf)
\qff \subset\qff
N_{\dff n}$\sss
for every\sss strictly increasing\sss
$\theta\dff \colon\dff
[\halfff n\dff]\qff \ttoo\qff [\dff k\trf]$\nnsp.\oss
If\sss $K$\sss has non-degenerate core,\oss then\sss
$(\dff \core K\trf)_{\dff n}
\off =\off
N_{\dff n}$\nsp,\oss
i.e.\qss the $\Delta$\dnsp-space $\core K$\sss 
is\dss the $\Delta$\dnsp-subspace of\sss $K$\sss
consisting of\dss non-de\-ge\-ner\-ate simplices.\oss
In\sss this case\sss the uniqueness part\sss of\qss 
Eilenberg--Zilber\dss lemma\sss implies\sss that\sss the map\sss
$\bm{\Delta}\dff \core K\qff \ttoo\qff K$\sss
is\dss injective and\sss hence\dss is\dss an\sss
isomorphism of\dss simplicial\sss sets.\oss

\mypar{Lemma.}{ndc-spaces}
\emph{Suppose\sss that\sss a simplicial\sss space\sss $K$\sss
has free degeneracies and\sss non-degenerate core.\oss
Then\sss
$\bm{\Delta}\dff \core K\qff \ttoo\qff K$\sss
is\dss an\sss isomorphism of\dss simplicial\sss spaces
and\sss the canonical\sss map\sss
$\num{\core K}_{\dff \Delta}
\qff \ttoo\qff
\num{\bm{\Delta}\dff \core K}
\off =\off
\num{K}$\sss
is\dss a homeomorphism.\oss}

\proof
Since\sss $K$\sss has non-degenerate core,\pss
$(\dff \bm{\Delta}\dff \core K\trf)_{\dff n}
\qff \ttoo\qff
K_{\dff n}$\sss
is\dss nothing else but\sss the map\qss (\ref{free-d}).\oss
This implies\sss the first\sss statement\sss of\dss the\sss lemma.\oss
In\sss turn,\oss this implies\sss that\sss
$\bm{\Delta}\dff \core K$\sss has free degeneracies
and\sss hence\sss
$\num{\core K}_{\dff \Delta}
\qff \ttoo\qff
\num{\bm{\Delta}\dff \core K}$\sss
is\dss a homeomorphism.\oss \eproof

\mysection{Topological\qss categories}{categories}

\myuppar{Topological\sss categories and\dss their nerves.}
All\sss simplicial\sss spaces in\sss this paper are arising\sss from\sss
topological\sss categories.\oss
\emph{Topological\sss categories}\pss were defined\sss by\trs Segal\qss \cite{s1}\qss
as a small\sss categories\sss such\sss that\sss their sets of\dss objects and
morphisms have\sss topologies and\dss the structure maps\qss
(the source and\dss target\sss of\dss morphisms,\oss
the identity\sss morphisms of\dss objects,\oss and\sss the composition)\qss
are continuous.\oss
Segal\dss attributed some of\trs the ideas of\pss \cite{s1}\qss 
to\dss Grothendieck.\oss

As\dss is\dss well\sss known,\oss each set\sss $[\halfff n\dff]$\sss
can\sss be consider as a category\sss having\sss $[\halfff n\dff]$\sss
as\sss the set\sss of\dss objects,\oss a single morphism\sss
$i\qff \ttoo\qff j$\dss if\trs $i\qff \leq\qff j$\nnsp,\oss
and\sss no morphisms\sss $i\qff \ttoo\qff j$\sss otherwise.\oss
From\sss this point\sss of\dss view\sss the non-decreasing maps\sss
$[\halfff m\dff]\qff \ttoo\qff [\halfff n\dff]$\sss 
are nothing else but\sss functors.\oss

A\sss topological\sss category $\mathcal{C}$ defines
a simplicial\sss space $\mathit{N}\trf \mathcal{C}$\dnsp,\oss
its\qss \emph{nerve}\pss in\sss the sense of\trs Segal\qss \cite{s1}.\oss
The $0$\dnsp-simplices of\dss $\mathit{N}\trf \mathcal{C}$\sss are\sss the objects
of\sss $\mathcal{C}$\nnsp,\oss
and,\oss in\sss general,\oss the $n$\dnsp-simplices are functors\dss
$\sigma\dff \colon\dff
[\halfff n\dff]\qff \ttoo\qff \mathcal{C}$\nnsp.\oss
The structure maps are defined\sss by\sss the composition of\dss functors,\oss
i.e.\qss if\trs
$\theta\dff \colon\dff
[\halfff m\dff]\qff \ttoo\qff [\halfff n\dff]$\sss
is\dss a non-decreasing\sss map,\oss then\dss
$\theta^{\fff *}\dff(\dff \sigma\trf)
\off =\off 
\sigma\dff \circ\trf \theta$\nnsp,\oss
where\sss $\theta$\sss in\sss the right\dss hand side
is\dss considered as a functor.\oss
Clearly,\oss a functor\dss
$[\halfff n\dff]\qff \ttoo\qff \mathcal{C}$\dss
is\dss determined\dss by\dss its values on\sss objects
and on morphisms\dss
$i\qff \ttoo\qff i\qff +\qff 1$\nnsp,\oss
where\dss $i\qff \in\qff [\halfff n\qff -\qff 1\dff]$\nnsp.\oss
Therefore $n$\dnsp-simplices of\dss $\mathit{N}\trf \mathcal{C}$
correspond\dss to sequences of\trs morphisms of\trs the form\vspace{1.5pt}
\begin{equation}
\label{simplex-category}
\quad
\begin{tikzcd}[column sep=large, row sep=boom]\dis
v_{\trf 0}
\arrow[r, "\dis p_{\trf 1}"]
&
v_{\dff 1}
\arrow[r, "\dis p_{\trf 2}"]
&
\ldots
\arrow[r, "\dis p_{\dff n}"]
&
v_{\dff n}\pff,
\end{tikzcd}
\end{equation}

\vspace{-12pt}\vspace{3pt}
where each $v_{\dff i}$\sss is\dss an object\sss of\dss $\mathcal{C}$\sss
and each\sss $p_{\dff i}$\sss is\dss a morphism\dss
$v_{\dff i\dff -\dff 1}\qff \ttoo\qff v_{\dff i}$\nnsp.\oss
Of\dss course,\oss the objects\sss $v_{\dff i}$\sss are determined\dss by\sss
the morphisms\sss $p_{\dff k}$\sss and\dss hence $n$\dnsp-simplices 
correspond\dss to sequences\sss
$(\trf p_{\dff 1}\dff,\off p_{\trf 2}\dff,\off \ldots\dff,\off p_{\dff n}\trf)$\sss
of\dss morphisms such\dss that\dss the composition\sss
$p_{\dff i\dff +\dff 1}\dff \circ\dff p_{\dff i}$\sss
is\dss defined\sss for each\sss $i$\sss between\sss $1$\sss and\sss $n\qff -\qff 1$\nnsp.\oss
So,\oss the set\sss of\sss $n$\dnsp-simplices\dss is\dss contained\sss
in\sss the $n$\dnsp-folds product\sss of\dss copies of\dss the space of\dss
morphisms and\sss hence has a natural\sss topology.\oss
A\sss trivial\sss verification shows\sss that\sss $\mathit{N}\trf \mathcal{C}$\sss
with\sss this\sss topology\dss is\dss a simplicial\sss space.\oss
We will\sss call\sss the $n$\dnsp-simplices
of\sss the nerve\sss $\mathit{N}\trf \mathcal{C}$\sss
simply\sss the\sss \emph{$n$\dnsp-simplices of\dss the category}\dss
$\mathcal{C}$\dnsp,\oss
and\sss will\sss denote\sss the space of\sss $n$\dnsp-simplices of\dss
$\mathcal{C}$\sss by\dss $\mathcal{C}_{\dff n}$\sss
and\sss the geometric realization\sss
$\num{\hnsp\mathit{N}\trf \mathcal{C}}$\sss by\dss
$\num{\mathcal{C}}$\nnsp.\oss

Basic examples are provided\sss by\sss the ordered sets $[\halfff n\dff]$\sss
equipped with\sss the discrete\sss topology and
considered as categories.\oss
The nerve\sss $\mathit{N} \trf[\halfff n\dff]$\sss
is\dss called\sss the\qss \emph{abstract\sss $n$\dnsp-simplex}.\oss
Its geometric realization\sss $\num{[\halfff n\dff]}$\sss
is\dss canonically\sss homeomorphic\sss to\sss the standard\sss
geometric $n$\dnsp-simplex $\Delta^n$\dnsp.\oss

\myuppar{Functors and\sss natural\dss transformations.}
A\qss \emph{continuous functor}\pss between\sss topological\sss categories\dss
is\dss a functor such\sss that\sss the corresponding maps of\dss the
spaces of\dss objects and\sss morphisms are continuous.\oss
Let\sss $\mathcal{C}\fff,\pff \mathcal{D}$\sss
be\sss topological\sss categories and\sss
$f\dff \colon\dff
\mathcal{C}\qff \ttoo\qff \mathcal{D}$\dss
be a continuous functor.\oss
Then\sss $f$\sss
defines a simplicial\dss map\dss
$\mathit{N}\fff f\dff \colon\dff
\mathit{N}\trf \mathcal{C}
\qff \ttoo\qff 
\mathit{N}\trf \mathcal{D}$\dnsp.\oss
In\sss fact,\oss the nerve\sss $\mathit{N}\trf \bullet$\sss
is\dss a functor\sss from\sss topological\sss categories\sss
to simplicial\sss spaces.\oss

Recall\sss that\sss a\sss homotopy\sss between\sss two simplicial\sss maps\sss
$a\fff,\qff b\dff \colon\dff K\qff \ttoo\qff L$\sss of\dss
simplicial\sss spaces\dss is\dss defined as a simplicial\sss map\sss
$K\dff \times\dff \mathit{N} \trf[\dff 1\dff]
\qff \ttoo\qff
L$\sss
such\sss that\sss its\sss restrictions\sss to\sss
$K\dff \times\dff 0$\sss and\sss $K\dff \times\dff 1$\sss
are\sss the maps $a$ and $b$\nnsp.\oss
For\sss nerves of\dss categories and simplicial\sss maps arising\sss
from\sss functors\sss
natural\sss transformations lead\sss to homotopies.\oss
Namely,\oss
given\dss two functors\dss
$f,\pff g\dff \colon\dff
\mathcal{C}\qff \ttoo\qff \mathcal{D}$\dnsp,\oss
a natural\dss transformation\dss
$f\qff \ttoo\qff g$\sss
defines a homotopy\dss between\sss
$\mathit{N}\fff f$\sss and\sss $\mathit{N}\dff g$\nnsp.\oss
Indeed,\oss a natural\dss transformation\dss
$t\dff \colon\dff f\qff \ttoo\qff g$\sss
can\sss be considered as\sss a\sss functor\dss
$\mathcal{C}\dff \times\dff [\dff 1\dff]
\qff \ttoo\qff
\mathcal{D}$\nnsp,\oss
where\sss $[\dff 1\dff]$\sss is\dss considered as a category.\oss
One can easily\sss see\sss that\dss the construction\dss
$\mathcal{C}\off \longmapsto\off \mathit{N}\trf \mathcal{C}$\dss
commutes with\dss the products,\oss
at\sss least\sss when all\dss involved\sss topological\sss spaces
are compactly\sss generated.\oss
Therefore\sss $t$\sss defines a simplicial\dss map\dss
$\mathit{N}\fff t\dff \colon\dff
\mathit{N}\trf \mathcal{C}\dff \times\dff \mathit{N} \trf[\dff 1\dff] 
\qff \ttoo\qff \mathit{N}\trf \mathcal{D}$\dnsp,\oss
i.e.\qss a\sss homotopy,\oss
and\sss this homotopy\dss is\dss a homotopy\dss between\sss
$\mathit{N}\fff f$\sss and\sss $\mathit{N}\dff g$\nnsp.\oss

A functor
$f\dff \colon\dff
\mathcal{C}\qff \ttoo\qff \mathcal{D}$\sss
defines a continuos map\sss
$\num{f}\dff \colon\dff
\num{\mathcal{C}}\qff \ttoo\qff \num{\mathcal{D}}$\nnsp.\oss
Since\sss $\num{[\dff 1\dff]}$\sss is\dss canonically\sss homeomorphic\sss to\sss
$[\dff 0\fff,\qff 1\dff]$\nnsp,\oss
a natural\sss transformation between\sss two such functors
defines a homotopy\sss
$\num{\mathcal{C}}\dff \times\dff [\dff 0\fff,\qff 1\dff]
\qff \ttoo\qff 
\num{\mathcal{D}}$\nnsp.\oss
This construction\sss turns out\sss to be a surprisingly efficient\sss
way\sss to construct\sss homotopies.\oss 
For example,\oss if\sss $\mathcal{C}$ and\sss $\mathcal{D}$ are equivalent\sss
topological\sss categories,\oss then\sss the spaces\sss 
$\num{\mathcal{C}}$ and\sss $\num{\mathcal{D}}$ are homotopy equivalent.\oss

\myuppar{Topological\sss categories with\sss free units.}
Clearly,\oss an $n$\dnsp-simplex\qss (\ref{simplex-category})\qss
of\dss a\sss topological\sss category\sss $\mathcal{C}$\sss
is\dss degenerate\sss if\trs and\dss only\trs if\dss
one of\dss the morphisms\sss
$p_{\dff i}\dff \colon\dff
v_{\dff i\dff -\dff 1}\qff \ttoo\qff v_{\dff i}$\sss
is\dss actually an\sss identity\sss morphism\qss
(and,\oss in\sss particular,\pss 
$v_{\dff i\dff -\dff 1}\off =\off v_{\dff i}$).\oss
Let\sss us\sss say\sss that\sss a\sss topological\sss category\sss $\mathcal{C}$\sss has\qss
\emph{free units}\pss if\dss the space of\dss the identity\sss morphisms of\dss objects of\dss
$\mathcal{C}$\sss is\dss open and closed\sss in\sss the space\sss
$\mor\trf \mathcal{C}$\sss of\dss morphisms of\dss $\mathcal{C}$\nnsp.\oss\vspace{-0.62pt}

Suppose\sss that\sss one can define a\qss \emph{dimension}\pss of\dss objects of\dss $\mathcal{C}$\sss
with\sss the following\sss two properties.\oss
First,\oss the dimension\dss is\dss
a non-negative integer or a pair of\dss non-negative integers such\sss that\sss
a morphism\sss $C\qff \ttoo\qff D$\sss of\dss $\mathcal{C}$\sss
is\dss an\sss identity\sss morphism\dss if\trs and\dss only\trs if\dss
the dimensions of\sss $C$\sss and\dss $D$\sss are equal.\oss
Second,\oss the dimension\dss is\dss a continuous function on\sss the space of\dss objects.\oss
Clearly,\oss these properties imply\sss that\sss $\mathcal{C}$\sss has free units.\oss
In\sss this case we say\sss that\sss $\mathcal{C}$\dss
has free units by\sss the\qss \emph{dimension\sss reasons}.\oss\vspace{-0.62pt}

\mypar{Lemma.}{free-d-categories}
\emph{If\pss $\mathcal{C}$ is\dss a\sss topological\sss category\sss
with\dss free units,\oss
then\sss the nerve\sss $\mathit{N}\trf \mathcal{C}$\sss 
is\dss a simplicial\sss space with\dss free degeneracies.\oss}\vspace{-0.62pt}

\proof
For an {\nsp}$n$\dnsp-simplex $\sigma$ corresponding\sss
to\sss the sequence\qss (\ref{simplex-category})\pss
let\dss 
$\theta_{\dff \sigma}
\dff \colon\dff
[\dff n\trf]\qff \ttoo\qff [\dff k\trf]$\dss
be\sss the unique non-decreasing surjective map\sss
such\sss that\sss 
$\theta\dff(\trf i\qff -\qff 1\trf)
\off =\off
\theta\dff(\trf i\trf)$\dss
if\trs and\dss only\trs if\dss the morphism\sss
$p_{\dff i}\dff \colon\dff
v_{\dff i\dff -\dff 1}\qff \ttoo\qff v_{\dff i}$\sss 
is\dss an identity\sss morphism.\oss
Clearly,\pss $\sigma\off =\off \theta_{\dff \sigma}^{\dff *}\dff(\trf \tau\trf)$\sss
for a unique $k$\dnsp-simplex $\tau$\nnsp,\oss
and\sss $\tau$\sss is\dss non-degenerate.\oss
Therefore\sss $\mathcal{C}_{\dff n}$\sss 
is\dss equal\sss to\sss the union of\dss disjoint\sss subsets corresponding\sss
to non-decreasing surjective maps\sss
$[\dff n\trf]\qff \ttoo\qff [\dff k\trf]$\nnsp.\oss
Since\sss $\mathcal{C}$\sss has free units,\oss
this disjoint\sss union of\dss subsets\dss is\dss actually a 
disjoint\sss union of\dss subspaces and\sss hence\sss the nerve\sss
$\mathit{N}\trf \mathcal{C}$\sss 
has\sss free degeneracies.\oss  \eproof

\mysection{Partially\qss ordered\qss  spaces}{pos-section}

\myuppar{Partially\sss ordered\sss topological\sss spaces.}
Most\sss of\dss our categories will\sss arise from\sss partially\sss ordered sets.\oss
Suppose\sss that\sss $S$\sss is\dss a set\sss
with a partial\sss order\sss $\leq$\nnsp.\oss
We may consider\sss $S$\sss as a category\sss having $S$ as its set\sss of\dss objects,\oss
and\sss having exactly\sss one morphism\sss $a\qff \ttoo\qff b$\sss
if\dss $a\qff \leq\qff b$\nnsp,\oss and none otherwise.\oss
Functors between such categories are simply\sss the non-decreasing maps,\oss
i.e.\qss maps\sss $f$\sss such\sss that\sss $a\qff \leq\qff b$\sss
implies\sss $f\trf(\dff a\trf)\qff \leq\qff f\trf(\dff b\trf)$\nnsp.\oss
We already\sss used a special\sss case of\dss this construction,\oss
namely\sss the case of\dss
the sets\sss $[\halfff n\dff]$\sss 
with\dss their natural\sss order.\oss 

A\qss \emph{partially\sss ordered\dss topological\sss space},\oss
or\sss simply an\qss \emph{ordered\dss topological\sss space}\pss is\dss a\sss
topological\sss space $S$\sss together with 
a partial\sss order\sss $\leq$\sss on $S$\nnsp.\oss
If\dss an ordered\sss topological\sss space $S$\sss is\dss
considered as a category,\oss then\sss its set\sss of\dss morphisms\dss
is\dss the subspace of\sss $S\dff \times\dff S$\sss consisting of\dss
pairs\sss $(\dff a\fff,\qff b\trf)$\sss such\sss that\sss $a\qff \leq\qff b$\sss
and\sss hence\dss is\dss a\sss topological\sss space in a natural\sss way.\oss
Clearly,\oss the structure maps of\dss the category $S$ are continuous
and\sss hence $S$\sss is\dss a\sss topological\sss category.\oss
As\sss it\sss turns out,\oss working\sss with such\sss topological\sss categories\dss
is\dss simpler\sss than with\sss general\sss ones,\oss and\sss
most\sss of\dss our\sss topological\sss categories
will\sss arise from ordered\sss topological\sss spaces.\oss

\myuppar{Partially\sss ordered spaces with\sss free equalities.}
Let\sss $S$\sss be a\sss topological\sss space\sss together with a
partial\sss order\sss $\leq$\nnsp.\oss
The order\sss $\leq$\sss is\dss said\sss to have\qss
\emph{free equalities},\oss or\sss simply\sss to be\qss 
\emph{free},\oss if\trs the\sss sets\vspace{-0.44pt}
\[
\quad
\left\{\pff 
(\trf x\fff,\qff y\trf)\qff \in\qff S\dff \times\dff S
\off \bigl|\off x\off =\off y
\pff\right\}\dff,
\]

\vspace{-36pt}\vspace{-0.44pt}
\[
\quad
\left\{\pff 
(\trf x\fff,\qff y\trf)\qff \in\qff S\dff \times\dff S
\off \bigl|\off x\off <\off y
\pff\right\}\dff,
\quad
\mbox{and}\quad
\]

\vspace{-36pt}\vspace{-0.44pt}
\[
\quad
\left\{\pff 
(\trf x\fff,\qff y\trf)\qff \in\qff S\dff \times\dff S
\off \bigl|\off y\off <\off x
\pff\right\}
\]

\vspace{-12pt}\vspace{-0.44pt}
are closed\sss in\sss $S\qff \times\qff S$\nnsp.\oss
It\dss is\dss worth\sss to mention\sss that\dss if\sss $S$\sss 
is\dss compactly generated and\sss the product\sss
$S\dff \times\dff S$\sss is\dss taken within\sss the category of\dss
compactly generated spaces,\oss then\sss the closedness of\dss the first\sss
set\dss is\dss equivalent\sss to $S$\sss being\qss \emph{weakly\dss Hausdorff}.\oss
But\sss we will\sss not\sss use\sss neither this fact,\oss
nor\sss weakly\dss Hausdorff\trs spaces.\oss
The property of\dss being free\dss is\dss fairly restrictive.\oss
For example,\oss the set\sss $\rrr$ with its usual\sss topology and order\sss $\leq$\sss
is\dss not\sss free.\oss 
But\sss $\rrr$ with\sss the discrete\sss topology and\sss the usual\sss order\sss $\leq$\sss
is\dss free,\oss
as\dss is\dss every partially ordered set\sss 
equipped with discrete\sss topology.\oss

\mypar{Lemma.}{free-orders}
\emph{If\qss a partially\dss ordered\sss space\sss $S$\sss has\qss free equalities,\oss
then\sss $S$\sss has free units as a\sss topological\sss category.\oss
In\sss particular,\oss its\sss nerve\dss
is\dss a simplicial\sss space with free degeneracies.\oss}

\proof
The space\sss
$\mor\trf S$\sss of\dss morphisms of\dss the\sss topological\sss category\sss $S$\sss 
is\dss the space of\dss pairs\sss
$(\trf x\fff,\qff y\trf)\qff \in\qff S\dff \times\dff S$\sss
such\sss that\sss $x\qff \leq\qff y$\nnsp.\oss
It\sss follows\sss that\dss if\sss $S$\sss is\dss has free equalities,\oss 
then $S$\sss has free units as a\sss topological\sss category.\oss
The second claim\sss follows from\trs Lemma\qss \ref{free-d-categories}.\oss  \eproof

\mypar{Lemma.}{pos} 
\emph{Let\sss $\mathcal{C}$\sss be a partially\sss ordered space\sss
considered as a category.\oss
Then\sss the nerve\sss $\mathit{N}\trf \mathcal{C}$\sss
is\dss a simplicial\sss space with non-degenerate core.\oss}

\proof
Let\dss $\leq$\dss be\sss the partial\sss order of\dss $\mathcal{C}$\dnsp.\oss
The $k$\dnsp-simplices of\dss $\mathcal{C}$\sss are\sss non-decreasing\sss
sequences\sss
$v_{\dff 0}\qff \leq\qff
v_{\dff 1}\qff \leq\qff
\ldots\qff \leq\qff
v_{\dff k}$\nsp.\oss
Such a simplex\dss is\dss degenerate\dss
if\trs and\dss only\trs if\dss
$v_{\dff i\dff -\dff 1}\off =\off v_{\dff i}$\sss for some $i$\nnsp.\oss
Therefore\sss the subspace $N_{\dff n}$ of\dss non-degenerate
$n$\dnsp-simplices consists of\dss strictly\sss increasing sequences\sss
$v_{\dff 0}\qff <\qff
v_{\dff 1}\qff <\qff
\ldots\qff <\qff
v_{\dff k}$\nsp.\oss 
Applying\sss to such a simplex\sss the structural\sss map $\theta^{\dff *}$
for a strictly\sss increasing\sss
$\theta\dff \colon\dff
[\halfff n\dff]\qff \ttoo\qff [\dff k\trf]$\sss
amounts\sss to removing from\sss the sequence\sss terms $v_{\dff i}$
with $i$ not\sss in\sss the image of\sss $\theta$\nnsp.\oss
It\sss follows\sss that\sss
$\theta^{\dff *}\dff(\trf N_{\dff k}\trf)
\qff \subset\qff
N_{\dff n}$\sss
for every\sss strictly increasing\sss
$\theta\dff \colon\dff
[\halfff n\dff]\qff \ttoo\qff [\dff k\trf]$\nnsp,\oss
i.e.\qss $\mathit{N}\trf \mathcal{C}$\sss
has non-degenerate core.\oss  \eproof

\mypar{Corollary.}{pos-is-nice}
\emph{If\qss a partially\dss ordered\sss space\sss $S$\sss has\qss free equalities,\oss
then\sss the canonical\dss map\dss
$\num{\nsp\core S}_{\dff \Delta}
\qff \ttoo\qff
\num{S}$\sss
is\dss a homeomorphism.\oss}

\proof
The nerve of\sss $S$\sss has free equalities by\trs
Lemma\qss  \ref{free-orders}\qss
and\dss non-degenerate core by\trs
Lem\-ma\qss \ref{pos}.\oss
It\dss remains\sss to apply\trs Lemma\qss \ref{ndc-spaces}.\oss  \eproof

\mysection{Topological\qss simplicial\qss complexes}{topological-simplicial-complexes}

\myuppar{Simplicial\sss complexes.}
Recall\sss that\sss a\qss \emph{simplicial\sss complex}\pss
is\dss a set\sss $S$\sss together with a collection of\dss
its finite subsets,\oss called\qss \emph{simplices},\oss
such\sss that\sss a subset\sss of\dss a simplex\dss is\dss also a simplex
and\sss $S$\sss is\dss the union of\dss simplices.\oss
If\sss $\sigma$\sss is\dss a simplex and\sss $\tau\qff \subset\qff \sigma$\nnsp,\oss
the simplex $\tau$\sss is\dss said\sss to be a\qss \emph{face}\pss of\sss $\sigma$\nnsp.\oss
A simplex $\sigma$\sss is\dss said\sss to be an\dss \emph{$n$\dnsp-simplex}\pss
if\sss $\sigma$\sss consists of\sss $n\qff +\qff 1$\sss elements.\oss
The\dss \emph{$n${\dnsp}th skeleton}\dss $\ssk_{\dff n}\dff S$\sss is\dss defined\sss
as\sss the simplicial\sss complex\sss having as simplices $k$\dnsp-simplices of\sss $S$\sss
with\sss $k\qff \leq\qff n$\nnsp.\oss
If\trs $S\fff,\qff T$\sss are simplicial\sss complexes,\oss
then a\qss \emph{simplicial\sss map}\qss from $S$\sss to $T$\sss
is\dss defined as a map\sss $S\qff \ttoo\qff T$\sss taking simplices\sss to simplices.\oss
Such maps may decrease\sss the dimension
of\dss simplices,\oss in\sss contrast\sss with simplicial\sss maps of\dss
simplicial\sss and $\Delta$\dnsp-sets and spaces.\oss

A simplicial\sss complex $S$\sss is\dss said\sss to be\qss
\emph{partially\sss ordered}\pss if\sss $S$\sss is\dss equipped with 
a partial\sss order\sss with respect\sss to which\sss every simplex\dss
is\dss linearly ordered.\oss
In\sss this case $S$ defines a $\Delta$\dnsp-set\sss $\Delta\fff S$\sss
and a simplicial\sss set\sss $\bm{\Delta}\dff S$\nnsp.\oss
Their $n$\dnsp-simplices are,\oss 
respectively,\oss increasing or\sss non-decreasing maps\sss
$\sigma\dff \colon\dff [\halfff n\dff]\qff \ttoo\qff S$\sss
such\sss that\sss  
$\image \sigma$\sss is\dss a simplex.\oss 
The structure maps are\sss
$\theta^{\dff *}\dff(\trf \sigma\trf)
\off =\off 
\sigma\dff\circ\dff \theta$\nnsp.\oss

\myuppar{Topological\sss simplicial\sss complexes.}
A\qss \emph{topological\sss simplicial\sss complex}\pss
is\dss a\sss topological\sss space $S$\sss together with a structure of\dss
a simplicial\sss complex on\sss $S$\nnsp.\oss
A\qss \emph{simplicial\sss map}\pss from $S$\sss to $T$\sss
is\dss a continuous map\sss $S\qff \ttoo\qff T$\sss
taking\sss simplices\sss to simplices.\oss
The\qss \emph{geometric realization}\qss $\bbnum{S}$\sss of\sss 
$S$\sss 
is\dss the set\sss of\dss maps\sss
$t\dff \colon\dff S\qff \ttoo\qff [\dff 0\fff,\qff 1\dff]$\sss
such\sss that\sss
$\{\trf v\qff \in\qff S \qff\mid\qff t\trf(\dff v\trf)\off \neq\off 0 \qff\}$\sss
is\dss a simplex and 
$\sum\nolimits_{\qff v\qff \in\qff S}\qff t\trf(\dff v\trf)
\off =\off
1$\nnsp.\oss
Equivalently,\oss 
$\bbnum{S}$\sss is\dss the set of\dss
are\sss 
weighted sums of\dss the form\vspace{1.5pt}
\begin{equation}
\label{weighted-sums}
\quad
t_{\dff 0}\dff v_{\dff 0}\pff +\pff
t_{\dff 1}\dff v_{\dff 1}\pff +\pff
\ldots\pff +\pff
t_{\dff n}\dff v_{\dff n}
\off,
\end{equation}

\vspace{-12pt}\vspace{1.5pt}
such\sss that\sss
$\left\{\qff 
v_{\dff 0}\dff,\off
v_{\dff 1}\dff,\off
\ldots\dff,\off
v_{\dff n}
\pff\right\}
$\sss
is\dss a simplex\sss
and\dss
$(\trf t_{\dff 0}\dff,\qff t_{\dff 1}\dff,\qff \ldots\dff,\qff t_{\dff n}\trf)
\qff \in\qff
\Delta^n$\nnsp.\oss
If\dss
$\sigma$\sss
is\dss the simplex\sss
$\{\trf v_{\dff 0}\dff,\qff
v_{\dff 1}\dff,\qff
\ldots\dff,\qff
v_{\dff n} \qff\}$\nnsp,\oss
we will\sss denote by\sss $\num{\sigma}$\sss the subspace of\dss 
weighted sums of\dss the form\qss (\ref{weighted-sums})\qss
and\sss by\sss $\inte \num{\sigma}$\sss
the subspace of\dss 
weighted sums\qss (\ref{weighted-sums})\qss
such\sss that\sss $t_{\dff i}\qff >\qff 0$\sss for every $i$\nnsp.\oss
Similarly,\pss $\inte \Delta^n$\sss is\dss the subspace of\dss points\sss
$(\trf t_{\dff 0}\dff,\qff t_{\dff 1}\dff,\qff \ldots\dff,\qff t_{\dff n}\trf)
\qff \in\qff
\Delta^n$\sss
such\sss that\sss $t_{\dff i}\qff >\qff 0$\sss for every $i$\nnsp.\oss

Of\dss course,\oss we need also a\sss topology on\sss the sets\sss $\bbnum{S}$\nnsp.\oss
Since we will\sss need only\sss this case,\oss
we will\sss assume\sss that\sss $S$\sss is\qss \emph{partially\sss ordered}\pss
in\sss the sense\sss that\sss $S$\sss is\dss 
equipped\sss with a partial\sss order\sss $\leq$\sss inducing\sss a\sss linear order on each simplex.\oss
Then we can assume\sss that\sss
$v_{\dff 0}\qff <\qff
v_{\dff 1}\qff <\qff
\ldots\qff <\qff
v_{\dff n}$\sss
in\sss sums\qss (\ref{weighted-sums}).\oss
Let\sss 
$\bbnum{S}_{\dff k}$\sss 
be\sss the set\sss of\dss maps $t$\sss such\sss that\sss
$t\trf(\dff v\trf)\off \neq\off 0$\sss
for\sss $\leq\qff k\qff +\qff 1$\sss points $v$\nnsp.\oss
In other\sss terms,\pss
$\bbnum{S}_{\dff k}$\sss is\dss the set\sss of\dss weighted sums\qss (\ref{weighted-sums})\qss
with\sss $n\qff \leq\qff k$\nnsp.\oss
Therefore\sss $\bbnum{S}_{\dff k}$\sss
can\sss be identified with\sss a quotient\sss of\dss 
the disjoint\sss union\vspace{1.5pt} 
\[
\quad
\coprod\nolimits_{\dff n\qff \leq\qff k}\off
S_{\fff n}\dff \times\dff \Delta^n
\off,
\]

\vspace{-12pt}\vspace{1.5pt}
where $S_{\fff n}\qff \subset\qff S^{\dff n\dff +\dff 1}$ is\dss
the set\sss of\sss $\dnsp n$\dnsp-tuples\sss 
$(\dff 
v_{\dff 0}\dff,\qff
v_{\dff 1}\dff,\qff
\ldots\dff,\qff
v_{\dff n}
\dff)$\sss
such\sss that\sss
$\{\dff 
v_{\dff 0}\dff,\qff
v_{\dff 1}\dff,\qff
\ldots\dff,\qff
v_{\dff n}
\dff\}$\sss
is\dss a simplex and
$v_{\dff 0}\qff <\qff
v_{\dff 1}\qff <\qff
\ldots\qff <\qff
v_{\dff n}$\nsp.\oss
We equip\sss $\bbnum{S}_{\dff k}$\sss with\sss the\sss topology of\dss this quotient.\oss
Clearly,\pss 
$\bbnum{S}$\sss 
is\dss equal\sss to\sss the union of\dss subspaces\sss
$\bbnum{S}_{\dff k}$\nsp,\oss
and we equip $\bbnum{S}$ with\sss the direct\sss limit\sss topology.\oss
If\sss $S$\sss is\dss not\sss partially\sss ordered,\pss
$S_{\fff n}$\sss should\sss be replaced\sss
by\sss the set\sss of\sss orderings of $n$\dnsp-simplices 
and\sss $S_{\fff n}\dff \times\dff \Delta^n$\sss by\sss the quotient\sss by\sss 
the diagonal\sss action of\dss the symmetric group $\Sigma_{\dff n\dff +\dff 1}$\nsp.\oss
We\sss omit\sss the details.\oss
If\dss 
$f\dff \colon\dff S\qff \ttoo\qff T$\sss
is\dss a simplicial\sss map,\oss
then\sss the formula\vspace{1.5pt}
\[
\quad
f\dff(\trf
t_{\dff 0}\dff v_{\dff 0}\pff +\pff
t_{\dff 1}\dff v_{\dff 1}\pff +\pff
\ldots\pff +\pff
t_{\dff n}\dff v_{\dff n}
\trf)
\off =\off
t_{\dff 0}\dff f\dff(\fff v_{\dff 0}\trf)\pff +\pff
t_{\dff 1}\dff f\dff(\fff v_{\dff 1}\trf)\pff +\pff
\ldots\pff +\pff
t_{\dff n}\dff f\dff(\fff v_{\dff n}\trf)
\pff.
\]

\vspace{-12pt}\vspace{1.5pt}
defines a map\sss
$\bbnum{f}\dff \colon\dff 
\bbnum{S}
\qff \ttoo\qff 
\bbnum{T}$\nnsp.\oss
Suppose\sss that\dss $f\dff \colon\dff S\qff \ttoo\qff T$\sss
is\dss a simplicial\sss map and\dss is\qss \emph{order-preserving}\pss
in\sss the sense\sss that\sss
$x\qff \leq\qff y$\sss implies\sss 
$f\dff(\dff x\trf)\qff \leq\qff f\dff(\dff y\trf)$\nnsp.\oss
Clearly,\oss in\sss this case\sss the map\sss
$\bbnum{f}\dff \colon\dff 
\bbnum{S}
\qff \ttoo\qff 
\bbnum{T}$\sss
is\dss continuous.\oss
In\sss fact,\oss this\dss is\dss true\qss
(and easy\sss to prove)\qss for arbitrary simplicial\sss maps
and\sss simplicial\sss complexes without\sss any orders.\oss

\myuppar{Skeletons and\sss the geometric realization.}
The\dss \emph{$n${\dnsp}th skeleton}\dss $\ssk_{\dff n}\dff S$\sss of\dss a\sss
topological\sss simplicial\sss complex\sss $S$\sss is\dss simply\sss the
$n${\dnsp}th skeleton of\sss $S$ as a simplicial\sss complex.\oss
The\sss topology on $S$ remains\sss the same.\oss
Suppose\sss that\sss the\sss topological\sss simplicial\sss complex\sss $S$\sss
is\dss partially ordered.\oss
Then for every\sss $n$\sss the space\sss $\bbnum{S}_{\dff n}$\sss is\dss
equal\sss to\sss the geometric realization\sss $\bbnum{\ssk_{\dff n}\dff S}$\sss
by\sss the definition of\dss the\sss latter.\oss
As in\sss the cases of\dss simplicial\sss spaces and $\Delta$\dnsp-set,\oss
the spaces\sss $\bbnum{\ssk_{\dff n}\dff S}$\sss can\sss be described\sss
in\sss terms of\dss push-out\sss diagrams.\oss
Namely,\oss let\sss $S_{\fff n}$\sss be,\oss as above,\oss
the space of\sss $n$\dnsp-tuples\sss 
$(\dff 
v_{\dff 0}\dff,\qff
v_{\dff 1}\dff,\qff
\ldots\dff,\qff
v_{\dff n}
\dff)$\sss
such\sss that\sss
$\{\dff 
v_{\dff 0}\dff,\qff
v_{\dff 1}\dff,\qff
\ldots\dff,\qff
v_{\dff n}
\dff\}$\sss
is\dss a simplex and
$v_{\dff 0}\qff <\qff
v_{\dff 1}\qff <\qff
\ldots\qff <\qff
v_{\dff n}$\nsp.\oss
Since $S$\sss is\dss partially ordered,\pss $S_{\fff n}$\sss
is\sss the space of $n$\dnsp-simplices of\sss $S$\nnsp.\oss
The diagram\vspace{3pt}
\begin{equation}
\label{push-out-tsc}
\quad
\begin{tikzcd}[column sep=boom, row sep=boomm]
S_{\fff n}\dff \times\dff \partial\dff \Delta^n
\arrow[r]
\arrow[d]
&
\protect{\bbnum{\nsp\ssk_{\dff n\dff -\dff 1}\dff S}}
\arrow[d, "\dis i\dff"']
\\
S_{\fff n}\dff \times\dff \Delta^n
\arrow[r, "\dis \pi"]
&
\protect{\bbnum{\nsp\ssk_{\dff n}\dff S}}\dff
\end{tikzcd}
\end{equation}

\vspace{-10.5pt}\vspace{3pt}
is\dss a\sss push-out\sss square of\dss topological\sss spaces.\oss
The proof\dss is\dss completely similar\sss to\sss the proof\dss for\sss
the diagrams\qss (\ref{push-out-delta})\qss and\qss (\ref{push-out-free-d}).\oss
By\sss the definition\sss $\bbnum{S}$\sss is\dss  
the direct\sss limit\sss of\dss the spaces  
$\bbnum{\nsp\ssk_{\dff n}\dff S}$\nnsp.\oss

\myuppar{Partially ordered spaces as\sss topological\sss simplicial\sss complexes.}
Suppose now\sss that\sss $S$\sss is\dss a\sss topological\sss space\sss
together with a partial\sss order\sss $\leq$\nnsp.\oss
Usually we consider $S$ either as a\sss topological\sss category
or\sss the associated simplicial\sss space.\oss
At\sss the same\sss time $S$\sss can\sss be considered\sss 
as a\sss topological\sss simplicial\sss complex\sss
having as simplices finite subsets of\sss $S$\sss linearly ordered\sss by\sss $\leq$\nnsp.\oss
By\sss the definition,\oss every simplex of\sss $S$ as a\sss
topological\sss simplicial\sss complex\dss is\dss linearly ordered.\oss
This allows\sss to define a $\Delta$\dnsp-space\sss $\Delta\fff S$ 
exactly an\sss in\sss the discrete case.\oss
Clearly,\oss the $n$\dnsp-simplices of\sss $\Delta\dff S$\sss 
and of\sss $S$ as a\sss topological\sss simplicial\sss complex\sss
are\sss the increasing sequences\sss
$v_{\dff 0}\qff <\qff
v_{\dff 1}\qff <\qff
\ldots\qff <\qff
v_{\dff n}$\sss
of\dss points\sss $v_{\dff i}\qff \in\qff S$\nnsp.\oss
Hence\sss the $\Delta$\dnsp-set\sss $\Delta\dff S$\sss 
and $S$ as a\sss topological\sss simplicial\sss complex\sss
carry\sss the same information and\dss there\dss is\dss hardly any difference
between\sss them.

One can also define a simplicial\sss space\sss $\bm{\Delta}\dff S$\sss
exactly an\sss in\sss the discrete case.\oss
Clearly,\oss the simplicial\sss space\sss $\bm{\Delta}\dff S$\sss
is\dss nothing else but\sss $S$\sss considered as a simplicial\sss space.\oss
Indeed,\pss $n$\dnsp-simplices of\sss $\bm{\Delta}\dff S$\sss
are\sss non-decreasing\sss sequences\sss
$v_{\dff 0}\qff \leq\qff
v_{\dff 1}\qff \leq\qff
\ldots\qff \leq\qff
v_{\dff k}$\sss
of\dss points\sss $v_{\dff i}\qff \in\qff S$\nnsp,\oss
as also $n$\dnsp-simplices of\sss $S$\sss
as a simplicial\sss space.\oss
By\trs Lemma\qss \ref{pos}\qss the simplicial\sss space $S$ has non-degenerate core.\oss
Moreover,\oss the proof\dss of\dss this\sss lemma\sss shows\sss that\sss
$\core S\off =\off \Delta\dff S$\nnsp.\oss
By\sss the remarks preceding\trs Lemma\qss \ref{ndc-spaces}\qss this implies\sss that\sss
$S$\sss as a\qss \emph{simplicial\sss set}\pss is\dss canonically\sss isomorphic\sss to\sss
$\bm{\Delta}\dff \core S
\off =\off
\bm{\Delta}\dff \Delta\dff S$\nnsp.\oss
But\sss $S$\sss as a\qss \emph{simplicial\sss space}\pss
may\sss be different\sss from\sss $\bm{\Delta}\dff \Delta\dff S$\nnsp.\oss\vspace{-0.3pt}

\mypar{Lemma.}{free-delta-realizations}
\emph{Let\dss $S$ be a partially\sss ordered space.\oss
Then\sss the geometric realization\dss $\bbnum{S}$ of\pss $S$ 
as a\sss topological\sss simplicial\sss complex\dss
is\dss canonically\dss homeomorphic\sss to\dss
$\num{\Delta\fff S}_{\dff \Delta}$\nsp.\oss}\vspace{-0.3pt}

\proof
As we saw,\oss the spaces of\sss $n$\dnsp-simplices of\sss $\Delta\dff S$\sss and\sss of\dss $S$\sss
as a\sss topological\sss simplicial\sss complex are\sss the same.\oss
Moreover,\oss the diagram\qss (\ref{push-out-tsc})\qss
can\sss be identified\sss with\sss the diagram\qss (\ref{push-out-delta})\qss
for\sss $D\off =\off \Delta\dff S$\nnsp.\oss
An\sss induction shows\sss that\sss
$\bbnum{\ssk_{\dff k}\dff S}$\sss
is\dss canonically\sss homeomorphic\sss to\sss
$\num{\nsp\ssk_{\dff k}\dff D}_{\dff \Delta}$\nsp.\oss
Alternatively,\oss one can\sss simply\sss observe\sss that\sss
$\bbnum{\ssk_{\dff k}\dff S}
\off =\off
\bbnum{S}_{\dff k}$\sss
and\sss
$\num{\nsp\ssk_{\dff k}\dff D}_{\dff \Delta}$\sss
are obtained\sss from\sss the same disjoint\sss union\sss
by\sss taking\sss the quotient\sss by\sss the same equivalence relations.\oss
It\sss remains\sss to pass\sss to\sss 
the direct\sss limit\sss over $k$\nnsp.\oss  \eproof\vspace{-0.3pt}

\mypar{Corollary.}{free-full-realizations}
\emph{If\qss a partially\dss ordered\sss space\dss $S$ has\qss free equalities,\oss
then\sss $\bm{\Delta}\dff S$\nnsp,\oss i.e.\qss
$S$\sss as a simplicial\sss space,\oss is\dss canonically\sss isomorphic\sss
to\sss $\bm{\Delta}\dff \Delta\dff S$\sss and\dss 
the\sss geometric realizations\dss $\num{S}$\sss and\qss $\bbnum{S}$\sss 
are\sss canonically\dss homeomorphic.\oss}\vspace{-0.3pt}

\proof
As we observed\dss right\sss before\trs Lemma\qss \ref{free-delta-realizations},\pss
$\core S
\off =\off 
\Delta\dff S$\nnsp.\oss
If\sss $S$\sss has free equalities,\oss
then\trs Lemma\qss \ref{ndc-spaces}\qss implies\sss that\sss
the canonical\sss map\sss
$\bm{\Delta}\dff \core S\qff \ttoo\qff S$\sss
is\dss an\sss isomorphism of\dss simplicial\sss spaces,\oss
and\sss therefore\sss
$\bm{\Delta}\dff \Delta\dff S
\off =\off
\bm{\Delta}\dff \core S$\sss
is\dss canonically\sss isomorphic\sss to $S$\nnsp.\oss
Also,\oss
Corollary\qss \ref{pos-is-nice}\qss implies\sss that\sss
$\num{S}$\sss is\dss canonically\dss homeomorphic\sss to\sss
$\num{\nsp\core S}_{\dff \Delta}
\off =\off
\num{\Delta\dff S}_{\dff \Delta}$\nsp,\oss
and\trs Lemma\qss \ref{free-delta-realizations}\qss implies\sss that\sss
$\num{\Delta\dff S}_{\dff \Delta}$\sss
is\dss canonically\dss homeomorphic\sss to\sss $\bbnum{S}$\nnsp.\oss
The corollary\sss follows.\oss  \eproof\vspace{-0.3pt}

\myuppar{Order-preserving\sss maps.}
Suppose\sss that\sss $S\fff,\qff T$\sss are 
partially\sss ordered\sss topological\sss simplicial\sss complexes
and\sss
$f\dff \colon\dff S\qff \ttoo\qff T$\dss
is\dss a simplicial\sss map.\oss
If\dss $f$\sss is\qss \emph{strictly\sss order-preserving}\pss
in\sss the sense\sss that\sss
$x\qff <\qff y$\sss implies\sss $f\dff(\dff x\trf)\qff <\qff f\dff(\dff y\trf)$\nnsp,\oss
then\sss $f$\sss induces a simplicial\sss map\sss
$\Delta\fff S\qff \ttoo\qff \Delta\fff T$\dnsp.\oss

Suppose now\sss that\sss $f$\dss is\dss only order-preserving,\oss
i.e.\qss
$x\qff \leq\qff y$\sss implies\sss 
$f\dff(\dff x\trf)\qff \leq\qff f\dff(\dff y\trf)$\nnsp.\oss
Then\sss $f$\sss induces a simplicial\sss map\sss
$\bm{\Delta}\dff S\qff \ttoo\qff \bm{\Delta}\dff T$\dnsp.\oss
By\sss restricting\sss this map\sss to\sss the $\Delta$\dnsp-subspace 
$\Delta\fff S$\sss of\dss $\bm{\Delta}\dff S$\sss
we\sss get\sss a simplicial\sss map of\sss $\Delta$\dnsp-spaces\sss
$f_{\dff *}\dff \colon\dff
\Delta\dff S\qff \ttoo\qff \bm{\Delta}\dff T$\nnsp.\oss\vspace{1pt}
More explicitly,\oss the map\sss
$f_{\dff *}\dff \colon\dff
\Delta\dff S\qff \ttoo\qff \bm{\Delta}\dff T$\sss
can be described as follows.\oss
The $n$\dnsp-simplices of\dss $\Delta\fff S$\sss
are\sss the increasing maps\sss
$\sigma\dff \colon\dff 
[\halfff n\dff]\qff \ttoo\qff S$\nnsp.\oss
The compositions\sss
$f\dff \circ\trf \sigma\dff \colon\dff 
[\halfff n\dff]\qff \ttoo\qff T$\sss
are non-decreasing\sss and\dss hence are {\nsp}$n$\dnsp-simplices of\sss
$\bm{\Delta}\fff T$\nnsp.\oss
A\sss easy check shows\sss that\sss 
$f_{\dff *}\dff(\trf \sigma\trf)
\off =\off
f\dff \circ\trf \sigma$\nnsp.\oss

Suppose now\sss that,\oss in addition,\qss $S\fff,\qff T$\sss
have free equalities.\oss
Let\sss us consider\sss the composition\vspace{0pt}
\[
\quad
\begin{tikzcd}[column sep=sboom, row sep=large]
\protect{\num{\Delta\fff S}_{\dff \Delta}}
\arrow[r, "\dis \protect{\num{f_{\dff *}}_{\dff \Delta}}"]
&
\protect{\num{\bm{\Delta}\dff T}_{\dff \Delta}}
\arrow[r]
&
\protect{\num{\bm{\Delta}\dff T}}
\arrow[r]
&
\protect{\num{\Delta\fff T}_{\dff \Delta}}\qff,
\end{tikzcd}
\]

\vspace{-12pt}\vspace{3pt}
where\sss the middle arrow\dss is\dss the natural\sss quotient\sss map,\oss 
and\sss the right\sss arrow\dss is\dss the canonical\sss homeomorphism\sss
$\num{\bm{\Delta}\dff T}
\qff \ttoo\qff
\num{\bm{\Delta}\dff \Delta\fff T}
\qff \ttoo\qff
\num{\Delta\fff T}_{\dff \Delta}$\qss
(see\dss Section\qss \ref{simplicial-spaces},\oss
remarks after\sss the diagram\qss (\ref{push-out-delta})\qss
for\sss the\sss last\sss arrow).\oss
In view of\trs Lemma\qss \ref{free-delta-realizations}\qss
the geometric realizations\sss
$\num{\Delta\fff S}_{\dff \Delta}$\sss
and\sss
$\num{\Delta\fff T}_{\dff \Delta}$\sss
can\sss be identified\sss
with\sss the geometric realizations\sss
$\bbnum{S}$\sss and\sss $\bbnum{T}$\sss
of\sss $S$\sss and\sss $T$\sss as\sss topological\sss simplicial\sss complexes.\oss
Hence\sss the above composition\dss is\dss a map\sss
$\bbnum{S}\qff \ttoo\qff \bbnum{T}$\nnsp.\oss
A direct\sss verification shows\sss that\sss this map\dss
is\dss nothing else but\sss the induced\sss map\vspace{1.5pt}
\[
\quad
\bbnum{f}\dff \colon\dff 
\bbnum{S}\qff \ttoo\qff \bbnum{T}
\]

\vspace{-12pt}\vspace{1.5pt}
defined above directly\sss in\sss terms of\dss the map of\dss simplicial\sss complexes
$f\dff \colon\dff S\qff \ttoo\qff T$\dnsp.\oss
In view of\trs Corollary\qss \ref{free-full-realizations}\qss
one can also consider\sss this map as a map\sss
$\num{S}\qff \ttoo\qff \num{T}$\dnsp.\oss
Another direct\sss verification shows\sss that\sss this map\dss
is\dss nothing else but\sss
$\num{f}\dff \colon\dff 
\num{S}\qff \ttoo\qff \num{T}$\nnsp.\oss

\mysection{Segal's\pss introductory\qss example}{segal-example}

\myuppar{The category of\dss vector spaces and\qss a\qss Grassmannian.}
Let $\vect$ be\sss the\sss topological\sss category\sss having
finitely dimensional\dss Hilbert\dss spaces 
as objects and\sss isometries as morphisms.\oss
Strictly speaking,\oss we should\sss choose as $\ob\dff \vect$ a set\sss of\dss such\dss Hilbert\dss spaces
containing at\sss least\sss one representative from each isomorphism class.\oss
The set\sss of\dss objects\dss is\dss equipped\sss with\sss the discrete\sss topology,\oss
but\sss for objects\sss $V\fff,\qff W$\sss the set\sss of\dss morphisms\sss
$V\qff \ttoo\qff W$\sss has\sss the usual\sss topology.\oss

Of\dss course,\oss different\sss choices of\dss $\ob\dff \vect$\sss
lead\sss to equivalent\sss topological\sss categories,\oss
and\sss the homotopy\sss type of\dss the classifying space\sss $\num{\vect}$\sss
does not\sss depend on\sss this choice.\oss
For example,\oss it\sss could\sss be\sss the set\sss 
of\dss the spaces $\ccc^{\dff n}$\dnsp,\qss $n\qff \in\qff \nnn$\nnsp.\oss
In\sss this case\sss $\num{\vect}$\sss is\dss the disjoint\sss union of\dss
classifying spaces\sss $B\dff U\fff(\trf \ccc^{\dff n}\trf)$\nnsp,\qss $n\qff \in\qff \nnn$\nnsp,\oss
where $U\fff(\trf \ccc^{\dff n}\trf)$\sss are\sss the unitary\sss groups,\oss
i.e.\qss groups of\dss isometries\sss $\ccc^{\dff n}\qff \ttoo\qff \ccc^{\dff n}$\dnsp.\oss
The classifying spaces are understood\sss in\sss the sense of\trs Segal\qss
\cite{s1},\oss Section\qss 3.\oss
They are different\sss from\trs Milnor's\dss classifying spaces,\oss
but\sss for nice groups such as $U\fff(\trf \ccc^{\dff n}\trf)$\sss
are homotopy equivalent\sss to\sss them.\oss
This determines\sss the homotopy\sss type of\dss $\num{\vect}$\nnsp.\oss

Now we would\sss like\sss to relate\sss $\num{\vect}$\nnsp,\oss
and\sss hence\sss the classifying spaces\sss $B\dff U\fff(\trf \ccc^{\dff n}\trf)$\nnsp,\oss
with a version of\trs Grassmann\trs manifolds.\oss
Let\sss us\sss fix\sss an\sss 
infinitely dimensional\sss separable\dss Hilbert\sss space $H$\sss
and\sss let\sss $G\trf(\trf H\trf)$\sss be\sss the space of\dss finitely dimensional\sss
subspaces of\dss $H$\sss with\sss the usual\sss topology.\oss
The space $G\trf(\trf H\trf)$\sss is\dss the disjoint\sss union of\dss its
connected components $G_{\dff n}\dff(\trf H\trf)$\nnsp,\oss
where\sss $G_{\dff n}\dff(\trf H\trf)$\sss consists of\dss
subspaces of\dss dimension $n$\nnsp.\oss
We can consider\sss $G\trf(\trf H\trf)$\sss as a\sss topological\sss category\sss
having\sss $G\trf(\trf H\trf)$\sss as\sss the space of\dss objects 
and only\sss the identity\sss morphisms.\oss
The\sss topology on\sss the set\sss of\dss morphisms\dss is\dss 
the\sss topology of\sss $G\trf(\trf H\trf)$\nnsp.\oss
Then\sss $\num{G\trf(\trf H\trf)}\off =\off G\trf(\trf H\trf)$\nnsp.\oss

\myuppar{An\sss intermediate category.}
Following\dss Segal\qss \cite{s4},\oss let\sss us\sss introduce a
category\sss $\vect\fff/H$\sss serving as an\sss intermediary\sss between\sss $\vect$\sss
and\sss $G\trf(\trf H\trf)$\nnsp.\oss
The objects of\sss $\vect\fff/H$\sss are pairs\sss $(\trf V\fff,\qss \alpha\trf)$\nnsp,\oss
where $V$\sss is\dss an object\sss of\sss $\vect$\sss and\sss
$\alpha\dff \colon\dff V\qff \ttoo\qff H$\sss is\dss an\sss isometric embedding.\oss
The\sss topology on\sss the set\sss of\dss objects results from\sss 
the discrete\sss topology on\sss $\ob\dff \vect$\sss
and\sss the usual\sss topology on\sss the spaces of\dss 
isometric embeddings\sss $V\qff \ttoo\qff H$\nnsp.\oss
Morphisms\sss 
$(\trf V\fff,\qss \alpha\trf)
\qff \ttoo\qff
(\trf V\fff'\fff,\qss \alpha'\trf)$\sss
are morphisms\sss $\beta\dff \colon\dff V\qff \ttoo\qff V\fff'$\sss 
such\sss that\sss $\alpha\off =\off \alpha'\dff \circ\qff \beta$\nnsp.\oss
So,\oss a morphism\sss
$(\trf V\fff,\qss \alpha\trf)
\qff \ttoo\qff
(\trf V\fff'\fff,\qss \alpha'\trf)$\sss
exists\dss if\trs and\dss only\trs if\dss $V$ and\sss $V\fff'$\sss are isomorphic\dss
Hilbert\sss spaces,\oss
and\dss if\dss they are,\oss it\dss is\dss unique.\oss
In\sss particular,\oss the set\sss of\dss morphisms can\sss be identified\sss with\sss
the set\sss of\dss pairs of\dss isomorphic objects,\oss
and we equip\sss it\sss with\sss the\sss topology\sss induced\sss from\sss
$\ob\dff \vect\qff \times\qff \ob\dff \vect$\nnsp.\oss
In contrast\sss with\sss $G\trf(\trf H\trf)$\sss the\sss topological\sss
category\sss $\vect\fff/H$\sss does not\sss have\sss free units.\oss
There are\sss two forgetting\sss functors,\vspace{0.88pt}
\[
\quad
\vect\fff/\fff H\qff \ttoo\qff \vect
\quad
\mbox{and}\quad
i\dff \colon\dff
\vect\fff/\fff H\qff \ttoo\qff G\trf(\trf H\trf)
\pff.
\]

\vspace{-12pt}\vspace{0.88pt}
The\sss functor\sss
$\vect\fff/\fff H\qff \ttoo\qff \vect$ forgets about\sss the isometric embeddings,\oss
i.e.\qss takes\sss $(\trf V\fff,\qss \alpha\trf)$\sss to\sss $V$\sss
and\sss takes morphisms\sss
$(\trf V\fff,\qss \alpha\trf)
\qff \ttoo\qff
(\trf V\fff'\fff,\qss \alpha'\trf)$\sss
to\sss the corresponding\sss morphisms\sss $V\qff \ttoo\qff V\fff'$\dnsp.\oss
The\sss functor $i$\sss takes embeddings\sss to\sss their\sss images,\oss
i.e.\qss takes\sss $(\trf V\fff,\qss \alpha\trf)$\sss to\sss $\alpha\trf(\trf V\trf)$\nnsp.\oss
For each\sss $W\qff \in\qff G\trf(\trf H\trf)$\sss we will\sss need also\sss the 
full\sss subcategory\sss of\dss $\vect\fff/\fff H$\sss having as objects pairs\sss
$(\trf V\fff,\qss \alpha\trf)$\sss such\sss that\sss $\alpha\trf(\trf V\trf)\off =\off W$\dnsp.\oss
We will\sss denote it\sss by\sss
$\vect\downarrow W$\dnsp.\oss

\mypar{Lemma.}{over-is-contractible}
\emph{The\sss geometric realization\sss
$\num{\textup{\vect}\downarrow W}$\sss 
is\dss contractible.\oss}

\proof
Given\sss two isometries\sss
$\alpha\dff \colon\dff V\qff \ttoo\qff W$\sss
and\sss
$\alpha'\dff \colon\dff V\fff'\qff \ttoo\qff W$\dnsp,\oss
there\dss is\dss a unique\sss isometry\sss
$\beta\dff \colon\dff V\qff \ttoo\qff V\fff'$\sss 
such\sss that\sss $\alpha\off =\off \alpha'\dff \circ\qff \beta$\nnsp.\oss
Therefore\sss for every ordered\sss pair of\dss objects of\dss $\vect\downarrow W$\sss
there\dss is\dss a unique morphism\sss from\sss the first\sss to\sss the second.\oss
It\sss follows\sss that\sss $\vect\downarrow W$\sss is\dss equivalent\sss to\sss
a category with only one object\sss and only one morphism.\oss
The classifying\sss space of\dss the\sss latter category\dss is\dss a one-point\sss space,\oss
and\sss hence\sss $\num{\vect\downarrow W}$\sss
is\dss contractible.\oss  \eproof

\mypar{Theorem.}{vect}
\emph{The maps\sss
$\num{\textup{\vect}}
\off \longleftarrow\off 
\num{\textup{\vect}\fff/H}
\qff \ttoo\qff
\num{G\trf(\trf H\trf)}
\off =\off
G\trf(\trf H\trf)$\sss
induced\sss by\sss these functors 
are homotopy equivalences.\oss}

\proof
Let\sss us consider\sss
$\num{\vect\fff/\fff H}
\qff \ttoo\qff
\num{\vect}$\nnsp.\oss
The space of\dss objects of\sss $\vect$\sss is\dss discrete,\oss
and\sss the space of\dss objects of\sss $\vect\fff/\fff H$\sss is\dss
the disjoint\sss union over\sss the objects $V$ of\sss $\vect$ of\sss
spaces of\dss isometric embeddings $V\qff \ttoo\qff H$\nnsp.\oss
As\dss is\dss well\sss known,\oss
for a finitely dimensional\dss Hilbert\sss space $V$\sss
the space of\dss isometric embeddings\sss $V\qff \ttoo\qff H$\sss
is\dss contractible.\oss
It\dss follows\sss that\sss $\vect\fff/\fff H\qff \ttoo\qff \vect$\sss 
induces a homotopy equivalence on objects.\oss
The space of\dss morphisms of\sss $\vect\fff/\fff H$\sss is\dss 
a\sss locally\sss trivial\dss bundle over\sss the space of\dss
morphisms of\dss $\vect$\nnsp,\oss
with\sss the fiber over a morphism\sss
$V\qff \ttoo\qff V\fff'$\sss of\dss $\vect$\sss being\sss the space 
of\dss isometric embeddings
$V\fff'\qff \ttoo\qff H$\nnsp.\oss
Since\sss the spaces of\dss isometric embeddings
$V\fff'\qff \ttoo\qff H$\sss are contractible,\oss
$\vect\fff/\fff H\qff \ttoo\qff \vect$\sss 
is\dss a homotopy equivalence on\sss morphisms.\oss
An $n$\dnsp-simplex\sss of\dss $\vect\fff/\fff H$\sss
can\sss be identified\sss with\sss an $n$\dnsp-simplex\sss
$V_{\dff 0}\qff \ttoo\qff
V_{\dff 1}\qff \ttoo\qff
\ldots\qff \ttoo\qff
V_{\fff n}$\sss
of\dss $\vect$\sss together\sss with an\sss isometric embedding\sss
$V_{\dff n}\qff \ttoo\qff H$\nnsp.\oss
Therefore\sss the space of $n$\dnsp-simplices of\dss $\vect\fff/\fff H$\sss
is\dss a\sss locally\sss trivial\dss bundle over\sss the space of
$n$\dnsp-simplices of $\vect$\nnsp,\oss
and\sss hence\sss the argument\sss for morphisms works for $n$\dnsp-simplices.

The\sss topological\sss category\sss $\vect$ does not\sss have free units,\oss
and\dss its nerve does not\sss have free degeneracies.\oss
But\sss it\dss is\dss easy\sss to see\sss that\dss its nerve\dss is\dss
both\sss good and\sss proper.\oss
It\dss is\dss good\sss because for every\sss object\sss $V$\sss
the identity\sss $\id_{\trf V}$\sss is\dss a neighborhood deformation\sss 
retract\sss of\dss the group of\dss isometries $V\qff \ttoo\qff V$\dnsp.\oss
It\dss is\dss proper\sss by essentially\sss the same reason\fff:\oss
the subspace which needs\sss to be a cofibration\dss is\dss 
a\dss CW-subcomplex of\dss the whole space.\oss
For\sss the category\sss $\vect\fff/\fff H$\sss the above description
of\sss $n$\dnsp-simplices shows\sss that\sss the space of\sss $n$\dnsp-simplices\dss
is\dss also a\sss locally\sss trivial\sss bundle over\sss the space of\dss
isometric embeddings\sss $V_{\dff n}\qff \ttoo\qff H$\sss
having as fibers\sss the $n$\dnsp-simplices of\dss $\vect$ with\sss the\sss last\sss
object\sss $V_{\fff n}$\nsp.\oss
Now one can easily see\sss that\sss  $\vect\fff/\fff H$\sss is\dss also good and\sss proper.\oss
By a classical\dss theorem of\qss Segal\qss (see\dss Section\qss \ref{simplicial-spaces})\qss
this,\oss together\sss with\sss the results of\dss the previous paragraph,\oss 
implies\sss that\sss
$\num{\vect\fff/\fff H\halfff}
\qff \ttoo\qff
\num{\vect}$\sss
is\dss a\sss homotopy equivalence.\oss  

Let\sss us consider\sss now\sss 
$\num{i}\dff \colon\dff
\num{\vect\fff/H}
\qff \ttoo\qff
\num{G\trf(\trf H\trf)}$\nnsp.\oss
An $n$\dnsp-simplex of\dss $G\trf(\trf H\trf)$\sss as a\sss topological\sss category\sss
has\sss the form\sss
$W\qff \ttoo\qff
W\qff \ttoo\qff
\ldots\qff \ttoo\qff
W$\dnsp,\oss
where\sss $W$\sss is\dss a\sss finitely dimensional\sss subspace of\sss $H$\sss
and each of\dss the {\nsp}$n$ arrows\dss is\dss the identity\sss morphism.\oss
An $n$\dnsp-simplex\vspace{1.5pt}\vspace{1.5pt}
\[
\quad
V_{\dff 0}\qff \ttoo\qff
V_{\dff 1}\qff \ttoo\qff
\ldots\qff \ttoo\qff
V_{\fff n}\qff \ttoo\qff
H
\]

\vspace{-12pt}\vspace{1.5pt}\vspace{1.5pt}
of\dss $\vect\fff/H$\sss
is\dss mapped\sss to\sss
$W\qff \ttoo\qff
W\qff \ttoo\qff
\ldots\qff \ttoo\qff
W$\sss
if\trs and\dss only\trs if\dss the image of\dss
$V_{\fff n}\qff \ttoo\qff H$\sss is\sss $W$\dnsp.\oss
Hence such $n$\dnsp-simplices can\sss be identified\sss with sequences of\dss
isometries of\dss the form\vspace{-0.25pt}
\[
\quad
V_{\dff 0}\qff \ttoo\qff
V_{\dff 1}\qff \ttoo\qff
\ldots\qff \ttoo\qff
V_{\fff n}\qff \ttoo\qff
W
\pff.
\]

\vspace{-12pt}\vspace{-0.25pt}
It\sss follows\sss that\sss the preimage of\dss
$W\qff \in\qff \num{G\trf(\trf H\trf)}\off =\off G\trf(\trf H\trf)$\sss
in\sss $\num{\vect\fff/H}$\sss is\dss equal\sss to 
$\num{\vect\downarrow W}$\nnsp.\oss 
Lemma\qss \ref{over-is-contractible}\qss implies\sss that\sss
$\num{\vect\downarrow W}$\sss is\dss contractible.\oss
It\sss follows\sss that\sss
$\num{\vect\fff/H}
\qff \ttoo\qff
\num{G\trf(\trf H\trf)}$\sss
is\dss a\sss locally\sss trivial\sss bundle with contractible fibers.\oss
Since $G\trf(\trf H\trf)$\sss is\dss paracompact,\oss
a classical\dss theorem\sss of\qss Dold\qss \cite{d}\qss implies\sss that\sss
$\num{\vect\fff/H}
\qff \ttoo\qff
\num{G\trf(\trf H\trf)}$\sss
is\dss a homotopy equivalence.\oss  \eproof\vspace{-0.25pt}

\mypar{Corollary.}{fixed-dim}
\emph{Let\trs $\textup{\vect}_{\dff n}$ be\sss the full\sss subcategory\sss
of\qss $\textup{\vect}$ having as objects vector spaces of\dss dimension $n$\nnsp.\oss
Then\dss $\num{\fff\textup{\vect}_{\dff n}}$ is\dss homotopy equivalent\sss to\dss
$G_{\dff n}\dff(\trf H\trf)$\sss
and\sss hence\sss $G_{\dff n}\dff(\trf H\trf)$\sss is\dss homotopy equivalent\sss to\sss
$B\dff U\fff(\trf \ccc^{\dff n}\trf)$\nnsp.\oss}\vspace{-0.25pt}\vspace{-0.125pt}

\proof
The first\sss claim\sss follows from\sss the fact\sss that\sss our\sss forgetting\sss
functors preserve\sss the dimension of\dss vector spaces involved.\oss
The second one follows\sss from\sss the fact\sss $\num{\vect}$\sss
is\dss homotopy equivalent\sss to\sss the disjoint\sss union of\dss classifying spaces\sss
$B\dff U\fff(\trf \ccc^{\dff n}\trf)$\nnsp.\oss  \eproof\vspace{-0.25pt}

\myuppar{Remark.}
The classifying space\sss $\num{\vect\downarrow W}$\sss
can\sss be identified\sss
with\sss the\sss total\sss space of\dss the universal\sss
$U\fff(\trf W\trf)$\dnsp-bundle over\sss 
$B\dff U\fff(\trf W\trf)$\nnsp,\oss
where\sss $U\fff(\trf W\trf)$\sss is\dss the unitary\sss group of\dss $W$\dnsp,\oss
i.e.\qss the group of\dss isometries\sss $W\qff \ttoo\qff W$\dnsp.\oss
See\qss \cite{s1},\oss Section\qss 3.\oss
This fact\sss can\sss be\sss thought\sss as\sss the\sss true reason of\dss
the contractibility of\dss $\num{\vect\downarrow W}$\dnsp.\oss\vspace{-0.25pt}

\mypar{Proposition.}{not-eq}
\emph{$i\dff \colon\dff
\textup{\vect}\fff/H
\qff \ttoo\qff
G\trf(\trf H\trf)$
is\dss not\sss an equivalence of\dss topological\sss categories.}\vspace{-0.25pt}

\proof
Suppose\sss that\sss
$k\dff \colon\dff
G\trf(\trf H\trf)\qff \ttoo\qff \vect\fff/H$\sss
is\dss a continuous functor such\sss that\sss there exists a\sss natural\sss
transformation\sss 
$i\dff \circ\dff k\qff \ttoo\qff \id$\nnsp.\oss
Let\sss us consider\sss the composition of\sss $k$\sss with\sss
the functor\sss
$\vect\fff/H\qff \ttoo\qff \vect$\nnsp.\oss
Since\sss the space of\dss objects of\dss $\vect$\sss has discrete\sss topology,\oss
this composition\dss is\dss constant\sss on\sss the components of\sss $G\trf(\trf H\trf)$\nnsp.\oss
Let\sss $K_{\dff n}$\sss be\sss the value of\sss 
$k$\sss on\sss the component\sss
$G_{\dff n}\dff(\trf H\trf)$\sss consisting of\dss
subspaces of\dss dimension $n$\nnsp.\oss
Then\sss $k$
provides a continuous map\sss
$V\off \longmapsto\off k\trf(\trf V\trf)$\sss 
from\sss $G_{\dff n}\dff(\trf H\trf)$\sss 
to objects of\dss $\vect\fff/H$\sss of\dss the form\sss 
$(\trf K_{\dff n}\dff,\qff \alpha_{\trf V}\trf)$\nnsp,\oss
and\sss the natural\sss transformation\sss
$i\dff \circ\dff k\qff \ttoo\qff \id$
provides morphisms\sss
$V\qff \ttoo\qff \alpha_{\trf V}\trf(\trf K_{\dff n}\trf)$\sss
of\dss the category\sss $G_{\dff n}\dff(\trf H\trf)$\nnsp.\oss
But\sss $G_{\dff n}\dff(\trf H\trf)$\sss has only\sss identity\sss morphisms,\oss
and\sss hence\sss $\alpha_{\trf V}\trf(\trf K_{\dff n}\trf)\off =\off V$\nnsp.\oss
The maps $\alpha_{\trf V}$ continuously\sss depend on\sss $V$\nnsp,\oss
and\sss since\sss $K_{\dff n}$\sss is\dss a fixed space,\oss
they define a\sss trivialization of\dss the canonical\sss vector bundle over\sss
$G_{\dff n}\dff(\trf H\trf)$\nnsp.\oss
But\sss this bundle\dss is\dss well\sss known\sss to be nontrivial.\oss
In\sss fact,\oss it\dss is\dss a universal\sss bundle.\oss
The proposition\sss follows.\oss  \eproof\vspace{-0.25pt}

\myuppar{Remark.}
Segal\qss \cite{s4}\qss claimed\sss that\sss the functor\sss
$i\dff \colon\dff
\vect\fff/H
\qff \ttoo\qff
G\trf(\trf H\trf)$\sss
is\dss an equivalence of\qss ({\fff}topo\-log\-i\-cal\fff)\qss categories.\oss
Would\sss this be\sss true,\oss
this would\sss immediately\sss imply\sss that\sss the induced\sss map
$\num{\vect\fff/H}
\qff \ttoo\qff
\num{G\trf(\trf H\trf)}$\sss
is\dss a homotopy equivalence.\oss
By\trs Proposition\qss \ref{not-eq}\qss
this claim\dss is\dss not\sss correct.\oss
Strictly speaking,\oss Segal\dss deals 
with\sss real\dss vector spaces and\sss $\rrr^{\dff \infty}$\dnsp,\oss
the union of\dss the sequence\sss
$\rrr^{\dff 0}\qff \subset\qff
\rrr^{\dff 1}\qff \subset\qff
\rrr^{\dff 2}\qff \subset\qff \ldots\qff $\nnsp,\oss
instead of\dss complex\sss vector spaces and\dss $H$\nnsp,\oss
but\sss the above arguments work equally\sss well\sss for\sss
$\ccc^{\dff \infty}$\sss in\sss the role of\dss $H$\nnsp,\oss
as also in\sss the real\sss case.\oss

\mysection{Coverings\qss and\qss simplicial\qss spaces}{coverings}

\myuppar{Simplicial\sss spaces associated\sss to coverings.}
Let\sss $X$\sss be a\sss topological\sss space
and\sss $U_{\dff a}\dff,\pff a\qff \in\qff \Sigma$\sss
be a covering of\sss $X$\nnsp.\oss
Let\sss $\Sigma^{\dff \fin}$\sss be\sss the set\sss of\dss finite non-empty\sss
subsets\sss of\sss $\Sigma$\nnsp.\oss
For\sss $\sigma\qff \in\qff \Sigma^{\dff \fin}$\dss let\vspace{-0.56pt}
\[
\quad
U_{\dff \sigma}
\off =\off
\bigcap\nolimits_{\pff a\qff \in\qff \sigma}\dff U_{\dff a}
\qff.
\]

\vspace{-12pt}
Following\qss Segal\qss \cite{s1},\oss
let\sss us consider\sss the following category\sss $X_{\dff U}$\nsp.\oss
Its objects are pairs\sss $(\trf x\fff,\qff \sigma\trf)$\sss such\sss that\sss
$\sigma\qff \in\qff \Sigma^{\dff \fin}$\sss and\sss
$x\qff \in\qff U_{\dff \sigma}$\nsp.\oss
If\qss $\tau\qff \supset\qff \sigma$ and\sss
$x\qff \in\qff U_{\dff \tau}$\nsp,\oss
then\sss there\dss is\dss a unique morphism\sss
$(\trf x\fff,\qff \tau\trf)
\qff \ttoo\qff
(\trf x\fff,\qff \sigma\trf)$\sss in\sss $X_{\dff U}$\nsp,\oss
and\sss there are no other morphisms.\oss
The space of\dss objects of\dss $X_{\dff U}$\sss is\dss 
the disjoint\sss union of\dss the subspaces\sss
$U_{\dff \sigma}\dff,\pff \sigma\qff \in\qff \Sigma^{\dff \fin}$\dnsp,\oss
and\sss the structure of\dss a category\dss results from\sss
reversing\sss the
order of\dss finite subsets of\sss $\Sigma$\sss by\sss inclusion.\oss
If\dss we consider $X$ as\sss the\sss topological\sss category\sss
having $X$ as\sss the space of\dss objects and only\sss identity\sss morphisms,\oss
then\sss the rule\sss
$(\trf x\fff,\qff \sigma\trf)
\off \longmapsto\off
x$\sss
defines a continuous functor\sss
$\pr\dff \colon\dff
X_{\dff U}\qff \ttoo\qff X$\nnsp,\oss
and\sss hence defines a map\sss
\[
\quad
\num{\pr}\dff \colon\dff
\num{X_{\dff U}}\qff \ttoo\qff \num{X}
\off =\off
X
\qff.
\]

\vspace{-15pt}
\mypar{Theorem.}{covering-categories}
\emph{If\qss the covering\sss
$U_{\dff a}\dff,\pff a\qff \in\qff \Sigma$\sss
is\dss numerable,\oss then\sss  
$\num{\pr}\dff \colon\dff
\num{X_{\dff U}}\qff \ttoo\qff X$\sss
is\dss a homotopy equivalence.\oss
Moreover,\oss 
there\dss exists\dss a homotopy\sss inverse\qss
$\psi\dff \colon\dff
X\qff \ttoo\qff \num{X_{\dff U}}$\dss
such\sss that\qss
$\num{\pr}\qff \circ\qff \psi
\off =\off
\id_{\qff X}$\dss
and\qss
$\psi\qff \circ\qff \num{\pr}$\sss
is\dss homotopic\sss to\sss the identity\sss
in\sss the class of\dss maps 
$f\dff \colon\dff
\num{X_{\dff U}}\qff \ttoo\qff \num{X_{\dff U}}$\dss
such\sss that\qss 
$\num{\pr}\dff \circ\dff f
\off =\off
\num{\pr}$\nnsp.\oss}

\proof
This\dss is\dss a classical\dss theorem of\qss Segal.\oss
See\qss \cite{s1},\oss Proposition\qss 4.1.\oss  \eproof

\myuppar{The case of\dss linearly\sss ordered coverings.}
Suppose\sss that\sss $\Sigma$\sss is\dss linearly ordered.\oss
Then\sss $\num{X_{\dff U}}$\sss admits 
a slightly\sss different\sss description.\oss
Let\sss us\sss consider,\oss
for each non-negative integer $n$\nnsp,\oss 
the set\sss $\Sigma^{\fff n}$ of\dss non-decreasing sequences\sss
$s
\off =\off
(\trf a_{\dff 0}\qff \leq\qff a_{\dff 1}\qff 
\leq\qff \ldots\qff 
\leq\qff a_{\dff n}\trf)$\sss
of\dss elements of\dss $\Sigma$\nnsp.\oss
Such a sequence can\sss be considered as a non-decreasing\sss map\sss
$[\halfff n\dff]\qff \ttoo\qff \Sigma$\nnsp.\oss
For\sss $s\qff \in\qff \Sigma^{\fff n}$\dss let\vspace{-0.56pt}
\[
\quad
U_{\dff s}
\off =\off
\bigcap\nolimits_{\pff 0\qff \leq\qff i\qff \leq\qff n}\dff U_{\dff a_{\dff i}}
\qff.
\]

\vspace{-12pt}
Let\sss $X_{\trf U}^{\dff \leq}$\sss be\sss the\sss simplicial\sss space
having as\sss $n$\dnsp-simplices\sss the pairs\sss 
$(\trf s\fff,\qff x\trf)$\sss such\sss that\sss
$s\qff \in\qff \Sigma^{\fff n}$\sss and\sss
$x\qff \in\qff U_{\dff s}$\nsp.\oss
Equivalently,\oss the space of\sss $n$\dnsp-simplices of\sss $X_{\trf U}^{\dff \leq}$\sss
is\dss the disjoint\sss union of\dss spaces\sss $U_{\dff s}$\sss
over\sss sequences $s\qff \in\qff \Sigma^{\fff n}$\dnsp.\oss
For a non-decreasing\sss
$\theta\dff \colon\dff
[\halfff m\dff]\qff \ttoo\qff [\halfff n\dff]$\sss
the map\sss $\theta^{\dff *}$\sss is\dss defined\sss by\sss
$\theta^{\dff *}\dff(\trf s\fff,\qff x\trf)
\off =\off
(\trf s\dff \circ\dff \theta\fff,\qff x\trf)$\nnsp.\oss
The rule\sss 
$(\trf s\fff,\qff x\trf)\off \longmapsto\off x$\sss
defines a simplicial\sss map\sss
\[
\quad\pr^{\dff \leq}\dff \colon\dff
X_{\trf U}^{\dff \leq}\qff \ttoo\qff X
\pff.
\]

\vspace{-15pt}
\mypar{Theorem.}{covering-ordered}
\emph{There\dss exists\dss a canonical\sss homeomorphism\dss
$h_{\dff U}\dff \colon\dff
\num{X_{\trf U}}
\qff \ttoo\qff
\num{X_{\trf U}^{\dff \leq}}$\sss
such\sss that\dss
$\num{\pr^{\dff \leq}}\qff \circ\pff h_{\dff U}
\off =\off
\num{\pr}$\nnsp.\oss}

\proof
See\dss Dugger\dss and\dss Isaksen\qss \cite{di},\oss
Proposition\qss 2.7\qss and\sss the remarks following it.\oss  \eproof

\mysection{Self-adjoint\qss Fredholm\qss operators}{spaces-operators}

\myuppar{Self-adjoint\dss Fredholm\sss operators and spectral\sss projections.}
For\sss the rest\sss of\dss the paper $H$ denotes a fixed
separable\sss infinitely\sss dimensional\dss Hilbert\sss space.\oss
We are interested\sss mainly\sss in\sss self-adjoint\dss Fredholm\sss operators
$A\dff \colon\dff H\dff \ttoo\dff H$\nnsp.\oss
The cases of\dss bounded\sss and of\dss closed densely defined operators 
are\sss the most\sss fundamental\sss ones.\oss
Most\sss of\dss arguments apply\sss
also\sss to closed densely defined operators with compact\sss resolvent,\oss
but\sss we\sss leave\sss this case\sss to a future occasion.\oss
As usual,\oss we denote by\sss 
$\sigma\dff(\trf A\trf)$\sss the spectrum of\sss $A$\nnsp.\oss
If\sss $A$\sss is\dss closed densely defined and self-adjoint\sss and 
$\rho\qff \subset\qff \rrr$\sss 
is\dss a\dss Borel\sss subset,\oss
we denote by $P_{\dff \rho}\dff(\trf A\trf)$ 
the spectral\sss projection of $A$
associated with $\rho$\nnsp.\oss
For $r\qff \in\qff \rrr$\sss let\vspace{3pt}
\[
\quad 
P_{\dff \geq\dff r}\dff(\trf A\trf)
\off =\off
P_{\dff [\dff r\fff,\dff \infty\dff)}\dff(\trf A\trf)
\qff,
\quad
P_{\dff \leq\dff r}\dff(\trf A\trf)
\off =\off
P_{\dff (\dff -\dff \infty\fff,\dff r\dff]}\dff(\trf A\trf)
\qff,\quad
\mbox{and}\quad
\]

\vspace{-33pt}
\[
\quad
P_{\dff <\dff r}\dff(\trf A\trf)
\off =\off
P_{\dff (\dff -\dff \infty\fff,\dff r\dff)}\dff(\trf A\trf)
\qff,
\quad
P_{\dff >\dff r}\dff(\trf A\trf)
\off =\off
P_{\dff (\dff r\fff,\dff \infty\dff)}\dff(\trf A\trf)
\qff.
\]

\vspace{-12pt}\vspace{3pt}
All\sss these operators are self-adjoint\sss projections.\oss
Clearly,\oss
$P_{\dff\geq\dff r}\dff(\trf A\trf)\qff +\qff P_{\fff<\dff r}\dff(\trf A\trf)
\off =\off
\id$\nnsp,\oss
where\dss $\id$\sss is\dss the identity operator
$H\dff \ttoo\dff H$\nnsp.\oss
We equip\sss the space of\dss self-adjoint\sss projections in\sss $H$\sss
with\sss the usual\dss topology defined\sss by\sss the norm.\oss\vspace{1pt}

\myuppar{The space of\dss self-adjoint\dss Fredholm\sss operators.}
By a well\sss known\sss reason\sss discovered\sss by\dss
Atiyah\dss and\dss Singer\qss \cite{as}\qss 
we consider\sss only
self-adjoint\dss Fredholm\sss operators 
which are not\sss essentially\sss positive or negative and denote by\sss 
$\hat{\mathcal{F}}$\sss the space of\dss such operators.\oss
There are\sss two versions of\dss the space $\hat{\mathcal{F}}$\dnsp.\oss
Namely,\oss one can require\sss that\sss operators in $\hat{\mathcal{F}}$ are bounded,\oss
or only\sss that\sss they are closed and densely defined.\oss 
Our results and\sss proof\sss work equally well\sss for\sss both versions.\oss

In\sss the case of\dss bounded operators\sss the usual\sss norm\sss topology\sss on
$\hat{\mathcal{F}}$\sss is\dss the most\sss natural\sss one.\oss
The choice of\dss topology\sss in\sss the other cases\dss
is\dss less obvious.\oss
We will\sss take a partially\sss axiomatic approach
and\sss begin\sss by\sss listing\sss the key desirable properties 
of\dss topologies on $\hat{\mathcal{F}}$\dnsp.\oss
One can\sss hardly\sss justify an attempt\sss at\sss a\sss full\sss
axiomatisation,\oss and our\sss list\dss is\dss not\sss intended\sss to be full.\oss
After\sss this we will\sss discuss some specific choices.\oss\vspace{1pt}

\myuppar{The essential\sss spectrum and stability.}
If\sss $A$ is\dss a self-adjoint\dss Fredholm\sss operator,\oss then 
$0\qff \in\qff \rrr$ does not\sss belong\sss to its essential\sss spectrum,\oss
i.e.\qss
there exist\sss $u\fff,\qff v\qff \in\qff \rrr$\sss
such\sss that\sss $u\qff <\qff 0\qff <\qff v$\nnsp,\oss
$u\fff,\qff v
\qff \not\in\qff \sigma\dff(\trf A\trf)$\nnsp,\oss
and\sss the image of\dss the projection\sss 
$P_{\dff (\dff u\fff,\dff v\dff)}\dff(\trf A\trf)$\sss
is\dss finitely\sss dimensional.\oss
Clearly,\oss 
$P_{\dff (\dff u\fff,\dff v\dff)}\dff(\trf A\trf)
\off =\off
P_{\dff [\dff u\fff,\dff v\dff]}\dff(\trf A\trf)$\nnsp.\oss\vspace{0.5pt}

Our\sss first\sss requirement\sss to\sss the\sss topology of\dss
$\hat{\mathcal{F}}$\sss is\dss the\qss \emph{stability}\pss
of\dss these properties in\sss the following sense.\oss
If\dss $A\fff,\qff u\fff,\qff v$\sss are
as\sss in\sss the previous paragraph,\oss
then\sss 
$u\fff,\qff v
\qff \not\in\qff \sigma\dff(\trf B\trf)$\sss
and\dss
$P_{\dff [\dff u\fff,\dff v\dff]}\dff(\trf B\trf)$\sss
is\dss finitely dimensional\sss projection for every\sss $B$\sss
belonging\sss to a neighborhood\sss $U$\sss of\sss $A$\sss
in\sss $\hat{\mathcal{F}}$\dnsp.\oss
Moreover,\oss the map\dss
$B
\off \longmapsto\off
P_{\dff [\dff u\fff,\dff v\dff]}\dff(\trf B\trf)$\sss
is\dss a continuous map from\sss $U$\sss to\sss the space of\dss
self-adjoint\sss projections.\oss
In\sss particular,\pss $[\dff u\fff,\qff v\trf]$\sss
is\dss disjoint\sss from\sss the essential\sss spectrum of\dss $B$\nnsp.\oss

\myuppar{Uniformly\sss invertible operators and contractibility.}
Suppose\sss that\sss 
$u\fff,\qff v\qff \in\qff \rrr$\sss
and\sss
$u\qff <\qff 0\qff <\qff v$\nnsp.\oss
Let\sss 
$\hat{\mathcal{F}}^{\dff \inv}\trf[\dff u\fff,\fff v\trf]
\qff \subset\pff
\hat{\mathcal{F}}$\sss
be\sss the subspace of\dss operators\sss $A\qff \in\qff \hat{\mathcal{F}}$\sss
such\sss that\sss 
$\sigma\dff(\dff A\dff)$\sss is\dss
disjoint\sss from\sss $[\dff u\fff,\qff v\trf]$\nnsp,\oss
or,\oss equivalently,\oss that\sss 
$\sigma\dff(\dff A\dff)
\qff \subset\qff 
(\dff -\qff \infty\fff,\qff u\trf)
\dff \cup\dff
(\trf v\fff,\qff \infty\dff)$\nnsp.\oss
One may say\sss that\sss such operators are\qss
\emph{uniformly\sss invertible}\pss
with bounds\sss $u\fff,\qff v$\nnsp.\oss\vspace{-0.62pt}

Our\sss second\sss requirement\sss to\sss the\sss topology of\dss the space\sss
$\hat{\mathcal{F}}$\sss is\dss the\qss \emph{contractibility}\pss
of\dss the subspaces\sss 
$\hat{\mathcal{F}}^{\dff \inv}\trf[\dff u\fff,\fff v\trf]$\nnsp.\oss
This\dss is\dss a more precise and slightly stronger version of\dss
the contractibility alluded\sss to by\trs Segal\qss \cite{s4}.\oss
See\sss the quote at\sss the beginning of\dss the present\sss paper.\oss\vspace{-0.62pt}

This contractibility\sss property\dss is\dss highly\sss non-trivial.\oss
When\sss it\dss is\dss known,\oss its\sss proof\dss depends either
on\dss Kuiper's\dss theorem\qss \cite{ku}\qss on\sss the contractibility
of\dss the unitary group of\dss a\dss Hilbert\sss space equipped\sss with\sss the norm\sss topology,\oss
or on an argument\sss of\trs Dixmier\sss and\dss Douady\qss \cite{dd}\qss
proving\sss the contractibility of\dss this unitary group equipped\sss with\sss
the strong operator\sss topology.\oss
In\sss the present\sss paper\sss we will\sss rely on\trs Kuiper's\dss theorem,\oss
but\sss not\sss on\trs Dixmier--Douady\trs argument.\oss
Kuiper's\dss theorem\dss usually\sss enters\sss the proofs not\sss directly,\oss
but\sss via a\sss theorem of\trs Atiyah\sss and\dss Singer\qss \cite{as},\oss
which we will\sss explain\sss now.\oss\vspace{-0.62pt}

\myuppar{The unrestricted\dss Grassmannian.}
The\qss \emph{unrestricted\dss Grassmannian}\qss
$\mathbf{G{\fff}r}$\sss is\dss the set\sss of\dss closed subspaces $K$ of\sss $H$\sss
such\sss that\sss both\sss the dimension and\sss the codimension of\sss $K$\sss
is\dss infinite.\oss
By identifying a closed subspace $K$ with\sss the orthogonal\sss projection\sss
$H\qff \ttoo\qff K$\sss we can\sss identify\sss $\mathbf{G{\fff}r}$\sss
with\sss a space of\dss bonded operators in\sss $H$\nnsp.\oss
We equip $\mathbf{G{\fff}r}$\sss by\sss the\sss topology\sss induced\sss from\sss
the norm\sss topology\sss by\sss this identification.\oss
A key\sss fact\dss is\dss the contractibility of\sss $\mathbf{G{\fff}r}$\dnsp.\oss
This\dss is\dss due\sss to\dss Atiyah\dss and\dss Singer\qss \cite{as}\qss
and\sss follows\sss from\trs Kuiper's\trs theorem about\sss
the contractibility of\dss the unitary groups of\sss infinitely
dimensional\dss Hilbert\dss spaces.\oss 
See\qss \cite{as},\oss the proof\dss of\qss Lemma\qss 3.6.\vspace{-0.62pt}

\myuppar{Half-line projections and\sss uniformly\sss positive operators.}
A natural\sss way\sss to prove\sss the contractibility\sss property\dss
is\dss to reduce\sss it\sss to\sss two other.\oss
The first\sss one\dss is\dss the
stability of\dss half-line projections\sss
in\sss the following sense.\oss
Let\sss $u\fff,\qff v$\sss be as above and\sss 
$A\qff \in\qff \hat{\mathcal{F}}^{\dff \inv}\trf[\dff u\fff,\fff v\trf]$\nnsp.\oss
The\qss \emph{stability of\dss half-line projections}\pss property\sss
requires\sss that\sss then\sss
$B\qff \in\qff \hat{\mathcal{F}}^{\dff \inv}\trf[\dff u\fff,\fff v\trf]$\sss
for every\sss $B$\sss belonging\sss to a neighborhood\sss $U$\sss of\sss $A$\sss
in\sss $\hat{\mathcal{F}}$\sss and\sss 
$B
\off \longmapsto\off
P_{\dff \geq\dff v}\dff(\trf B\trf)$\sss
is\dss a continuous map from\sss $U$\sss to\sss the space of\dss
self-adjoint\sss projections.\oss
Clearly,\oss then\sss 
$B
\off \longmapsto\off
P_{\dff \leq\dff u}\dff(\trf B\trf)$\sss
is\dss also continuous.\oss\vspace{-0.62pt}

Suppose\sss that\sss $\varepsilon\qff \geq\qff 0$\nnsp.\oss
Let\sss $\hat{\mathcal{F}}_{\qff >\qff \varepsilon}$\sss
be\sss the space of\dss self-adjoint\sss operators
$A\dff \colon\dff H\dff \ttoo\dff H$
such\sss that
$\sigma\dff(\dff A\dff)
\qff \subset\qff 
(\dff \varepsilon\fff,\qff \infty\dff)$\nnsp.\oss
Clearly,\oss such operators are\dss Fredholm.\oss
One may say\sss that\sss they are\qss
\emph{uniformly\sss positive}\pss
with bound\sss $\varepsilon$\nnsp.\oss
The second\sss property\dss is\dss the\qss \emph{contractibility}\pss
of\dss the spaces\sss $\hat{\mathcal{F}}_{\qff >\qff \varepsilon}$\nsp.\oss\vspace{-0.62pt}

\myuppar{Topologies on\sss $\hat{\mathcal{F}}$\dnsp.}
In\sss the case of\dss bounded operators\sss the space $\hat{\mathcal{F}}$\sss
is\dss considered\sss with\sss its norm\sss topology.\oss
In\sss this case\sss the first\sss stability\sss property\sss and\sss the stability\sss
of\dss half-line projections are well\sss known.\oss
In\sss the case of\dss closed and densely defined operators\sss 
the first\sss choice\dss is\dss the\sss topology of\dss 
convergence in\sss the uniform resolvent\sss sense.\oss
We refer\sss to\qss \cite{rs},\oss Section\qss VIII.7\qss for\sss its definition 
and\sss basic properties.\oss
The first\sss stability\sss property\dss follows,\qss for example,\oss from\qss
\cite{rs},\oss Theorem\qss VIII.23\qss (as stated,\oss this\sss theorem\dss is\dss
concerned only\sss with\sss sequences,\oss
but\sss the arguments are completely\sss general\fff).\oss

But\sss the stability of\dss half-line projections does not\sss hold,\oss
as a classical\sss example of\trs Rellich\dss shows.\oss
See,\oss for example,\oss \cite{k},\oss 
Chapter\qss V,\oss Example\qss 4.13.\oss
One can deal\sss with\sss this problem by\sss refining\sss 
the\sss topology of\sss $\hat{\mathcal{F}}$\sss 
simply\sss by declaring\sss the maps\sss
$B
\off \longmapsto\off
P_{\dff \geq\dff v}\dff(\trf B\trf)$\sss 
to be continuous.\oss
Alternatively,\oss let\sss
$\chi\dff \colon\dff \rrr\qff \ttoo\qff \rrr$ 
be a strictly\sss increasing continuous function such\sss that\vspace{3pt}
\[
\quad
\lim\nolimits_{\dff t\qff \to\qff \infty}\qff \chi\dff(\dff t\trf)
\off =\off 
1\quad
\mbox{and}\quad
\lim\nolimits_{\dff t\qff \to\qff -\dff\infty}\qff \chi\dff(\dff t\trf)
\off =\off 
-\qff 1
\qff.
\]

\vspace{-12pt}\vspace{3pt}
Let\sss us\sss refine\sss the\sss topology\sss of\dss
$\hat{\mathcal{F}}$\sss by\sss requiring\sss the map\sss
$A\qff \longmapsto\qff \chi\dff(\dff A\dff)$\sss to be continuous.\oss
This refined\sss topology\dss
is\dss independent\sss on\sss the choice of\dss $\chi$\nnsp.\oss
See\qss \cite{rs},\oss Theorem\qss VIII.20.\oss
It\dss is\dss known\sss as\sss the\qss
\emph{Riesz\dss topology}\pss and\trs is\dss often\sss defined\sss in\sss
terms of\dss a particular choice of\sss $\chi$\nnsp,\oss
namely,\qss
\[
\quad
\chi\dff \colon\dff
t\off \longmapsto\off 
t\left/\sqrt{t^{\dff 2}\qff +\qff 1}\right.
\off.
\]

\vspace{-12pt} 
A routine argument\sss based on\sss the functional\sss calculus\qss
(see\qss \cite{rs},\oss Theorem\qss VIII.20)\qss
shows\sss that\sss the the maps\sss
$B
\off \longmapsto\off
P_{\dff \geq\dff v}\dff(\trf B\trf)$\sss 
are continuous in\dss Riesz\dss topology.\oss
Also,\oss it\dss is\dss well\sss known\sss that\sss on\sss the subspace
of\dss bounded operators\dss Riesz\dss topology and\sss the norm\sss topology agree.\oss
From\sss now on\sss we will\sss equip\sss $\hat{\mathcal{F}}$\sss
with\sss the\dss Riesz\dss topology\sss for both classes of\dss operators.\oss

\mypar{Proposition.}{positive-contractible}
\emph{In\sss both cases of\qss bounded and of\qss closed densely defined operators\sss
the spaces\sss $\hat{\mathcal{F}}_{\qff >\qff \varepsilon}$\sss are contractible.\oss}

\proof
In\sss the case of\dss bounded operators,\oss
if\dss $r\qff >\qff \varepsilon$\nnsp,\oss
then $r\fff \id\qff \in\qff \hat{\mathcal{F}}_{\qff >\dff \varepsilon}$\dss 
and\sss the\sss linear\sss homotopy\sss 
$(\trf A\fff,\qff t\trf)
\off \longmapsto\off
(\trf 1\qff -\qff t\trf)\trf A
\qff +\qff
t\fff r\fff \id$\nnsp,\qss
$0\qff \leq\qff t\qff \leq\qff 1$\nnsp,\oss
contracts\sss $\hat{\mathcal{F}}_{\qff >\qff \varepsilon}$\sss
to $r\fff \id$\nnsp.\oss
Let\sss us consider\sss the case of\dss closed\sss densely\sss defined operators.\oss
Let\sss $\mathcal{Z}$\sss be\sss the space of\dss 
bound\-ed self-adjoint\sss operators\sss $A$\sss
such\sss that\sss $\norm{A}\qff \leq\qff 1$\sss and\sss the kernel\sss 
$\kernel\dff \id_{\trf H}\qff -\qff A^{\dff 2}$\sss is\sss equal\sss to\sss $0$\nnsp,\oss
equipped\sss with\sss the norm\sss topology.\oss
Then\sss for\sss the above choice of\sss $\chi$\sss the map\sss 
$A\off \longmapsto\off \chi\trf(\trf A\trf)$\sss defines a homeomorphism\sss
between\sss the space of\dss closed densely defined self-adjoint\sss operators
with\sss the\dss Riesz\dss topology and\sss $\mathcal{Z}$\nnsp.\oss
See\qss \cite{sch},\oss Section\qss 7.3.\oss
Clearly,\oss this homeomorphism\sss takes\sss 
$\hat{\mathcal{F}}_{\qff >\qff \varepsilon}$\sss
to\sss
$\hat{\mathcal{F}}_{\qff >\qff \chi\trf(\dff \varepsilon\dff)}
\qff \cap\qff
\mathcal{Z}$\nnsp.\oss
Let\sss us\sss choose\sss
$c\qff \in\qff 
(\trf \chi\trf(\dff \varepsilon\dff)\fff,\qff 1\trf)$\sss
and consider\sss families\sss
$A_{\dff t}
\off =\off
t\dff A\qff +\qff (\trf 1\qff -\qff t\trf)\dff c\dff \id_{\trf H}$\nsp,\qss
$0\qff \leq\qff t\qff \leq\qff 1$\nnsp.\oss
If\dss 
$A
\qff \in\qff 
\hat{\mathcal{F}}_{\qff >\qff \chi\trf(\dff \varepsilon\dff)}$\nsp,\oss
then\sss
$A_{\dff t}
\qff \in\qff 
\hat{\mathcal{F}}_{\qff >\qff \chi\trf(\dff \varepsilon\dff)}$\sss
for every\sss $t$\nnsp.\oss
An easy check\sss shows\sss that\trs if,\oss in\sss addition,\pss
$A\qff \in\qff \mathcal{Z}$\nnsp,\oss
then\sss
$A_{\dff t}\qff \in\qff \mathcal{Z}$\sss
for every\sss $t$\nnsp.\oss
It\sss follows\sss that\sss these families define a homotopy contracting\sss
$\hat{\mathcal{F}}_{\qff >\qff \chi\trf(\dff \varepsilon\dff)}
\qff \cap\qff
\mathcal{Z}$\sss
to\sss $c\dff \id_{\trf H}$\nsp.\oss
Therefore\sss this intersection\dss is\dss contractible,\oss
and\sss hence\sss $\hat{\mathcal{F}}_{\qff >\qff \varepsilon}$\sss
is\dss also contractible.\oss  \eproof

\mypar{Proposition.}{invertible-contractible}
\emph{In\sss both cases of\qss bounded and of\qss closed densely defined operators\sss
the spaces\sss 
$\hat{\mathcal{F}}^{\dff \inv}\trf[\dff u\fff,\fff v\trf]$\sss 
are contractible.\oss}

\proof
The rule\sss
$A
\off \longmapsto\off 
\image\dff P_{\qff \geq\qff v}\dff(\dff A\dff)$\sss
defines a continuous map\sss
$\hat{\mathcal{F}}^{\dff \inv}\trf[\dff u\fff,\fff v\trf]
\qff \ttoo\qff
\mathbf{G{\fff}r}$\nnsp.\oss
It\dss is\dss a\sss locally\sss trivial\sss bundle with\sss the fibers homeomorphic\sss to\sss
$\hat{\mathcal{F}}_{\qff >\dff -\dff u}
\dff \times\qff 
\hat{\mathcal{F}}_{\qff >\dff v}$\nsp.\oss
As we mentioned above,\oss the base\sss $\mathbf{G{\fff}r}$\sss of\dss this bundle\dss
is\dss contractible.\oss
By\dss Proposition\qss \ref{positive-contractible}\qss the fibers are also contractible.\oss
Since\sss the base\sss $\mathbf{G{\fff}r}$\dnsp,\oss being a metric space,\oss is\dss paracompact,\oss
this implies\sss that\sss the\sss total\sss space\sss
$\hat{\mathcal{F}}^{\dff \inv}\trf[\dff u\fff,\fff v\trf]$\dss 
is\dss also contractible.\oss  \eproof

\mysection{Classifying\qss spaces\qss 
for\qss self-adjoint\qss Fredholm\qss operators}{classifying-spaces-saf}

\myuppar{Enhanced operators.}
An\qss \emph{enhanced\qss (self-adjoint\trs Fredholm)\qss operator}\pss is\dss
defined as a pair\sss $(\trf A\dff,\qff \varepsilon\trf)$\nnsp,\oss
where $A\qff \in\qff \hat{\mathcal{F}}$ and $\varepsilon\qff \in\qff \rrr$\sss
are such\sss that\sss $\varepsilon\qff >\qff 0$\nnsp,\oss the\sss interval\sss
$[\dff -\qff \varepsilon\fff,\qff \varepsilon\trf]$\sss is\dss disjoint\sss from\sss the essential\sss
spectrum of\sss $A$\nnsp,\oss
and\sss $-\qff \varepsilon\fff,\qff \varepsilon\qff \not\in\qff \sigma\dff(\dff A\dff)$\nnsp.\oss
Let\sss $\hat{\mathcal{E}}$ be\sss the set\sss of\dss enhanced operators
equipped\sss with\sss the\sss topology defined\sss by\sss
the\sss topology of\sss $\hat{\mathcal{F}}$\sss
and\sss the discrete\sss topology on\sss the\qss
\emph{controlling\sss parameters}\qss $\varepsilon\qff \in\qff \rrr$\nnsp.\oss
The space\sss $\hat{\mathcal{E}}$\sss 
is\dss ordered\sss by\sss the relation\sss
\[
\quad
(\trf A\dff,\qff \varepsilon\trf)
\off \leq\off
(\trf A'\dff,\qff \varepsilon'\trf)
\quad
\mbox{if}\quad 
A\off =\off A'
\quad
\mbox{and}\quad 
\varepsilon\qff \leq\qff \varepsilon'
\pff.
\]

\vspace{-12pt}
This order defines a structure of\dss a\sss topological\sss
category on $\hat{\mathcal{E}}$\sss having\sss $\hat{\mathcal{E}}$\sss as\sss the space of\dss objects
and a single morphism\sss
$(\trf A\dff,\qff \varepsilon\trf)
\qff \ttoo\qff
(\trf A\dff,\qff \varepsilon'\trf)$\dss
if\trs $A\off =\off A'$\sss and\sss
$\varepsilon\qff \leq\qff \varepsilon'$\nnsp.\oss

We can also consider $\hat{\mathcal{F}}$ as a\sss topological\sss category\sss
having $\hat{\mathcal{F}}$ 
as\sss the space of\dss objects\qss and only\sss identity\sss morphisms.\oss
Then\sss the classifying space $\num{\hat{\mathcal{F}}}$\sss
is\dss equal\sss to $\hat{\mathcal{F}}$ considered as a\sss topological\sss space.\oss
The obvious forgetting\sss functor\sss
$\hat{\varphi}\dff \colon\dff \hat{\mathcal{E}}\qff \ttoo\qff \hat{\mathcal{F}}$\sss
induces a map of\dss classifying spaces
$\num{\hat{\varphi}}\dff \colon\dff 
\num{\hat{\mathcal{E}}}
\qff \ttoo\qff 
\num{\hat{\mathcal{F}}}
\off =\off
\hat{\mathcal{F}}$\dnsp.\oss

\mypar{Theorem.}{forgetting-enhancement}
\emph{The map\sss
$\num{\hat{\varphi}}\dff \colon\dff 
\num{\hat{\mathcal{E}}}\qff \ttoo\qff \num{\hat{\mathcal{F}}}
\off =\off 
\hat{\mathcal{F}}$\sss
is\dss a\sss homotopy\sss equivalence.\oss}

\proof
Let\sss $\Sigma\off =\off \rrr_{\trf >\trf 0}$\nsp.\oss
For every\sss $a\qff >\qff 0$\sss
let\sss $U_{\dff a}\qff \subset\qff \hat{\mathcal{F}}$\sss
be\sss the set\sss of\sss operators\sss $A\qff \in\qff \hat{\mathcal{F}}$\sss
such\sss that\sss $(\trf A\fff,\qff a\trf)$\sss
is\dss an enhanced operator.\oss
Then\sss $U_{\dff a}\dff,\qff a\qff \in\qff \Sigma$\sss
is\dss an open covering of\sss $\hat{\mathcal{F}}$\dnsp.\oss
Let\sss us consider\sss 
the category\sss $\hat{\mathcal{F}}_{\trf U}$
and\sss the functor\sss
$\pr\dff \colon\dff \hat{\mathcal{F}}_{\trf U}\qff \ttoo\qff \hat{\mathcal{F}}$\nsp.\oss
Using\sss the standard order on\sss $\Sigma\off =\off \rrr_{\trf >\trf 0}$\sss
we can also construct\sss the simplicial\sss space\sss
$\hat{\mathcal{F}}_{\trf U}^{\dff \leq}$\nsp.\oss

By\sss the definition,\pss $\hat{\mathcal{E}}$\sss
is\dss the disjoint\sss union of\dss sets\sss
$U_{\dff a}\dff,\qff a\qff \in\qff \Sigma$\nnsp.\oss
The order on\sss $\hat{\mathcal{E}}$\sss is\dss induced\sss by\sss
the usual\sss order on\sss $\Sigma\off =\off \rrr_{\trf >\trf 0}$\nsp.\oss
An $n$\dnsp-simplex of\dss $\hat{\mathcal{E}}$\sss is\dss determined\sss
by an operator $A\qff \in\qff \hat{\mathcal{F}}$ 
and a non-decreasing\sss sequence\sss
$s
\off =\off
(\trf a_{\dff 0}\qff \leq\qff a_{\dff 1}\qff \leq\qff \ldots\qff \leq\qff a_{\dff n}\trf)$\sss 
of\dss positive numbers such\sss that\sss
$(\trf A\fff,\qff a_{\dff i}\trf)$\sss is\dss an enhanced operator for every\sss $i$\nnsp.\oss
Since\sss $\rrr_{\trf >\trf 0}$\sss is\dss equipped\sss with\sss the discrete\sss topology,\oss
the space of\sss $n$\dnsp-simplices of\sss $\hat{\mathcal{E}}$\sss is\dss
the disjoint\sss union of\dss the spaces
\[
\quad
U_{\dff s}
\off =\off
\bigcap\nolimits_{\pff 0\qff \leq\qff i\qff \leq\qff n}\dff U_{\dff a_{\dff i}}
\qff
\]

\vspace{-12pt}
over all\sss sequences $s\qff \in\qff \Sigma^{\fff n}$\dnsp.\oss
Therefore\sss we can\sss identify\sss the spaces of\dss $n$\dnsp-simplices of\sss
$\hat{\mathcal{E}}$\sss and\sss
$\hat{\mathcal{F}}_{\trf U}^{\dff \leq}$\nsp.\oss
Comparing\sss the definitions shows\sss that\sss this\sss identification\sss
respects\sss the structure maps
$\theta^{\dff *}$\dnsp.\oss
It\sss follows\sss that\sss the simplicial\sss spaces\sss
$\hat{\mathcal{E}}$\sss and\sss
$\hat{\mathcal{F}}_{\trf U}^{\dff \leq}$\dss
are\sss canonically\sss isomorphic.\oss
One may even say\sss that\sss they are equal.\oss
Comparing\sss the definitions shows\sss that\sss the maps\vspace{1.5pt}
\[
\quad
\hat{\varphi}\dff \colon\dff \hat{\mathcal{E}}\qff \ttoo\qff \hat{\mathcal{F}}
\quad
\mbox{and}\dff\quad
\pr^{\dff \leq}\dff \colon\dff
\hat{\mathcal{F}}_{\trf U}^{\dff \leq}\qff \ttoo\qff \hat{\mathcal{F}}
\]

\vspace{-12pt}\vspace{1.5pt}
are\sss also equal.\oss
Therefore\sss it\dss is\dss sufficient\sss to prove\sss that\sss the map\sss 
$\num{\pr^{\dff \leq}}\qff \colon\qff
\num{\hat{\mathcal{F}}_{\trf U}^{\dff \leq}}\qff \ttoo\qff \hat{\mathcal{F}}$\sss
is\dss a\sss homotopy equivalence.\oss
Since\sss $\hat{\mathcal{F}}$\sss is\dss a metric space and,\oss
in\dss particular,\oss is\dss paracompact,\oss
the\sss covering\sss $U_{\dff a}\dff,\qff a\qff \in\qff \Sigma$\sss is\dss numerable.\oss
By\trs Theorem\qss \ref{covering-categories}\qss the map\sss
$\num{\pr}\dff \colon\dff
\num{\hat{\mathcal{F}}_{\trf U}}\qff \ttoo\qff \hat{\mathcal{F}}$\sss
is\dss a homotopy equivalence.\oss
It\sss remains\sss to apply\trs Theorem\qss \ref{covering-ordered}.\oss  \eproof

\myuppar{Enhanced operator\dss models.}
An\qss \emph{enhanced\sss operator\dss model}\pss is\dss defined as a\sss triple\sss
$(\trf V\fff,\pff F\fff,\pff \varepsilon\trf)$\nnsp,\oss
where\sss $V$\sss is\dss a finitely dimensional\sss subspace of\sss $H$\nnsp,\oss
$F\dff \colon\dff V\qff \ttoo\qff V$\sss is\dss a self-adjoint\sss operator,\oss
and $\varepsilon\qff \in\qff \rrr$\sss is\dss 
such\sss that\sss $\varepsilon\qff >\qff 0$\sss and\dss 
$\sigma\dff(\trf F\trf)
\off \subset\off
(\dff -\qff \varepsilon\fff,\pff \varepsilon\trf)$\nnsp.\oss
Let\sss $\mathcal{E}\hat{\mathcal{O}}$\sss be\sss the set\sss of\dss enhanced operator\sss models\dss
equipped\sss with\sss topology\sss defined\sss by\sss the obvious\sss
topology on\sss pairs $(\trf V\fff,\pff F\trf)$\sss and\sss the 
discrete\sss topology on\sss the parameters\sss 
$\varepsilon\qff \in\qff \rrr$\nnsp.\oss
The space\sss $\mathcal{E}\hat{\mathcal{O}}$\sss 
is\dss ordered\sss by\sss the relation\vspace{1.5pt}\vspace{0.75pt}
\[
\quad
(\trf V\fff,\pff F\fff,\pff \varepsilon\trf)
\off \leq\off
(\trf V\fff'\fff,\pff F\fff'\fff,\pff \varepsilon'\trf)
\]

\vspace{-36pt}\vspace{0.75pt}
\[
\quad
\mbox{if}\fff\quad
V\qff \subset\qff V\fff'\dff,\quad 
F\off =\off F\fff'\qff \bigl|\halfff\qff V\qff,\quad  
\varepsilon
\qff \leq\qff
\varepsilon'\qff,\off\off  
\]

\vspace{-12pt}\vspace{1.5pt}\vspace{0.75pt}
and\sss all\sss eigenvectors of\sss $F\fff'$ with eigenvalues\sss in\sss
$[\dff -\qff \varepsilon\fff,\qff \varepsilon\trf]$ are,\oss in\sss fact,\oss
eigenvectors of\sss $F$\dnsp.\oss 
In\sss particular,\qss
${}-\qff \varepsilon\fff,\pff \varepsilon$ are not\sss eigenvalues of\sss $F\fff'$\dnsp.\oss

As in\sss the case of\sss $\hat{\mathcal{E}}$\dnsp,\oss
we can consider\sss this order as\sss a structure of\dss a\sss topological\sss
category on\sss $\mathcal{E}\hat{\mathcal{O}}$\dnsp.\oss
There\dss is\dss an obvious functor\sss
$\hat{\psi}\dff \colon\dff 
\hat{\mathcal{E}}
\qff \ttoo\qff 
\mathcal{E}\hat{\mathcal{O}}$\sss
taking an enhanced operator\sss
$(\trf A\dff,\qff \varepsilon\trf)$\sss
to\qss \emph{its\sss enhanced\sss operator\sss model}\qss
$(\trf V\fff,\pff F\fff,\pff \varepsilon\trf)$\nnsp,\oss
where\vspace{1.5pt}
\[
\quad
V
\off =\off 
\image\dff P_{\dff [\dff -\qff \varepsilon\fff,\qff \varepsilon\trf]}\dff(\trf A\trf)
\off =\off 
\image\dff P_{\dff (\dff -\qff \varepsilon\fff,\qff \varepsilon\trf)}\dff(\trf A\trf)
\]

\vspace{-12pt}\vspace{1.5pt}
and\sss the operator\sss 
$F\dff \colon\dff V\qff \ttoo\qff V$\sss
is\trs induced\sss by $A$\nnsp.\oss

\mypar{Theorem.}{to-enhanced-models}
\emph{The map\sss
$\num{\hat{\psi}}\dff \colon\dff
\num{\hat{\mathcal{E}}}
\qff \ttoo\qff 
\num{\mathcal{E}\hat{\mathcal{O}}}$\sss
is\dss a\sss homotopy\sss equivalence.\oss}

\proof
An $n$\dnsp-simplex of\dss the category\sss $\hat{\mathcal{E}}$\sss
can\sss be identified\sss with an operator $A\qff \in\qff \hat{\mathcal{F}}$
together\sss with\sss a non-decreasing\sss sequence\sss
$a_{\dff 0}\qff \leq\qff a_{\dff 1}\qff \leq\qff \ldots\qff \leq\qff a_{\dff n}$\sss 
of\dss positive numbers such\sss that\sss
$(\trf A\fff,\qff a_{\dff i}\trf)$\sss is\dss 
an enhanced operator for every\sss $i$\nnsp.\oss

Similarly,\oss an $n$\dnsp-simplex of\dss the category\sss $\mathcal{E}\hat{\mathcal{O}}$\sss
can\sss be identified\sss with a finitely\sss dimensional\sss subspace\sss 
$V_{\fff n}\qff \subset\qff H$\sss
together\sss with\sss a self-adjoint\sss operator
$F_{\fff n}\dff \colon\dff
V_{\fff n}\qff \ttoo\qff V_{\fff n}$\sss
and\sss a non-de\-creas\-ing\sss sequence\sss
$a_{\dff 0}\qff \leq\qff a_{\dff 1}\qff \leq\qff \ldots\qff \leq\qff a_{\dff n}$\sss
of\dss positive numbers such\sss that\sss
$\sigma\dff(\dff F_{\fff n}\dff)
\off \subset\off
(\dff -\qff a_{\dff n}\fff,\pff a_{\dff n}\trf)$\sss
and\sss the numbers\sss $-\qff a_{\dff i}\dff,\qff a_{\dff i}$\sss 
are not\sss eigenvalues of\sss $F_{\fff n}$\nsp.\oss
The subspaces $V_{\fff i}$\sss with $i\qff <\qff n$\sss can\sss be recovered as\sss
the images\sss
$\image\dff P_{\dff [\dff -\qff a_{\dff i}\fff,\qff a_{\dff i}\trf]}\dff(\trf F_{\fff n}\trf)$\nnsp.\oss

The map\sss
$\hat{\psi}_{\dff n}\dff \colon\dff
\hat{\mathcal{E}}_{\dff n}
\qff \ttoo\qff
\mathcal{E}\hat{\mathcal{O}}_{\dff n}$\sss
induced\sss by $\hat{\psi}$\sss
keeps\sss the sequence\sss
$a_{\dff 0}\qff \leq\qff a_{\dff 1}\qff \leq\qff \ldots\qff \leq\qff a_{\dff n}$\sss
and\sss takes $A$\sss to\sss the subspace\sss
$V_{\fff n}$\sss and\sss the operator\sss $F_{\fff n}\dff \colon\dff
V_{\fff n}\qff \ttoo\qff V_{\fff n}$\sss induced\sss by\sss $A$\nnsp,\oss
where\vspace{1.5pt}
\begin{equation}
\label{vn}
\quad
V_{\fff n}
\off =\off 
\image\dff P_{\dff [\dff -\qff a_{\dff n}\fff,\qff a_{\dff n}\trf]}\dff(\dff A\dff)
\pff.
\end{equation}

\vspace{-12pt}\vspace{1.5pt}
In\sss particular,\oss the preimage under\sss $\hat{\psi}_{\dff n}$ of\sss
an $n$\dnsp-simplex such as in\sss the previous paragraph
can\sss be identified\sss with\sss the space 
of\dss operators $A\qff \in\qff \hat{\mathcal{F}}$\sss such\sss that\sss
$A\qff \bigl|\halfff\qff V_{\fff n}
\off =\off
F_{\fff n}$
and\qss (\ref{vn})\qss holds.\oss

By\sss restricting\sss such operators\sss to\sss the orthogonal\sss
complement\sss $H\dff \ominus\dff V_{\fff n}$ of\dss $V_{\fff n}$\sss in\sss $H$\sss
this space can\sss be identified\sss with\sss the space 
of\dss self-adjoint\dss Fredholm\sss operators\sss $B$ in\sss
$H\dff \ominus\dff V_{\fff n}$\sss  such\sss that\sss
$\sigma\trf(\trf B\trf)
\dff \cap\trf 
[\dff -\qff a_{\dff n}\fff,\qff a_{\dff n}\trf]
\off =\off
\varnothing$\nnsp,\oss
i.e.\qss with\sss the space\sss
$\hat{\mathcal{F}}^{\dff \inv}\trf[\dff -\qff a_{\dff n}\dff,\dff a_{\dff n}\trf]$\sss
with\sss $H\dff \ominus\dff V_{\fff n}$\sss in\sss the role of\sss $H$\nnsp.

By\trs Proposition\qss \ref{invertible-contractible}\qss the space\dss
$\hat{\mathcal{F}}^{\dff \inv}\trf[\dff -\qff a_{\dff n}\dff,\dff a_{\dff n}\trf]$\sss
is\dss contractible.\oss
On\sss the other hand,\oss
by\sss choosing isomorphisms\sss $H\dff \ominus\dff V\qff \ttoo\qff H$
continuously depending on $V$ for $V$ close\sss to a given $V_{\fff n}$
one can see\sss that\sss the map\sss 
$\hat{\mathcal{E}}_{\dff n}
\qff \ttoo\qff
\mathcal{E}\hat{\mathcal{O}}_{\dff n}$\sss
is\dss a\sss locally\dss trivial\sss bundle.\oss
Since its fibers are contractible and\sss its base\sss $\mathcal{E}\hat{\mathcal{O}}_{\dff n}$\sss
is\dss paracompact,\oss the map\sss 
$\hat{\mathcal{E}}_{\dff n}
\qff \ttoo\qff
\mathcal{E}\hat{\mathcal{O}}_{\dff n}$\sss
is\dss a\sss homotopy\sss equivalence by a classical\dss theorem of\trs
Dold\qss (see \cite{d},\oss Corollary\qss 3.2).\oss
It\sss remains\sss to check\sss that\sss 
$\hat{\mathcal{E}}$ and\sss $\mathcal{E}\hat{\mathcal{O}}$
are good enough\sss for\sss this property\sss to imply\sss the conclusion of\dss the\sss theorem.\oss

Since we are using discrete\sss topology for $\varepsilon$\nnsp,\oss
both categories\sss $\hat{\mathcal{E}}$ and\sss $\mathcal{E}\hat{\mathcal{O}}$\sss
have free units.\oss
Therefore\trs Lemma\qss \ref{free-d-categories}\qss implies\sss that\sss 
$\mathit{N}\dff \hat{\mathcal{E}}$ and\sss $\mathit{N}\dff \mathcal{E}\hat{\mathcal{O}}$\sss
have\sss free degeneracies.\oss
Now\trs Proposition\qss \ref{level-heq}\qss
implies\sss that\sss
$\num{\hat{\psi}}\dff \colon\dff
\num{\hat{\mathcal{E}}}
\qff \ttoo\qff 
\num{\mathcal{E}\hat{\mathcal{O}}}$\sss
is\dss a\sss homotopy\sss equivalence.\oss  \eproof

\myuppar{Remark.}
The\sss last\sss step of\dss the proof\dss can\sss be carried out\sss
slightly\sss differently.\oss
Since\sss the simplicial\sss spaces\sss
$\mathit{N}\dff \hat{\mathcal{E}}$\sss 
and\sss 
$\mathit{N}\dff \mathcal{E}\hat{\mathcal{O}}$\sss
have free degeneracies,\oss they
are also good and\sss proper.\oss
By a\dss theorem of\trs Segal\qss
(see\dss Section\qss \ref{simplicial-spaces})\qss
this implies\sss that\sss
$\num{\hat{\psi}}$\sss
is\dss a\sss homotopy\sss equivalence.\oss

\myuppar{Operator\dss models.}
An\qss \emph{operator\dss model}\oss is\dss defined\sss as a\sss pair\sss
$(\trf V\fff,\pff F\trf)$\nnsp,\oss
where\sss $V$\sss is\dss a finitely dimensional\sss subspace of\sss $H$\sss
and\sss
$F\dff \colon\dff V\qff \ttoo\qff V$\sss 
is\dss a self-adjoint\sss operator.\oss
Let\sss $\hat{\mathcal{O}}$\sss be\sss the set\sss of\dss operator\sss models\dss
equipped\sss with\sss the obvious\sss topology  
and ordered\sss by\sss the relation\dss $\leq$\nnsp,\oss where\vspace{1.5pt}
\[
\quad
(\trf V\fff,\pff F \trf)
\off \leq\off
(\trf V\fff'\fff,\pff F\fff' \trf)
\quad
\mbox{if}\fff\quad
V\qff \subset\qff V\fff'\dff,\quad 
F\off =\off F\fff'\qff \bigl|\halfff\qff V\qff,\quad  
\]

\vspace{-12pt}\vspace{1.5pt}
and\sss $\num{\lambda}\qff <\qff \num{\lambda'}$\sss
for every eigenvalue $\lambda$ of\dss $F$\sss
and\sss every\sss eigenvalue $\lambda'$ of\dss the restriction of\dss $F\fff'$\sss
to\sss the orthogonal\sss complement\sss $V\fff'\dff \ominus\dff V$\dnsp.\oss
If\trs $V\off =\off V\fff'$\dnsp,\oss
then\sss
$(\trf V\fff,\pff F\trf)
\off \leq\off
(\trf V\fff'\fff,\pff F\fff'\trf)$\sss
if\trs and\dss only\trs if\dss also\sss $F\fff'\off =\off F$\dnsp.\oss
There\dss is\dss an obvious forgetting functor\sss
$\hat{o}\dff \colon\dff 
\mathcal{E}\hat{\mathcal{O}}
\qff \ttoo\qff 
\hat{\mathcal{O}}$\nnsp.\oss

\myuppar{{\dnsp}$\mathcal{E}\hat{\mathcal{O}}$ and\sss $\hat{\mathcal{O}}$
as\sss topological\sss simplicial\sss complexes.}
Both\sss $\mathcal{E}\hat{\mathcal{O}}$ and\sss $\hat{\mathcal{O}}$
are partially ordered\sss topological\sss spaces.\oss
Clearly,\oss if\vspace{1.2pt}
\[
\quad
(\trf V\fff,\pff F\fff,\pff \varepsilon\trf)
\off <\off
(\trf V\fff'\fff,\pff F\fff'\fff,\pff \varepsilon'\trf)
\qff,
\]

\vspace{-12pt}\vspace{1.2pt}
then either\sss
$\dim V\off <\off \dim V\fff'$\dnsp,\oss or\sss
$\varepsilon\qff <\qff \varepsilon'$\nnsp.\oss
Similarly,\oss if\dss\vspace{1.2pt}
\[
\quad
(\trf V\fff,\pff F \trf)
\off <\off
(\trf V\fff'\fff,\pff F\fff' \trf)
\qff,
\]

\vspace{-12pt}\vspace{1.2pt}
then\sss $\dim\dff V\off <\off \dim\dff V\fff'$\dnsp.\oss
Therefore\sss
$\mathcal{E}\hat{\mathcal{O}}$ and\sss $\hat{\mathcal{O}}$
have free equalities by dimension\sss reasons.\oss

As explained\sss in\dss Section\qss \ref{topological-simplicial-complexes},\oss
the partial\sss orders allow\sss to consider\sss
$\mathcal{E}\hat{\mathcal{O}}$ and\sss $\hat{\mathcal{O}}$\sss
as\sss topological\sss simplicial\sss complexes\sss
and\sss to define\sss
the $\Delta$\dnsp-spaces\sss
$\Delta\dff \mathcal{E}\hat{\mathcal{O}}$ 
and\sss $\Delta\dff \hat{\mathcal{O}}$\dnsp,\oss
as also\sss the simplicial\sss spaces\sss
$\bm{\Delta}\dff \mathcal{E}\hat{\mathcal{O}}$ 
and\sss $\bm{\Delta}\dff \hat{\mathcal{O}}$\dnsp.\oss
The\sss latter are nothing else but\sss the nerves of\dss 
$\mathcal{E}\hat{\mathcal{O}}$ and\sss $\hat{\mathcal{O}}$\sss
considered as\sss topological\sss categories.\oss
Since\sss the partially\sss ordered spaces\sss
$\mathcal{E}\hat{\mathcal{O}}$ and\sss $\hat{\mathcal{O}}$\sss
have free equalities,\oss
Corollary\qss \ref{free-full-realizations}\qss implies\sss that\sss the geometric realizations
of\sss each of\dss them as a\sss topological\sss simplicial\sss space
and as a category are canonically\sss isomorphic.\oss
The map\sss
$\hat{o}\dff \colon\dff 
\mathcal{E}\hat{\mathcal{O}}
\qff \ttoo\qff 
\hat{\mathcal{O}}$\sss
can be considered\sss both as a functor and as an order-preserving map.\oss
As explained at\sss the end of\trs Section\qss \ref{topological-simplicial-complexes},\oss
the resulting\sss maps\sss
\[
\quad
\num{\hat{o}}\dff \colon\dff 
\num{\mathcal{E}\hat{\mathcal{O}}}\qff \ttoo\qff \num{\hat{\mathcal{O}}}
\quad
\mbox{and}\quad
\bbnum{\hat{o}}\dff \colon\dff 
\bbnum{\mathcal{E}\hat{\mathcal{O}}}\qff \ttoo\qff \bbnum{\hat{\mathcal{O}}}
\pff
\]

\vspace{-12pt}
are\sss the same.\oss
Hence\sss we can\sss treat\sss
$\mathcal{E}\hat{\mathcal{O}}$
and\sss
$\hat{\mathcal{O}}$\sss
as\sss topological\sss simplicial\sss complexes.\oss

\mypar{Theorem.}{to-models}
\emph{The map\sss
$\num{\hat{o}}\dff \colon\dff
\num{\mathcal{E}\hat{\mathcal{O}}}
\qff \ttoo\qff 
\num{\hat{\mathcal{O}}}$\sss
is\dss a\sss homotopy\sss equivalence.\oss}

\proof
Let\sss $\mu\off =\off \num{\hat{o}}\off =\off \bbnum{\hat{o}}$\nnsp.\oss
Considering\sss $\mathcal{E}\hat{\mathcal{O}}$
and\sss
$\hat{\mathcal{O}}$\sss
as\sss topological\sss simplicial\sss complexes,\oss
we will\sss use an\sss induction\sss by $n$ in order\sss to prove\sss that\sss the maps
\[
\quad
\mu^{\dff -\dff 1}\qff \bbnum{\ssk_{\dff n}\dff \hat{\mathcal{O}}}
\off \ttoo\qff
\bbnum{\ssk_{\dff n}\dff \hat{\mathcal{O}}}
\]

\vspace{-12pt}
induced\sss by $\mu$\sss are homotopy equivalences
and\sss then pass\sss to\sss the direct\sss limit.\oss
In order\sss to do\sss this,\oss we need a more explicit\sss description of\dss
simplices of\sss $\hat{\mathcal{O}}$\sss and\sss $\mathcal{E}\hat{\mathcal{O}}$\dnsp.\oss

An $n$\dnsp-simplex $\sigma$ of\dss the\sss topological\sss simplicial\sss complex\sss $\hat{\mathcal{O}}$\sss
can\sss be identified\sss with a strictly\sss increasing sequence of\sss elements of\sss
$\hat{\mathcal{O}}$\sss and\sss hence with a flag\vspace{1.25pt}
\begin{equation}
\label{flag}
\quad
V_{\dff 0}\off \subset\off
V_{\dff 1}\off \subset\off
\ldots\off \subset\off
V_{\fff n}
\end{equation}

\vspace{-12pt}\vspace{1.25pt}
of\dss finitely dimensional\sss subspaces of\dss $H$\sss
together\sss with\sss a self-adjoint\sss operator
$F\dff \colon\dff
V_{\fff n}\qff \ttoo\qff V_{\fff n}$\sss
such\sss that\sss 
$V_{\dff i}\off \neq\off V_{\dff i\dff +\dff 1}$\dss 
and\sss
$\num{\lambda}\qff <\qff \num{\lambda'}$\sss
for every eigenvalue $\lambda$ of\dss the restriction\sss 
$F\qff \bigl|\halfff\qff V_{\dff i}$\sss
and\sss every\sss eigenvalue $\lambda'$ of\dss 
the restriction of\dss $F$\sss to\sss
$V_{\dff i\dff +\dff 1}\qff \ominus\pff V_{\dff i}$\nsp,\oss
where\sss $0\qff \leq\qff i\qff \leq\qff n\qff -\qff 1$\nnsp.\oss

We already\sss used\sss in\sss the proof\dss of\trs Theorem\qss \ref{to-enhanced-models}\qss
the fact\sss that\sss 
an $l$\dnsp-simplex $\tau$\sss of\dss
$\mathcal{E}\hat{\mathcal{O}}$\sss
as a\sss topological\sss category\sss 
can\sss be identified\sss with a finitely\sss dimensional\sss subspace\sss 
$W_{\fff l}\qff \subset\qff H$\sss
together\sss with\sss a self-adjoint\sss operator
$F\dff \colon\dff
W_{\fff l}\qff \ttoo\qff W_{\fff l}$\sss
and\sss a sequence\sss
$a_{\dff 0}\qff \leq\qff a_{\dff 1}\qff \leq\qff \ldots\qff \leq\qff a_{\dff l}$\sss
of\dss positive numbers such\sss that\sss
$\sigma\dff(\trf F\trf)
\off \subset\off
(\dff -\qff a_{\dff l}\fff,\pff a_{\dff l}\trf)$\sss
and\sss the numbers\sss $-\qff a_{\dff i}\dff,\qff a_{\dff i}$\sss 
are not\sss eigenvalues of\sss $F$\nnsp.\oss
The $l$\dnsp-simplices $\tau$ of\sss 
$\mathcal{E}\hat{\mathcal{O}}$\sss
as a simplicial\sss complex\sss
correspond\sss to\sss non-degenerate $\tau$\nnsp,\oss
i.e.\qss to $\tau$ such\sss that\sss
$a_{\dff i}\qff <\qff a_{\dff i\dff +\dff 1}$\dss 
for every\sss $i\qff \leq\qff l\qff -\qff 1$\nnsp.\oss

Let\sss us\sss begin\sss the induction\sss by\sss considering\sss the map\vspace{1.25pt}
\begin{equation}
\label{zero-skeleton}
\quad
\mu^{\dff -\dff 1}\qff \bbnum{\ssk_{\trf 0}\dff \hat{\mathcal{O}}}
\off \ttoo\qff
\bbnum{\ssk_{\trf 0}\dff \hat{\mathcal{O}}}
\pff.
\end{equation}

\vspace{-12pt}\vspace{1.25pt}
The space $\bbnum{\ssk_{\trf 0}\dff \hat{\mathcal{O}}}$\sss is\dss
nothing else but\sss the space of\dss objects of\dss $\hat{\mathcal{O}}$\dnsp,\oss
i.e.\qss the space of\dss operator models.\oss
Let\sss $(\trf V_{\dff 0}\fff,\pff F_{\dff 0}\trf)$\sss be an operator model.\oss
The preimage\sss 
$\mu^{\dff -\dff 1}\qff (\trf V_{\dff 0}\fff,\pff F_{\dff 0}\trf)$\sss
is\dss equal\sss to\sss the geometric realization of\dss the\sss
topological\sss simplicial\sss complex\sss
having as $n$\dnsp-sim\-plices\sss the sim\-plices $\tau$ as above such\sss that\sss
$W_{\fff l}\off =\off V_{\dff 0}$\nsp,\qss $F\off =\off F_{\dff 0}$\nsp,\oss
and\sss 
$\sigma\dff(\trf F\trf)$\sss
is\dss contained\sss in\sss
$(\dff -\qff a_{\dff i}\fff,\pff a_{\dff i}\trf)$\sss
for every $i$\nnsp.\oss
The\sss last\sss condition\dss is\dss equivalent\sss to\sss
$a_{\dff i}\qff \in\qff (\trf u_{\dff 0}\fff,\qff \infty\trf)$\nnsp,\oss
where $u_{\dff 0}\off =\off \max\qff \num{\lambda}$ over\sss the eigenvalues
$\lambda$ of\sss $F_{\dff 0}$\nnsp.\oss
It\sss follows\sss that\sss
$\mu^{\dff -\dff 1}\qff (\trf V_{\dff 0}\fff,\pff F_{\dff 0}\trf)$\sss
can\sss be identified\sss with\sss the geometric realization of\dss
the\sss topological\sss simplicial\sss complex\sss having as $n$\dnsp-simplices\sss
increasing sequences\sss
$a_{\dff 0}\qff <\qff a_{\dff 1}\qff <\qff \ldots\qff <\qff a_{\dff n}$\sss
of\dss numbers\sss $>\qff u_{\dff 0}$\nsp.\oss
Clearly,\oss the geometric realization of\dss this complex\dss is\dss contractible.\oss
On\sss the other\sss hand,\oss the number\sss 
$u_{\dff 0}$ continuously
depends on\sss the operator model\sss
$(\trf V_{\dff 0}\fff,\pff F_{\dff 0}\trf)$\nnsp.\oss
This\sss implies\sss that\qss 
(\ref{zero-skeleton})\qss
is\dss the projection of\dss 
a\sss locally\sss trivial\sss bundle.\oss
Since\sss its\sss base\dss is\dss paracompact\sss
and\sss its\sss fibers are contractible,\oss 
(\ref{zero-skeleton})\qss is\dss a homotopy equivalence.\oss

Preparing\sss for\sss the step of\dss the induction,\oss
let\sss us\sss fix an $n$\dnsp-simplex $\sigma$ of\dss 
the simplicial\sss complex\sss $\hat{\mathcal{O}}$\sss as above.\oss
Let {\nsp}$\tau$ as above be an 
$l$\dnsp-simplex of\sss $\mathcal{E}\hat{\mathcal{O}}$
as a simplicial\sss complex.\oss
Then $\hat{o}\trf(\trf \tau\trf)$ is\sss
a simplex of\sss $\hat{\mathcal{O}}$
corresponding\sss to\sss the same operator\sss
$F\dff \colon\dff
W_{\fff l}\qff \ttoo\qff W_{\fff l}$\sss
and\sss the flag of\dss subspaces\vspace{3pt}
\[
\quad
W_{\dff i}
\off =\off
\image\dff P_{\dff [\dff -\qff a_{\dff i}\fff,\qff a_{\dff i}\trf]}\dff(\trf F\trf)
\]

\vspace{-12pt}\vspace{3pt}
with\sss the possible repetitions removed.\oss
Hence\sss $\hat{o}\trf(\trf \tau\trf)\off =\off \sigma$\sss
if\trs and\dss only\trs if\trs the corresponding operators $F$ are equal\sss
and\sss removing\sss the repetitions\sss
turns\sss the flag 
\[
\quad
\quad
W_{\dff 0}\off \subset\off
W_{\dff 1}\off \subset\off
\ldots\off \subset\off
W_{\fff l}
\]

\vspace{-12pt}
into\sss the flag\qss (\ref{flag}).\oss
The second condition\dss holds\sss if\dss and\sss only\trs if\dss
$W_{\dff i}\off =\off V_{\trf \theta\trf(\trf i\trf)}$
for some surjective non-decreasing map
$\theta\dff \colon\dff
[\dff l\trf]\qff \ttoo\qff [\halfff n\dff]$
and every $i\qff \in\qff [\dff l\trf]$\nnsp.\oss
This happens\dss if\dss and\dss only\trs if\dss
each $a_{\dff i}$
belongs\sss to a non-empty\sss interval\sss between some
eigenvalues of\dss $F$\sss
and each of\dss these intervals contains some $a_{\dff i}$\nsp.\oss 
In\sss more details,\oss
let\vspace{3pt}  
\[
\quad
u\trf(\trf k\trf)
\off =\off
\max\qff \num{\lambda}
\quad
\mbox{and}\quad
l\trf(\trf k\trf)
\off =\off
\min\qff \num{\lambda'}
\qff,
\]

\vspace{-12pt}\vspace{3pt}
where {\hnsp}$\lambda${\hnsp} runs over\sss the eigenvalues of\sss
$F\qff \bigl|\halfff\qff V_{\fff k}$\sss
and\sss $\lambda'$ over\sss the eigenvalues of\dss
$F\qff 
\bigl|\halfff\qff 
V_{\fff k\dff +\dff 1}\dff \ominus\dff V_{\fff k}$\nsp.\oss
If\sss $k\off =\off n$\nnsp,\oss
then\sss $u\trf(\trf k\trf)$\sss is\dss defined\sss in\sss the same way,\oss
but\sss $l\trf(\trf k\trf)\off =\off \infty$\sss by definition.\oss
Then\sss $u\trf(\trf k\trf)\pff <\pff l\trf(\trf k\trf)$ for every 
{\nsp}$k$\nnsp.\oss
Clearly,\qss
$\hat{o}\trf(\trf \tau\trf)\off =\off \sigma$\sss
if\dss and\dss only\dss if\dss the operators $F$ are equal\sss and\vspace{3pt}\vspace{0.3pt}
\begin{equation}
\label{intervals}
\quad
a_{\dff i}
\off \in\pff 
\bigl(\trf u\trf(\trf k\trf)\fff,\qff l\trf(\trf k\trf) \trf\bigr)
\dff,
\quad
\mbox{where}\quad
k
\off =\off
\theta\trf(\trf i\trf)
\pff
\end{equation}

\vspace{-12pt}\vspace{3pt}\vspace{0.3pt}
for some surjective non-decreasing map
$\theta\dff \colon\dff
[\dff l\trf]\qff \ttoo\qff [\halfff n\dff]$
and every\sss
$i\qff \in\qff [\dff l\trf]$\nnsp.\oss
Let\sss us define discrete simplicial\sss complex\sss
$S\dff(\trf \sigma\fff,\qff F\trf)$\sss as\sss the union\vspace{3pt}\vspace{0.3pt}
\[
\quad
S\dff(\trf \sigma\fff,\qff F\trf)
\off =\off
\bigcup\nolimits_{\off 0\qff \leq\qff k\qff \leq\qff n}\off
\bigl(\trf u\trf(\trf k\trf)\fff,\qff l\trf(\trf k\trf)\trf\bigr)
\]

\vspace{-12pt}\vspace{3pt}\vspace{0.3pt}
with all\sss finite subsets of\sss
$S\dff(\trf \sigma\fff,\qff F\trf)$\sss
being simplices.\oss
For\sss $a\qff \in\qff S\dff(\trf \sigma\fff,\qff F\trf)$\sss
let\sss\vspace{3pt}\vspace{0.3pt} 
\[
\quad
O\trf(\dff a\trf)
\off =\off
\bigl(\qff 
V_{\dff k\trf(\dff a\trf)}\fff,\pff F\qff \bigl|\halfff\qff V_{\fff k\trf(\dff a\trf)}\fff,\pff a
\qff\bigr)
\off \in\off
\mathcal{E}\hat{\mathcal{O}}
\pff,
\]

\vspace{-12pt}\vspace{3pt}\vspace{0.3pt}
where\sss $k\trf(\dff a\trf)$\sss
is\dss the unique\sss
$k\qff \in\qff [\halfff n\dff]$\sss
such\sss that\sss
$a\off \in\pff 
(\trf u\trf(\trf k\trf)\fff,\qff l\trf(\trf k\trf)\trf)$\nnsp.\oss
The rules\sss
$a\off \longmapsto\off O\trf(\dff a\trf)$\sss
and\sss
$a\off \longmapsto\off k\trf(\dff a\trf)$\sss
define simplicial\sss maps\sss
$O\dff \colon\dff
S\dff(\trf \sigma\fff,\qff F\trf)
\qff \ttoo\qff 
\mathcal{E}\hat{\mathcal{O}}$\sss
and\sss
$k\dff \colon\dff
S\dff(\trf \sigma\fff,\qff F\trf)
\qff \ttoo\qff
[\halfff n\dff]$\nnsp.\oss

The partial\sss order of\sss $\hat{\mathcal{O}}$\sss
leads\sss to an\sss isomorphism\sss 
$[\halfff n\dff]\qff \ttoo\qff \sigma$
such\sss that\sss the square\vspace{-3pt}
\[
\quad
\begin{tikzcd}[column sep=boom, row sep=boomm]
S\dff(\trf \sigma\fff,\qff F\trf)
\arrow[r, "\dis O"]
\arrow[d, "\dis k\dff"']
&
\mathcal{E}\hat{\mathcal{O}}
\arrow[d, "\dis \hat{o}"]
\\
\protect{[\halfff n\dff]}
\arrow[r]
&
\hat{\mathcal{O}}\dff.
\end{tikzcd}
\]

\vspace{-12pt}\vspace{-1.5pt}
is\dss commutative.\oss
By\sss passing\sss to geometric realizations 
we get\sss the commutative diagram\vspace{-3pt}
\[
\quad
\begin{tikzcd}[column sep=boom, row sep=boomm]
\protect{\num{S\dff(\trf \sigma\fff,\qff F\trf)}}
\arrow[r, "\dis \protect{\num{O}}"]
\arrow[d, "\dis \kappa\dff"']
&
\protect{\num{\mathcal{E}\hat{\mathcal{O}}}}
\arrow[d, "\dis \mu"]
\\
\Delta^n
\arrow[r]
&
\protect{\num{\hat{\mathcal{O}}}}\dff,
\end{tikzcd}
\]

\vspace{-12pt}\vspace{-1.5pt}
where $\kappa\off =\off \num{k}$\nnsp.\oss
Clearly,\oss if\sss $\tau$\sss is\dss simplex of\dss $\mathcal{E}\hat{\mathcal{O}}$\nnsp,\oss
then
$\hat{o}\trf(\trf \tau\trf)\off =\off \sigma$\sss
if\dss and\dss only\dss if\dss
$\tau\off =\off O\trf(\trf \alpha\trf)$\sss
for a simplex $\alpha$ of\dss
$S\dff(\trf \sigma\fff,\qff F\trf)$
such\sss that\sss
$k\trf(\trf \alpha\trf)\off =\off [\halfff n\dff]$\nnsp,\oss
and\dss if\trs such $\alpha$ exists,\oss
it\dss is\dss unique.\fff\footnotemark\oss
It\sss follows\sss that\sss $\num{O}$\sss induces 
a canonical\sss homeomorphism\vspace{1.5pt}
\[
\quad
h_{\dff \sigma}\qff \colon\qff
\kappa^{\dff -\dff 1}\qff \bigl(\qff \inte \Delta^n\qff\bigr)
\off \ttoo\off
\mu^{\dff -\dff 1}\qff \bigl(\qff \inte \num{\sigma}\qff\bigr)
\pff.
\]

\vspace{-10.5pt}
Moreover,\oss if\dss $U_{\dff \sigma}$\sss
is\sss the simplex $\num{\sigma}$\sss with\sss its center removed,\oss
then\sss there\dss is\dss a deformation\sss retraction of\dss
$\mu^{\dff -\dff 1}\trf(\trf U_{\dff \sigma}\trf)$\sss
onto\sss
$\mu^{\dff -\dff 1}\trf(\trf \partial\trf \num{\sigma}\trf)$\sss
covering\sss the radial\sss deformation\sss retraction of\dss
$U_{\dff \sigma}$\sss onto\sss the  boundary\sss $\partial\trf \num{\sigma}$\nnsp.\oss
The simplicial\sss map\sss $k$\sss and\sss its geometric realization
$\kappa$ are standard\sss maps
depending\sss only\sss on\sss the intervals in\qss (\ref{intervals}).\oss
In\sss more details,\pss $S\dff(\trf \sigma\fff,\qff F\trf)$\sss is\dss the join
of\sss $n\qff +\qff 1$ infinite\qss ``simplices''\qss having intervals\sss
$(\trf u\trf(\trf k\trf)\fff,\qff l\trf(\trf k\trf) \trf)$ as\sss their sets of\dss vertices\qss
(and\sss finite sets of\dss vertices as simplices),\oss
and $k$\sss is\sss the canonical\sss map of\dss this join\sss to $[\halfff n\dff]$\nnsp.\oss
This shows,\oss in\sss particular,\oss that\sss
$\kappa\off =\off \num{k}$\sss
is\dss a\sss locally\sss trivial\sss bundle with contractible fibers
over\sss $\inte \Delta^n$\dnsp.\oss

\footnotetext{It\dss is\dss worth\sss 
to point\sss out\sss that\sss there are
simplices $\tau$ of\dss $\mathcal{E}\hat{\mathcal{O}}$\sss
such\sss that\sss
$\hat{o}\trf(\trf \tau\trf)\off \subset\off \sigma$\nnsp,\oss
but\sss $\tau$\sss is\dss not\sss in\sss the image of\sss $O$\nnsp.\oss
The reason\dss is\dss that\sss removing\sss the $k${\nnsp}th vertex from $\sigma$
adds\sss the interval\sss $[\trf l\trf(\trf k\trf)\fff,\qff u\trf(\trf k\trf)\trf]$\sss
to\sss the allowed\sss values of\dss controlling\sss parameters $a_{\dff i}$\nsp.\oss} 

Let\sss us\sss allow $\sigma$ vary among $n$\dnsp-simplices of\dss 
the\sss topological\sss simplicial\sss complex\sss $\hat{\mathcal{O}}$\dnsp.\oss
Clearly,\oss the union of\dss the interiors $\inte \sigma$\sss is\dss
equal\sss to\sss 
$\bbnum{\ssk_{\dff n}\dff \hat{\mathcal{O}}}
\off \smallsetminus\off
\bbnum{\ssk_{\dff n\dff -\dff 1}\dff \hat{\mathcal{O}}}$\nnsp.\oss 
Also,\oss 
the intervals in\qss (\ref{intervals})\qss 
continuously depend on\sss $\sigma$\nnsp,\oss
as also homeomorphisms $h_{\dff \sigma}$\nsp.\oss
It\sss follows\sss that\sss the map\vspace{1.5pt}
\[
\quad
\mu^{\dff -\dff 1}\qff \bigl(\qff 
\bbnum{\ssk_{\dff n}\dff \hat{\mathcal{O}}}
\off \smallsetminus\off
\bbnum{\ssk_{\dff n\dff -\dff 1}\dff \hat{\mathcal{O}}}
\qff\bigr)
\off \ttoo\off
\bbnum{\ssk_{\dff n}\dff \hat{\mathcal{O}}}
\off \smallsetminus\off
\bbnum{\ssk_{\dff n\dff -\dff 1}\dff \hat{\mathcal{O}}}
\]

\vspace{-10.5pt}
induced\sss by\sss
$\mu$\sss is\dss a\sss locally\sss trivial\sss bundle
with contractible fibers and\sss hence\dss is\dss a homotopy equivalence.\oss
Clearly,\oss the union\sss\vspace{1.5pt}
\[
\quad
\mathcal{U}_{\dff n\dff -\dff 1}
\off =\off
\cup_{\dff \sigma}\qff U_{\dff \sigma}
\]

\vspace{-10.5pt}
is\dss an open\sss neighborhood of\dss
$\bbnum{\ssk_{\dff n\dff -\dff 1}\dff \hat{\mathcal{O}}}$\sss
in\sss
$\bbnum{\ssk_{\dff n}\dff \hat{\mathcal{O}}}$\sss
and\vspace{1.5pt}
\[
\quad
\mu^{\dff -\dff 1}\qff 
\bigl(\qff 
\bbnum{\ssk_{\dff n\dff -\dff 1}\dff \hat{\mathcal{O}}}
\qff\bigr)
\]

\vspace{-10.5pt}
is\dss a deformation\sss retract\sss of\dss 
$\mu^{\dff -\dff 1}\qff (\trf \mathcal{U}_{\dff n\dff -\dff 1}\trf)$\nnsp.\oss
Suppose now\sss that\sss the map\vspace{1.5pt}
\[
\quad
\mu^{\dff -\dff 1}\qff \bigl(\qff 
\bbnum{\ssk_{\dff n\dff -\dff 1}\dff \hat{\mathcal{O}}}
\qff\bigr)
\qff \ttoo\qff
\bbnum{\ssk_{\dff n\dff -\dff 1}\dff \hat{\mathcal{O}}}
\]

\vspace{-10.5pt}
induced\sss by\sss $\mu$\sss
is\dss a\sss homotopy equivalence.\oss
Then\sss the map\dss 
$\mu^{\dff -\dff 1}\qff (\qff 
\mathcal{U}_{\dff n\dff -\dff 1}
\qff)
\qff \ttoo\qff
\mathcal{U}_{\dff n\dff -\dff 1}$\dss
induced\sss by\sss $\mu$\sss
is\dss also a\sss homotopy equivalence.\oss
We see\sss that\sss the map $\mu$\sss is\dss
a homotopy equivalence over\sss the open subsets\sss
$\bbnum{\ssk_{\dff n}\dff \hat{\mathcal{O}}}
\off \smallsetminus\off
\bbnum{\ssk_{\dff n\dff -\dff 1}\dff \hat{\mathcal{O}}}$\sss
and\dss $\mathcal{U}_{\dff n\dff -\dff 1}$\sss
of\trs
$\bbnum{\ssk_{\dff n}\dff \hat{\mathcal{O}}}$\nnsp.\oss
Also,\pss $\mu$\sss is\dss a\sss locally\sss trivial\sss bundle
with contractible fibers over\sss the intersection of\dss these subsets,\oss
and\sss hence\dss is\dss a homotopy equivalence
over\sss this\sss intersection.\oss
The glueing\sss theorem\sss for\sss homotopy equivalences\sss
implies\sss that\sss $\mu$\sss is\dss a homotopy equivalence over\sss
$\bbnum{\ssk_{\dff n}\dff \hat{\mathcal{O}}}$\nnsp,\oss
i.e.\qss that\sss the map\vspace{1.5pt}
\[
\quad
\mu^{\dff -\dff 1}\qff \bbnum{\ssk_{\dff n}\dff \hat{\mathcal{O}}}
\off \ttoo\qff
\bbnum{\ssk_{\dff n}\dff \hat{\mathcal{O}}}
\]

\vspace{-12pt}\vspace{1.5pt}
induced\sss by $\mu$\sss is\dss a\sss homotopy equivalence.\oss
See,\oss for example,\oss tom\dss Dieck\qss \cite{td1},\oss
Theorem\qss 1\qss for a much more general\sss result.\oss
This completes\sss the induction step.\oss
Passing\sss to\sss the direct\sss limit\sss completes\sss the proof.\oss
Cf.\qss the proof\dss of\qss Proposition\qss \ref{level-heq}.\oss \eproof

\myuppar{Remarks.}
Clearly,\oss the main\sss part\sss of\dss the proof\dss of\qss Theorem\qss \ref{to-models}\qss
is\dss the induction step.\oss
It\dss follows\sss approximately\sss the same outline as\sss the induction step\sss in\dss 
Atiyah--Singer\dss proof\qss \cite{as}\qss of\dss their\dss Proposition\qss 3.5.\oss
At\sss the same\sss time\sss it\dss is\dss similar\sss to\sss the lemma about\sss
quasi-fibrations used\sss by\trs Quillen\qss \cite{q}\qss
in\sss his proof\dss of\trs his\trs Theorem\qss B.\oss
In\sss fact,\oss the proof\dss of\qss Theorem\qss \ref{to-models}\qss
shows\sss that\sss the map\sss 
$\mu\off =\off \num{\hat{o}}\off =\off \bbnum{\hat{o}}$\dss 
is\dss a quasi-fibration.\oss

\myuppar{Subspace models.}
Operator\sss models\sss $(\trf V,\qff F\trf)$\sss
still\sss contain extra\sss information\sss not\sss affecting\sss the homotopy\sss type
of\dss the classifying space\sss $\num{\hat{\mathcal{O}}}$\nnsp.\oss
Namely,\pss 
the space of\sss self-adjoint\sss operators\sss 
$F\dff \colon\dff V\qff \ttoo\qff V$\sss
is\sss trivially\sss contractible.\oss
It\sss follows\sss that\sss the space of\sss operator\sss
models,\oss which\dss is\dss the space of\sss objects of\dss the\sss
topological\sss category\sss $\hat{\mathcal{O}}$\nnsp,\oss 
is\dss homotopy\sss equivalent\sss to\sss the space of\dss
finitely dimensional\sss subspaces of\dss $H$\nnsp.\oss
But\sss the space of\dss morphisms\dss is\qss \emph{not}\pss
homotopy\sss equivalent\sss to\sss the space of\dss pairs\sss $V\fff,\qff V\fff'$\sss
of\dss finitely dimensional\sss subspaces such\sss that\sss 
$V\qff \subset\qff V\fff'$\dnsp.\oss

Indeed,\oss a\sss morphism $f$ of\sss $\hat{\mathcal{O}}$ 
corresponds\sss to an\sss inequality\sss
$(\trf V,\qff F\trf)
\off \leq\off
(\trf V\fff',\qff F\fff'\trf)$\nnsp.\oss
By\sss the definition,\oss
there exists $\varepsilon\qff >\qff 0$\sss
such\sss that\sss every\sss eigenvector of\sss $F\fff'$ 
with\sss the eigenvalue\sss $\lambda$ such\sss that\sss 
$\num{\lambda}\qff <\qff \varepsilon$\sss 
is,\oss in\sss fact,\oss
an eigenvector of\sss $F$\dnsp,\oss
and every other\dss is\dss not.\oss
Therefore admissible deformations of\dss the morphism $f$ do not\sss allow\sss to\sss
turn an eigenvalue of\dss $F\fff'$\sss to\sss the right\sss of\sss $\sigma\dff(\trf F\trf)$
into an eigenvalue\sss to\sss the\sss left\sss of\sss $\sigma\dff(\trf F\trf)$\nnsp.\oss

In fact,\qss $f$ defines a splitting of\dss the orthogonal\sss
complement\sss $V\fff'\dff \ominus\dff V$\sss 
into\sss a direct\sss sum\vspace{1.5pt}
\[
\quad
V\fff'\dff \ominus\dff V
\off =\off
V\fff'_{\fff -}
\qff \oplus\qff
V\fff'_{\fff +}
\qff
\]

\vspace{-12pt}\vspace{1.5pt}
such\sss that\sss $F\fff'$\sss is\dss negative on\sss $V\fff'_{\fff -}$
and\sss positive on\sss $V\fff'_{\fff +}$\nsp.\oss
Either of\dss these subspaces could\sss be $0$\nsp,\oss
and\dss if\trs both are $0$\nnsp,\oss
then\sss $V\fff'\off =\off V$\sss and\sss $f$\sss is\dss an\sss identity morphism.\oss
The space of\dss splittings\sss of\dss the complement\sss 
$V\fff'\dff \ominus\dff V$\sss into\sss the direct\sss sum of\dss two subspaces\dss is\dss 
homeomorphic\sss to\sss 
the union of\dss several\sss classical\dss Grassmann\dss manifolds
and\sss hence\dss is\dss not\sss contractible.\oss

The above observations\sss partially\sss motivate\sss the following definitions.\oss
A\qss \emph{subspace\sss model}\pss is\dss simply a\sss finitely\sss dimensional\sss
subspace\sss $V$\sss of\dss $H$\nnsp.\oss
The set\sss of\dss subspace models\dss is\dss equipped\sss with\sss the usual\sss
topology on\sss the set\sss of\dss subspaces.\oss
A\qss \emph{morphism\sss of\dss subspace models}\dss $V\qff \ttoo\qff V\fff'$\sss 
is\dss an ordered\sss pair\sss $U_{\dff -}\dff,\pff U_{\dff +}$\dss
of\dss subspaces of\dss $V\fff'$\sss defining an orthogonal\sss decomposition\vspace{1.375pt}
\[
\quad
V\fff'
\off =\off
U_{\dff -}\qff \oplus\qff
V\qff \oplus\qff U_{\dff +}
\qff.
\]

\vspace{-12pt}\vspace{1.375pt}
The subspaces\sss $U_{\dff -}$ and\sss $U_{\dff +}$\sss are called\sss the\qss
\emph{negative}\qss and\qss \emph{positive parts}\oss of\dss the morphism\sss in question.\oss
The composition of\dss morphisms\dss is\dss defined\sss by\sss taking\sss the sum of\dss
negative parts and\sss the sum of\dss positive parts\sss to get,\oss
respectively,\oss the negative and\sss
the positive parts of\dss the composition.\oss 
The\sss topology on\sss the set\sss of\dss
morphisms\dss is\dss defined\sss in\sss the obvious manner.\oss
This defines a\sss topological\sss category\sss $\hat{\mathcal{S}}$\sss
having\sss subspace models as\sss objects.\oss

Let\sss us\sss assign\sss to an operator\sss model\sss 
$(\trf V,\qff F\trf)$\sss 
the subspace model\sss
$V$\dnsp,\oss
and\sss to a morphism
\[
\quad
(\trf V,\qff F\trf)
\qff \ttoo\qff
(\trf V\fff',\qff F\trf)
\qff
\]

\vspace{-12pt}
of\sss operator models\sss the morphism\sss
of\dss subspace models corresponding\sss 
to\sss the decomposition\vspace{1.5pt}
\[
\quad
V\fff'
\off =\off
U_{\dff -}\qff \oplus\qff
V\qff \oplus\qff U_{\dff +}
\pff,
\quad
\]

\vspace{-19pt}
where\vspace{-9pt}
\[
\quad
U_{\dff -}
\off =\off
\image\dff P_{\dff (\trf -\qff \infty\fff,\qff 0\trf)}\dff
\left(\trf F\fff'\qff \bigl|\halfff\qff V\fff'\dff \ominus\dff V\trf\right)
\pff,
\]

\vspace{-36pt}
\[
\quad
U_{\dff +}
\off =\off
\image\dff P_{\dff (\trf 0\fff,\qff \infty\trf)}\dff
\left(\trf F\fff'\qff \bigl|\halfff\qff V\fff'\dff \ominus\dff V\trf\right)
\qff,
\]

\vspace{-12pt}\vspace{1.5pt}
i.e.\qss having\sss
$U_{\dff -}$\sss
and\dss
$U_{\dff +}$\sss
as\sss its\sss negative and\sss positive parts respectively.\oss
Clearly,\oss these rules define\sss a\sss forgetting\sss functor\sss
$\hat{\omega}\dff \colon\dff
\hat{\mathcal{O}}
\qff \ttoo\qff
\hat{\mathcal{S}}$\nnsp.\oss

\mypar{Lemma.}{s-free-degeneracies}
\emph{The nerve of\dss $\hat{\mathcal{S}}$\sss
is\dss a simplicial\sss space with\sss free degeneracies.\oss}

\proof
While\sss the\sss category\sss $\hat{\mathcal{S}}$\sss
is\dss not\sss associated\sss with a\sss partial\sss order,\oss
it\sss has free units\sss by\sss the same reasons
as,\oss for example,\oss the\sss topological\sss category\sss $\hat{\mathcal{O}}$\dnsp.\oss
Indeed,\oss if\dss 
$f\dff \colon\dff
V\qff \ttoo\qff V\fff'$\sss 
is\dss a morphism of\sss $\hat{\mathcal{S}}$\sss
corresponding\sss to\sss an orthogonal\sss decomposition\sss
$V\fff'
\off =\off
U_{\dff -}\qff \oplus\qff
V\qff \oplus\qff U_{\dff +}$\nsp,\oss
then\sss $f$\sss is\dss an\sss identity\sss morphism\dss
if\trs and\dss only\trs if\dss
$\dim\dff V\fff'\off =\off \dim\dff V$\dnsp,\oss
and\sss the dimension\dss is\dss a continuous function on\sss the space of\dss objects.\oss
By\trs Lemma\qss \ref{free-d-categories}\qss the nerve of\sss $\hat{\mathcal{S}}$\sss
is\dss a simplicial\sss space with\sss free degeneracies.\oss  \eproof

\mypar{Theorem.}{forgetting-operators}
\emph{The map\sss
$\num{\hat{\omega}}\dff \colon\dff
\num{\hat{\mathcal{O}}}
\qff \ttoo\qff 
\num{\hat{\mathcal{S}}}$\sss
is\dss a\sss homotopy\sss equivalence.\oss}

\proof
Clearly,\oss the space of\sss $n$\dnsp-simplices of\sss $\hat{\mathcal{O}}$\sss
is\dss a\sss locally\sss trivial\sss bundle over\sss the space of\sss
$n$\dnsp-simplices of\sss $\hat{\mathcal{S}}$\nnsp.\oss
Since\sss the space of\dss self-adjoint\sss operators 
a finitely dimensional\dss Hilbert\dss space with eigenvalues in\sss
a given\sss interval\dss is\dss contractible,\oss
its fibers are contractible.\oss
Hence\sss the\sss theorem of\qss Dold\qss \cite{d}\qss used\sss in\sss the proof\dss of\trs
Theorem\qss \ref{to-enhanced-models}\qss implies\sss that\sss the projection\dss
is\dss a homotopy equivalence.\oss
The nerves of\sss $\hat{\mathcal{O}}$\sss and\sss  $\hat{\mathcal{S}}$\sss 
have free degeneracies
by\trs Lemmas\qss \ref{free-orders}\qss
and\qss \ref{s-free-degeneracies}\qss respectively.\oss
Therefore\sss it\dss remains\sss to apply\trs
Proposition\qss \ref{level-heq}.\oss  \eproof\vspace{0.125pt}

\myuppar{Vector space models.}
A\qss \emph{vector space\sss model}\pss is\dss a\sss finitely dimensional\sss
vector space\sss $V$\sss together\sss with 
a\sss Hermitian\sss scalar\sss product\sss on $V$\dnsp,\oss
i.e.\qss a finitely\sss dimensional\dss Hilbert\sss space.\oss
The set\sss of\dss vector space models\qss
(as usual,\oss we need\sss to\sss limit\sss ourselves 
by\sss vector space models belonging\sss
to a fixed\sss universal\sss set\halfff)\qss
is\dss considered as a discrete set.\oss
A\qss \emph{morphism\sss of\dss vector\sss space models}\dss $V\qff \ttoo\qff V\fff'$\sss 
is\dss a\sss triple 
$(\trf U_{\dff -}\dff,\qff U_{\dff +}\dff,\qff f\qff)$\nnsp,\oss
where\sss
$f\dff \colon\dff
V\qff \ttoo\qff V\fff'$\sss
is\dss an\sss isometric embedding\sss and\dss
$U_{\dff -}\dff,\pff U_{\dff +}$\dss
are subspaces of\dss $V\fff'$
defining an orthogonal\sss decomposition\vspace{3pt}
\[
\quad
V\fff'
\off =\off
U_{\dff -}\qff \oplus\qff
f\dff(\trf V\trf)\qff \oplus\qff U_{\dff +}
\qff.
\]

\vspace{-12pt}\vspace{3pt}
The composition of\dss morphisms\dss is\dss defined\sss in\sss the same
way as\sss for subspace models,\oss
and\sss the\sss topology on\sss the set\sss of\dss
morphisms\dss is\dss defined\sss in\sss the obvious manner.\oss
This defines a\sss topological\sss category\sss $Q$\sss
having\sss vector space models as\sss objects.\oss
In contract\sss with\sss the\sss topology on\sss the set\sss of\dss objects,\oss
the\sss topology on\sss the set\sss of\dss morphisms\dss is\dss not\sss discrete.\oss
This definition\dss is\dss a slightly\sss rephrased adaptation
of\qss Quillen's\dss $Q$\dnsp-construction\qss \cite{q}\qss
due\sss to\dss Segal\qss \cite{s4}.\oss\vspace{0.125pt}

\myuppar{The intermediate category\sss $Q/\fff H$\nnsp.}
The obvious functor\sss 
$\hat{\mathcal{S}}\qff \ttoo\qff Q$\nnsp,\oss
which assigns\sss to a subspace\sss $V\qff \subset\qff H$\sss
the space $V$ considered as an object\sss of\sss $Q$\nnsp,\oss
is\qss \emph{not\sss continuous}\pss because\sss the space of\dss objects of\sss $Q$\sss
is\dss discrete.\oss
But\trs Segal\qss \cite{s4}\qss defined an\sss intermediate category\sss
$Q/\fff H$\sss
such\sss that\sss there are continuous functors\sss
$\hat{\mathcal{S}}\off \longleftarrow\off Q/\fff H\qff \ttoo\qff Q$\nnsp.\oss
The objects of\dss $Q/\fff H$\sss are pairs\sss 
$(\trf V\fff,\qff h\trf)$\sss 
such\sss that\sss
$V$\sss is\dss an object\sss of\dss $Q$\sss 
and\sss $h$\sss is\dss a\sss isometric embedding\sss $V\qff \ttoo\qff H$\nnsp.\oss
A morphism\sss
$(\trf V\fff,\qss h\trf)
\qff \ttoo\qff
(\trf V\fff'\fff,\qss h\fff'\trf)$\sss
is\dss defined as morphism\sss
$(\trf U_{\dff -}\dff,\qff U_{\dff +}\dff,\qff f\qff)$
of\dss vector space models\sss
$V\qff \ttoo\qff V\fff'$\sss such\sss that\sss 
$h
\off =\off
h\fff'\dff \circ\dff f$\nnsp.\oss
The composition\dss is\dss defined\sss in\sss the obvious manner.\oss

Let\sss us assign\sss to an object\sss $(\trf V\fff,\qff h\trf)$\sss
the image\sss $h\dff(\trf V\trf)$\nnsp,\oss and\sss to a morphism\sss
$(\trf V\fff,\qff h\trf)
\qff \ttoo\qff
(\trf V\fff'\fff,\qff h\fff'\trf)$\sss
defined\sss by\sss
$(\trf U_{\dff -}\dff,\qff U_{\dff +}\dff,\qss f\qff)$\sss
the morphism of\sss $\hat{\mathcal{S}}$\sss 
corresponding\sss to\sss the decomposition\vspace{3pt}
\[
\quad
h\fff'\dff(\trf V\fff'\trf)
\off =\off
h\fff'\dff(\trf U_{\dff -}\trf)\qff \oplus\qff
h\dff(\trf V\trf)\qff \oplus\qff h\fff'\dff(\trf U_{\dff +}\trf)
\qff.
\]

\vspace{-12pt}\vspace{3pt}
These rules define\sss the promised\sss functor\sss 
$i\dff \colon\dff Q/\fff H\qff \ttoo\qff \hat{\mathcal{S}}$\dnsp.\oss
And\dss ignoring\sss the\sss isometric embeddings defines 
a forgetting\sss functor\sss
$Q/\fff H\qff \ttoo\qff Q$\nnsp.\oss

\mypar{Theorem.}{intermediary}
\emph{The maps\sss
$\num{Q}
\off \longleftarrow\off 
\num{Q/\fff H}
\qff \ttoo\qff
\num{\hat{\mathcal{S}}}$\sss
induced\sss by\sss these functors are homotopy equivalences.\oss}

\proof
Let\sss us consider\sss
$\num{Q/\fff H}
\qff \ttoo\qff
\num{Q}$\nnsp.\oss
This part\sss of\dss proof\dss is\dss similar\sss to\sss the corresponding\sss part\sss
of\dss the proof\dss of\qss Theorem\qss \ref{vect}.\oss
The space of\dss objects of\sss $Q$\sss is\dss discrete,\oss
and\sss the space of\dss objects of\sss $Q/\fff H$\sss is\dss
the disjoint\sss union over\sss the objects $V$ of\sss $Q$ of\sss
spaces of\dss isometric embeddings $V\qff \ttoo\qff H$\nnsp.\oss
As\dss is\dss well\sss known,\oss
for a finitely dimensional\sss vector space $V$ with a scalar product\sss
the space of\dss isometric embeddings\sss $V\qff \ttoo\qff H$\sss
is\dss contractible.\oss
It\dss follows\sss that\sss $Q/\fff H\qff \ttoo\qff Q$\sss 
induces a homotopy equivalence on objects.\oss
The space of\dss morphisms of\sss $Q/\fff H$\sss is\dss 
a\sss locally\sss trivial\dss bundle over\sss the space of\dss
morphisms of\dss $Q$\nnsp,\oss
with\sss the fiber over a morphism\sss
$V\qff \ttoo\qff V\fff'$\sss of\dss $Q$\sss being\sss the space 
of\dss isometric embeddings
$V\fff'\qff \ttoo\qff H$\nnsp.\oss
Since\sss the spaces of\dss isometric embeddings
$V\fff'\qff \ttoo\qff H$\sss are contractible,\oss
$Q/\fff H\qff \ttoo\qff Q$\sss 
induces a homotopy equivalence on\sss morphisms also.\oss
The same argument\sss applies\sss to $n$\dnsp-simplices.\oss
The categories $Q$ and\sss $Q/\fff H$\sss do not\sss have free units,\oss
but\sss they are good and\sss proper\sss by\sss the same reasons
as\sss the categories\sss $\vect$ and\sss $\vect\fff/\fff H$\sss
from\dss Section\qss \ref{segal-example}.\oss
See\sss the proof\dss of\qss Theorem\qss \ref{vect}.\oss
Therefore\sss the same\sss theorem of\qss Segal\dss implies\sss that\sss
$\num{Q/\fff H\halfff}
\qff \ttoo\qff
\num{Q}$\sss
is\dss a\sss homotopy equivalence.\oss  

Let\sss us consider\sss the map
$\iota
\off =\off
\num{i}\dff \colon\dff
\num{Q/\fff H}
\qff \ttoo\qff
\num{\hat{\mathcal{S}}}$\nnsp.\oss
This part\sss of\dss the proof\dss follows\sss the same outline as\sss the proof\dss
of\trs Theorem\qss \ref{to-models}.\oss
Let\sss us\sss begin\sss with some preliminary\sss remarks about\sss $\hat{\mathcal{S}}$\nsp\dnsp.\oss

By\trs Lemma\qss \ref{s-free-degeneracies}\qss the simplicial\sss space\sss
$\hat{\mathcal{S}}$\sss has\sss free degeneracies.\oss
The arguments in\sss the proof\dss of\trs Lemma\qss \ref{s-free-degeneracies}\qss
also imply\sss that\sss the composition of\dss non-identity\sss morphisms 
of\dss the category\sss
$\hat{\mathcal{S}}$\sss cannot\sss be an\sss identity\sss morphism.\oss
It\sss follows\sss that $\hat{\mathcal{S}}$\sss has non-degenerate core.\oss
Now\trs Lemma\qss \ref{ndc-spaces}\qss implies\sss that\sss the canonical\sss map\sss
$\bm{\Delta}\dff \core \hat{\mathcal{S}}
\qff \ttoo\qff
\hat{\mathcal{S}}$\dss
is\dss an\sss isomorphism and\sss the canonical\sss map\sss
$\num{\core \hat{\mathcal{S}}}_{\dff \Delta}
\qff \ttoo\qff 
\num{\hat{\mathcal{S}}}$\sss
is\dss a\sss homeomorphism.\oss
Therefore we can use\sss
$\num{\core \hat{\mathcal{S}}}_{\dff \Delta}$\sss
instead\sss of\dss $\num{\hat{\mathcal{S}}}$\sss
and\sss treat\sss
$Q/\fff H
\qff \ttoo\qff
\hat{\mathcal{S}}$\sss
as a simplicial\sss map\sss
$Q/\fff H
\qff \ttoo\qff
\bm{\Delta}\dff \core \hat{\mathcal{S}}$\nsp\dnsp.\oss

Similarly\sss to\sss the proof\dss of\trs Theorem\qss \ref{to-models},\oss
we are going\sss to use an\sss induction by\sss the skeletons of\dss the $\Delta$\dnsp-space
$\core \hat{\mathcal{S}}$\nsp\dnsp.\oss
In order\sss to simplify\sss notations,\oss we will\sss treat\sss $\hat{\mathcal{S}}$
as a $\Delta$\dnsp-space and abbreviate\sss $\core \hat{\mathcal{S}}$\sss
to $\hat{\mathcal{S}}$\nsp\dnsp.\oss
Let\sss us\sss begin\sss the induction\sss by considering\sss the map\vspace{1.5pt}
\begin{equation}
\label{zero-skeleton-q}
\quad
\iota^{\dff -\dff 1}\pff \num{\ssk_{\trf 0}\dff \hat{\mathcal{S}}}
\off \ttoo\qff
\num{\ssk_{\trf 0}\dff \hat{\mathcal{S}}}
\pff.
\end{equation}

\vspace{-12pt}\vspace{1.5pt}
The space\sss $\num{\ssk_{\trf 0}\dff \hat{\mathcal{S}}}$\sss 
is\dss equal\sss to\sss the space of\dss objects of\sss $\hat{\mathcal{S}}$\dnsp,\oss
i.e.\qss the space\sss $G\trf(\trf H\trf)$ 
from\dss Section\qss \ref{segal-example}.\oss
If\sss $V\qff \in\pff G\trf(\trf H\trf)$\nnsp,\oss
then\sss $\iota^{\dff -\dff 1}\trf(\trf V\trf)$\sss
is\dss the preimage of\dss $V$\sss under\sss the map\sss
$\num{i}\dff \colon\dff \num{\vect}\qff \ttoo\qff G\trf(\trf H\trf)$\sss
from\dss Section\qss \ref{segal-example}.\oss
As we saw\sss in\sss the proof\dss of\trs Theorem\qss \ref{vect},\oss
this preimage\dss is\dss contractible.\oss
Hence\sss the map\qss (\ref{zero-skeleton-q})\qss is\dss a homotopy equivalence.\oss

Next,\oss let\sss us\sss fix a non-degenerate 
$n$\dnsp-simplex $\sigma$ of\sss
$\hat{\mathcal{S}}$\nsp\dnsp.\oss
It\sss has\sss the form\sss\vspace{1.5pt}
\[
\quad
W_{\dff 0}\qff \ttoo\qff
W_{\dff 1}\qff \ttoo\qff
\ldots\qff \ttoo\qff
W_{\dff n}
\pff,
\]

\vspace{-12pt}\vspace{1.5pt}
where\sss $W_{\fff i}$\sss are\sss objects of\sss $\hat{\mathcal{S}}$\sss 
such\sss that\sss $W_{\fff i}\off \neq\off W_{\fff i\dff -\dff 1}$\sss
for every\sss $i\off =\off 1\fff,\qff 2\fff,\qff \ldots\fff,\qff n$\nnsp,\oss
and arrows are morphisms.\oss
The simplex $\sigma$ can\sss be considered as a subcategory of\sss 
$\hat{\mathcal{S}}$\nsp\dnsp.\oss
Let\sss $\vect\downarrow \sigma$\sss be\sss the preimage of\dss this 
subcategory\sss under\sss the functor\sss 
$i\dff \colon\dff Q/\fff H\qff \ttoo\qff \hat{\mathcal{S}}$\dnsp.\oss
Then\sss the set\sss of\dss objects of\sss
$\vect\downarrow \sigma$\sss 
is\dss the disjoint\sss union of\dss the sets
$\ob\qff \vect\downarrow W_{\fff i}$,\qss 
$0\qff \leq\qff i\qff \leq\qff n$\nnsp,\oss
where\sss $\vect\downarrow W_{\fff i}$\sss are\sss
the\sss topological\sss categories from\dss Section\qss \ref{segal-example}.\oss 
The category\sss $\vect\downarrow \sigma$\sss
has\sss two\sss types of\dss morphisms.\oss
First,\oss morphisms
of\dss categories\sss $\vect\downarrow W_{\fff i}$\sss
are morphisms of\dss
$\vect\downarrow \sigma$\nnsp.\oss
Second,\oss for every\sss pair of\dss objects\sss
$V\qff \ttoo\qff W_{\fff i}$\sss
and\sss
$V\fff'\qff \ttoo\qff W_{\fff k}$\sss
of\dss the categories\sss
$\vect\downarrow W_{\fff i}$\sss
and\sss
$\vect\downarrow W_{\fff k}$\sss
respectively\sss
such\sss that\sss $i\qff <\qff k$\dss
there\dss is\dss a unique morphism\sss from\sss
$V\qff \ttoo\qff W_{\fff i}$\dss
to\sss
$V\fff'\qff \ttoo\qff W_{\fff k}$\nsp.\oss
There\dss is\dss a commutative diagram\vspace{-1.5pt}\vspace{1pt}
\[
\quad
\begin{tikzcd}[column sep=boomm, row sep=boomm]
\vect\downarrow \sigma
\arrow[r]
\arrow[d, "\dis k\dff"']
&
Q/\fff H
\arrow[d, "\dis \dff i"]
\\
\sigma
\arrow[r]
&
\hat{\mathcal{S}}\dff,
\end{tikzcd}
\]

\vspace{-12pt}\vspace{1pt}
where $k$\sss is\dss induced\sss by $i$\sss
and\sss the horizontal\sss arrows are\sss the inclusions.\oss
Let\sss $\kappa\off =\off \num{k}$\nnsp.\oss
By\sss passing\sss to geometric realizations 
we get\sss the commutative diagram\vspace{-1.5pt}\vspace{1pt}
\[
\quad
\begin{tikzcd}[column sep=boom, row sep=boomm]
\protect{\num{\vect\downarrow \sigma}}
\arrow[r]
\arrow[d, "\dis \kappa\dff"']
&
\protect{\num{Q/\fff H}}
\arrow[d, "\dis \qff \iota"]
\\
\protect{\num{\sigma}}
\arrow[r]
&
\protect{\num{\hat{\mathcal{S}}}}\dff,
\end{tikzcd}
\]

\vspace{-12pt}\vspace{1pt}
Let\sss $\tau$\sss be
an $l$\dnsp-simplex of\sss $Q/\fff H$\nnsp.\oss
It\dss is\dss uniquely\sss determined\sss by a sequence\vspace{1.5pt}\vspace{1pt}
\[
\quad
V_{\dff 0}\qff \ttoo\qff
V_{\dff 1}\qff \ttoo\qff
\ldots\qff \ttoo\qff
V_{\dff l}
\]

\vspace{-12pt}\vspace{1.5pt}\vspace{1pt}
of\dss morphisms of\sss $Q$\sss together\sss with an\sss
isometric embedding\sss
$h_{\dff l}\dff \colon\dff
V_{\dff l}\qff \ttoo\qff H$\nnsp.\oss
The corresponding\sss isometric embeddings\sss
$h_{\dff k}\dff \colon\dff
V_{\fff k}\qff \ttoo\qff H$\sss
are\sss equal\dss to\sss 
the compositions of\sss $h_{\dff l}$\sss with\sss 
the isometric embeddings corresponding\sss to\sss
morphisms $V_{\fff k}\qff \ttoo\qff V_{\dff l}$\nsp.\oss
The map $\iota$ takes\sss the simplex $\num{\tau}$ 
to\sss the simplex $\num{\sigma}$
if\trs and\dss only\trs if\dss
$i\trf(\trf \tau\trf)
\off =\off 
\theta^{\dff *}\dff(\trf \sigma\trf)$\sss
for some surjective non-decreasing\sss map\sss
$\theta\dff \colon\dff
[\dff l\qff]\qff \ttoo\qff [\halfff n\dff]$\nnsp.\oss
The\sss latter condition holds\sss implies\sss that\sss
$\tau$\sss is\dss a simplex of\dss $\vect\downarrow \sigma$\nnsp.\oss
It\sss follows\sss that\sss the inclusion\sss 
$\num{\vect\downarrow \sigma}
\qff \ttoo\qff
\num{Q/\fff H}$\sss 
induces a canonical\sss homeomorphism\vspace{1.5pt}\vspace{1pt}
\[
\quad
h_{\dff \sigma}\qff \colon\qff
\kappa^{\dff -\dff 1}\qff \bigl(\qff \inte \num{\sigma}\qff\bigr)
\off \ttoo\off
\iota^{\dff -\dff 1}\qff \bigl(\qff \inte \num{\sigma}\qff\bigr)
\pff.
\]

\vspace{-10.5pt}\vspace{1pt}
The simplicial\sss map\sss $k$\sss and\sss its geometric realization
$\kappa$ are standard\sss maps
depending\sss only\sss on\sss the dimension of\dss subspaces\sss $W_{\fff i}$.\oss
In\sss more details,\oss the geometric realization\sss 
$\num{\vect\downarrow \sigma}$\sss 
is\dss the join
of\sss $n\qff +\qff 1$\sss geometric realizations\sss
$\num{\vect\downarrow W_{\fff i}}$,\oss
and $k$\sss is\sss the canonical\sss map of\dss this join\sss to $\num{\sigma}$\nnsp.\oss
By\trs Lemma\qss \ref{over-is-contractible}\qss the spaces\sss
$\num{\vect\downarrow W_{\fff i}}$\sss are contractible.\oss
This implies,\oss in\sss particular,\oss that\sss
$\kappa\off =\off \num{k}$\sss
is\dss a\sss locally\sss trivial\sss bundle with contractible fibers
over\sss $\inte \num{\sigma}$\nnsp.\oss\vspace{0.525pt}

As in\sss the proof\dss of\trs Theorem\qss \ref{to-models},\oss
we now allow $\sigma$ vary among $n$\dnsp-simplices of\dss 
the\sss $\Delta$\dnsp-space\sss $\hat{\mathcal{S}}$\dnsp,\oss
or,\oss more precisely,\oss of\dss the core $\core \hat{\mathcal{S}}$\dnsp.\oss
The union of\dss the interiors $\inte \sigma$\sss is\dss
equal\sss to\sss\vspace{1.5pt}\vspace{0.525pt} 
\[
\quad
\num{\ssk_{\dff n}\dff \hat{\mathcal{S}}}
\off \smallsetminus\off
\num{\ssk_{\dff n\dff -\dff 1}\dff \hat{\mathcal{S}}}
\pff,
\]

\vspace{-12pt}\vspace{1.5pt} \vspace{0.525pt}
and\sss the rest\sss of\dss the proof\dss is\dss similar\sss to\sss 
the\sss last\sss part\sss of\dss the proof\dss of\trs 
Theorem\qss \ref{to-models}.\oss  \eproof\vspace{0.525pt}

\mypar{Theorem.}{operators-categories}
\emph{The maps}\vspace{1.5pt}
\[
\quad
\hat{\mathcal{F}}
\off \longleftarrow\off
\num{\hat{\mathcal{E}}}
\qff \ttoo\qff
\num{\mathcal{E}\hat{\mathcal{O}}}
\qff \ttoo\qff
\num{\hat{\mathcal{O}}}
\qff \ttoo\qff
\num{\hat{\mathcal{S}}}
\]

\vspace{-12pt}\vspace{1.5pt}\vspace{0.525pt}
\emph{are\sss homotopy\sss equivalences,\oss
and\qss $\num{\hat{\mathcal{S}}}$\sss
is\dss canonically\dss homotopy\sss equivalent\dss 
to\trs $\num{Q}$\nnsp.\oss}\vspace{0.525pt}

\proof
Theorems\qss \ref{forgetting-enhancement},\oss
\ref{to-enhanced-models},\oss 
\ref{to-models},\oss 
\ref{forgetting-operators}\pss
together\sss imply\sss the first\sss statement\sss of\dss the\sss theorem.\oss
The second statement\sss follows 
from\sss Theorem\qss \ref{intermediary}.\oss  \eproof\vspace{0.525pt}

\myuppar{Remarks.}
The\sss second statement\sss of\qss Theorem\qss \ref{operators-categories}\qss 
is\dss due\sss to\dss Segal\qss \cite{s4},\oss
who only outlined\sss the main\sss ideas behind\dss the proof.\oss
The above proof\dss is\dss inspired\sss by\dss Segal's\dss outline,\oss
but\sss uses different\dss tools.\oss
In\sss particular,\oss it\sss avoids using\dss Proposition\qss 2.7\qss
from\qss \cite{s2},\oss
referred\sss by\dss Segal\dss as\sss the\sss key\sss tool\sss for making
his outline rigorous.\oss
The use of\dss coverings in\sss the proof\dss of\trs
Theorem\qss \ref{forgetting-enhancement}\qss seems\sss to be new,\oss
as also\sss the use of\dss the discrete\sss topology on\sss the set\sss
of\dss controlling\sss parameters $\varepsilon$\nnsp.\oss
Segal\qss \cite{s2}\qss uses\sss the usual\sss topology of\dss $\rrr$\sss
for similar purposes.\oss\vspace{0.525pt}

In\dss Segal's\dss outline\sss the role of\dss categories\sss
$\hat{\mathcal{E}}$\dnsp,\qss
$\mathcal{E}\hat{\mathcal{O}}$\dnsp,\qss
$\hat{\mathcal{O}}$\dnsp,\oss and\sss 
$\hat{\mathcal{S}}$\sss
is\dss played\sss by\sss
an abstract\sss version of\sss 
$\hat{\mathcal{O}}$\nnsp.\oss
This version\sss has as objects pairs\sss $(\trf V\fff,\qff F\trf)$\nnsp,\oss
where $V$\sss is\dss a finitely dimensional\dss Hilbert\sss space and\sss
$F\dff \colon\dff V\qff \ttoo\qff V$\sss is\dss a self-adjoint\sss operator.\oss
The morphisms\sss
$(\trf V\fff,\qff F\trf)\qff \ttoo\qff (\trf V\fff'\fff,\qff F\fff'\trf)$
are isometric embeddings\sss
$\iota\dff \colon\dff V\qff \ttoo\qff V\fff'$\sss
such\sss that\sss
$\iota\qff \circ\qff F
\off =\off
F\fff'\dff \circ\pff \iota$
and\dss the map\sss
$\kernel F\qff \ttoo\qff \kernel F\fff'$\sss 
induced\sss by\dss $\iota$\dss is\dss an\sss isomorphism.\oss
See\sss the\sss last\sss paragraph of\qss \cite{s4}.\oss\vspace{0.525pt}

\myuppar{Small\sss versions of\dss the category $Q$\nnsp.}
Instead of\dss using\sss the set\sss of\qss ``all''\qss finitely dimensional\sss vector spaces
with a scalar product\sss in\sss the definition of\sss $Q$\sss
one can use any set\sss
containing\sss a representative from each\sss isomorphism class.\oss
The inclusion of\dss the resulting version $Q'$ of\sss $Q$\sss
into\sss the original\sss one is\dss an equivalence of\dss categories,\oss
and\sss hence\sss the inclusion $\num{Q'}\qff \ttoo\qff \num{Q}$\sss is\dss
homotopy equivalence.\oss 
The smallest\sss version\sss has as objects\sss 
the standard vector spaces $\ccc^{\dff n}$\dnsp\dnsp.\vspace{0.525pt}

\myuppar{Finite-unitary\sss groups.}
For {\nsp}$n\qff \in\qff \nnn$\dss
let\sss $U^{\fff n}$\sss
be\sss the group of\dss unitary operators\sss
$H\qff \ttoo\qff H$\sss
equal\sss to\sss the identity on some closed\sss subspace of\sss $H$\sss
of\dss codimension $n$ with\sss the norm\sss topology.\oss
Clearly,\pss 
$U^{\trf 0}\off \subset\off
U^{\fff 1}\off \subset\off
U^{\dff 2}\off \subset\off
\dots
\off$.\oss
Let\sss $U^{\dff \mathrm{fin}}$ be\sss the union of\dss the spaces $U^{\fff n}$
with\sss the direct\sss limit\sss topology,\oss
and\dss let\dss
$-\qff U^{\dff \mathrm{fin}}
\off =\off
\{\pff -\qff u\qff \mid\qff u\qff \in\qff U^{\dff \mathrm{fin}} \pff\}$\nnsp.\oss

\mypar{Theorem.}{harris-h}
\emph{There\dss is\dss a canonical\dss homeomorphism\dss
$h\dff \colon\dff
\num{\hat{\mathcal{S}}}\qff \ttoo\qff -\qff U^{\dff \mathrm{fin}}${\nsp}.\oss}

\proof
An $n$\dnsp-simplex $\sigma$ of\sss $\hat{\mathcal{S}}$\sss
is\dss defined\sss by a sequence\sss\vspace{1.5pt}
\[
\quad
V_{\dff 0}\qff \ttoo\qff
V_{\dff 1}\qff \ttoo\qff
\ldots\qff \ttoo\qff
V_{\dff n}
\pff,
\]

\vspace{-12pt}\vspace{1.5pt}
of\dss morphisms of\sss $\hat{\mathcal{S}}$\nsp.\oss
Let\sss\vspace{1.5pt}
\[
\quad
V_{\dff i}
\off =\off 
U_{\dff -}^{\dff i}\dff \oplus\dff V_{\dff i\dff -\dff 1}\dff \oplus\dff U_{\dff +}^{\dff i}
\]

\vspace{-12pt}\vspace{1.5pt}
be\sss the decomposition defining\sss the morphism\sss
$V_{\dff i\dff -\dff 1}\qff \ttoo\qff V_{\dff i}$\nsp,\oss
and\sss let\dss
$P_{\dff -}^{\dff i}$,\pss
$P_{\dff +}^{\dff i}$,\pss
$P_{\dff i}$\dss
be\sss the orthogonal\sss projections of\sss $H$\sss onto\sss
the subspaces\dss
$U_{\dff -}^{\dff i}$,\pss
$U_{\dff +}^{\dff i}$,\pss
$V_{\dff i}$\dss
respectively.\oss

Let\sss us\sss redefine\sss the  standard $n$\dnsp-dimensional\sss simplex\sss
$\Delta^n$\sss as\sss the space of\dss points\sss\vspace{1.5pt}
\[
\quad
u
\off =\off
(\trf u_{\dff 1}\dff,\qff \ldots\dff,\qff u_{\dff n}\trf)
\qff \in\qff \rrr^{\dff n}
\]

\vspace{-12pt}\vspace{1.5pt}
such\sss that\sss
$0\qff \leq\qff
u_{\dff 1}\qff \leq\qff u_{\trf 2}\qff \leq\qff \ldots\qff \leq\qff u_{\dff n}
\qff \leq\qff 1$\nnsp.\oss
For an $n$\dnsp-simplex $\sigma$ as above and\sss 
$u\qff \in\qff \Delta^n$\sss 
let\dss $S\trf(\trf \sigma\fff,\qff u \trf)$\sss
be\sss the operator\sss $H\qff \ttoo\qff H$\sss 
equal\sss to\vspace{3pt}\vspace{1.5pt}
\begin{equation}
\label{simplex-operator}
\quad
-\off
u_{\dff n}\qff P_{\dff -}^{\dff n}
\off -\off
\ldots
\off -\off
u_{\dff 2}\qff P_{\dff -}^{\dff 2}
\off -\off
u_{\dff 1}\qff P_{\dff -}^{\dff 1}
\off +\off
u_{\dff 1}\qff P_{\dff +}^{\dff 1}
\off +\off
u_{\dff 2}\qff P_{\dff +}^{\dff 2}
\off +\off
\ldots
\off +\off
u_{\dff n}\qff P_{\dff +}^{\dff n}
\end{equation}

\vspace{-12pt}\vspace{3pt}\vspace{1.5pt}
on\sss $V_{\dff n}$\sss
and\sss to\sss $\id_{\trf H}$\dss on\sss $H\dff \ominus\off V_{\dff n}$\nsp.\oss
Let\dss\vspace{1.5pt}\vspace{1.5pt}
\[
\quad 
U\trf(\trf \sigma\fff,\qff u \trf)
\off =\off
\exp\qff\bigl(\trf
\pi\fff i\pff
S\trf(\trf \sigma\fff,\qff u \trf)
\trf\bigr)
\qff.
\]

\vspace{-12pt}\vspace{1.5pt}\vspace{1.5pt}
Then\sss $S\trf(\trf \sigma\fff,\qff u \trf)$\sss
is\dss a self-adjoint\sss operator\sss in\sss $H$\nnsp,\oss
and\sss $U\trf(\trf \sigma\fff,\qff u \trf)$\sss
is\dss a unitary\sss operator equal\sss to\sss 
$-\qff \id_{\trf H}$\dss on\sss $H\dff \ominus\off V_{\dff n}$\nsp.\oss
The maps\dss
$\hat{\mathcal{S}}_{\dff n}\dff \times\dff \Delta^n
\qff \ttoo\qff
-\qff U^{\dff \mathrm{fin}}$\dss
defined\sss by\sss
$(\trf \sigma\fff,\qff u \trf)
\off \longmapsto\off
U\trf(\trf \sigma\fff,\qff u \trf)$\sss
agree with\sss the equivalence relation defining\sss
$\num{\hat{\mathcal{S}}}$\sss
and\dss hence induce a map\sss\vspace{1.5pt}
\[
\quad 
h\dff \colon\dff
\num{\hat{\mathcal{S}}}
\qff \ttoo\qff 
-\qff U^{\dff \mathrm{fin}}
\pff.
\]

\vspace{-12pt}\vspace{1.5pt}
The spectral\sss theorem\sss implies\sss
that $h$\sss is\dss a\sss bijection,\oss and since
$\num{\hat{\mathcal{S}}}$ and\sss $-\qff U^{\dff \mathrm{fin}}$ 
have\sss the direct\sss limit\sss
topology,\pss $h$ is\dss a\sss homeomorphism.\oss  \eproof

\mypar{Corollary.}{subspace-unitary}
\emph{There\dss is\dss a canonical\dss homotopy\sss equivalence\dss
$\hat{\mathcal{F}}\qff \ttoo\qff U^{\dff \mathrm{fin}}${\nsp}.\oss}

\proof
This immediately\sss follows\sss from\trs Theorems\qss \ref{harris-h}\qss
and\qss \ref{operators-categories}.\oss  \eproof

\myuppar{Remark.}
Theorem\qss \ref{harris-h}\qss and\dss its proof\dss
are an adaptation\sss to\dss Hilbert\dss spaces,\oss
the group\sss $U\ffin$\dnsp,\oss and\sss the category\sss 
$\hat{\mathcal{S}}$\sss
of\dss a\sss theorem of\qss Harris\qss \cite{h}.\oss
See\qss \cite{h},\oss Theorem\dss in\dss Section\qss 2.\oss
We will\sss need also\sss the\sss theorem of\trs Harris\dss
itself\dss and\sss will\sss return\dss to\dss it\sss
at\dss the end of\trs Section\qss \ref{categories-grassmannians}.\oss

\mysection{Polarizations\qss and\qss splittings}{polarizations-splittings}

\myuppar{Polarizations\sss of\trs Hilbert\sss spaces.}
A\qss \emph{polarization}\qss of\dss a separable infinitely dimensional\dss
Hilbert\sss space\sss $K$\sss is\dss 
a presentation of\dss $K$\sss as\sss an orthogonal\sss direct\sss sum\sss
$K\off =\off K_{\dff -}\dff \oplus\dff K_{\dff +}$
of\dss an ordered\sss pair of\dss two infinitely dimensional\dss subspaces.\oss
The set\sss of\dss all\sss polarizations\dss is\dss equipped\sss with\sss
topology\sss defined\sss by\sss the norm\sss topology\sss on\sss the spaces
of\dss orthogonal\sss projections\sss
$K\qff \ttoo\qff K_{\dff -}$\nsp,\qss
$K\qff \ttoo\qff K_{\dff +}$\nsp.\oss
We will\sss call\sss this\sss topology\sss simply\sss the\qss
\emph{norm\dss topology}.\oss
By a well\sss known\sss theorem of\trs Atiyah\sss and\dss Singer\qss \cite{as}\qss
the space of\dss polarizations of\dss $K$\sss is\dss contractible.\oss

\myuppar{Polarized\sss enhanced operator\sss models.}
A\qss \emph{polarized\sss enhanced\sss operator\sss model}\pss is\dss a $5$\dnsp-tuple
\[
\quad
(\trf V\fff,\pff F\fff,\pff \varepsilon\fff,\pff H_{\dff -}\fff,\qff H_{\dff +}\trf)
\]

\vspace{-12pt}
such\sss that $(\trf V\fff,\pff F\fff,\pff \varepsilon\trf)$\sss
is\dss an enhanced operator\sss model\sss and
$H_{\dff +}\dff,\pff H_{\dff -}$ is\dss 
a\sss polarization of\sss $H\dff \ominus\dff V$\dnsp,\oss 
i.e.\qss
$H
\off =\off
H_{\dff -}\dff \oplus\qff V\qff \oplus\qff H_{\dff +}$\nsp.\oss
Let\sss $\mathcal{PE}\hat{\mathcal{O}}$\sss be\sss the space 
of\dss polarized\sss enhanced\sss operator\sss models
with\sss topology defined\sss by\sss the\sss topology\sss
on\sss the space of\dss enhanced operator models\sss
$(\trf V\fff,\pff F\fff,\pff \varepsilon\trf)$\sss and\sss
the norm\sss topology on\sss polarizations.\oss
The space $\mathcal{PE}\hat{\mathcal{O}}$\sss
is\dss ordered\sss by\sss the relation\sss $\leq$\nnsp,\oss where\vspace{0.9pt}
\[
\quad
(\trf V\fff,\pff F\fff,\pff \varepsilon\fff,\pff H_{\dff -}\fff,\pff H_{\dff +}\trf)
\off \leq\off
(\trf V\fff'\fff,\pff F\fff'\fff,\pff \varepsilon'\fff,\pff H\fff'_{\dff -}\fff,\pff H\fff'_{\fff +}\trf)
\quad
\]

\vspace{-39pt}\vspace{0.9pt}
\[
\quad
\mbox{if}\dff\quad
(\trf V\fff,\pff F\fff,\pff \varepsilon\trf)
\off \leq\off
(\trf V\fff'\fff,\pff F\fff'\fff,\pff \varepsilon'\trf)
\quad
\mbox{and}
\]

\vspace{-39pt}\vspace{0.9pt}
\[
\quad
H_{\dff -}
\off =\off
H\fff'_{\fff -}
\qff \oplus\pff
\image\dff P_{\dff [\dff -\qff \varepsilon'\fff,\qff -\qff \varepsilon\trf]}\dff(\trf F\fff'\trf)
\qff,
\qquad
H_{\dff +}
\off =\off
H\fff'_{\dff +}
\qff \oplus\pff
\image\dff P_{\dff [\dff \varepsilon\fff,\qff \varepsilon'\trf]}\dff(\trf F\fff'\trf)
\pff.
\]

\vspace{-12pt}\vspace{0.9pt}
The\sss last\sss condition can\sss be reformulated\sss in\sss terms of\dss
vector space models.\oss
The inequality\sss
$(\trf V\fff,\pff F\fff,\pff \varepsilon\trf)
\off \leq\off
(\trf V\fff'\fff,\pff F\fff'\fff,\pff \varepsilon'\trf)$\sss
means\sss that\sss there\dss is\dss a morphism\dss
$(\trf V\fff,\pff F\fff,\pff \varepsilon\trf)
\qff \ttoo\qff
(\trf V\fff'\fff,\pff F\fff'\fff,\pff \varepsilon'\trf)$\dss
of\dss enhanced operator models.\oss
This\sss morphism\sss leads\sss to a morphism\dss
$V\qff \ttoo\qff V\fff'$\sss of\dss subspace models.\oss
Let\sss
$V\fff'
\off =\off
U_{\dff -}\qff \oplus\qff
V\qff \oplus\qff U_{\dff +}$\sss
be\sss the corresponding\sss orthogonal\sss decomposition.\oss
In\sss terms of\dss this decomposition\sss the\sss last\sss
condition means\sss that\vspace{0pt}
\[
\quad
H_{\dff -}
\off =\off
H\fff'_{\fff -}
\qff \oplus\qff
U_{\dff -}
\quad
\mbox{and}\quad
H_{\dff +}
\off =\off
H\fff'_{\fff +}
\qff \oplus\qff
U_{\dff +}
\pff.
\]

\vspace{-12pt}
In\sss particular,\oss
if\dss the morphism\sss
$V\qff \ttoo\qff V\fff'$\sss
is\dss given,\oss
then\sss polarizations\sss
$H\dff \ominus\dff V
\off =\off
H_{\dff -}\dff \oplus\dff H_{\dff +}$\sss
and\sss
$H\dff \ominus\dff V\fff'
\off =\off
H\fff'_{\fff -}\dff \oplus\dff H\fff'_{\dff +}$\sss
determine each other.\oss
The above order\sss $\leq$\sss defines a structure of\dss
a\sss topological\sss category on\sss $\mathcal{PE}\hat{\mathcal{O}}$\dnsp.\oss
There\dss is\dss an obvious forgetting\sss functor\sss
$\hat{\pi}\dff \colon\dff 
\mathcal{PE}\hat{\mathcal{O}}
\qff \ttoo\qff 
\mathcal{E}\hat{\mathcal{O}}$\nnsp.\oss
There\dss is\dss also a functor\sss
$\mathcal{P}\dff\hat{\psi}\dff \colon\dff 
\hat{\mathcal{E}}
\qff \ttoo\qff 
\mathcal{PE}\hat{\mathcal{O}}$\dss
taking 
$(\trf A\fff,\qff \varepsilon\trf)$\sss
to\qss\vspace{0.9pt} 
\[
\quad
(\trf V\fff,\pff F\fff,\pff \varepsilon\fff,\pff H_{\dff -}\fff,\qff H_{\dff +}\trf)
\qff,
\qquad
\mbox{where}
\]

\vspace{-39pt}\vspace{0.9pt}
\[
\quad
V
\off =\off\dff 
\image\dff P_{\dff [\dff -\qff \varepsilon\fff,\qff \varepsilon\trf]}\dff(\trf A\trf)
\qff,\quad
H_{\dff -}
\off =\off\dff 
\image\dff P_{\dff (\dff -\qff \infty\fff,\qff -\qff \varepsilon\trf)}\dff(\trf A\trf)
\qff,\quad
H_{\dff +}
\off =\off\dff 
\image\dff P_{\dff (\dff \varepsilon\fff,\qff \infty\trf)}\dff(\trf A\trf)
\qff,
\]

\vspace{-12pt}\vspace{0.9pt}
and\sss the operator\sss 
$F\dff \colon\dff V\qff \ttoo\qff V$\sss
is\dss induced\sss by $A$\nnsp.\oss
In\sss particular,\pss
$(\trf V\fff,\qff F\fff,\qff \varepsilon\trf)$\sss
is\dss an\sss enhanced\sss operator\sss model\sss of\dss the enhanced operator\sss
$(\trf A\fff,\qff \varepsilon\trf)$\nnsp.\oss

\mypar{Theorem.}{to-p-models}
\emph{The map\sss
$\num{\mathcal{P}\dff\hat{\psi}}\dff \colon\dff
\num{\hat{\mathcal{E}}}
\qff \ttoo\qff 
\num{\mathcal{PE}\hat{\mathcal{O}}}$\sss
is\dss a\sss homotopy\sss equivalence.\oss}

\proof
The proof\dss is\dss also similar\sss to\sss the proof\dss of\qss Theorem\qss \ref{to-enhanced-models}.\oss
Recall\sss that\sss an $n$\dnsp-simplex of\dss the category\sss $\hat{\mathcal{E}}$\sss
can\sss be identified\sss with an operator $A\qff \in\qff \hat{\mathcal{F}}$
together\sss with\sss a non-decreasing\sss sequence\sss
$a_{\dff 0}\qff \leq\qff a_{\dff 1}\qff \leq\qff \ldots\qff \leq\qff a_{\dff n}$\sss 
of\dss positive numbers such\sss that\sss
$(\trf A\fff,\qff a_{\dff i}\trf)$\sss is\dss 
an enhanced operator for every\sss $i$\nnsp.\oss
Similarly,\oss
an $n$\dnsp-simplex $\sigma$ of\dss $\mathcal{PE}\hat{\mathcal{O}}$\sss
can\sss be identified\sss with a finitely\sss dimensional\sss subspace\sss 
$V_{\fff n}\qff \subset\qff H$\sss
together\sss with\sss a self-adjoint\sss operator
$F_{\fff n}\dff \colon\dff
V_{\fff n}\qff \ttoo\qff V_{\fff n}$\nsp,\oss
a non-de\-creas\-ing\sss sequence\sss
$a_{\dff 0}\qff \leq\qff a_{\dff 1}\qff \leq\qff \ldots\qff \leq\qff a_{\dff n}$\sss
of\dss positive numbers as\sss in\sss the proof\dss of\trs 
Theorem\qss \ref{forgetting-polarization},\oss
and a polarization\sss
$H\dff \ominus\dff V_{\fff n}
\off =\off
H_{\dff -}\dff \oplus\qff H_{\dff +}$\sss
of\dss $H\dff \ominus\dff V_{\fff n}$\nsp.\oss
The preimage of\sss $\sigma$\sss in\sss the space of\sss
$n$\dnsp-simplices of\sss $\hat{\mathcal{E}}$\sss
can\sss be identified\sss with\sss the space 
of\dss operators $A\qff \in\qff \hat{\mathcal{F}}$\sss such\sss that\sss
$V_{\fff n}
\off =\off 
\image\dff P_{\dff [\dff -\qff a_{\dff n}\fff,\qff a_{\dff n}\trf]}\dff(\trf A\trf)$\nnsp,\oss
the operator\sss
$V_{\fff n}\qff \ttoo\qff V_{\fff n}$\sss
induced\sss by\sss $A$\sss is\dss equal\sss to\sss $F_{\fff n}$\nsp,\oss
and\vspace{1.5pt}
\[
\quad
H_{\dff -}
\off =\off\dff 
\image\dff P_{\dff (\dff -\qff \infty\fff,\qff -\qff a_{\dff n}\trf)}\dff(\trf A\trf)
\qff,\quad
H_{\dff +}
\off =\off\dff 
\image\dff P_{\dff (\dff a_{\dff n}\fff,\qff \infty\trf)}\dff(\trf A\trf)
\qff.
\]

\vspace{-10.5pt}
The space of\dss such operators\dss is\dss homeomorphic\sss to\sss
$\left(\qff \hat{\mathcal{F}}_{\qff >\qff a_{\dff n}}\qff\right)^{\dff 2}$\dnsp.\oss
See\dss Section\qss \ref{spaces-operators}.\oss
By\dss Proposition\qss \ref{positive-contractible}\qss this space\dss
is\dss contractible.\oss
As usual,\oss this implies\sss that\sss $\mathcal{P}\dff\hat{\psi}$\sss
induces a homotopy equivalence of\dss the spaces of\sss $n$\dnsp-simplices
for every $n$\nnsp,\oss
and\sss this,\oss in\sss turn,\oss implies\sss that\sss
the geometric realization\sss
$\num{\mathcal{P}\dff\hat{\psi}}\dff \colon\dff
\num{\hat{\mathcal{E}}}
\qff \ttoo\qff 
\num{\mathcal{PE}\hat{\mathcal{O}}}$\sss
is\dss a\sss homotopy\sss equivalence.\oss  \eproof

\mypar{Theorem.}{forgetting-polarization}
\emph{The map\sss
$\num{\hat{\pi}}\dff \colon\dff
\num{\mathcal{PE}\hat{\mathcal{O}}}
\qff \ttoo\qff 
\num{\mathcal{E}\hat{\mathcal{O}}}$\sss
is\dss a\sss homotopy\sss equivalence.\oss}

\proof
This\sss proof\dss is\dss also similar\sss to\sss the proof\dss of\qss
Theorem\qss \ref{to-enhanced-models}.\oss
As we pointed out\sss in\sss that\sss proof,\oss
an $n$\dnsp-simplex of\dss the category\sss $\mathcal{E}\hat{\mathcal{O}}$\sss
can\sss be identified\sss with a finitely\sss dimensional\sss subspace\sss 
$V_{\fff n}\qff \subset\qff H$\sss
together\sss with\sss a self-adjoint\sss operator
$F_{\fff n}\dff \colon\dff
V_{\fff n}\qff \ttoo\qff V_{\fff n}$\sss
and\sss a non-de\-creas\-ing\sss sequence\sss
$a_{\dff 0}\qff \leq\qff a_{\dff 1}\qff \leq\qff \ldots\qff \leq\qff a_{\dff n}$\sss
of\dss positive numbers such\sss that\sss
$\sigma\dff(\dff F\dff)
\off \subset\off
(\dff -\qff a_{\dff n}\fff,\pff a_{\dff n}\trf)$\sss
and\sss the numbers\sss $-\qff a_{\dff i}\dff,\qff a_{\dff i}$\sss 
are not\sss eigenvalues of\sss $F_{\fff n}$\nsp.\oss
Similarly,\oss an $n$\dnsp-simplex of\dss the category\sss 
$\mathcal{PE}\hat{\mathcal{O}}$\sss
can\sss be identified\sss with such data\sss together\sss
with a polarization of\dss $H\dff \ominus\dff V_{\fff n}$\nsp.\oss
Applying $\hat{\pi}$ simply\sss forgets\sss the polarization.\oss 
Since\sss the spaces of\dss polarizations are contractible,\pss
$\hat{\pi}$\sss induces a homotopy equivalence between\sss the spaces
of\sss $n$\dnsp-simplices for every $n$\nnsp.\oss
Since both categories\sss $\mathcal{PE}\hat{\mathcal{O}}$\sss
and\sss $\mathcal{E}\hat{\mathcal{O}}$\sss are,\oss
obviously,\oss categories with\sss free units,\oss
this implies\sss that\sss the geometric realization\sss
$\num{\hat{\pi}}\dff \colon\dff
\num{\mathcal{PE}\hat{\mathcal{O}}}
\qff \ttoo\qff 
\num{\mathcal{E}\hat{\mathcal{O}}}$\sss
is\dss a\sss homotopy\sss equivalence.\oss  \eproof

\myuppar{Polarized\sss operator\sss models.}
A\qss \emph{polarized\sss operator\sss model}\pss is\dss a\sss quadruple\sss
$(\trf V\fff,\pff F\fff,\pff H_{\dff -}\dff,\pff H_{\dff +}\trf)$\sss
such\sss that\sss $(\trf V\fff,\pff F \trf)$\sss is\dss an operator\sss model\sss
and\dss 
$H\dff \ominus\dff V
\off =\off
H_{\dff -}\dff \oplus\dff H_{\dff +}$\sss 
is\dss  a\sss polarization of\dss $H\dff \ominus\dff V$\dnsp.\oss
The space\sss $\mathcal{P}\fff\hat{\mathcal{O}}$\sss 
of\dss polarized\sss subspace\sss model\dss is\dss ordered\sss 
by\sss the relation\sss $\leq$\nnsp,\oss where\vspace{1.5pt}\vspace{1.0pt}\vspace{-0.125pt}
\[
\quad
(\trf V\fff,\off F\fff,\off H_{\dff -}\dff,\off H_{\dff +}\trf)
\off \leq\off
(\trf V\fff'\fff,\off F\fff'\fff,\off H\fff'_{\dff -}\dff,\off H\fff'_{\dff +}\trf)
\quad
\]

\vspace{-37.5pt}\vspace{1.0pt}
\[
\quad
\mbox{if}\dff\quad
(\trf V\fff,\off F \trf)
\off \leq\off
(\trf V\fff'\fff,\off F\fff' \trf)
\quad
\mbox{and}\quad
H\fff'_{\dff -}\qff \subset\pff H_{\dff -}
\qff,\quad
H\fff'_{\dff +}\qff \subset\pff H_{\dff +}
\pff.
\]

\vspace{-10.5pt}\vspace{1.0pt}\vspace{-0.125pt}
As before,\oss we consider\sss this order 
as\sss a structure of\dss a\sss topological\sss
category on\sss $\mathcal{P}\fff\hat{\mathcal{O}}$\dnsp.\oss
There are obvious forgetting\sss functors\sss
$\hat{\pi}\dff \colon\dff
\mathcal{P}\fff\hat{\mathcal{O}}
\qff \ttoo\qff
\hat{\mathcal{O}}$\sss
and\qss
$\mathcal{P}\hat{o}\dff \colon\dff
\mathcal{PE}\hat{\mathcal{O}}
\qff \ttoo\qff
\mathcal{P}\fff\hat{\mathcal{O}}$\dnsp.\oss

\myuppar{Polarized\sss subspace\sss models.}
A\qss \emph{polarized\sss subspace\sss model}\pss is\dss a\sss triple
$(\trf V\fff,\pff H_{\dff -}\dff,\pff H_{\dff +}\trf)$\nnsp,\oss
where\sss $V$\sss is\dss a\sss finitely dimensional\sss subspace of\sss $H$\sss
and\dss 
$H\dff \ominus\dff V
\off =\off
H_{\dff -}\dff \oplus\dff H_{\dff +}$\sss 
is\dss  a\sss polarization of\dss $H\dff \ominus\dff V$\dnsp.\oss
The space\sss $\mathcal{P}{\nsp}\hat{\mathcal{S}}$\sss 
of\dss polarized\sss subspace\sss models\dss is\dss ordered\sss 
by\sss the relation\sss $\leq$\nnsp,\oss where\vspace{0pt}
\[
\quad
(\trf V\fff,\off H_{\dff -}\dff,\off H_{\dff +}\trf)
\off \leq\off
(\trf V\fff'\fff,\off H\fff'_{\dff -}\dff,\off H\fff'_{\dff +}\trf)
\quad
\]

\vspace{-39pt}
\[
\quad
\mbox{if}\dff\quad
V\off \subset\off V\fff'
\quad
\mbox{and}\quad
H\fff'_{\dff -}\qff \subset\pff H_{\dff -}
\qff,\quad
H\fff'_{\dff +}\qff \subset\pff H_{\dff +}
\pff.
\]

\vspace{-12pt}
As before,\oss we consider\sss this order on\sss the\sss topological\sss space\sss 
$\mathcal{P}{\nsp}\hat{\mathcal{S}}$\sss 
as\sss a structure of\dss a\sss topological\sss
category on\sss $\mathcal{P}{\nsp}\hat{\mathcal{S}}$\dnsp.\oss
The morphisms of\dss $\mathcal{P}{\nsp}\hat{\mathcal{S}}$\sss 
have\sss the form\vspace{0pt}
\[
\quad
(\trf V\fff,\off H_{\dff -}\dff,\off H_{\dff +}\dff)
\qff \ttoo\qff
(\trf U_{\dff -}\qff \oplus\qff
V\qff \oplus\qff U_{\dff +}\dff,\off 
H_{\dff -}\dff \ominus\dff U_{\dff -}\dff,\off 
H_{\dff +}\dff \ominus\dff U_{\dff +}
\trf)
\qff,
\]

\vspace{-12pt}
where\sss $U_{\dff -}\dff,\pff U_{\dff +}$\sss
are finitely dimensional\sss subspaces of\dss  
$H_{\dff -}\dff,\pff H_{\dff +}$\sss respectively.\oss
There\dss is\dss a\sss forgetting\sss functor\sss
$\mathcal{P}\hat{\omega}\dff \colon\dff
\mathcal{P}\fff\hat{\mathcal{O}}
\qff \ttoo\qff
\mathcal{P}{\nsp}\hat{\mathcal{S}}$\sss
defined\sss in\sss the same way as\sss 
$\hat{\omega}\dff \colon\dff
\hat{\mathcal{O}}
\qff \ttoo\qff
\hat{\mathcal{S}}$\nsp\dnsp.\oss
Also,\oss there\dss is\dss a\qss ``forgetting''\qss functor\sss
$\hat{\pi}\dff \colon\dff
\mathcal{P}{\nsp}\hat{\mathcal{S}}
\qff \ttoo\qff 
\hat{\mathcal{S}}$\sss
taking a polarized subspace model\sss
$(\trf V\fff,\pff H_{\dff -}\dff,\pff H_{\dff +}\trf)$\sss
to\sss the subspace model\sss $V$\dnsp,\oss
and\sss a\sss morphism\sss\vspace{0pt}
\[
\quad
(\trf V\fff,\off H_{\dff -}\dff,\off H_{\dff +}\trf)
\off \ttoo\off
(\trf V\fff'\fff,\off H\fff'_{\dff -}\dff,\off H\fff'_{\dff +}\trf)
\quad
\]

\vspace{-12pt}
of\dss  
$\mathcal{P}{\nsp}\hat{\mathcal{S}}$\dss
to\sss the morphism\sss
$V\qff \ttoo\qff V\fff'$\sss of\dss 
$\hat{\mathcal{S}}$\sss defined\sss by\sss 
the orthogonal\sss decomposition\vspace{0pt}
\[
\quad
V\fff'
\off =\off
U_{\dff -}\qff \oplus\qff
V\qff \oplus\qff U_{\dff +}
\pff,
\]

\vspace{-12pt}
where\sss
$U_{\dff -}
\off =\off 
H_{\dff -}\dff \ominus\qff H\fff'_{\fff -}$
and\dss
$U_{\dff +}
\off =\off 
H_{\dff +}\dff \ominus\qff H\fff'_{\fff +}$.\oss

\myuppar{Forgetting\sss functors.}
Let\sss us\sss consider\sss the following diagram of\dss 
forgetting\sss functors.\vspace{-1pt}
\[
\quad
\begin{tikzcd}[column sep=sboom, row sep=sboom]
\hat{\mathcal{E}}
\arrow[r, "\dis \mathcal{P}\dff\hat{\psi}"]
\arrow[d, "\dis =\dff"']
&
\mathcal{PE}\hat{\mathcal{O}}
\arrow[r, "\dis \mathcal{P}\hat{o}\vphantom{\mu}"]
\arrow[d, "\dis \hat{\pi}\dff"']
&
\mathcal{P}\hat{\mathcal{O}}
\arrow[r, "\dis \mathcal{P}\hat{\omega}"]
\arrow[d, "\dis \hat{\pi}\dff"']
&
\mathcal{P}{\nsp}\hat{\mathcal{S}}
\arrow[d, "\dis \hat{\pi}\dff"']
\\
\hat{\mathcal{E}}
\arrow[r, "\dis \hat{\psi}"]
&
\mathcal{E}\hat{\mathcal{O}}
\arrow[r, "\dis \hat{o}\vphantom{\mu}"]
&
\hat{\mathcal{O}}
\arrow[r, "\dis \hat{\omega}"]
&
\hat{\mathcal{S}}
\end{tikzcd}
\]

\vspace{-12pt}
\mypar{Theorem.}{forgetting-equivalence}
\emph{All\trs forgetting\dss functors\sss in\sss
the above diagram\sss induce homotopy equivalences of\dss
geometric\sss realizations.\oss}

\proof
For\sss the\sss lower horizontal\sss arrows\sss this was proved\sss
in\dss Section\qss \ref{classifying-spaces-saf},\oss
for\sss $\mathcal{P}\dff\hat{\psi}$\sss in\trs Theorem\qss \ref{to-p-models},\oss
and\sss for\dss
$\hat{\pi}\dff \colon\dff 
\mathcal{PE}\hat{\mathcal{O}}
\qff \ttoo\qff 
\mathcal{E}\hat{\mathcal{O}}$\dss
in\dss Theorem\qss \ref{forgetting-polarization}.\oss
The proof\dss for\sss the\sss two other arrows $\hat{\pi}$\sss
is\dss completely\sss similar\sss to\sss the proof\dss of\trs
Theorem\qss \ref{forgetting-polarization}.\oss
Now\sss the commutativity of\dss the diagram\sss implies\sss
that\sss $\mathcal{P}\hat{o}$\sss and\sss $\mathcal{P}\hat{\omega}$\sss
also induce\sss homotopy equivalences of\dss geometric realizations.\oss
One can also prove\sss this by adapting\sss the proofs of\trs
Theorems\qss \ref{to-models}\qss and\qss \ref{forgetting-operators}.\oss  \eproof

\myuppar{Split\dss polarized subspace models.}
Let\sss $\mathcal{P}$\sss
be\sss the full\sss subcategory of\sss
$\mathcal{P}{\nsp}\hat{\mathcal{S}}$\sss
having as objects\sss the objects of\sss
$\mathcal{P}{\nsp}\hat{\mathcal{S}}$\sss
of\dss the form\sss
$(\dff 0\fff,\qff H_{\dff -}\fff,\qff H_{\dff +}\dff)$\nnsp.\oss
This subcategory\sss has only\sss identity\sss morphisms,\oss
and\dss its space of\dss objects\dss is\dss nothing else but\sss the
space of\dss polarizations of $H$ and\sss hence\dss is\dss contractible.\oss
It\dss follows\sss that\sss 
$\num{\mathcal{P}}$\sss
is\dss contractible.\oss
A\qss \emph{split\dss polarized\sss subspace model},\oss
or\sss simply\sss a\qss \emph{split\dss model},\oss 
is\dss a morphism\sss
$N\qff \ttoo\qff M$\sss 
of\trs $\mathcal{P}{\nsp}\hat{\mathcal{S}}$
such\sss that\sss $N$\sss is\dss an object\sss of\trs
$\mathcal{P}{\nsp}$\dnsp.\oss
A split\sss model\sss can\sss be identified\sss with an object\sss
$M\off =\off(\trf V\fff,\qff H_{\dff -}\dff,\pff H_{\dff +}\trf)$ of\sss
$\mathcal{P}{\nsp}\hat{\mathcal{S}}$\sss
together\sss with a decomposition\sss
$V\off =\off V_{\dff -}\dff \oplus\trf V_{\dff +}$\nsp.\oss
We will\sss call\sss such a decomposition\sss a\qss \emph{splitting}\pss
of\dss $V$\dnsp.\oss
Under\sss this identification\sss
$N$\sss corresponds\sss to\sss 
$(\dff 0\fff,\qff V_{\dff -}\dff \oplus\dff H_{\dff -}\dff,\pff 
V_{\dff +}\dff \oplus\dff H_{\dff +}\trf)$
and\sss
$N\qff \ttoo\qff M$\sss
is\dss the unique morphism\sss from\sss $N$\sss to\sss $M$\nnsp.\oss 
A morphism of\dss split\sss models\dss is\dss a commutative diagram 
of\dss the form\vspace{-3pt}\vspace{-0.375pt}
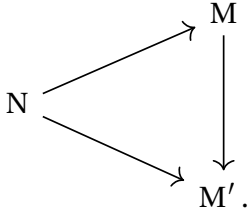
\[
\quad
\begin{tikzcd}[column sep=boom, row sep=hqboom]
&
M
\arrow[dd]
\\
N
\arrow[ru]
\arrow[rd]
&
\\
&
M\fff'\qff.
\end{tikzcd}
\]

\vspace{-12pt}\vspace{-3pt}\vspace{-0.375pt}
Therefore split\sss models are\sss the objects of\dss a\sss topological\sss category,\oss
which we denote by
$s\dff \hat{\mathcal{S}}$\nnsp.\oss
There\dss is\dss a canonical\dss forgetting functor\sss
$\phi\dff \colon\dff
s\dff \hat{\mathcal{S}}
\qff \ttoo\qff
\mathcal{P}{\nsp}\hat{\mathcal{S}}$\sss
taking\sss
$N\qff \ttoo\qff M$\sss
to\sss $M$\nnsp.\oss\vspace{-0.5pt}

\mypar{Theorem.}{split-contractible}
\emph{The classifying\sss space\sss 
$\num{s\dff \hat{\mathcal{S}}}$\sss
is\dss contractible.\oss}\vspace{-0.5pt}

\proof
Let\sss
$\rho\dff \colon\dff
s\dff \hat{\mathcal{S}}
\qff \ttoo\qff
\mathcal{P}$\sss
be\sss the functor\sss taking\sss
$N\qff \ttoo\qff M$\sss
to\sss $N$ 
and every\sss morphism of\dss
$s\dff\hat{\mathcal{S}}$\sss
to an\sss identity\sss morphism of\dss
$\mathcal{P}$\dnsp.\oss
Let\sss
$i\dff \colon\dff
\mathcal{P}{\nsp}
\qff \ttoo\qff
s\dff \hat{\mathcal{S}}$\sss
be\sss the functor\sss taking an object\sss
$N$ of\dss
$\mathcal{P}$\sss
to\sss the split\sss model\sss
$N\qff \ttoo\qff N$\nnsp.\oss
Clearly,\pss $\rho\dff \circ\dff i$\sss
is\dss the identity\sss functor.\oss
For every split\sss model\sss
$N\qff \ttoo\qff M$\sss 
there\dss is\dss a canonical\dss morphism\vspace{-3pt}\vspace{-0.375pt}
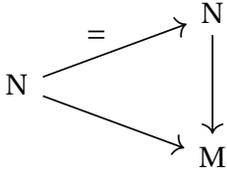
\[
\quad
\begin{tikzcd}[column sep=boom, row sep=shqboom]
&
N
\arrow[dd]
\\
N
\arrow[ru, "\dis ="]
\arrow[rd]
&
\\
&
M
\end{tikzcd}
\]

\vspace{-12pt}\vspace{-3pt}\vspace{-0.375pt}
of\dss $s\dff \hat{\mathcal{S}}$\dnsp.\oss
These morphisms define a natural\sss transformation\sss
from $i\dff \circ\dff \rho$\sss to\sss the identity\sss
functor of\dss $s\dff \hat{\mathcal{S}}$\dnsp,\oss
and\sss hence\sss
$\num{i}\dff \circ\dff \num{\rho}
\off =\off
\num{i\dff \circ\dff \rho}$\sss
is\dss homotopic\sss to\sss the identity.\oss
Since\sss
$\num{\rho}\dff \circ\dff \num{i}
\off =\off
\num{\rho\dff \circ\dff i}$\sss
is\dss equal\sss to\sss the identity,\oss
it\sss follows\sss that\sss
$\num{s\dff \hat{\mathcal{S}}}$\sss
is\dss homotopy equivalent\sss to 
$\num{\mathcal{P}}$\nnsp.\oss
Since\sss the\sss latter space\dss is\dss contractible,\oss
this proves\sss the\sss theorem.\oss  \eproof

\myuppar{Polarizations and\sss partial\sss orders.}
The main\sss technical\sss advantage of\dss working\sss with\sss
polarized\sss models\dss is\dss that\sss the categories of\dss models\sss
$\mathcal{PE}\hat{\mathcal{O}}$\dnsp,\qss 
$\mathcal{P}\hat{\mathcal{O}}$\dnsp,\oss
and\sss 
$\mathcal{P}{\nsp}\hat{\mathcal{S}}$\dss
are all\sss defined\sss in\sss terms of\dss partial\sss orders,\oss
in\sss contrast\sss with\sss the category\sss $\hat{\mathcal{S}}$\dnsp,\oss
which\dss is\dss not.\oss
All\sss these partial\sss orders obviously\sss have free equalities
in\sss the sense of\trs Section\qss \ref{topological-simplicial-complexes},\oss
and\sss hence we can consider\sss
$\mathcal{PE}\hat{\mathcal{O}}$\dnsp,\qss 
$\mathcal{P}\hat{\mathcal{O}}$\dnsp,\oss
and\sss 
$\mathcal{P}{\nsp}\hat{\mathcal{S}}$\dss
as\sss topological\sss simplicial\sss complexes.\oss

\myuppar{Finite-polarized\sss operators.}
For $n\qff \in\qff \nnn$\dss let\sss
$\hat{F}^{\dff n}$\sss be\sss the space of\dss
self-adjoint\sss operators\sss $A$\sss in\sss $H$\sss
such\sss that\sss $\norm{A}\off =\off 1$\nnsp,\oss
the essential\sss spectrum of\sss $A$\sss consists of\dss two points\sss
$-\qff 1\fff,\qff 1$\nnsp,\oss
and\sss the rank of\dss the spectral\sss projection\sss
$P_{\trf (\dff -\qff 1\fff,\qff 1\trf)}\trf(\trf A\trf)$\sss
is\dss $\leq\qff n$\nnsp.\oss
Clearly,\pss 
$\hat{F}^{\trf 0}\off \subset\off
\hat{F}^{\fff 1}\off \subset\off
\hat{F}^{\dff 2}\off \subset\off
\dots
\off$.\oss
Let\sss $\hat{F}\ffin$\sss be\sss the union\sss of\dss spaces\sss
$\hat{F}^{\dff n}$
with\sss the direct\sss limit\sss topology.\oss
While\sss the\sss topology of\dss $\hat{F}\ffin$\sss
is\dss different\sss from\sss the\sss one\sss 
induced\sss from\sss $\hat{\mathcal{F}}$\dnsp,\oss
the inclusion\sss
$\hat{F}\ffin\qff \ttoo\qff\hat{\mathcal{F}}$\sss
is\dss continuous.\oss

The space\sss $\hat{F}\ffin$\sss is\dss closely\sss related\sss to\sss $-\qff U\ffin$\dnsp.\oss
In\sss particular,\oss
the\sss formula\sss
$A\off \longmapsto\off \exp\qff(\trf \pi\fff i\qff A\trf)$\sss
defines a canonical\sss map\sss
$\exp\ffin\dff \colon\dff \hat{F}\ffin\qff \ttoo\qff -\qff U\ffin$\dnsp.\oss

\mypar{Theorem.}{polarized-homeo}
\emph{There\dss is\dss a canonical\dss homeomorphism\dss
$\mathcal{P}h\dff \colon\dff
\num{\mathcal{P}{\nsp}\hat{\mathcal{S}}}
\qff \ttoo\qff 
\hat{F}\ffin${\nsp}.\oss}

\proof
The proof\dss follows\sss the proof\dss of\trs Theorem\qss \ref{harris-h},\oss
and we will\sss use\sss the notations from\sss the\sss latter.\oss
In\sss particular,\oss we will\sss use\sss the version of\dss the standard\sss
simplices\sss $\Delta^n$\sss from\sss that\sss proof.\oss
An $n$\dnsp-simplex\sss $\overline{\sigma}$\sss of\dss $\mathcal{P}{\nsp}\hat{\mathcal{S}}$\sss
can\sss be identified\sss with a pair\sss consisting of\dss an $n$\dnsp-simplex $\sigma$
of\sss $\hat{\mathcal{S}}$\sss defined\sss by a 
sequence of\dss morphisms\sss\vspace{1.5pt}
\[
\quad
V_{\dff 0}\qff \ttoo\qff
V_{\dff 1}\qff \ttoo\qff
\ldots\qff \ttoo\qff
V_{\dff n}
\]

\vspace{-12pt}\vspace{1.5pt}\vspace{-0.375pt}
of\dss $\hat{\mathcal{S}}$
and\sss a\sss splitting\dss 
$H\dff \ominus\dff V_{\dff n}
\off =\off
H_{\dff -}\dff \oplus\dff H_{\dff +}$\nsp.\oss
Let\dss
$P_{\dff -}^{\dff i}$,\pss
$P_{\dff +}^{\dff i}$,\pss
$P_{\dff i}$\dss
be\sss the same orthogonal\sss projections as\sss in\sss the proof\dss 
of\trs Theorem\qss \ref{harris-h},\oss
and\sss let\sss $P_{\dff -}\dff,\pff P_{\dff +}$\sss be\sss
the orthogonal\sss projections of\dss $H$\sss onto\sss
$H_{\dff -}\dff,\pff H_{\dff +}$\sss respectively.\oss
For\sss $u\qff \in\qff \Delta^n$\sss let\sss
$S\trf(\trf \sigma\fff,\qff u \trf)$\sss
be\sss the operator\qss (\ref{simplex-operator}),\oss 
and\sss let\vspace{2.5pt}
\begin{equation}
\label{simplex-pol-operator}
\quad
\overline{S}\qff(\qff \overline{\sigma}\dff,\qff u\trf)
\off =\off
-\qff P_{\dff -}
\off +\off
S\trf(\trf \sigma\fff,\qff u \trf)
\off +\off
P_{\dff +}
\pff.
\end{equation}

\vspace{-12pt}\vspace{2.5pt}
Then\sss 
$\overline{S}\qff(\qff \overline{\sigma}\dff,\qff u\trf)
\qff \in\qff 
\hat{F}\ffin$\dnsp.\oss
The resulting\sss maps\dss
$\overline{S}\trf \colon\trf
\mathcal{P}{\nsp}\hat{\mathcal{S}}_{\dff n}\dff \times\dff \Delta^n
\qff \ttoo\qff
\hat{F}\ffin$\dss
agree with\sss the equivalence relation defining\sss
$\num{\hat{\mathcal{S}}}$\sss
and\dss hence induce a map\sss\vspace{1.5pt}
\[
\quad 
\mathcal{P}h\dff \colon\dff
\num{\mathcal{P}{\nsp}\hat{\mathcal{S}}}
\qff \ttoo\qff 
\hat{F}\ffin
\pff.
\]

\vspace{-12pt}\vspace{1.5pt}
The spectral\sss theorem\sss implies\sss
that $\mathcal{P}h$\sss is\dss a\sss bijection,\oss and since
$\num{\mathcal{P}{\nsp}\hat{\mathcal{S}}}$ and\sss $\hat{F}\ffin$ 
have\sss the direct\sss limit\sss
topology,\pss $\mathcal{P}h$ is\dss a\sss homeomorphism.\oss \eproof

\mypar{Theorem.}{exp}
\emph{The\sss following\sss diagram\dss is\dss commutative.}\vspace{-2.5pt}
\[
\quad
\begin{tikzcd}[column sep=sboom, row sep=sboom]
\protect{\num{\mathcal{P}{\nsp}\hat{\mathcal{S}}}}
\arrow[d, "\dis \mathcal{P}h\dff"']
\arrow[r, "\dis \protect{\num{\hat{\pi}}}"]
&
\protect{\num{\hat{\mathcal{S}}}}
\arrow[d, "\dis \dff h"]
\\
\hat{F}\ffin
\arrow[r, "\dis \exp\ffin"]
&
-\qff U\ffin
\end{tikzcd}
\]

\vspace{-12pt}\vspace{2.5pt}
\emph{Here\sss $\mathcal{P}h$\sss and\trs $h$\sss are\sss homeomorphisms\sss from\trs
Theorems\qss \ref{polarized-homeo}\qss and\pss \ref{harris-h}.\oss}

\proof
We will\sss use\sss the notations from\sss the proof\dss of\trs
Theorem\qss \ref{polarized-homeo}.\oss
Clearly,\oss the operator\sss\vspace{0.75pt}
\[
\quad
\exp\qff\bigl(\trf -\qff \pi\fff i\qff P_{\dff -}\qff +\qff \pi\fff i\qff P_{\dff +}\trf\bigr)
\off =\off
\exp\qff\bigl(\trf -\qff \pi\fff i\qff P_{\dff -}\trf\bigr)
\qff
\exp\qff\bigl(\trf \pi\fff i\qff P_{\dff +}\trf\bigr)
\]

\vspace{-12pt}\vspace{0.75pt}
is\dss equal\sss to\sss $-\qff \id_{\trf H}$\sss on\sss $H\dff \ominus\dff V_{\dff n}$\nsp,\oss
and\sss hence\vspace{0.75pt}
\[
\quad
\exp\qff\bigl(\pff \overline{S}\qff(\qff \overline{\sigma}\dff,\qff u\trf)\trf\bigr)
\off =\off
\exp\qff\bigl(\trf S\trf(\trf \sigma\fff,\qff u \trf)\trf\bigr)
\pff.
\]

\vspace{-12pt}\vspace{0.75pt}
Now\sss the\sss theorem\sss follows from\sss the definitions of\dss the maps involved.\oss  \eproof

\myuppar{Comparing\dss the spaces\sss 
$\hat{\mathcal{F}}$\sss and\sss $\hat{F}\ffin$\dnsp.}
By combining\trs Theorems\qss \ref{operators-categories},\pss
\ref{forgetting-equivalence},\pss and\qss \ref{polarized-homeo}\qss
we see\sss that\sss there exists a canonical\sss
homotopy\sss equivalence between\sss 
$\hat{\mathcal{F}}$\sss and\dss $\hat{F}\ffin$\dnsp.\oss
At\sss the same\sss time\sss 
$\hat{F}\ffin\qff \subset\pff \hat{\mathcal{F}}$\sss
and\sss the inclusion\sss map\sss
$i\dff \colon\dff
\hat{F}\ffin
\qff \ttoo\qff
\hat{\mathcal{F}}$\sss
is\dss continuous.\oss
This\sss leads\sss to\sss the unavoidable question\trs if\trs the inclusion\sss map\sss
$i$\sss is\dss a homotopy equivalence.\oss
The answer\dss is\dss positive,\oss but\dss for\sss the proof\dss one needs\sss to
introduce a\qss ``finite''\qss analogue of\dss the category\sss $\hat{\mathcal{E}}$\dnsp.\oss

Let\sss $\hat{\mathcal{E}}\ffin$ be\sss the\sss topological\sss category\sss
having as objects enhanced operators\sss $(\trf A\fff,\qff \varepsilon\trf)$\sss
such\sss that\sss $A\qff \in\qff \hat{F}\ffin$\sss
and\sss with\sss morphisms defined exactly as for\sss $\hat{\mathcal{E}}$\dnsp.\oss
The\sss topology\dss is\dss defined\sss by\sss the direct\sss limit\sss
topology of\sss $\hat{F}\ffin$\sss and\sss the discrete\sss topology on\sss
the controlling\sss parameters $\varepsilon$\nnsp.\oss
As a category,\oss $\hat{\mathcal{E}}\ffin$\sss is\dss a\sss full\sss subcategory of\dss
$\hat{\mathcal{E}}$\dnsp,\oss but\sss the\sss topology\dss is\dss not\sss the one induced\sss
from\sss $\hat{\mathcal{E}}$\dnsp.\oss
Still,\oss the inclusion\sss
$\hat{\mathcal{E}}\ffin
\qff \ttoo\qff 
\hat{\mathcal{E}}$\sss
is\dss a continuous functor.\oss
Let\sss
\[
\quad 
\mathcal{P}\dff\hat{\psi}\ffin\dff \colon\dff 
\hat{\mathcal{E}}\ffin
\qff \ttoo\qff 
\mathcal{PE}\hat{\mathcal{O}}
\]

\vspace{-12pt}
be\sss the composition of\dss this inclusion\sss with\sss 
$\mathcal{P}\dff\hat{\psi}\dff \colon\dff 
\hat{\mathcal{E}}
\qff \ttoo\qff 
\mathcal{PE}\hat{\mathcal{O}}$\dnsp,\oss
and\sss
\[
\quad 
\hat{\varphi}\ffin\dff \colon\dff
\hat{\mathcal{E}}\ffin
\qff \ttoo\qff 
\hat{F}\ffin
\]

\vspace{-12pt}\vspace{-0.125pt}
be\sss the obvious\sss forgetting\sss functor.\oss

\mypar{Theorem.}{finite-to-enhanced-models}
\emph{The\dss maps}\oss\vspace{0.5pt}
\[
\quad 
\num{\hat{\varphi}\ffin}\qff \colon\qff
\num{\hat{\mathcal{E}}\ffin}
\qff \ttoo\qff 
\num{\hat{F}\ffin}
\off =\dff\off
\hat{F}\ffin
\pff,\qquad
\num{\mathcal{P}\dff\hat{\psi}\ffin}\qff \colon\qff
\num{\hat{\mathcal{E}}\ffin}
\qff \ttoo\qff 
\num{\mathcal{P}\mathcal{E}\hat{\mathcal{O}}}
\pff,
\]

\vspace{-12pt}\vspace{0.5pt}
\emph{and\dss the\sss inclusion\sss
$\num{\hat{\mathcal{E}}\ffin}
\qff \ttoo\qff 
\num{\hat{\mathcal{E}}}$\sss
are\sss homotopy\sss equivalences.\oss}

\proof
The\sss proof\trs for\sss $\num{\hat{\varphi}\ffin}$\sss
is\dss the same as\sss the proof\dss of\dss
Theorem\qss \ref{forgetting-enhancement}.\oss
As of\dss $\num{\mathcal{P}\dff\hat{\psi}\ffin}$\nnsp,\oss
the proofs\sss of\trs Theorems\qss \ref{to-enhanced-models}\qss and\qss \ref{forgetting-equivalence}\qss
will\dss apply\sss once we prove\sss that\dss the spaces\vspace{1.5pt}
\[
\quad 
\hat{F}^{\dff \mathrm{fin}\qff \inv}\qff[\dff -\qff \varepsilon\fff,\fff \varepsilon\trf]
\off =\off\dff
\hat{F}\ffin
\qff \cap\pff
\hat{\mathcal{F}}^{\dff \inv}\trf[\dff -\qff \varepsilon\fff,\fff \varepsilon\trf]
\]

\vspace{-12pt}\vspace{1.5pt}
are contractible for every $\varepsilon\qff \in\qff (\trf 0\fff,\qff 1\dff)$\nnsp.\oss
The spectral\sss deformation contracting\sss
$[\dff -\qff 1\fff,\qff -\qff \varepsilon\trf]
\qff \cup\qff
[\trf \varepsilon\fff,\qff 1\trf]$\sss
to\sss $\{\trf -\qff 1\fff,\qff 1\trf\}$\sss
deforms\sss
$\hat{F}^{\dff \mathrm{fin}\qff \inv}\qff[\dff -\qff \varepsilon\fff,\fff \varepsilon\trf]$\sss
to\sss the subspace of\dss self-adjoint\sss operators $A$ such\sss that\sss
$\sigma\trf(\trf A\trf)
\off =\off 
\{\trf -\qff 1\fff,\qff 1\trf\}$\sss
and\sss both $-\qff 1$ and\sss $1$\sss are eigenvalues of\dss infinite multiplicity.\oss
Such operators\sss $A$\sss are determined\sss by\sss decompositions of\sss $H$\sss
into\sss the sum of\dss two eigenspaces,\oss
and\sss hence\sss the space of\dss such operators\dss is\dss homeomorphic\sss to\sss
the space of\dss polarizations of\sss $H$\nnsp.\oss
Since\sss the\sss latter\dss is\dss contractible,\oss
the spaces\sss
$\hat{F}^{\dff \mathrm{fin}\qff \inv}\qff[\dff -\qff \varepsilon\fff,\fff \varepsilon\trf]$\sss
are also contractible.\oss
The claim\sss about\sss $\num{\mathcal{P}\dff\hat{\psi}\ffin}$\sss follows.\oss
Now\sss the claim about\sss the inclusion\sss follows from\trs
Theorem\qss \ref{forgetting-equivalence}.\oss  \eproof

\mypar{Theorem.}{enhanced-and-finite}
\emph{The composition}\vspace{1.5pt}\vspace{-0.427pt}
\[
\quad
\hat{\eta}
\dff \colon\dff
\num{\hat{\mathcal{E}}\ffin}
\qff \ttoo\qff
\num{\mathcal{P}{\nsp}\hat{\mathcal{S}}}
\qff \ttoo\qff 
\hat{F}\ffin
\pff,
\]

\vspace{-12pt}\vspace{1.5pt}\vspace{-0.427pt}
\emph{where\sss the second\sss map\dss is\dss the homeomorphism\sss $\mathcal{P}h$\nnsp,\oss
is\dss homotopic\sss to\sss 
$\num{\hat{\varphi}\ffin}\qff \colon\qff
\num{\hat{\mathcal{E}}\ffin}
\qff \ttoo\qff 
\hat{F}\ffin$\nsp\dnsp.}

\proof
Let\sss $(\trf A\fff,\qff \varepsilon\trf)$\sss be an object\sss of\sss
$\hat{\mathcal{E}}\ffin$\dnsp.\oss
Then\sss $A\qff \in\pff \hat{F}\ffin$\sss
and\sss $\varepsilon\qff <\qff 1$\nnsp.\oss
In\sss the notations of\dss proof\dss of\trs 
Theorem\qss \ref{polarized-homeo}\qss
the operator\sss $A$\sss has\sss the form\vspace{1.5pt}\vspace{-0.427pt}
\[
\quad
-\qff P_{\dff -}
\off -\off
u_{\dff n}\qff P_{\dff -}^{\dff n}
\off -\off
\ldots
\off -\off
u_{\dff 2}\qff P_{\dff -}^{\dff 2}
\off -\off
u_{\dff 1}\qff P_{\dff -}^{\dff 1}
\off +\off
u_{\dff 1}\qff P_{\dff +}^{\dff 1}
\off +\off
u_{\dff 1}\qff P_{\dff +}^{\dff 2}
\off +\off
\ldots
\off +\off
u_{\dff n}\qff P_{\dff +}^{\dff n}
\off +\off
P_{\dff +}
\]

\vspace{-12pt}\vspace{1.5pt}\vspace{-0.427pt}
with\sss $u_{\dff i}\off \neq\off \varepsilon$\sss for every\sss $i$\nnsp.\oss
Let\sss $a$\sss be\sss equal\dss to\sss the maximal\dss $i$\sss such\sss that\sss
$u_{\dff i}\qff <\qff \varepsilon$\sss and\sss to\sss $0$\sss
if\trs there\dss is\dss no such $i$\nnsp.\oss
In\sss the notations of\dss proofs of\trs 
Theorems\qss \ref{polarized-homeo}\qss and\qss \ref{harris-h}\qss
the image of\dss $(\trf A\fff,\qff \varepsilon\trf)$\sss
in\sss $\mathcal{P}{\nsp}\hat{\mathcal{S}}$\dss
is\dss equal\dss to\sss the\sss triple\sss
$(\trf V\fff,\pff K_{\dff -}\dff,\pff K_{\dff +}\trf)$\nnsp,\oss
where\vspace{0pt}\vspace{-0.427pt}
\[
\quad
V\off =\off V_{\fff a}
\pff,
\]

\vspace{-39pt}\vspace{-0.427pt}
\[
\quad
K_{\dff -}
\off =\off
H_{\dff -}
\qff \oplus\qff
U_{\dff -}^{\dff n}
\qff \oplus\qff
\ldots
\qff \oplus\qff
U_{\dff -}^{\dff a\dff +\dff 1}
\pff,\quad
\mbox{and}
\]

\vspace{-36pt}\vspace{-0.427pt}
\[
\quad
K_{\dff +}
\off =\off
U_{\dff +}^{\dff a\dff +\dff 1}
\qff \oplus\qff
\ldots
\qff \oplus\qff
U_{\dff -}^{\dff n}
\qff \oplus\qff
H_{\dff +}
\pff.
\] 

\vspace{-12pt}\vspace{1.5pt}\vspace{-0.427pt}
This\sss object\sss of\sss 
$\mathcal{P}{\nsp}\hat{\mathcal{S}}$\dss 
can\sss be considered as a $0$\dnsp-simplex of\sss 
$\mathcal{P}{\nsp}\hat{\mathcal{S}}$\dss
and\sss hence as a point\sss of\dss 
$\num{\mathcal{P}{\nsp}\hat{\mathcal{S}}}$\nnsp.\oss
The homeomorphism\dss
$\num{\mathcal{P}{\nsp}\hat{\mathcal{S}}}
\qff \ttoo\qff 
\hat{F}\ffin$\dss
takes\sss this\sss point\sss to\sss the operator\vspace{1.5pt}\vspace{-0.427pt}
\[
\quad
\hat{\eta}\trf (\trf A\fff,\qff \varepsilon\trf)
\off =\off
-\qff P_{\dff -}
\off -\off
P_{\dff -}^{\dff n}
\off -\off
\ldots
\off -\off
P_{\dff -}^{\dff a\dff +\dff 1}
\off +\off
P_{\dff +}^{\dff a\dff +\dff 1}
\off +\off
\ldots
\off +\off
P_{\dff +}^{\dff n}
\off +\off
P_{\dff +}
\pff.
\]

\vspace{-12pt}\vspace{1.5pt}\vspace{-0.427pt}
The map\sss
$(\trf A\fff,\qff \varepsilon\trf)
\off \longmapsto\off
\hat{\eta}\trf (\trf A\fff,\qff \varepsilon\trf)$\sss
replaces each eigenvalue of\sss $A$\sss by\sss $-\qff 1$\nnsp,\qss $0$\nnsp,\oss
or\sss $1$\sss and\sss preserves\sss the signs of\dss 
eigenvalues not\sss replaced\sss by\sss $0$\nnsp.\oss
Also,\pss $\hat{\eta}\trf (\trf A\fff,\qff \varepsilon\trf)$\sss
is\dss equal\sss to\sss $A$\sss on\sss the subspace\sss
$H_{\dff -}\dff \oplus\trf H_{\dff +}$\sss 
of\dss finite codimension\sss in $H$\nnsp.\oss
This implies,\oss in\sss particular,\oss
that\sss the\sss linear\sss path connecting\sss $A$\sss
with\sss $\hat{\eta}\trf (\trf A\fff,\qff \varepsilon\trf)$\sss 
is\dss contained\sss in\sss $\hat{F}\ffin$\dnsp.\oss

These\sss linear\sss paths define a homotopy 
of\dss the space of\dss objects of\dss
$\hat{\mathcal{E}}\ffin$\dnsp.\oss
In order\sss to extend\sss this homotopy\sss to\sss 
$\num{\hat{\mathcal{E}}\ffin}$\nnsp,\oss
let\sss us\sss consider\sss $\hat{\mathcal{E}}\ffin$ as 
a\sss topological\sss simplicial\sss complex.\oss
Clearly,\oss it\dss is\dss defined\sss by 
an ordered space with\sss free equalities.\oss
Hence\trs Corollary\qss \ref{free-full-realizations}\qss 
implies\sss that\sss
$\num{\hat{\mathcal{E}}\ffin}
\off =\off 
\bbnum{\hat{\mathcal{E}}\ffin}$\nnsp,\oss
and we can\sss represent\sss points of\dss 
$\num{\hat{\mathcal{E}}\ffin}$\sss
by\sss weighted sums of\dss the form\vspace{1.5pt}\vspace{-0.427pt}
\[
\quad
t_{\trf 0}\trf (\trf A\fff,\qff \varepsilon_{\trf 0}\trf)
\pff +\pff
t_{\trf 1}\trf (\trf A\fff,\qff \varepsilon_{\dff 1}\trf)
\pff +\pff
\ldots
\pff +\pff
t_{\dff n}\trf (\trf A\fff,\qff \varepsilon_{\dff n}\trf)
\pff,
\]

\vspace{-12pt}\vspace{1.5pt}\vspace{-0.427pt}
where\sss $t_{\dff i}\qff \geq\qff 0$\dss for\sss every $i$\nnsp,\oss
$t_{\trf 0}\qff +\qff
t_{\trf 1}\qff +\qff
\ldots\qff +\qff
t_{\dff n}
\off =\off
1$\nnsp,\oss
and\dss
$\varepsilon_{\trf 0}\qff <\qff
\varepsilon_{\dff 1}\qff <\qff
\ldots\qff <\qff
\varepsilon_{\dff n}$\nsp.\oss
The map\sss $\hat{\eta}$\sss takes\sss
the point\sss represented\sss by\sss this sum\dss 
to\sss the operator\vspace{1.5pt}\vspace{0.35pt}
\[
\quad
t_{\trf 0}\trf \hat{\eta}\trf (\trf A\fff,\qff \varepsilon_{\trf 0}\trf)
\pff +\pff
t_{\trf 1}\trf \hat{\eta}\trf (\trf A\fff,\qff \varepsilon_{\dff 1}\trf)
\pff +\pff
\ldots
\pff +\pff
t_{\dff n}\trf \hat{\eta}\trf (\trf A\fff,\qff \varepsilon_{\dff n}\trf)
\pff,
\]

\vspace{-12pt}\vspace{1.5pt}\vspace{0.35pt}
which\dss is\dss equal\dss 
to\sss $A$\sss on\sss 
$H_{\dff -}\dff \oplus\trf H_{\dff +}$\nsp.\oss
Therefore\sss the\sss linear\sss path connecting\sss it\sss with\sss $A$\sss
is\dss contained\sss in\sss $\hat{F}\ffin$\dnsp.\oss
These\sss linear\sss paths 
define a homotopy\sss between\sss $\hat{\eta}$\sss
and\sss $\num{\hat{\varphi}\ffin}$\nnsp.\oss  \eproof

\mypar{Theorem.}{inclusion-is-he}
\emph{The composition}\vspace{1.5pt}\vspace{0.35pt}
\[
\quad
i\trf \colon\trf
\hat{F}\ffin
\qff \ttoo\qff 
\hat{\mathcal{F}}
\qff \ttoo\qff
\num{\mathcal{P}{\nsp}\hat{\mathcal{S}}}
\qff \ttoo\qff 
\hat{F}\ffin
\pff,
\]

\vspace{-12pt}\vspace{1.5pt}\vspace{0.35pt}
\emph{where\sss the\sss first\sss map\dss is\dss the inclusion,\oss 
the second\sss map\dss is\dss the canonical\dss homotopy equivalence,\oss 
and\sss the\sss last\sss one\dss is\dss the homeomorphism\sss $\mathcal{P}h$\nnsp,\oss 
is\trs homotopic\sss to\sss the\sss identity.\oss}

\proof
The map\sss
$\hat{\mathcal{F}}
\qff \ttoo\qff
\num{\mathcal{P}{\nsp}\hat{\mathcal{S}}}$\sss
results from\sss the diagram\sss
$\hat{\mathcal{F}}
\off \longleftarrow\off
\num{\hat{\mathcal{E}}}
\qff \ttoo\qff
\num{\mathcal{P}{\nsp}\hat{\mathcal{S}}}$\dss
by\sss inverting\sss the\sss left\sss arrow up\sss to homotopy.\oss
The corresponding\qss ``finite''\qss diagram\dss is\dss\vspace{1.5pt}\vspace{0.35pt}
\[
\quad
\hat{F}\ffin
\off \longleftarrow\off
\num{\hat{\mathcal{E}}\ffin}
\qff \ttoo\qff
\num{\mathcal{P}{\nsp}\hat{\mathcal{S}}}
\pff,
\]

\vspace{-12pt}\vspace{1.5pt}\vspace{0.35pt}
where\sss the\sss left\sss arrow\dss is\dss $\num{\hat{\varphi}\ffin}$\nnsp.\oss
Since\sss $\num{\hat{\varphi}\ffin}$\sss is\dss a homotopy equivalence by\trs
Theorem\qss \ref{finite-to-enhanced-models},\oss 
this\sss arrow\dss is\dss invertible up\sss to homotopy and\dss
this diagram\dss leads\sss to a map\sss
$\hat{F}\ffin\qff \ttoo\qff \num{\mathcal{P}{\nsp}\hat{\mathcal{S}}}$\dss
homotopic\sss to\sss the composition\dss
$\hat{F}\ffin
\qff \ttoo\qff 
\hat{\mathcal{F}}
\qff \ttoo\qff
\num{\mathcal{P}{\nsp}\hat{\mathcal{S}}}$\sss
and such\sss that\sss its composition\sss with\sss
$\num{\hat{\varphi}\ffin}\dff \colon\dff
\num{\hat{\mathcal{E}}\ffin}
\qff \ttoo\qff
\hat{F}\ffin$\sss
is\dss homotopic\sss to\sss the canonical\sss map\sss
$\num{\hat{\mathcal{E}}\ffin}
\qff \ttoo\qff
\num{\mathcal{P}{\nsp}\hat{\mathcal{S}}}$\nnsp.\oss
It\sss follows\sss that\sss $i\trf \circ\trf \num{\hat{\varphi}\ffin}$\sss
is\dss homotopic\sss to\sss the map\sss $\hat{\eta}$\sss from\trs
Theorem\qss \ref{enhanced-and-finite}\qss and\sss hence\sss to\sss
$\num{\hat{\varphi}\ffin}$\nnsp.\oss
Since\sss $\num{\hat{\varphi}\ffin}$\sss is\dss a homotopy equivalence,\oss
it\sss follows\sss that\sss $i$\sss is\dss homotopic\sss to\sss
the identity.\oss  \eproof

\mypar{Corollary.}{subspace-polarized}
\emph{The\sss inclusion\dss
$\hat{F}\ffin
\qff \ttoo\qff 
\hat{\mathcal{F}}$\sss
is\dss a\sss homotopy\sss equivalence.\oss}

\proof
Since\sss 
$\hat{\mathcal{F}}
\qff \ttoo\qff
\num{\mathcal{P}{\nsp}\hat{\mathcal{S}}}$\sss
is\dss a homotopy equivalence,\oss and\sss
$\num{\mathcal{P}{\nsp}\hat{\mathcal{S}}}
\qff \ttoo\qff 
\hat{F}\ffin$\sss
is\dss even a homeomorphism,\oss
Theorem\qss \ref{inclusion-is-he}\qss
shows\sss that\sss the inclusion\sss
$\hat{F}\ffin
\qff \ttoo\qff 
\hat{\mathcal{F}}$\sss
is\dss equivalent\sss to\sss
the identity\sss
$\id\dff \colon\dff
\hat{F}\ffin
\qff \ttoo\qff 
\hat{F}\ffin$\sss
in\sss the homotopy category.\oss  \eproof

\myuppar{Compactly-polarized\sss operators.}
Let\sss $\hat{F}\comp$\sss be\sss the subspace of\sss $\hat{\mathcal{F}}$\sss
consisting of\dss operators\sss
$A\qff \in\qff \hat{\mathcal{F}}$\sss such\sss that\sss
$\norm{A}\off =\off 1$\sss and\sss
the essential\sss spectrum of\sss $A$\sss is\dss equal\sss to\sss
$\{\qff -\qff 1\fff,\qff 1 \qff\}$\nnsp.\oss
In\sss particular,\oss the\sss topology\sss of\sss $\hat{F}\comp$\sss 
is\dss the norm\sss topology.\oss
The spaces\sss $\hat{F}\ffin$\sss and\sss $\hat{F}\comp$\sss
are related\sss in\sss the same way\sss as\sss the space of\dss
operators of\dss finite rank and\sss the space of\dss compact\sss operators.\oss
In\sss fact,\oss it\dss is\dss easy\sss to see\sss that\sss $\hat{F}\comp$\sss
is\dss the norm closure of\dss $\hat{F}\ffin$\dnsp.\oss
Replacing self-adjoint\sss operators $A$ by\sss skew-adjoint\sss operators\sss 
$i\dff A$\nnsp,\oss where\sss $i\off =\off \sqrt{-\qff 1}$\nnsp,\oss 
turns\sss the space\sss $\hat{F}\comp$\sss 
into\sss the space\sss $\hat{F}_{\fff *}$\sss playing a\sss key\sss role\sss in\sss
the work\sss of\qss Atiyah\dss and\trs Singer\qss \cite{as}.\oss
The\sss topology of\sss $\hat{F}\ffin$\sss is\dss different\sss
from\sss the\sss topology\sss induced\sss from\sss $\hat{F}\comp$\dnsp,\oss
but\sss the inclusion\sss 
$\hat{F}\ffin\qff \ttoo\qff \hat{F}\comp$\sss
is\dss continuous.\oss
Moreover,\oss this inclusion\dss is\dss a homotopy equivalence.\oss
This\sss follows\sss from\dss Corollary\qss \ref{subspace-polarized}\qss
and\sss the fact\sss that\sss $\hat{F}\comp$\sss
is\dss a deformation\sss retract\sss of\sss $\hat{\mathcal{F}}$\dnsp,\oss
proved\sss in\qss \cite{as},\oss Section\qss 2.

\mysection{Restricted\pss Grassmannians}{restricted-grassmannians}

\myuppar{Restricted\dss Grassmannians.}
Suppose\sss that\sss a polarization\dss 
$H\off =\off K_{\dff -}\dff \oplus\dff K_{\dff +}$\sss
of\dss $H$\sss is\dss fixed.\oss
Restricted\dss Grassmannians\sss are spaces of\dss closed\sss subspaces of\sss $H$\sss
which are\qss ``comparable''\qss in size with\sss $K_{\dff -}$\sss in a sense\sss to be specified,\oss
which depends on\sss the problem at\sss hand.\oss
Several\sss versions are discussed\dss by\trs 
Pressley\sss and\dss Segal\qss \cite{ps},\oss Chapter\qss 7.\oss\vspace{-0.1pt}

We will\sss need\sss only\sss two versions of\dss this notion.\oss
The first\sss version\dss is\dss
the space\sss $\gr$\sss of\dss subspaces\sss $K$\dss
\emph{commensurable}\pss with\sss $K_{\dff -}$\nsp,\oss
i.e.\qss such\sss that\sss the intersection\sss
$K\dff \cap\dff K_{\dff -}$\sss is\dss closed and\sss has finite codimension\sss
in\sss both $K$ and\sss $K_{\dff -}$\nsp.\oss
The\sss topology of\dss $\gr$\sss is\dss induced\sss by\sss the 
norm\sss topology\sss of\dss orthogonal\sss projections\sss
$H\qff \ttoo\qff K$\nnsp.\oss
The\dss Grassmannian\sss $\gr$\sss is\dss mentioned\sss
under\sss the name\sss
$\gr_{\dff 1}\dff(\trf H\trf)$\sss 
by\dss Pressley\sss and\dss Segal\qss \cite{ps},\oss Section\qss 7.2,\oss
but\trs is\dss not\sss discussed\sss in any details.\oss

Having\sss applications\sss to differential\sss operators in\sss mind,\oss
we will\sss define\sss the second version\sss in a more general\sss context\dss
than\sss the usual\sss one.\oss
Suppose\sss that\sss a presentation of\dss $H$\sss 
as a direct\sss sum\vspace{3pt}\vspace{-0.1pt}
\begin{equation}
\label{spectral-decomposition}
\quad
H
\off =\off
\bigoplus\nolimits_{\dff n\qff \in\qff \zzz}\qff H_{\dff n}
\end{equation}

\vspace{-12pt}\vspace{3pt}\vspace{-0.1pt}
of\dss finitely dimensional\sss subspaces $H_{\dff n}$\nsp,\qss 
$n\qff \in\qff \zzz$\qss
is\trs fixed.\oss
For\sss $a\fff,\qff b\qff \in\qff \zzz$\nnsp,\qss 
$a\qff \leq\qff b$\dss let\vspace{3pt}\vspace{-0.1pt}
\[
\quad
H_{\dff \geq\dff a}
\off =\off
\bigoplus\nolimits_{\dff n\qff \geq\qff a}\qff H_{\dff n}
\qff,\quad
H_{\dff \leq\dff a}
\off =\off
\bigoplus\nolimits_{\dff n\qff \leq\qff a}\qff H_{\dff n}
\qff,\quad
\]

\vspace{-33pt}\vspace{-0.1pt}
\[
\quad
\mbox{and}\quad
H_{\dff [\dff a\fff,\dff b\dff]}
\off =\off
\bigoplus\nolimits_{\dff a\qff \leq\qff n\qff \leq\qff b}\qff H_{\dff n}
\qff,\quad
\]

\vspace{-12pt}\vspace{3pt}\vspace{-0.1pt}
The subspaces\sss
$H_{\dff >\dff a}$\sss and\sss $H_{\dff <\dff a}$\sss are defined similarly.\oss
We will\sss assume\sss that\sss
$K_{\dff -}\off =\off H_{\dff \leq\dff 0}$\nsp.\oss
Let\sss us\sss say\sss that\sss a subspace\sss $K\qff \subset\qff H$\sss is\qss
\emph{admissible}\pss with respect\sss to\sss the decomposition\qss 
(\ref{spectral-decomposition})\pss
if\qss\vspace{3pt}\vspace{-0.1pt}
\begin{equation}
\label{ab-subspaces}
\quad
H_{\dff \leq\dff a}
\off \subset\off
K
\off \subset\off
H_{\dff \leq\dff b}
\end{equation}

\vspace{-12pt}\vspace{3pt}\vspace{-0.1pt}
for some\sss
$a\fff,\qff b\qff \in\qff \zzz$\nnsp,\qss $a\qff \leq\qff b$\nnsp.\oss
Clearly,\oss every admissible subspace\dss is\dss closed.\oss
Let\sss $\gr\trf(\dff a\fff,\qff b\trf)$\sss be\sss the space of\dss
subspaces\sss $K$\sss such\sss that\qss (\ref{ab-subspaces})\sss holds,\oss
and\sss let\sss
$\gr\trf(\dff \infty\dff)$\sss be\sss the set\sss of\dss all\sss
subspaces admissible with respect\sss to\qss 
(\ref{spectral-decomposition}).\oss
Clearly,\pss $\gr\trf(\dff \infty\dff)$\sss
is\dss equal\sss to\sss the union of\dss 
$\gr\trf(\dff a\fff,\qff b\trf)$\sss
over all\sss $a\fff,\qff b$\sss as above,\oss
and we equip\sss $\gr\trf(\dff \infty\dff)$\sss
with\sss the corresponding direct\sss limit\sss topology.\oss

Up\sss to homeomorphism\sss $\gr\trf(\dff \infty\dff)$\sss
is\dss independent\sss 
from\sss the choice of\dss the decomposition\qss
(\ref{spectral-decomposition}),\oss
and\sss usually\dss is\dss defined\sss in\sss terms of\dss
a\sss basis of\dss $H$\nnsp.\oss
The assumption\sss
$K_{\dff -}\off =\off H_{\dff \leq\dff 0}$\sss
implies\sss that\sss
$\gr\trf(\dff \infty\dff)
\qff \subset\qff
\gr$\nnsp.\oss
The inclusion\sss map\sss
$\gr\trf(\dff \infty\dff)
\qff \ttoo\qff
\gr$\sss
is\dss continuous,\oss
but\sss the\sss topology of\dss the space\sss
$\gr\trf(\dff \infty\dff)$\sss
is\dss different\sss from\sss the\sss topology\sss
induced\sss from\sss $\gr$\nnsp.\oss

The\sss goal\sss of\dss this section\dss is\dss to prove\sss
that\sss the inclusion\sss
$\gr\trf(\dff \infty\dff)\qff \ttoo\qff \gr$\sss
is\dss a homotopy equivalence.\oss
While\sss the\dss Grassmannian\sss $\gr\trf(\dff \infty\dff)$\sss
is\dss a\sss well\sss known space,\oss
the\dss Grassmannian\sss $\gr$\sss is\trs less so.\oss
Still,\oss this result\sss should\sss be known,\oss 
but\sss the author\sss
failed\sss to find\sss it\sss  
in\sss the\sss literature.\oss

\mypar{Lemma.}{filtration-of-gr}
\emph{Let\qss $U_{\fff n}$\nsp,\dss $C_{\dff k}$
and\dss $D_{\dff m}$\sss be\sss the sets\sss
of\dss subspaces $K\qff \in\qff \gr$
such\sss that}
\[
\quad
K\trf \cap\qff H_{\dff >\dff n}
\off =\off
0\qff,\qff\quad
K\pff \subset\off H_{\dff \leq\dff k}
\qff,\quad
\mbox{\emph{and}}\qff\quad 
K
\off \supset\pff
H_{\dff \leq\dff m}
\]

\vspace{-12pt}
\emph{respectively.\oss
Then\sss $\gr$\sss is\dss equal\sss to\sss the union of\qss sets\dss $U_{\fff n}$\nsp,\oss
the sets\sss $U_{\fff n}$\sss are open\sss in\sss $\gr$\nnsp,\oss
the sets\sss $C_{\dff k}$\dss and\trs $D_{\dff m}$\sss
are closed\sss in\sss $\gr$\nnsp,\oss
and\dss $C_{\dff n}$\sss is\dss a deformation\sss 
retract\sss of\qss $U_{\dff n}$\sss
by a deformation\sss preserving\sss subsets\trs 
$U_{\fff n}\dff \cap\qff C_{\dff k}$\sss 
and\qss
$U_{\fff n}\dff \cap\qff D_{\dff m}$\sss
for every\sss $k\fff,\qff m$\sss
and\sss the subset\qss
$U_{\fff n}\dff \cap\qff \gr\trf(\dff \infty\dff)$\nnsp.\oss}

\proof
Let\sss us\sss begin\sss with a simple observation.\oss
Let\sss $E\qff \subset\qff H$\sss be a finitely dimensional\sss subspace.\oss
We claim\sss that\sss $E\trf \cap\trf H_{\dff >\dff n}\off =\off 0$\sss 
for sufficiently\sss large$n$\nnsp.\oss
If\dss this\dss is\dss not\sss the case,\oss 
then for every\sss natural\sss $n$\sss
there exists a unit\sss vector $v_{\dff n}$\sss belonging\sss to $E$\sss and\sss
to $H_{\dff >\dff n}$\nsp.\oss
Since $E$\sss is\dss finitely dimensional,\oss the sequence of\dss these vectors
contains a subsequence $v_{\dff i}$ converging\sss to some unit\sss vector $v$\nnsp.\oss
For every\sss natural\sss number\sss $m$\sss the vectors $v_{\dff i}$\sss belong\sss
to $H_{\dff > m}$\sss for sufficiently\sss large $i$ and\sss hence $v$ also
belongs to $H_{\dff > m}$\nsp.\oss
It\sss follows\sss that\sss $v\qff \in\qff H_{\dff > m}$\sss for every natural\sss $m$\nnsp.\oss
But\sss $v$\sss is\dss a unit\sss vector,\oss and\sss the intersection of\dss the subspaces
$H_{\dff > m}$ over all\sss natural\sss $m$\sss is\dss the zero subspace.\oss
The contradiction\sss proves our claim.\oss

Let\sss $K\qff \in\qff \gr$\nnsp.\oss
Then\sss $K$\sss 
is\dss commensurable with\sss $H_{\dff \leq\dff 0}$\sss
and\sss hence\sss $K\dff \cap\trf H_{\dff >\dff 0}$\sss is\dss a finitely 
dimensional\sss subspace.\oss
If\dss $n$\sss is\dss a sufficiently\sss large,\oss
then\sss by\sss the previous paragraph\vspace{0.625pt}
\[
\quad
K\trf \cap\trf H_{\dff >\dff n}
\off =\off
\bigl(\trf K\dff \cap\trf H_{\dff >\dff 0}\trf\bigr)
\trf \cap\trf 
H_{\dff >\dff n}
\off =\off
0\qff,
\]

\vspace{-12pt}\vspace{0.625pt}
i.e.\trs $K\qff \in\qff U_{\fff n}$\nsp.\oss
This proves\sss that\sss $\gr$\sss is\dss equal\sss 
to\sss the union of\dss sets\sss $U_{\fff n}$\nsp.\oss

If\dss $K\qff \in\qff U_{\fff n}$\nsp,\oss
then\sss the orthogonal\sss projection\sss 
$K\qff \ttoo\qff H_{\dff \leq\dff n}$\sss 
is\dss injective.\oss
We claim\sss that\sss the orthogonal\sss projection\dss 
$K\fff'
\qff \ttoo\qff
H_{\dff \leq\dff n}$\sss 
is\dss injective for every\sss $K\fff'$\sss
sufficiently close\sss to $K$\nnsp.\oss 
Indeed,\oss since $K$\sss is\dss commensurable with\sss
$H_{\dff \leq\dff 0}$\nsp,\oss
the image of\sss $K$\sss in\sss $H_{\dff \leq\dff n}$\sss
has finite codimension\sss in\sss $H_{\dff \leq\dff n}$\nsp.\oss
Let\sss $F$\sss be\sss the orthogonal\sss complement\sss of\dss 
this image in\sss $H_{\dff \leq\dff n}$\nsp.\oss
Then\sss the map\sss
$F\dff \oplus\dff K\qff \ttoo\qff H_{\dff \leq n}$\sss
equal\sss to\sss the identity\sss on\sss $F$\sss and\sss to\sss the projection
on\sss $K$\sss is\dss an\sss isomorphism.\oss
The open\sss mapping\sss theorem\sss implies\sss that\sss similar maps\sss
$F\dff \oplus\dff K\fff'\qff \ttoo\qff H_{\dff \leq\dff n}$\sss
are\sss isomorphisms for\sss $K\fff'$\sss sufficiently close\sss to $K$\nnsp.\oss
In\sss particular,\oss the projection\dss 
$K\fff'
\qff \ttoo\qff
H_{\dff \leq\dff n}$\sss 
is\dss injective for such\sss $K\fff'$\dnsp.\oss
This proves our claim.\oss
Clearly,\oss if\dss the projection\sss
$K\fff'
\qff \ttoo\qff
H_{\dff \leq\dff n}$\sss 
is\dss injective,\oss
then\sss
$K\fff'\trf \cap\trf H_{\dff >\dff n}
\off =\off
0$\sss
and\sss hence\sss
$K\qff \in\qff U_{\fff n}$\nsp.\oss
This proves\sss that\sss $U_{\fff n}$\sss is\dss open.\oss

Obviously,\oss the sets $C_{\dff n}$\sss and\sss $D_{\dff n}$ are closed.\oss
For\sss $K\qff \in\qff U_{\fff n}$\sss let\sss $\pi\dff K$\sss
be\sss the image of\dss the orthogonal\sss projection\sss
$K\qff \ttoo\qff H_{\dff \leq\dff n}$\nsp.\oss
Then\sss $\pi\dff K\qff \in\qff C_{\dff n}$\sss
and\sss $K$\sss is\dss the graph\sss of\dss a\sss linear map\sss
$l_{\qff K}\dff \colon\dff 
\pi\dff K
\qff \ttoo\qff 
H_{\dff >\dff n}$\nsp.\oss
The graphs of\dss the maps\sss $t\dff l_{\qff K}$\nsp,\pss
$t\qff \in\qff [\dff 0\fff,\qff 1\dff]$\sss
form a continuous path connecting\sss $K$\sss with\sss $\pi\dff K$\nnsp.\oss
If\dss $K\qff \in\qff C_{\dff n}$\nsp,\oss 
then\sss $l_{\qff K}\off =\off 0$\sss and\sss
this path\dss is\dss a constant\sss path.\oss
Clearly,\oss these paths continuously\sss depend\sss on\sss $K$\sss
and\sss define a deformation of\dss the identity\sss map of\dss $U_{\fff n}$\sss
to\sss the map\sss 
$K\off \longmapsto\off \pi\dff K$\nnsp.\oss
In\sss particular,\pss $C_{\dff n}$\sss
is\dss a deformation\sss retract\sss of\dss $U_{\fff n}$\nsp.\oss

If\dss $k\qff \leq\qff n$\nnsp,\oss
then\sss $C_{\dff k}\qff \subset\qff C_{\dff n}$\sss
and\sss hence\sss the above deformation\sss preserves\sss
$U_{\fff n}\dff \cap\qff C_{\dff k}$\nsp.\oss
If\dss $k\qff >\qff n$\sss
and\dss
$K\pff \subset\off H_{\dff \leq\dff k}$\nsp,\oss
then\sss the image of\sss $l_{\qff K}$\sss
is\dss contained\sss in\sss $H_{\trf [\dff n\dff +\dff 1\fff,\dff k\trf]}$\sss
and\sss hence\sss the property\sss
$K\pff \subset\off H_{\dff \leq\dff k}$\sss
is\dss preserved during\sss the deformation,\oss
and\sss hence\sss 
$U_{\fff n}\dff \cap\qff C_{\dff k}$\sss
is\dss preserved.\oss

If\dss $m\qff >\qff n$\nnsp,\oss
then\sss
$U_{\fff n}\dff \cap\qff D_{\dff m}
\off =\off
\varnothing$\nnsp.\oss
If\dss $m\qff \leq\qff n$\sss and\dss
$K\off \supset\pff H_{\dff \leq\dff m}$\nsp,\oss
then\sss $\pi\dff K\off \supset\pff H_{\dff \geq\dff m}$\sss
and\sss $l_{\qff K}$\sss is\dss equal\sss to $0$ on\sss 
$H_{\dff \geq\dff m}$\nsp.\oss
It\sss follows\sss that\sss the property\sss
$K\off \supset\pff H_{\dff \geq\dff m}$\sss
is\dss preserved during\sss the deformation,\oss
and\sss hence\sss 
$U_{\fff n}\dff \cap\qff D_{\dff m}$\sss
is\dss preserved.\oss

Since\sss
$\gr\trf(\dff \infty\dff)$\sss is\dss the union of\dss
intersections\sss
$D_{\fff m}\dff \cap\qff C_{\dff k}$\sss
over all\sss $m\qff \leq\qff k$\nnsp,\oss
the results of\dss the\sss two previous paragraph\sss
imply\sss that\sss
$U_{\fff n}\dff \cap\qff \gr\trf(\dff \infty\dff)$\sss
is\dss preserved during\sss the deformation.\oss  \eproof

\mypar{Lemma.}{gr-is-cw}
\emph{The space\sss $\gr$ has\sss the homotopy\sss type of\trs a\sss
CW-complex.}

\proof
The map\sss
$K\qff \longmapsto\qff H_{\dff \leq\dff n}\dff \ominus\trf K$\sss
is\dss a homeomorphism\sss between\sss
$C_{\dff n}$\sss and\sss the space of\dss
finitely dimensional\sss subspaces of\sss $H_{\dff \leq\dff n}$\nsp.\oss
The\sss latter\dss is\dss known\sss to have\sss the homotopy\sss type of\dss
a\dss CW-complex.\oss
By\trs Lemma\qss \ref{filtration-of-gr}\qss this implies\sss that\sss
$U_{\fff n}$ has\sss the homotopy\sss type of\dss
a\dss CW-complex.\oss
Since $U_{\fff m}\qff \subset\qff U_{\fff n}$\sss for\sss
$m\qff <\qff n$\nnsp,\oss
every\sss finite intersection of\dss sets $U_{\fff n}$\sss is\dss
equal\sss to one of\dss them.\oss
Therefore\sss $U_{\fff n}\dff,\off n\qff \in\qff \zzz$\sss is\dss a covering of\dss $\gr$\sss
such\sss that\sss every finite intersections of\dss its\sss elements
has\sss the homotopy\sss type of\dss a\dss CW-complex.\oss
Since\sss $U_{\fff n}$ are open\sss by\trs
Lemma\qss \ref{filtration-of-gr}\qss and\sss
$\gr$\nnsp,\oss being a metric space,\oss is\dss paracompact,\oss
this covering\dss is\dss numerable.\oss
It\sss follows\sss that\sss $\gr$\sss
has\sss the homotopy\sss type of\dss a\dss CW-complex.\oss
See\trs tom\dss Dieck\qss \cite{td1},\oss Theorem\qss 4.\oss  \eproof

\mypar{Theorem.}{two-grassmannians}
\emph{The\sss inclusion\trs
$\gr\trf(\dff \infty\dff)\qff \ttoo\qff \gr$\sss
is\dss a\sss homotopy\sss equivalence.\oss}

\proof
It\dss is\dss well\sss known\sss that\sss the space\sss
$\gr\trf(\dff \infty\dff)$\sss has\sss the homotopy\sss type of\dss a\dss CW-complex.\oss
By\trs Lemma\qss \ref{gr-is-cw}\qss the space $\gr$ also
has\sss the homotopy\sss type of\dss a\dss CW-complex.\oss
Therefore\sss it\dss is\dss sufficient\sss to prove\sss that\sss
the inclusion\sss
$\gr\trf(\dff \infty\dff)\qff \ttoo\qff \gr$\sss
is\dss a\sss weak\sss homotopy\sss equivalence.\oss
Let\sss $X$\sss be a compact\sss space,\oss
and\sss let\sss 
$K\dff \colon\dff
x\off \longmapsto\off K\dff(\dff x\trf)$\sss
be a continuous map\sss $X\qff \ttoo\qff \gr$\nnsp.\oss
Since $\gr$\sss is\dss equal\sss to\sss the union of\dss
the increasing sequence of\dss open subsets\sss $U_{\fff m}$
from\trs Lemma\qss \ref{filtration-of-gr},\oss
the image of\dss $K$\sss is\dss contained\sss in\sss
$U_{\fff k}$\sss for some $k$\nnsp,\oss
and\trs Lemma\qss \ref{filtration-of-gr}\qss
implies\sss that\sss $K$\sss is\dss homotopic\sss to a map\sss
with\sss the image contained\sss in\sss $C_{\dff k}$\nsp.\oss
Therefore we can assume\sss that\sss $K$ maps $X$\sss to $C_{\dff k}$\nsp,\oss
i.e.\qss that\sss
$K\dff(\dff x\trf)\pff \subset\off H_{\dff \leq\dff k}$
for every\sss $x\qff \in\qff X$\nnsp.\oss

Let\sss 
$L\dff(\dff x\trf)
\off =\off
H\qff \ominus\qff
K\dff(\dff x\trf)$\nnsp.\oss
Then\sss
$L\dff(\dff x\trf)
\pff \supset\off
H_{\dff >\dff k}
\off =\off
H_{\dff \geq\dff k\dff +\dff 1}$
for every\sss $x\qff \in\qff X$\nnsp.\oss
Let\sss us\sss argue as\sss in\sss the previous paragraph,\oss
but\sss with\sss the roles of\dss positive 
and\sss negative numbers interchanged.\oss
Then\trs Lemma\qss \ref{filtration-of-gr}\qss
implies\sss that\sss $L$\sss
is\dss homotopic\sss to a map\sss
$L'\dff \colon\dff
X\qff \ttoo\qff \gr$\sss
such\sss that\sss
$L'\dff(\dff x\trf)
\pff \subset\off
H_{\dff \geq\dff n}$\sss
for some $n$\nnsp.\oss
Moreover,\oss
we can assume\sss that\sss
$L'\dff(\dff x\trf)
\pff \supset\off
H_{\dff \geq\dff k\dff +\dff 1}$
for every\sss $x\qff \in\qff X$\nnsp.\oss
Let\sss
$K\fff'\dff(\dff x\trf)
\off =\off
H\qff \ominus\qff
L'\dff(\dff x\trf)$\nnsp.\oss
Then\sss $K\fff'$\sss is\dss homotopic\sss to $K$
and\sss\vspace{3pt}
\[
\quad
H_{\dff <\dff n}
\pff \subset\off
K\fff'\dff(\dff x\trf)
\pff \subset\off
H_{\dff \leq\dff k}
\]

\vspace{-9pt}
for every\sss $x\qff \in\qff X$\nnsp.\oss
It\sss follows\sss that\sss $K$\sss is\dss homotopic\sss to a map\sss
with\sss the image in\sss
$\gr\trf(\dff \infty\dff)$\nnsp.\oss
The homotopies provided\sss by\trs Lemma\qss \ref{filtration-of-gr},\oss
as also\sss the operation\sss $K\off \longmapsto\off H\dff \ominus\dff K$\nnsp,\oss
preserve\sss the subspace\sss $\gr\trf(\dff \infty\dff)$\nnsp.\oss
It\sss follows\sss that\sss every map\sss of\dss pairs\sss
$(\trf X\fff,\qff Y\trf)
\qff \ttoo\qff
(\trf \gr\fff,\qff \gr\trf(\dff \infty\dff)\trf)$\sss
with compact\sss $X$\sss is\dss homotopic\sss to a map\sss
with\sss the image in\sss $\gr\trf(\dff \infty\dff)$\nnsp.\oss
This\dss is\dss obviously a stronger\sss property\sss than\sss 
$\gr\trf(\dff \infty\dff)\qff \ttoo\qff \gr$\sss
being\sss a\sss weak\sss homotopy\sss equivalence.\oss
The\sss theorem\sss follows.\oss  \eproof

\mysection{Categories\qss related\qss to\qss restricted\qss Grassmannians}{categories-grassmannians}

\myuppar{Categories\sss related\sss to restricted\dss Grassmannians.}
Let\sss
$P\off =\off
(\trf V\dff,\pff H_{\dff -}\dff,\pff H_{\dff +}\trf)$\sss
be a polarized subspace model.\oss
Recall\sss that\sss the category
$P\nsp\downarrow\fff \mathcal{P}{\nsp}\hat{\mathcal{S}}$ 
of\qss \emph{objects under}\trs $P$\sss
has as objects morphisms of\sss $\mathcal{P}{\nsp}\hat{\mathcal{S}}$
of\dss the form\sss $P\qff \ttoo\qff M$\nnsp,\oss
and as morphism commutative diagrams of\dss the form\vspace{-7pt}
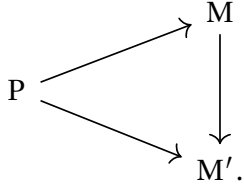
\[
\quad
\begin{tikzcd}[column sep=boom, row sep=sma]
&
M
\arrow[dd]
\\
P
\arrow[ru]
\arrow[rd]
&
\\
&
M\fff'\dff.
\end{tikzcd}
\]

\vspace{-16pt}
An object\sss of\sss 
$P\nsp\downarrow\fff \mathcal{P}{\nsp}\hat{\mathcal{S}}$\dnsp,\oss
i.e.\qss a morphism\sss $P\qff \ttoo\qff M$\sss has\sss the form
\[
\quad
(\trf V\fff,\qff H_{\dff -}\fff,\qff H_{\dff +}\dff)
\qff \ttoo\qff
(\trf W\fff,\qff K_{\dff -}\fff,\qff K_{\dff +}\dff)
\qff,
\]

\vspace{-12pt}
where $W$\sss is\dss a subspace admitting a decomposition\sss
$W
\off =\off
U_{\dff -}\qff \oplus\qff
V\qff \oplus\qff U_{\dff +}$\sss
for some subspaces\sss 
$U_{\dff -}\qff \subset\pff H_{\dff -}$\sss
and\dss
$U_{\dff +}\qff \subset\pff H_{\dff +}$\nsp.\oss
We will\sss call\sss such subspaces $W$\trs 
\emph{adapted}\pss to $P$\nnsp,\oss
or simply\qss \emph{adapted}.\oss
The subspaces $U_{\dff -}\dff,\pff U_{\dff +}$\sss
are uniquely\sss determined\sss by $W$\nnsp,\oss
and\sss the subspaces\dss 
$K_{\dff -}\fff,\off K_{\dff +}$\dss 
can\sss be recovered as\dss
$K_{\dff -}\off =\off H_{\dff -}\dff \ominus\qff U_{\dff -}$\nsp,\oss
$K_{\dff +}\off =\off H_{\dff +}\dff \ominus\qff U_{\dff +}$\nsp.\oss
Hence objects of\sss
$P\nsp\downarrow\fff \mathcal{P}{\nsp}\hat{\mathcal{S}}$\sss
can\sss be identified with adapted subspaces.\oss
In\sss these\sss terms,\oss
if\dss $W\fff,\pff W\fff'$\sss are adapted subspaces,\oss
then a morphism\sss $W\qff \ttoo\qff W\fff'$\sss exists\sss
if\trs and\dss only\trs if\dss
$W\qff \subset\pff W\fff'$\dnsp,\oss
and\sss in\sss this case
\begin{equation}
\label{w-morphism}
\quad
W\fff'
\off =\off
T_{\dff -}\qff \oplus\qff
W\qff \oplus\qff T_{\dff +}
\end{equation}

\vspace{-12pt}
for some subspaces\sss
$T_{\dff -}\qff \subset\qff H_{\dff -}$
and\sss
$T_{\dff +}\qff \subset\qff H_{\dff +}$\sss
uniquely determined\sss by\sss $W\fff,\pff W\fff'$\dnsp.\oss
If\dss a morphism\sss $W\qff \ttoo\qff W\fff'$\sss exists,\oss
it\dss is\dss unique.\oss

Let\sss $\mathcal{G}\dff (\trf P\trf)$\sss
be\sss the category\sss
having as objects diagrams in\sss 
$\mathcal{P}{\nsp}\hat{\mathcal{S}}$\sss
of\dss the form\sss
$P\qff \ttoo\qff M\off \longleftarrow\off N$\nnsp,\oss
where\sss $N$\sss is\dss an object\sss of\sss
$\mathcal{P}$\dnsp,\oss
and\sss as morphisms commutative diagrams 
of\dss the form\vspace{-7pt}
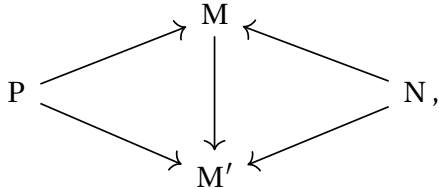
\[
\quad
\begin{tikzcd}[column sep=boom, row sep=sma]
&
M
\arrow[dd]
&
\\
P
\arrow[ru]
\arrow[rd]
&
&
N\qff,
\arrow[lu]
\arrow[ld]
\\
&
M\fff'
&
\end{tikzcd}
\]

\vspace{-16pt}
which can\sss be identified\sss with morphisms\sss
$M\qff \ttoo\qff M\fff'$\nnsp.\oss
Together\sss with\sss $\mathcal{P}{\nsp}\hat{\mathcal{S}}$\nnsp,\oss
it\dss
is\dss a\sss topological\sss category.\oss
In\sss terms of\dss adapted subspaces\sss the objects of\sss
$\mathcal{G}\dff (\trf P\trf)$ can be identified\sss with\sss
pairs consisting of\dss an
adapted subspace $W$\sss together with a splitting\sss
$W\off =\off W_{\dff -}\dff \oplus\trf W_{\dff +}$\nsp.\oss
A morphism\sss from\sss $W$\sss with\sss the splitting\sss
$W\off =\off W_{\dff -}\dff \oplus\trf W_{\dff +}$\sss
to another adapted subspace\sss $W\fff'$\sss with some splitting\sss
exists\sss if\trs and\dss only\trs if\dss
$W\qff \subset\pff W\fff'$\dnsp,\oss
and\sss in\sss this case\sss the splitting of\sss $W\fff'$\sss
should\sss be\vspace{1.5pt}
\begin{equation}
\label{induced-splitting}
\quad
W\fff'
\off =\off
\bigl(\trf T_{\dff -}\dff \oplus\dff W_{\dff -}\trf\bigr)
\qff \oplus\qff 
\bigl(\trf W_{\dff +}\dff \oplus\dff T_{\dff +}\trf\bigr)
\pff,
\end{equation}

\vspace{-12pt}\vspace{1.5pt}
where\sss $T_{\dff -}\fff,\pff T_{\dff +}$ are the subspaces from\qss (\ref{w-morphism}).\oss
Hence morphisms of\sss $\mathcal{G}\dff (\trf P\trf)$\sss can\sss be identified\sss
with\sss inclusions\sss $W\qff \subset\pff W\fff'$\sss together with a splitting\sss
of\dss the smaller one.\oss
One can consider objects of\dss $\mathcal{G}\dff (\trf P\trf)$\sss
as objects of\sss 
$P\nsp\downarrow\fff \mathcal{P}{\nsp}\hat{\mathcal{S}}$\sss
equipped\sss with a splitting of\dss the corresponding adapted subspace.\oss
Similarly,\oss one can consider\sss morphisms of\sss $\mathcal{G}\dff (\trf P\trf)$\sss
as morphisms of\dss 
$P\nsp\downarrow\fff \mathcal{P}{\nsp}\hat{\mathcal{S}}$\sss
equipped\sss with a splitting.\oss
Clearly,\pss 
$P\nsp\downarrow\fff \mathcal{P}{\nsp}\hat{\mathcal{S}}$\sss
is\dss the category defined\sss by\sss the inclusion 
partial\sss order on\sss the adapted subspaces.\oss
Hence\sss the category\sss
$\mathcal{G}\dff (\trf P\trf)$\sss
is\dss also defined\sss by a partial\sss order,\oss
the explicit\sss description of\dss which\sss 
we\sss leave\sss to\sss the reader.\oss

The category\sss $\mathcal{G}\dff (\trf P\trf)$\sss
should\sss be\sss thought\sss as a model\sss of\dss
the\dss Grassmannian\sss $\gr$\nnsp.\oss
In order\sss to define\sss the corresponding\sss model\sss of\dss
the\dss Grassmannian $\gr\trf(\dff \infty\dff)$\nnsp,\oss
let\sss us assume\sss that\sss
\[
\quad
V\off =\off H_{\trf 0}\qff,\quad
H_{\dff -}\off =\off H_{\trf <\dff 0}\qff,\quad
\mbox{and}\quad 
H_{\dff +}\off =\off H_{\trf >\dff 0}
\]

\vspace{-12pt} 
for an orthogonal\sss decomposition\qss (\ref{spectral-decomposition}).\oss
Ignoring\sss the\sss topological\sss structure for a moment,\oss
let\sss
$\mathcal{G}_{\qff \infty}\trf (\trf P\trf)$\sss
be\sss the full\sss subcategory
of\dss $\mathcal{G}\dff (\trf P\trf)$\sss having as objects adapted
subspaces $W$ such\sss that\sss
$W
\off \subset\off
H_{\trf [\trf -\dff n\fff,\dff n\trf]}$\sss 
for some $n$\nnsp,\oss together with a splitting.\oss
Let\sss
$\mathcal{G}_{\qff \nnn}\trf (\trf P\trf)$\sss
be\sss the full\sss subcategory
of\dss $\mathcal{G}\dff (\trf P\trf)$\sss having as objects 
subspaces 
$H_{\trf [\trf -\dff n\fff,\dff n\trf]}$\nnsp,\qss 
$n\qff \in\qff \nnn$\nnsp,\oss
together with a splitting,\oss
By\sss technical\sss reasons we will\sss need also\sss for each\sss $n\qff \in\qff \nnn$\sss
the full\sss subcategories\sss
$\mathcal{G}_{\dff n}\trf (\trf P\trf)$\sss
and\sss
$\mathcal{G}_{\qff \nnn,\qff n}\trf (\trf P\trf)$\sss
of\dss $\mathcal{G}\dff (\trf P\trf)$\sss having as objects adapted
subspaces $W$ such\sss that\sss
$W
\off \subset\off
H_{\trf [\trf -\dff n\fff,\dff n\trf]}$\sss
and\sss
$W
\off =\off
H_{\trf [\trf -\dff n\fff,\dff n\trf]}$\sss
respectively,\oss together with a splitting.\oss
Of\dss course,\oss the category\sss 
$\mathcal{G}_{\qff \nnn,\qff n}\trf (\trf P\trf)$\sss
has only\sss identity morphisms
and\sss its space of\dss objects\dss is\dss the space of\dss
splittings of\dss $H_{\trf [\trf -\dff n\fff,\dff n\trf]}$\nsp.\oss

Being\sss full\sss subcategories of\dss
$\mathcal{G}\dff (\trf P\trf)$\nnsp,\oss 
all\dss these categories\sss
are defined\sss by\sss partial\sss orders.\oss
In\sss particular,\oss to\sss turn\sss these categories into\sss topological\sss categories
one needs only\sss introduce\sss topologies into\sss their sets of\dss objects.\oss
In\sss the cases of\dss the categories\sss 
$\mathcal{G}_{\dff n}\trf (\trf P\trf)$\sss
and\sss
$\mathcal{G}_{\qff \nnn,\qff n}\trf (\trf P\trf)$\sss
there\dss is\dss only one natural\sss way\sss to do\sss this.\oss
The sets of\dss objects of\dss 
$\mathcal{G}_{\qff \infty}\trf (\trf P\trf)$\sss
and\sss
$\mathcal{G}_{\qff \nnn}\trf (\trf P\trf)$\sss
are\sss the unions over\sss $n\qff \in\qff \nnn$ of\dss spaces of\dss objects of\dss 
$\mathcal{G}_{\dff n}\trf (\trf P\trf)$\sss
and\sss
$\mathcal{G}_{\qff \nnn,\qff n}\trf (\trf P\trf)$\sss
respectively,\oss
and we equip\sss these sets with\sss the direct\sss limit\sss topologies
of\dss these unions.\oss

\mypar{Lemma.}{using-compactness}
\emph{Let\sss $n$ be a non-negative\sss integer and\dss let\dss
$(\trf X\fff,\qff Y\trf)$ be a pair of\dss compact\sss 
spaces\sss
$X\qff \supset\qff Y$\sss
having\sss the homotopy extension\sss property.\oss
Every\sss map of\qss pairs\dss}
\[
\quad
f\dff \colon\dff
\bigl(\trf X\fff,\qff Y\trf\bigr)
\qff \ttoo\qff
\bigl(\qff 
\ob\dff \mathcal{G}\dff (\trf P\trf)\fff,\pff 
\ob\dff \mathcal{G}_{\dff n}\dff (\trf P\trf)
\qff\bigr)
\]

\vspace{-12pt}
\emph{is\dss homotopic\sss to a map with\sss the image\sss in\dss
$\ob\dff \mathcal{G}_{\dff m}\dff (\trf P\trf)$\sss for some $m$ by\sss a\sss homotopy\sss
fixed on\dss $Y$\dnsp.}

\proof
Let\sss us consider\sss the case when\sss 
$Y\off =\off \varnothing$\sss
first.\oss
Let\sss $x\off \longmapsto\off W\dff(\dff x\trf)$\sss
be a map\sss from\sss $X$\sss to\sss the space of\dss adapted subspaces.\oss
Then\sss
$W\dff(\dff x\trf)
\off =\off
U_{\dff -}\dff(\dff x\trf)\qff \oplus\qff
V\qff \oplus\qff U_{\dff +}\dff(\dff x\trf)$\sss
for some maps\sss
$x\off \longmapsto\off U_{\dff -}\dff(\dff x\trf)$\nnsp,\qss
$x\off \longmapsto\off U_{\dff +}\dff(\dff x\trf)$\nnsp.\oss
These maps define vector\sss bundles on\sss $X$\nnsp.\oss
Since $X$\sss is\dss compact,\oss these bundles
are direct\sss summands of\dss trivial\dss bundles.\oss
Moreover,\oss there exist\sss maps\sss
$x\off \longmapsto\off C_{\dff -}\dff(\dff x\trf)$\nnsp,\qss
$x\off \longmapsto\off C_{\dff +}\dff(\dff x\trf)$\sss
such\sss that\sss the bundles corresponding\sss to\sss the maps\vspace{0pt}
\[
\quad
x
\off \longmapsto\off 
C_{\dff -}\dff(\dff x\trf)\dff \oplus\dff U_{\dff -}\dff(\dff x\trf)
\quad
\mbox{and}\quad
x
\off \longmapsto\off 
U_{\dff +}\dff(\dff x\trf)\dff \oplus\dff C_{\dff +}\dff(\dff x\trf)
\]

\vspace{-12pt}
are\sss trivial.\oss
It\dss follows\sss that\sss the subspaces 
$C_{\dff -}\dff(\dff x\trf)\dff \oplus\dff U_{\dff -}\dff(\dff x\trf)$ 
can\sss be equipped\sss with bases continuously depending on $x$\nnsp,\oss
i.e.\qss are spanned\sss by $k$\dnsp-frames in\sss $H_{\trf <\dff 0}$
continuously depending on $x$\sss 
where $k$\sss is\dss the dimension of\dss
$C_{\dff -}\dff(\dff x\trf)\dff \oplus\dff U_{\dff -}\dff(\dff x\trf)$\nnsp.\oss
Since\sss the space of\sss $k$\dnsp-frames in a\dss Hilbert\dss space\dss
is\dss contractible,\oss one can\sss deform\sss this family of\sss $k$\dnsp-frames\sss
to a constant\sss frame contained\sss in\sss
$H_{\dff [\dff -\dff m\fff,\qff -\dff 1\dff]}$\sss for some $m$\nnsp.\oss
Such a deformation can\sss be covered\sss by a deformation of\dss subspaces\sss
$U_{\dff -}\dff(\dff x\trf)$\nnsp,\oss
and\sss hence\sss the map\sss
$x\off \longmapsto\off U_{\dff -}\dff(\dff x\trf)$\sss
is\dss homotopic\sss to a map\sss
$x\off \longmapsto\off U\fff'_{\fff -}\dff(\dff x\trf)$\sss
such\sss that\sss 
$U\fff'_{\fff -}\dff(\dff x\trf)
\qff \subset\pff
H_{\dff [\dff -\dff m\fff,\qff -\dff 1\dff]}$\dss
for all $x$\nnsp.\oss
The same arguments apply\sss to\sss
$x\off \longmapsto\off U_{\dff +}\dff(\dff x\trf)$\nnsp.\oss
It\sss follows\sss that\sss the map\sss
$x\off \longmapsto\off W\dff(\dff x\trf)$\sss
is\dss homotopic\sss to a map\sss
$x\off \longmapsto\off W\fff'\dff(\dff x\trf)$\sss
such\sss that\sss\vspace{1.5pt}
\[
\quad
W\fff'\dff(\dff x\trf)
\off \subset\off
H_{\dff [\dff -\dff m\fff,\qff -\dff 1\dff]}
\qff \oplus\qff
V
\qff \oplus\qff
H_{\dff [\dff 1\fff,\qff m\dff]}
\off =\off
H_{\trf [\trf -\dff m\fff,\dff m\trf]}
\]

\vspace{-12pt}\vspace{1.5pt}
for some $m$ and all $x\qff \in\qff X$\nnsp.\oss
If\dss the spaces\sss $W\dff(\dff x\trf)$\sss are equipped\sss
with splittings continuously depending on $x$\sss
such a homotopy can\sss be covered\sss by a homotopy of\dss
splittings.\oss
It\sss follows\sss that\sss $f$
is\dss homotopic\sss to a map\sss with\sss the image in\sss 
$\ob\dff \mathcal{G}_{\dff m}\dff (\trf P\trf)$\nnsp.\oss

Let\sss us\sss consider now\sss the general\sss case.\oss
Suppose\sss that\sss\vspace{1.5pt}
\[
\quad
U_{\dff -}\dff(\trf y\trf)\qff \oplus\qff
V\qff \oplus\qff U_{\dff +}\dff(\trf y\trf)
\off =\off
W\dff(\trf y\trf)
\off \subset\off
H_{\trf [\trf -\dff n\fff,\dff n\trf]}
\]

\vspace{-39pt}
\[
\quad
\phantom{U_{\dff -}\dff(\trf y\trf)\qff \oplus\qff
V\qff \oplus\qff U_{\dff +}\dff(\trf y\trf)
\off =\off
W\dff(\trf y\trf)
\off }
=\off
H_{\dff [\dff -\dff n\fff,\qff -\dff 1\dff]}
\qff \oplus\qff
V
\qff \oplus\qff
H_{\dff [\dff 1\fff,\qff n\dff]}
\]

\vspace{-12pt}
for\sss some $n$ and every\sss $y\qff \in\qff Y$\dnsp,\oss
and\sss let\sss\vspace{1.5pt}
\[
\quad
D_{\dff -}\dff(\trf y\trf)
\off =\off
H_{\dff [\dff -\dff n\fff,\qff -\dff 1\dff]}\dff \ominus\dff U_{\dff -}\dff(\trf y\trf)
\quad
\mbox{and}\dff\quad
D_{\dff +}\dff(\trf y\trf)
\off =\off
H_{\dff [\dff 1\fff,\qff n\dff]}\dff \ominus\dff U_{\dff +}\dff(\trf y\trf)
\]

\vspace{-12pt}\vspace{1.5pt}
for\sss $y\qff \in\qff Y$\dnsp.\oss
Then\sss the direct\sss sum of\dss the bundle on\sss $Y$\sss defined\sss by\sss the map\sss
$y\off \longmapsto\off D_{\dff -}\dff(\trf y\trf)$\sss
with\sss the bundle defined\sss by\sss
$y\off \longmapsto\off U_{\dff -}\dff(\trf y\trf)$\sss
is\dss a\sss trivial\sss bundle.\oss
At\sss the same\sss time\sss the direct\sss sum of\dss 
the bundle on\sss $Y$\sss defined\sss by\sss the map\sss
$y\off \longmapsto\off C_{\dff -}\dff(\trf y\trf)$\sss
with\sss the bundle defined\sss by\sss
$y\off \longmapsto\off U_{\dff -}\dff(\trf y\trf)$\sss
is\dss also a\sss trivial\sss bundle.\oss
It\sss follows\sss that\sss the bundles defined\sss by\sss
$y\off \longmapsto\off C_{\dff -}\dff(\trf y\trf)$\sss
and\sss
$y\off \longmapsto\off D_{\dff -}\dff(\trf y\trf)$\sss
are stably\sss isomorphic,\oss
and\sss hence are\sss isomorphic after 
adding\sss to\sss them\sss trivial\sss bundles.\oss 
Since one can add\sss a\sss trivial\sss bundle\sss to\sss
the second one simply\sss by\sss replacing\sss $n$\sss
by some\sss $m\qff \geq\qff n$\nnsp,\oss
we can assume\sss that\sss these bundles on\sss $Y$\sss
are actually\sss isomorphic.\oss
Then\sss the maps\sss
$y\off \longmapsto\off C_{\dff -}\dff(\trf y\trf)\fff,\off D_{\dff -}\dff(\trf y\trf)$\nnsp,\oss
being\sss the classifying maps of\dss these bundles,\oss
are homotopic.\oss
Moreover,\oss they are homotopic\sss in\sss the class of\dss maps\sss
$y\off \longmapsto\off E\dff(\trf y\trf)\qff \subset\pff H_{\trf <\trf 0}$\sss
such\sss that\sss
$E\dff(\trf y\trf)$\sss
is\dss orthogonal\sss to $U_{\dff -}\dff(\trf y\trf)$\sss
for every\sss $y\qff \in\qff Y$\dnsp.\oss
In order\sss to\sss justify\sss the\sss last\sss claim,\oss
it\dss is\dss sufficient\sss to\sss trivialize\sss the bundle over\sss $Y$\sss
with\sss the fibers\sss 
$H_{\trf <\trf 0}\dff \ominus\dff U_{\dff -}\dff(\trf y\trf)$\sss
using\trs Kuiper's\trs theorem.\oss
By\sss using\trs Kuiper's\trs theorem once more,\oss
we can extend\sss this homotopy\sss to a homotopy\sss of\dss the map\sss
$x\off \longmapsto\off C_{\dff -}\dff(\dff x\trf)$\sss 
in\sss the class of\dss maps\sss
$x\off \longmapsto\off F\dff(\dff x\trf)\qff \subset\pff H_{\trf <\trf 0}$\sss
such\sss that\sss
$F\dff(\trf x\trf)$\sss
is\dss orthogonal\sss to $U_{\dff -}\dff(\dff x\trf)$\sss
for every\sss $x\qff \in\qff X$\nnsp.\oss
Therefore we can assume\sss that\sss
$C_{\dff -}\dff(\trf y\trf)
\off =\off
D_{\dff -}\dff(\trf y\trf)$\sss
for every\sss $y\qff \in\qff Y$\dnsp.\oss

Similarly,\oss we can assume\sss that\dss
$C_{\dff +}\dff(\trf y\trf)
\off =\off
D_{\dff +}\dff(\trf y\trf)$\dss
for every\sss $y\qff \in\qff Y$\dnsp.\oss
Then\sss the subspaces 
$C_{\dff -}\dff(\dff x\trf)\dff \oplus\dff U_{\dff -}\dff(\dff x\trf)$ 
and\sss
$C_{\dff +}\dff(\dff x\trf)\dff \oplus\dff U_{\dff +}\dff(\dff x\trf)$
can\sss be equipped\sss with bases continuously depending on $x$
and equal\sss to\sss some fixed\sss $k$\dnsp-frames contained\sss in\dss
$H_{\dff [\dff -\dff m\fff,\qff -\dff 1\dff]}$\sss
and\sss
$H_{\dff [\dff 1\fff,\qff -\dff m\dff]}$\sss
respectively\sss for every\sss $x\qff \in\qff Y$\dnsp.\oss
Using\sss the contractibility\sss of\dss the spaces of\dss $k$\dnsp-frames
we can\sss deform\sss these families of\dss frames\sss to constant\sss
families by\sss deformations fixed\sss on\sss $Y$\dnsp.\oss
These deformations can\sss be covered\sss by a deformation of\dss subspaces\sss
$U_{\dff -}\dff(\dff x\trf)$
and\sss $U_{\dff +}\dff(\dff x\trf)$\sss
fixed on\sss $Y$\dnsp,\oss
and\sss therefore\sss the maps\sss
$x\off \longmapsto\off U_{\dff -}\dff(\dff x\trf)\fff,\pff 
U_{\dff +}\dff(\dff x\trf)$\sss
are homotopic\sss to maps\sss
$x\off \longmapsto\off U\fff'_{\fff -}\dff(\dff x\trf)\fff,\pff 
U\fff'_{\fff +}\dff(\dff x\trf)$\sss
such\sss that\sss 
$U\fff'_{\fff -}\dff(\dff x\trf)
\qff \subset\pff
H_{\dff [\dff -\dff m\fff,\qff -\dff 1\dff]}$\dss
and\sss
$U\fff'_{\fff +}\dff(\dff x\trf)
\qff \subset\pff
H_{\dff [\dff 1\fff,\qff m\dff]}$\dss
for every $x$\nnsp.\oss
Moreover,\oss these deformations can\sss be assumed\sss to be fixed on\sss $Y$\dnsp.\oss
The\sss remaining\sss part\sss of\dss the proof\dss
is\dss concerned\sss with splittings and\dss is\dss exactly\sss
the same as in\sss the case\sss $Y\off =\off \varnothing$\nnsp.\oss  \eproof

\mypar{Lemma.}{free-degeneracies}
\emph{The\sss partially\sss ordered\sss spaces\qss
$\mathcal{G}_{\qff \nnn}\trf (\trf P\trf)$\nsp,\qss
$\mathcal{G}_{\qff \infty}\trf (\trf P\trf)$
and\trs
$\mathcal{G}\dff (\trf P\trf)$\sss
are free,\oss i.e.\qss
have free equalities,\oss
and\dss as simplicial\sss spaces\sss they\sss have\dss free degeneracies.\oss}

\proof
The objects $O$ of\dss $\mathcal{G}\dff (\trf P\trf)$\sss
are adapted subspaces $W$\sss together with a splitting of\dss $W$\nnsp,\oss
and\dss the relations\sss
$O\off =\off O\fff'$\nnsp,\qss
$O\off <\off O\fff'$\nnsp,\oss
and\dss
$O\off >\off O\fff'$\sss
between objects of\dss $\mathcal{G}\dff (\trf P\trf)$\sss
are equivalent\dss to\sss the relations\sss
$W\off =\off W\fff'$\nnsp,\pss
$W\off \subsetneq\off W\fff'$\nnsp,\oss
and\trs
$W\off \supsetneq\off W\fff'$\sss
respectively\sss between\sss the corresponding subspaces.\oss
The\sss latter\sss imply\sss that\dss
$\dim\dff W\off =\off \dim\dff W\fff'$\nnsp,\qss
$\dim\dff W\off <\off\dim\dff  W\fff'$\nnsp,\oss
and\dss
$\dim\dff W\off >\off \dim\dff W\fff'$\sss
respectively.\oss
These relations between\sss dimension define closed subsets\sss of\dss
the set\sss of\dss pairs\sss $W\fff,\pff W\fff'$\sss
comparable with\sss respect\sss to\sss the inclusion.\oss
It\sss follows\sss that\sss 
the partially ordered space\sss $\mathcal{G}\dff (\trf P\trf)$\sss
has\sss free equality\sss by\sss dimension reasons.\oss
The same argument\sss applies\sss to\sss 
$\mathcal{G}_{\qff \nnn}\trf (\trf P\trf)$\sss
and\sss
$\mathcal{G}_{\qff \infty}\trf (\trf P\trf)$\nsp.\oss
The second statement\sss now\sss follows\sss from\trs
Lemma\qss \ref{free-orders}.\oss  \eproof

\myuppar{The simplices of\dss $\mathcal{G}\dff (\trf P\trf)$\nnsp.}
By\sss the definition,\oss
an $n$\dnsp-simplex of\dss
$\mathcal{G}\dff (\trf P\trf)$\sss
is\dss a sequence\vspace{-0.5pt}
\[
\quad
P\qff \ttoo\qff
M_{\trf 0}\qff \ttoo\qff
M_{\dff 1}\qff \ttoo\qff
\ldots\qff \ttoo\qff
M_{\dff n}
\]

\vspace{-12pt}\vspace{-0.5pt}
of\dss morphisms in\sss the category\sss $\mathcal{P}{\nsp}\hat{\mathcal{S}}$\sss
together with a morphism\sss
$N\qff \ttoo\qff M_{\trf 0}$\nsp,\oss
where\sss $N$\sss is\dss an object\sss of\sss
$\mathcal{P}$\dnsp.\oss
Equivalently,\oss an $n$\dnsp-simplex\dss is\dss determined\sss by\sss
a non-decreasing\sss sequence\vspace{-0.5pt}
\[
\quad
W_{\fff 0}\off \subset\off
W_{\fff 1}\off \subset\off
\ldots\off \subset\off
W_{\fff n}
\]

\vspace{-12pt}\vspace{-0.5pt}
of\dss adapted subspaces\sss
together\sss with a splitting\sss of\dss $W_{\dff 0}$\nsp.\oss\vspace{-0.125pt}

\mypar{Lemma.}{infty-inclusion}
\emph{The natural\sss map\dss
$\num{\mathcal{G}_{\qff \infty}\trf (\trf P\trf)}
\qff \ttoo\qff
\num{\mathcal{G}\dff (\trf P\trf)}$\sss
is\dss a\sss weak\sss homotopy\sss equivalence.\oss}

\proof
Let\sss us\sss apply\trs
Lemma\qss \ref{using-compactness}\qss
to\sss the pairs\sss $(\trf S^{\fff n},\qff *\trf)$
and\sss $(\trf D^{\fff n\dff +\dff 1},\qff S^{\fff n}\trf)$\nnsp,\oss
where $S^{\fff n}$ and\sss $D^{\fff n\dff +\dff 1}$ are\sss the standard\sss
spheres and discs,\oss and $*\qff \in\qff S^{\fff n}$\dnsp.\oss
These special\sss cases\sss
imply\sss that\sss the inclusion\sss
$\mathcal{G}_{\qff \infty}\trf (\trf P\trf)
\qff \ttoo\qff
\mathcal{G}\dff (\trf P\trf)$\sss
induces weak\sss homotopy equivalence on\sss the spaces of\dss objects,\oss 
i.e.\qss on\sss the spaces of\sss $0$\dnsp-simplices 
of\dss these categories.\oss
Let\sss us consider\sss the spaces of\sss $n$\dnsp-simplices for an arbitrary $n$\nnsp.\oss
Let\sss $X$\sss be a compact\sss space.\oss
A map from\sss $X$\sss to\sss the space of\sss 
$n$\dnsp-simplices can\sss be considered as a continuous 
map\sss assigning\sss to points\sss
$x\qff \in\qff X$\sss chains\vspace{-0.5pt}
\[
\quad
W_{\fff 0}\trf(\dff x\trf)\off \subset\off
W_{\fff 1}\trf(\dff x\trf)\off \subset\off
\ldots\off \subset\off
W_{\fff n}\trf(\dff x\trf)
\]

\vspace{-12pt}\vspace{-0.5pt}
of\dss adapted subspaces
together\sss with a splitting\sss of\sss $W_{\fff 0}\trf(\dff x\trf)$
continuously depending on $x$\nnsp.\oss
A homotopy of\dss the map
$x\off \longmapsto\off W_{\fff n}\trf(\dff x\trf)$
can\sss be covered\sss by a homotopy of\dss the whole chain of\dss subspaces
$W_{\fff i}\trf(\dff x\trf)$ together with splittings of\dss subspaces
$W_{\fff 0}\trf(\dff x\trf)$\nnsp.\oss
Cf.\qss the proof\dss of\trs Lemma\qss \ref{using-compactness}.\oss
Therefore\trs Lemma\qss \ref{using-compactness}\qss 
implies\sss that\sss the inclusion\sss
$\mathcal{G}_{\qff \infty}\trf (\trf P\trf)
\qff \ttoo\qff
\mathcal{G}\dff (\trf P\trf)$\sss
induces weak\sss homotopy equivalences 
between\sss the spaces of\sss $n$\dnsp-simplices 
of\dss these categories.\oss
By\trs Lemma\qss \ref{free-degeneracies}\qss
the simplicial\sss spaces\sss
$\mathcal{G}_{\qff \infty}\trf (\trf P\trf)$
and\trs
$\mathcal{G}\dff (\trf P\trf)$\sss
have free degeneracies.\oss
Therefore\sss the\sss lemma\sss
follows from\trs Proposition\qss \ref{level-heq}.\oss  \eproof

\mypar{Lemma.}{infty-n-inclusion}
\emph{The inclusion\sss
$\num{\mathcal{G}_{\qff \nnn}\trf (\trf P\trf)}
\qff \ttoo\qff
\num{\mathcal{G}_{\qff \infty}\trf (\trf P\trf)}$\dss
is\dss a\sss weak\sss homotopy\sss equivalence.\oss}

\proof
We claim\sss that\sss
$\num{\mathcal{G}_{\qff \nnn,\qff n}\trf (\trf P\trf)}$\sss
is\dss a deformation\sss retract\sss of\dss
$\num{\mathcal{G}_{\dff n}\trf (\trf P\trf)}$\nnsp.\oss
Indeed,\oss if\dss $W$\sss is\dss an adapted subspace such\sss that\sss
$W\off \subset\off
H_{\trf [\trf -\dff n\fff,\dff n\trf]}$\nsp,\oss
then\sss there\dss is\dss a unique morphism\sss\vspace{-0.354pt}
\[
\quad
W\qff \ttoo\qff H_{\trf [\trf -\dff n\fff,\dff n\trf]}
\]

\vspace{-12pt}\vspace{-0.354pt} 
of\dss the category\sss
$P\nsp\downarrow\fff \mathcal{P}{\nsp}\hat{\mathcal{S}}$\nnsp.\oss
Moreover,\oss if\dss $W$\sss is\dss equipped\sss with a splitting,\oss
then\sss this morphism defines a splitting of\dss 
$H_{\trf [\trf -\dff n\fff,\dff n\trf]}$\nsp.\oss
In other\sss words,\oss every object\sss of\dss 
$\mathcal{G}_{\dff n}\trf (\trf P\trf)$
defines an object\sss of\dss 
$\mathcal{G}_{\qff \nnn,\qff n}\trf (\trf P\trf)$\sss
together with a morphism\sss from\sss the first\sss one\sss to\sss the second.\oss
If\dss there\dss is\dss a morphism\sss from one object\sss of\dss
$\mathcal{G}_{\dff n}\trf (\trf P\trf)$\sss to another,\oss
then\sss the objects of\dss
$\mathcal{G}_{\qff \nnn,\qff n}\trf (\trf P\trf)$\sss
assigned\sss to\sss them are equal.\oss
Therefore we can extend our construction\sss to morphisms\sss by\sss
assigning\sss to every\sss morphism of\dss
$\mathcal{G}_{\dff n}\trf (\trf P\trf)$\sss
some identity\sss morphism of\dss
$\mathcal{G}_{\qff \nnn,\qff n}\trf (\trf P\trf)$\sss
and get\sss a functor\vspace{-0.354pt}
\[
\quad
t\dff \colon\dff
\mathcal{G}_{\dff n}\trf (\trf P\trf)
\qff \ttoo\qff
\mathcal{G}_{\qff \nnn,\qff n}\trf (\trf P\trf)
\]

\vspace{-12pt}\vspace{-0.354pt}
together with a natural\sss transformation from\dss
the identity functor of\dss 
$\mathcal{G}_{\dff n}\trf (\trf P\trf)$\sss
to\sss $i\dff \circ\dff t$\nnsp,\oss
where\sss
$i\dff \colon\dff
\mathcal{G}_{\qff \nnn,\qff n}\trf (\trf P\trf)
\qff \ttoo\qff
\mathcal{G}_{\dff n}\trf (\trf P\trf)$\sss
is\dss the inclusion.\oss
At\sss the same\sss time\sss $t\dff \circ\dff i$\sss is\dss the identity functor.\oss
It\sss follows\sss that\sss
$\num{\mathcal{G}_{\qff \nnn,\qff n}\trf (\trf P\trf)}$\sss
is\dss a deformation\sss retract\sss of\dss
$\num{\mathcal{G}_{\dff n}\trf (\trf P\trf)}$\nsp,\oss
with\sss $\num{t\fff}$\sss being\sss the retraction.\oss
This proves our claim.\oss

Next,\oss we claim\sss that\sss for every compact\sss space $X$ and every 
continuous map\sss\vspace{-0.354pt}
\[
\quad
f\dff \colon\dff
X\qff \ttoo\qff
\num{\mathcal{G}_{\qff \infty}\trf (\trf P\trf)}
\]

\vspace{-12pt}\vspace{-0.354pt}
the image\sss $\image f$\dss is\dss contained\sss in\sss
$\num{\mathcal{G}_{\dff n}\trf (\trf P\trf)}$\sss
for some $n$\nnsp.\oss
By\trs Lemma\qss \ref{free-degeneracies}\qss
the partially ordered space\sss $\mathcal{G}_{\qff \infty}\trf (\trf P\trf)$\sss
is\dss free.\oss
This implies\sss that\sss the partially ordered space
$\mathcal{G}_{\dff n}\trf (\trf P\trf)$\sss
is\dss also free.\oss
Now\trs Corollary\qss \ref{free-full-realizations}\qss implies\sss that\sss
we can consider\sss 
$\mathcal{G}_{\qff \infty}\trf (\trf P\trf)$\sss
and\sss
$\mathcal{G}_{\dff n}\trf (\trf P\trf)$\sss
as\sss topological\sss simplicial\sss complexes 
without\sss affecting\sss the geometric realizations.\oss
Since $X$\sss is\dss compact,\oss a standard argument\sss about\sss
the direct\sss limit\sss topology\sss implies\sss that\sss 
$\image f$\sss is\dss contained\sss in\sss 
$\bbnum{\ssk_{\dff k}\dff \mathcal{G}_{\qff \infty}\trf (\trf P\trf)}$\sss
for some $k$\nnsp.\oss
See\dss Section\qss \ref{topological-simplicial-complexes}\qss
for\sss the definitions of\qss $\ssk_{\dff k}\dff \bullet$\sss and\dss $\bbnum{\bullet}$\nsp.\oss
These definitions,\oss together\sss with\sss the fact\sss that\sss
we are using\sss the direct\sss limit\sss topology on\sss objects of\dss
$\mathcal{G}_{\qff \infty}\trf (\trf P\trf)$\nnsp,\oss
imply\sss that\sss
$\bbnum{\ssk_{\dff k}\dff \mathcal{G}_{\qff \infty}\trf (\trf P\trf)}$\sss
is\dss the direct\sss limit\sss of\dss the subspaces\sss
$\bbnum{\ssk_{\dff k}\dff \mathcal{G}_{\dff n}\trf (\trf P\trf)}$\nsp.\oss
Therefore\sss the same standard argument\sss implies\sss that\sss the image\sss 
$\image f$\dss is\dss contained\sss in\sss
$\bbnum{\ssk_{\dff k}\dff \mathcal{G}_{\dff n}\trf (\trf P\trf)}$\sss
for some $n$\nnsp.\oss
Since\sss
$\bbnum{\ssk_{\dff k}\dff \mathcal{G}_{\dff n}\trf (\trf P\trf)}
\off \subset\off
\bbnum{\mathcal{G}_{\dff n}\trf (\trf P\trf)}
\off =\off
\num{\mathcal{G}_{\dff n}\trf (\trf P\trf)}$\nnsp,\oss
this proves our claim.\oss

It\sss follows\sss that\sss
every element\sss of\dss the homotopy group\sss
$\pi_{\dff k}\trf(\qff \num{\mathcal{G}_{\qff \infty}\trf (\trf P\trf)}\dff,\qff *\trf)$
can\sss be represented\sss by a continuos map\sss
$S^{k}\qff \ttoo\qff
\num{\mathcal{G}_{\dff n}\trf (\trf P\trf)}$\sss
for some $n$\nnsp.\oss
The base point\sss $*$\sss is\dss assumed\sss to be in\sss 
$\num{\mathcal{G}_{\trf 0}\trf (\trf P\trf)}$\nnsp.\oss
By\sss the first\sss claim above
such a map\dss is\dss homotopic\sss to a map with\sss the image contained\sss in
$\num{\mathcal{G}_{\qff \nnn,\qff n}\trf (\trf P\trf)}
\off \subset\off
\num{\mathcal{G}_{\qff  \nnn}\trf (\trf P\trf)}$\nsp.\oss
Moreover,\oss during\sss this homotopy\sss
the base point\sss $*$\sss moves along\sss the $1$\dnsp-simplex
defined\sss by\sss the morphism\sss
$H_{\trf 0}
\qff \ttoo\qff 
H_{\trf [\trf -\dff n\fff,\dff n\trf]}$\dnsp.\oss
It\dss follows\sss that\sss the homomorphisms of\dss
the homotopy\sss groups induced\sss by\sss the inclusion\sss
\[
\quad
\num{\mathcal{G}_{\qff  \nnn}\trf (\trf P\trf)}
\qff \ttoo\qff
\num{\mathcal{G}_{\qff \infty}\trf (\trf P\trf)}
\]

\vspace{-12pt}
are surjective.\oss
A similar argument\sss shows\sss that\sss they are injective.\oss
Therefore\sss this\sss inclusion\sss induces\sss isomorphisms of\dss
the homotopy groups and\dss is\dss a\sss weak\sss homotopy equivalence.\oss  \eproof

\mypar{Lemma.}{g-are-cw}
\emph{The spaces\sss
$\num{\mathcal{G}_{\qff \nnn}\trf (\trf P\trf)}$\nsp,\qss
$\num{\mathcal{G}_{\qff \infty}\trf (\trf P\trf)}$
and\trs
$\num{\mathcal{G}\dff (\trf P\trf)}$ 
are homotopy equivalent\sss to\sss CW-complexes.\oss}

\proof
By\trs Lemma\qss \ref{free-degeneracies}\qss
the partially ordered spaces\sss
$\mathcal{G}_{\qff \nnn}\trf (\trf P\trf)$\nsp,\qss
$\mathcal{G}_{\qff \infty}\trf (\trf P\trf)$
and\trs
$\mathcal{G}\dff (\trf P\trf)$\sss
have free equalities.\oss
Corollary\qss \ref{free-full-realizations}\qss
and\trs Lemma\qss \ref{free-delta-realizations}\qss 
imply\sss that\sss we can consider\sss 
them as\sss topological\sss simplicial\sss complexes or\sss $\Delta$\dnsp-spaces 
without\sss affecting\sss geometric realizations.\oss
Clearly,\oss the spaces of\sss $n$\dnsp-simplices of\dss these\sss
$\Delta$\dnsp-spaces are
homotopy equivalent\sss to\sss CW-complexes.\oss
Now\sss the construction of\dss the geometric realizations
by a sequence of\dss push-outs\qss
(see\trs Section\qss \ref{simplicial-spaces})\sss
shows\sss that\sss their geometric realizations are
homotopy equivalent\sss to\sss CW-complexes.\oss  \eproof

\mypar{Theorem.}{framing}
\emph{The canonical\sss maps\dss
$\num{\mathcal{G}_{\qff \nnn}\trf (\trf P\trf)}
\qff \ttoo\qff
\num{\mathcal{G}_{\qff \infty}\trf (\trf P\trf)}
\qff \ttoo\qff
\num{\mathcal{G}\dff (\trf P\trf)}$\dss
are\sss homotopy\sss equivalences.\oss}

\proof
By\trs Lemmas\qss \ref{infty-inclusion}\qss and\qss \ref{infty-n-inclusion}\qss
imply\sss that\sss 
these canonical\sss maps
are weak\sss homotopy equivalences.\oss 
Therefore\sss the\sss theorem\sss follows\sss from\trs
Lemma\qss \ref{g-are-cw}.\oss  \eproof

\myuppar{The spaces $\num{\mathcal{G}_{\qff \nnn}\trf (\trf P\trf)}$
and\dss $\gr\trf(\dff \infty\dff)$\nnsp.}
Let\sss us\sss consider\sss the set\sss $\nnn$\sss of\dss natural\sss numbers
as\sss the category defined\sss by\sss the usual\sss order on\sss $\nnn$\nnsp.\oss
So,\oss if\dss $n\qff \leq\qff m$\nnsp,\oss
then\sss there\dss is\dss a unique morphism\sss 
$n\qff \ttoo\qff m$\nnsp,\oss
if\dss $n\qff >\qff m$\nnsp,\oss
then\sss there are no such morphisms.\oss
For each\sss $n\qff \in\qff \nnn$\dss let\sss
$\mathbb{G}\trf(\dff n\trf)$\sss be\sss the space of\dss splittings
of\dss the subspace\sss $H_{\trf [\trf -\dff n\fff,\dff n\trf]}$\nnsp.\oss
or,\oss what\dss is\dss the same,\oss
the space of\dss objects of\sss
$\mathcal{G}_{\qff \nnn,\qff n}\trf (\trf P\trf)$\nnsp.\oss
The space $\mathbb{G}\trf(\dff n\trf)$ can\sss be identified\sss with\sss 
the\dss Grassmannian\sss of\dss subspaces of\dss 
$H_{\trf [\trf -\dff n\fff,\dff n\trf]}$\sss 
of\dss arbitrary dimension.\oss
If\dss $n\qff \leq\qff m$\nnsp,\oss
then\sss there\dss is\dss a unique morphism\sss
\[
\quad
H_{\trf [\trf -\dff n\fff,\dff n\trf]}
\off \ttoo\off
H_{\trf [\trf -\dff m\fff,\dff m\trf]}
\]

\vspace{-12pt}
of\sss the category
$\mathcal{P}{\nsp}\hat{\mathcal{S}}$\nnsp,\oss
and\dss this morphism defines a map\sss
$\mathbb{G}\trf(\dff n\trf)
\qff \ttoo\qff 
\mathbb{G}\trf(\dff m\trf)$\sss
of\dss the spaces of\dss splittings.\oss
In\sss more details,\oss
this map\sss takes a splitting\sss
$H_{\trf [\trf -\dff n\fff,\dff n\trf]}
\off =\off
W_{\dff -}\dff \oplus\dff W_{\dff +}$\sss
to\sss the splitting\vspace{1.5pt}
\[
\quad
H_{\trf [\trf -\dff m\fff,\dff m\trf]}
\off =\off
\bigl(\qff 
H_{\trf [\trf -\dff m\fff,\dff -\dff n\trf]}
\qff \oplus\qff W_{\dff -} 
\trf\bigr)
\off \oplus\off
\bigl(\trf 
W_{\dff +}
\qff \oplus\qff 
H_{\trf [\trf n\fff,\dff n\trf]} 
\qff\bigr)
\pff.
\]

\vspace{-10.5pt}
Compare\qss (\ref{induced-splitting}).\oss
The maps\sss
$\mathbb{G}\trf(\dff n\trf)
\qff \ttoo\qff 
\mathbb{G}\trf(\dff m\trf)$\sss
define a functor\sss $\mathbb{G}$\sss 
from\sss the category\sss $\nnn$\sss
to\sss the category of\dss topological\sss spaces,\oss
carrying\sss the same\sss information as\sss 
$\mathcal{G}_{\qff  \nnn}\trf (\trf P\trf)$\nnsp.\oss
Indeed,\oss objects of\dss
$\mathcal{G}_{\qff  \nnn}\trf (\trf P\trf)$
can\sss be considered as pairs\sss $(\dff n\fff,\qff x\trf)$\sss
such\sss that\sss $n\qff \in\qff \nnn$\sss and\sss
$x\qff \in\qff \mathbb{G}\trf(\dff n\trf)$\nnsp.\oss
A morphism\sss
$(\dff n\fff,\qff x\trf)
\qff \ttoo\qff
(\dff m\fff,\qff y\trf)$\sss
exists\sss if\trs and\dss only\trs if\dss
$n\qff \leq\qff m$\sss
and\sss $y$\sss is\dss the image of\dss $x$\sss under\sss the map\sss
$\mathbb{G}\trf(\dff n\trf)
\qff \ttoo\qff 
\mathbb{G}\trf(\dff m\trf)$\nnsp.\oss
Moreover,\oss
if\dss such a morphism exists,\oss then\sss it\dss is\dss unique.\oss
Essentially\sss by\sss the definition,\oss this means\sss that\sss
$\num{\mathcal{G}_{\qff \nnn}\trf (\trf P\trf)}$\sss
is\dss the\qss \emph{homotopy\sss colimit}\qss
of\dss the functor\sss $\mathbb{G}$\nnsp.\oss

We can\sss treat\sss the maps\sss
$\mathbb{G}\trf(\dff n\trf)
\qff \ttoo\qff 
\mathbb{G}\trf(\dff m\trf)$
as inclusions
and define\sss an\sss infinitely dimensional\trs
Grassmannian\sss $\mathbb{G}\trf(\dff \infty\dff)$\sss
as\sss the union of\sss $\mathbb{G}\trf(\dff n\trf)$
equipped\sss with\sss the direct\sss limit\sss topology.\oss
Let\sss
$g\dff \colon\dff
\mathbb{G}\trf(\dff \infty\trf)
\qff \ttoo\qff
\gr\trf(\dff \infty\trf)$\sss
be\sss the map\sss taking a splitting\sss
$H_{\trf [\trf -\dff n\fff,\dff n\trf]}
\off =\off
W_{\dff -}\dff \oplus\dff W_{\dff +}$\sss
to\sss the subspace\sss
$H_{\trf \leq\dff -\dff n}
\dff \oplus\dff
W_{\dff -}$\nsp,\oss
which\dss is\dss clearly admissible.\oss
It\dss is\dss easy\sss to see\sss that\sss $g$\sss
is\dss a\sss bijection.\oss
Moreover,\oss since both\sss
$\mathbb{G}\trf(\dff \infty\trf)$\sss
and\sss
$\gr\trf(\dff \infty\trf)$\sss
have direct\sss limit\sss topologies,\qss
$g$\sss is\dss a homeomorphism.\oss
We will\sss use\sss the map $g$\sss to identify\sss the spaces\sss
$\mathbb{G}\trf(\dff \infty\trf)$\sss 
and\sss
$\gr\trf(\dff \infty\trf)$\nnsp.\oss

Let\sss us consider\sss $\mathbb{G}\trf(\dff \infty\dff)$
as\sss a\sss topological\sss category\sss
having\sss $\mathbb{G}\trf(\dff \infty\dff)$ as\sss the space of\dss objects
and only\sss the identity\sss morphism.\oss
The classifying space of\dss this category\dss is\dss the\dss
space\sss $\mathbb{G}\trf(\dff \infty\dff)$\sss itself.\oss
Let\dss us\sss represent\sss the objects of\trs
$\mathcal{G}_{\qff \nnn}\trf (\trf P\trf)$\sss
by\sss pairs $(\dff n\fff,\qff x\trf)$\sss as above.\oss
By assigning\sss to an object\sss $(\dff n\fff,\qff x\trf)$\sss
the image of\sss $x$\sss in\sss $\mathbb{G}\trf(\dff \infty\dff)$\sss
and\sss to every morphism of\trs
$\mathcal{G}_{\qff \nnn}\trf (\trf P\trf)$\sss
an\sss identity\sss morphism of\dss 
$\mathbb{G}\trf(\dff \infty\dff)$\sss
we get\sss a\sss functor\sss
$\mathcal{G}_{\qff \nnn}\trf (\trf P\trf)
\qff \ttoo\qff
\mathbb{G}\trf(\dff \infty\dff)$
and\sss a continuous map\vspace{1.5pt}
\[
\quad
\num{\mathcal{G}_{\qff \nnn}\trf (\trf P\trf)}
\off \ttoo\off
\num{\mathbb{G}\trf(\dff \infty\dff)}
\off =\off
\mathbb{G}\trf(\dff \infty\dff)
\off =\off
\gr\trf(\dff \infty\dff)
\pff.
\]

\vspace{-12pt}\vspace{1.5pt}
\mypar{Theorem.}{grassmannian}
\emph{The map\dss
$\num{\mathcal{G}_{\qff \nnn}\trf (\trf P\trf)}
\qff \ttoo\qff 
\mathbb{G}\trf(\dff \infty\dff)
\off =\off
\gr\trf(\dff \infty\dff)$\dss
is\dss a\sss homotopy\sss equivalence.\oss}

\proof
For each\sss $m\qff \in\qff \nnn$\dss let\sss
$\mathcal{G}_{\qff \nnn,\qff \leq\qff m}\trf (\trf P\trf)$\sss
be\sss the full\sss subcategory of\dss $\mathcal{G}\dff (\trf P\trf)$\sss 
having as objects adapted
subspaces $W$ such\sss that\sss
$W
\off =\off
H_{\trf [\trf -\dff n\fff,\dff n\trf]}$\sss
for some\qss $n\qff \leq\qff m$\nnsp,\oss together with a splitting.\oss
The maps\sss
$\mathbb{G}\trf(\dff n\trf)
\qff \ttoo\qff 
\mathbb{G}\trf(\dff m\trf)$\sss
define a functor\sss\vspace{1.5pt}
\[
\quad
p_{\fff m}\dff \colon\dff
\mathcal{G}_{\qff \nnn,\qff \leq\qff m}\trf (\trf P\trf)
\qff \ttoo\qff
\mathcal{G}_{\qff \nnn,\qff m}\trf (\trf P\trf)
\]

\vspace{-12pt}\vspace{1.5pt}
and\sss a natural\sss transformation\sss from\sss
the\sss identity\sss functor of\sss 
$\mathcal{G}_{\qff \nnn,\qff \leq\qff m}\trf (\trf P\trf)$\sss
to $p_{\fff m}$\nsp.\oss
It\sss follows\sss that\sss $p_{\fff m}$\sss is\dss an equivalence of\dss categories
and\sss hence\sss the geometric realization\vspace{2.5pt}
\[
\quad
\num{p_{\fff m}}\dff \colon\dff
\num{\mathcal{G}_{\qff \nnn,\qff \leq\qff m}\trf (\trf P\trf)}
\qff \ttoo\qff
\num{\mathcal{G}_{\qff \nnn,\qff m}\trf (\trf P\trf)}
\off =\off
\mathbb{G}\trf(\dff m\trf)
\]

\vspace{-12pt}\vspace{2.5pt} 
is\dss a homotopy equivalence.\oss
Let\sss $n\qff \leq\qff m$\nnsp.\oss
Then\sss the square\vspace{1.5pt}
\[
\quad
\begin{tikzcd}[column sep=boom, row sep=rboom]
\protect{\num{\mathcal{G}_{\qff \nnn,\qff \leq\qff n}\trf (\trf P\trf)}}
\arrow[d, "\dis \protect{\num{p_{\fff n}}}\qff"']
\arrow[r]
&
\protect{\num{\mathcal{G}_{\qff \nnn,\qff \leq\qff m}\trf (\trf P\trf)}}
\arrow[d, "\dis \protect{\num{p_{\fff m}}}\qff"']
\\
\mathbb{G}\trf(\dff n\trf)
\arrow[r]
&
\mathbb{G}\trf(\dff m\trf)
\end{tikzcd}
\]

\vspace{-12pt}\vspace{0pt}
where\sss the horizontal\sss arrows are\sss inclusions,\oss
is\dss commutative.\oss
Clearly,\pss
$\num{\mathcal{G}_{\qff \nnn}\trf (\trf P\trf)}$\sss
is\dss the direct\dss limit\sss of\dss the subspaces\sss
$\num{\mathcal{G}_{\qff \nnn,\qff \leq\qff n}\trf (\trf P\trf)}$\nnsp.\oss
Since\sss
$\mathbb{G}\trf(\dff \infty\dff)
\off =\off
\gr\trf(\dff \infty\dff)$\sss 
is\dss the direct\dss limit\sss of\dss the subspaces\sss
$\mathbb{G}\trf(\dff n\trf)$\sss 
and\sss the maps\sss $\num{p_{\fff n}}$\sss
are homotopy equivalences,\oss
it\sss follows\sss that\sss\vspace{1.5pt}
\[
\quad
\num{\mathcal{G}_{\qff \nnn}\trf (\trf P\trf)}
\qff \ttoo\qff 
\mathbb{G}\trf(\dff \infty\dff)
\]

\vspace{-12pt}\vspace{1.5pt}
is\dss a homotopy equivalence.\oss
See\qss \cite{td1},\oss Lemma\qss 6.\oss  \eproof

\myuppar{Remark.}
The original\dss proof\dss of\trs Theorem\qss \ref{grassmannian}\qss
was based on\sss results of\trs Dugger\qss \cite{du},\oss
namely,\oss on\dss Theorem\qss 22.2\qss from\sss \cite{du}.\oss
The above proof\trs is\dss more elementary.\oss

\mypar{Theorem.}{s-sigma-homotopy-type}
\emph{There\sss is\dss a canonical\sss homotopy\sss
equivalence\sss 
$\num{\mathcal{G}\dff (\trf P\trf)}
\qff \ttoo\qff 
\mathbb{G}\trf(\dff \infty\dff)
\off =\off
\gr\trf(\dff \infty\dff)$\hnsp\dnsp.}

\proof
Theorem\qss \ref{framing}\qss together\sss with\trs
Theorem\qss \ref{grassmannian}\qss imply\sss that
$\num{\mathcal{G}\dff (\trf P\trf)}$\sss 
is\dss 
homotopy equivalent\sss to 
$\mathbb{G}\trf(\dff \infty\dff)$\hnsp\dnsp.\oss
While\sss the homotopy equivalence of\trs
Theorem\qss \ref{grassmannian}\qss does not\sss depend on any choices,\pss
the category\sss $\mathcal{G}_{\qff \nnn}\trf (\trf P\trf)$\sss
and\dss construction of\dss the homotopy equivalence\sss\vspace{1.5pt}
\[
\quad
\num{\mathcal{G}_{\qff \nnn}\trf (\trf P\trf)}
\qff \ttoo\qff
\num{\mathcal{G}\dff (\trf P\trf)}
\]

\vspace{-12pt}\vspace{1.5pt}
depend on\sss the orthogonal\sss decomposition\qss (\ref{spectral-decomposition}).\oss 
For\sss our current\sss purposes we can assume\sss that\sss
all\sss summand\sss $H_{\dff n}$\nsp,\oss except\sss $H_{\trf 0}$\nsp,\oss
are one-dimensional\sss subspaces generated\sss by\sss
vectors in orthonormal\sss bases of\dss subspaces\sss
$H_{\dff -}\off =\off H_{\trf <\dff 0}$\sss
and\sss 
$H_{\dff +}\off =\off H_{\trf >\dff 0}$\nsp.\oss
Then\sss the category 
$\mathcal{G}_{\qff \nnn}\trf (\trf P\trf)$
does not\sss depend on\sss the decomposition\qss (\ref{spectral-decomposition})\qss
up\sss to canonical\sss isomorphisms,\sss
but\sss different\sss choices of\dss bases\sss 
lead\sss to different\sss functors\sss
$\mathcal{G}_{\qff \nnn}\trf (\trf P\trf)
\qff \ttoo\qff
\mathcal{G}\dff (\trf P\trf)$\nnsp.\oss
At\sss the same\sss time,\oss 
if\dss two choices of\dss bases can\sss be connected\sss by
a continuous family,\oss 
then\sss the corresponding maps\sss
$\num{\mathcal{G}_{\qff \nnn}\trf (\trf P\trf)}
\qff \ttoo\qff
\num{\mathcal{G}\dff (\trf P\trf)}$\sss
are homotopic.\oss
Since\sss the space of\dss bases of\dss a\dss Hilbert\sss space
can\sss be identified\sss with\sss its\sss general\sss linear\sss group
and\sss hence\dss is\dss path connected,\oss
it\dss follows\sss that\sss up\sss to homotopy\sss the map\sss
$\num{\mathcal{G}_{\qff \nnn}\trf (\trf P\trf)}
\qff \ttoo\qff
\num{\mathcal{G}\dff (\trf P\trf)}$\sss
is\dss independent\sss on\sss the choices of\dss bases.\oss  \eproof

\myuppar{Remarks.}
The\dss Grassmannian\sss $\gr\trf(\dff \infty\dff)$\sss
appears\sss in\sss several\sss proofs of\qss 
Bott\trs periodicity\dss going\sss back\sss to\dss
McDuff\qss \cite{mc}\qss and\dss Quillen's\qss \cite{q}.\oss
In\dss Giffen's\dss proof\qss \cite{g}\qss
it\dss is\dss denoted\sss by 
$G\dff(\dff \infty\fff,\qff \infty\fff,\qff \ccc\trf)$
and appears\sss in a context\sss very close\sss to ours one.\oss
It\dss is\dss not\sss hard\sss to see\sss that\sss
$\gr\trf(\dff \infty\dff)$\sss
is\dss canonically\sss homotopy equivalent\dss to\sss
$\zzz\dff \times\dff B\fff U$\dnsp,\oss
the classifying space of\dss the $K$\dnsp-theory.\oss

\myuppar{The group\sss $U\dff(\dff \infty\dff)$\sss
and\sss category\sss $\hat{\mathcal{S}}_{\dff \infty}$.}
The methods of\dss this section work\sss in also in other situation.\oss
As an example,\oss in\sss the rest\sss of\dss this section\sss
we will\sss use\sss these methods\sss to prove\sss that\sss the group\sss
$U\ffin$\sss is\dss homotopy equivalent\sss to\sss $U\dff(\dff \infty\dff)$\nnsp,\oss
a more classical\dss infinite-dimensional\sss unitary\sss group.\oss
The group\sss $U\dff(\dff \infty\dff)$\sss is\dss an analogue of\sss
$\gr\trf(\dff \infty\dff)$\sss and\dss is\dss defined as follows.\oss
Suppose now\sss that\sss a decomposition of\sss $H$\sss in\sss the form\qss
(\ref{spectral-decomposition})\qss is\dss fixed.\oss
Let\sss us say\sss that\sss a unitary operator 
$u\dff \colon\dff H\qff \ttoo\qff H$\sss
is\qss \emph{admissible}\pss if\dss
for some\sss $n\qff \in\qff \nnn$\sss 
the operator\sss $u$\sss is\dss 
equal\sss to\sss the identity\sss on\vspace{1.5pt}
\[
\quad
H_{\qff \leq\dff n}\qff \oplus\qff H_{\qff \geq\dff n}
\pff.
\]

\vspace{-12pt}\vspace{1.5pt}
Let\sss $U\trf(\dff \infty\dff)$\sss be\sss the group of\dss all\sss
admissible unitary operators.\oss
Obviously,\pss
$U\trf(\dff \infty\dff)\qff \subset\pff U^{\dff \mathrm{fin}}$\dnsp.\oss
We equip\sss $U\trf(\dff \infty\dff)$\sss with\sss the\sss topology\sss
induced\sss from\sss $U^{\dff \mathrm{fin}}$\dnsp,\oss
or,\oss equivalently,\oss with\sss the direct\sss limit\sss topology\sss
defined\sss by\sss the subgroups of\dss operators equal\sss to\sss the identity on\sss
$H_{\qff \leq\dff n}\qff \oplus\qff H_{\qff \geq\dff n}$\sss for each\sss fixed $n$\nnsp.\oss
Let\dss
$-\qff U\trf(\dff \infty\dff)
\off =\off
\{\pff -\qff u\qff \mid\qff u\qff \in\qff U\dff(\dff \infty\dff) \pff\}$\nnsp.\oss
In order\sss to prove\sss that\sss the inclusion\sss
$U\dff(\dff \infty\dff)\qff \ttoo\pff U^{\dff \mathrm{fin}}$\sss
is\dss a homotopy equivalence,\oss
we need a category\sss $\hat{\mathcal{S}}_{\dff \infty}$\sss 
similar\sss to\sss the categories\sss $\mathcal{G}_{\qff \infty}\trf (\trf P\trf)$\nnsp.\oss
Let\sss us\sss define\sss the category\sss $\hat{\mathcal{S}}_{\dff \infty}$\sss 
as\sss the full\sss subcategory of\sss $\hat{\mathcal{S}}$\sss
having as objects finitely dimensional\sss subspaces\sss
$V\qff \subset\pff H$\sss such\sss that\sss
$V\qff \subset\pff H_{\trf [\trf -\qff n\fff,\qff n\trf]}$\sss for some $n$\nnsp.\oss
Then,\oss in\sss particular,\pss
$\num{\hat{\mathcal{S}}_{\dff \infty}}
\off \subset\off
\num{\hat{\mathcal{S}}}$\nnsp.\oss
One can easily check\sss that\sss the homeomorphism\sss
$\num{\hat{\mathcal{S}}}
\qff \ttoo\qff
-\qff U\ffin$\sss
of\trs Theorem\qss \ref{harris-h}\qss induces
a homeomorphism\sss
$\num{\hat{\mathcal{S}}_{\dff \infty}}
\qff \ttoo\qff
-\qff U\dff(\dff \infty\dff)$\nnsp.\oss
This\dss is\dss essentially\sss the\sss 
theorem of\qss Harris\qss \cite{h},\oss
mentioned at\sss the end of\qss Section\qss \ref{classifying-spaces-saf}.\oss
The square\vspace{-1.5pt}
\[
\quad
\begin{tikzcd}[column sep=sboom, row sep=sboom]
\protect{\num{\hat{\mathcal{S}}_{\dff \infty}}}
\arrow[d]
\arrow[r]
&
-\qff U\trf(\dff \infty\dff)
\arrow[d]
\\
\protect{\num{\hat{\mathcal{S}}}}
\arrow[r]
&
-\qff U\ffin
\end{tikzcd}
\]

\vspace{-9pt}
is\dss commutative and shows\sss that\sss in order\sss to prove\sss
that\sss $U\dff(\dff \infty\dff)\qff \ttoo\pff U^{\dff \mathrm{fin}}$\sss
is\dss a homotopy equivalence,\oss
it\dss is\dss sufficient\sss to prove\sss the following\sss theorem.\oss

\mypar{Theorem.}{s-framing}
\emph{The inclusion\dss
$\num{\hat{\mathcal{S}}_{\dff \infty}}
\qff \ttoo\qff
\num{\hat{\mathcal{S}}}$\dss
is\dss a\sss homotopy equivalence.\oss}

\proof
The proof\dss is\dss similar\sss to\sss the proof\dss of\trs
Theorem\qss \ref{framing}\qss and\dss is\dss actually simpler.\oss
The proof\dss starts with an analogue of\trs
Lemma\qss \ref{using-compactness}.\oss
One needs\sss to prove\sss that\sss a map\sss
$x\off \longmapsto\off V\dff(\dff x\trf)$\sss
from a compact\sss space\sss to\sss 
$\ob\dff \hat{\mathcal{S}}$\sss
can\sss be deformed\sss to a map with\sss the image in\sss
$\ob\dff \hat{\mathcal{S}}_{\dff \infty}$,\oss
as also\sss to prove a relative version of\dss this claim.\oss
The situation\dss is\dss simpler\sss than\sss in\trs
Lemma\qss \ref{using-compactness},\oss and\sss the same\sss method works.\oss
The category\sss $\hat{\mathcal{S}}_{\dff \infty}$\sss has free units
as a subcategory of\sss $\hat{\mathcal{S}}$\dnsp,\oss
and\sss hence\sss $\hat{\mathcal{S}}_{\dff \infty}$\sss
and\sss $\hat{\mathcal{S}}$\sss have free degeneracies as simplicial\sss spaces.\oss
This\dss is\dss an analogue of\qss Lemma\qss \ref{free-degeneracies}.\oss
Arguing as in\sss the proof\dss of\trs Lemma\qss \ref{infty-inclusion},\oss
we see\sss that\sss the inclusion\sss
$\num{\hat{\mathcal{S}}_{\dff \infty}}
\qff \ttoo\qff
\num{\hat{\mathcal{S}}}$\dss
is\dss a\sss weak\sss homotopy equivalence.\oss
Finally,\oss arguing as in\sss the proof\dss of\trs Theorem\qss \ref{framing},\oss
we see\sss that\sss $\num{\hat{\mathcal{S}}_{\dff \infty}}$\sss
and\sss $\num{\hat{\mathcal{S}}}$\sss
are homotopy equivalent\sss to\dss CW-complexes.\oss
The\sss theorem\sss follows.\oss  \eproof

\mypar{Corollary.}{classical-subspace-unitary}
\emph{The inclusion\sss 
$U\trf(\dff \infty\dff)\qff \ttoo\qff U^{\dff \mathrm{fin}}$
is\dss a homotopy equivalence.\oss
There\dss is\dss a canonical\dss  
homotopy\sss equivalence\sss
$\hat{\mathcal{F}}\qff \ttoo\qff U\trf(\dff \infty\trf)${\nsp}.\oss}\vspace{0.5pt}

\proof
As was explained above,\oss
the first\sss claim\sss follows\sss from\trs Theorem\qss \ref{s-framing}.\qss
The second claim\sss follows\sss from\dss the first\sss one and\dss
Corollary\qss \ref{subspace-unitary}.\oss  \eproof

\mysection{Grassmannian\qss bundle\qss and\qss Grassmannian\qss quasi-fibration}{grassmannian-fibrations}

\myuppar{The\dss Grassmannian\dss bundle.}
The\dss Grassmannian $\gr$ depends on\sss the choice of\dss
a polarization\sss
$H\off =\off K_{\dff -}\dff \oplus\dff K_{\dff +}$\nsp,\oss
or,\oss what\dss is\dss the same,\oss
on\sss the choice of\dss a closed\sss infinitely\sss dimensional\sss
subspace\sss $K_{\dff -}\qff \subset\qff H$\sss
of\dss infinite codimension.\oss
In\sss this section\sss we need\sss to\sss take\sss this dependence
into account.\oss
For a polarized subspace model\sss
$P\off =\off
(\trf V\dff,\pff H_{\dff -}\dff,\pff H_{\dff +}\trf)$\dss
let\sss
$\gr\trf(\trf P\trf)$\sss
be\sss the\dss Grassmannian $\gr$ associated\sss 
with\sss the subspace\sss $H_{\dff -}$\nsp.\oss
Let\sss us\sss consider\sss the space of\dss pairs\sss
$(\trf P\fff,\qff K\trf)$\sss
such\sss that\sss $P$\sss is\dss a polarized subspace model\sss
and\sss $K\qff \in\qff \gr\trf(\trf P\trf)$\nnsp.\oss
Clearly,\oss the projection\sss
$(\trf P\fff,\qff K\trf)
\off \longmapsto\off
P$\sss
from\sss this space of\dss pairs\sss to\sss the space of\dss
polarized subspace models,\oss
i.e.\qss to\sss the space of\dss objects of\sss 
$\mathcal{P}{\nsp}\hat{\mathcal{S}}$\dnsp,\oss 
is\dss a\sss locally\sss trivial\sss bundle
with\sss the fiber\sss $\gr$\nnsp.\oss

We would\sss like\sss to extend\sss this bundle\sss to a\sss
locally\sss trivial\sss bundle over
$\num{\mathcal{P}{\nsp}\hat{\mathcal{S}}}$\nnsp.\oss
It\dss is\dss convenient\sss to consider\sss
$\mathcal{P}{\nsp}\hat{\mathcal{S}}$\sss
as a\sss topological\sss simplicial\sss complex associated\sss
with\sss the partial\sss order defining\sss the category\sss 
$\mathcal{P}{\nsp}\hat{\mathcal{S}}$\dnsp.\oss
The same arguments as in\sss the proof\dss of\trs 
Lemma\qss \ref{free-degeneracies}\qss
show\sss that\sss $\mathcal{P}{\nsp}\hat{\mathcal{S}}$\sss
is\dss free as a\sss topological\sss partially ordered space.\oss
Hence\trs Corollary\qss \ref{free-full-realizations}\qss
implies\sss that\sss the geometric realizations\sss
$\num{\mathcal{P}{\nsp}\hat{\mathcal{S}}}$\sss
and\sss
$\bbnum{\mathcal{P}{\nsp}\hat{\mathcal{S}}}$\sss 
of\dss
$\mathcal{P}{\nsp}\hat{\mathcal{S}}$\sss
as a\sss topological\sss category and 
as a\sss topological\sss simplicial\sss complex
are\sss the same.\oss
The points of\dss $\bbnum{\mathcal{P}{\nsp}\hat{\mathcal{S}}}$\sss 
are represented\sss by\sss weighted sums\vspace{1.0pt}
\begin{equation}
\label{p-weighted-sum}
\quad
t_{\trf 0}\trf P_{\dff 0}\pff +\pff
t_{\dff 1}\trf P_{\dff 1}\pff +\pff
\ldots\pff +\pff
t_{\dff n}\trf P_{\dff n}
\pff,
\end{equation}

\vspace{-12pt}\vspace{1.0pt}
where\sss
$(\trf t_{\dff 0}\dff,\qff t_{\dff 1}\dff,\qff \ldots\dff,\qff t_{\dff n}\trf)
\qff \in\qff
\Delta^n$\dss
and\dss 
$P_{\dff 0}\fff,\qff
P_{\dff 1}\fff,\qff
\ldots\fff,\qff
P_{\dff n}$\sss 
are objects of\dss 
$\mathcal{P}{\nsp}\hat{\mathcal{S}}$
such\sss that\vspace{1.0pt}
\begin{equation}
\label{p-sequence}
\quad
P_{\dff 0}\off <\off
P_{\dff 1}\off <\off
\ldots\off <\off
P_{\dff n}
\pff.
\end{equation}

\vspace{-12pt}\vspace{1.0pt}
We claim\sss that\sss the\dss Grassmannians\sss
$\gr\trf(\trf P_{\dff i}\trf)$\sss
are equal.\oss
Indeed,\oss
if\trs $P\off \leq\off P\fff'$\nnsp,\oss where\vspace{1.0pt}
\[
\quad
P\off =\off
(\trf V\fff,\off H_{\dff -}\dff,\off H_{\dff +}\trf)
\quad
\mbox{and}\quad
P\fff'\off =\off
(\trf V\fff'\fff,\off H\fff'_{\dff -}\dff,\off H\fff'_{\dff +}\trf)
\quad
\]

\vspace{-12pt}\vspace{1.0pt}
are objects of\dss  
$\mathcal{P}{\nsp}\hat{\mathcal{S}}$\dnsp,\oss
then\sss $H\fff'_{\dff -}\qff \subset\qff H_{\dff -}$\sss
and\sss
$H_{\dff -}\dff \ominus\trf H\fff'_{\dff -}$\sss
is\dss finitely dimensional,\oss
and\sss hence\sss
$\gr\trf(\trf P\trf)
\off =\off
\gr\trf(\trf P\fff'\trf)$\nnsp.\oss
This proves\sss the claim.\oss

Given\sss  
$x\qff \in\qff 
\bbnum{\mathcal{P}{\nsp}\hat{\mathcal{S}}}$\nnsp,\oss
let\sss us represent\sss $x$\sss by a weighted sum\qss 
(\ref{p-weighted-sum})\qss and set\sss
$\gr\trf(\dff x\trf)
\off =\off
\gr\trf(\trf P_{\dff i}\trf)$\sss
for some $i$\nnsp.\oss
The above claim\sss implies\sss that\sss
$\gr\trf(\dff x\trf)$\sss
is\dss independent\sss on\sss the choice of\sss $i$ and,\oss 
moreover,\oss is\dss independent\sss on\sss the choice of\dss
the weighted sum\sss representing $x$\sss
because\sss the choice of\dss the\sss weighted sum\dss is\dss limited\sss to
adding or\sss removing\sss objects\sss $P_{\fff i}$\sss
with\sss the coefficient\sss $t_{\dff i}\off =\off 0$\nnsp.\oss
Let\sss $\mathbf{G}$\sss be\sss the space of\dss pairs\sss
$(\dff x\fff,\qff K\trf)$\sss 
such\sss that\sss
$x\qff \in\qff \bbnum{\mathcal{P}{\nsp}\hat{\mathcal{S}}}$\sss
and\dss $K\qff \in\qff \gr\trf(\dff x\trf)$\nnsp,\oss
and\sss let\vspace{-0.125pt}
\[
\quad
\bm{\pi}\dff \colon\dff
\mathbf{G}
\qff \ttoo\qff 
\bbnum{\mathcal{P}{\nsp}\hat{\mathcal{S}}}
\off =\off
\num{\mathcal{P}{\nsp}\hat{\mathcal{S}}}
\]

\vspace{-12pt}
be\sss the projection\sss
$(\dff x\fff,\qff K\trf)
\off \longmapsto\off
x$\nnsp.\oss
Lemma\qss \ref{free-delta-realizations}\qss allows\sss to identify\dss 
$\bbnum{\mathcal{P}{\nsp}\hat{\mathcal{S}}}$\sss
and\sss
$\num{\Delta\qff \mathcal{P}{\nsp}\hat{\mathcal{S}}}_{\dff \Delta}$.\oss
The construction of\trs
$\num{\Delta\qff \mathcal{P}{\nsp}\hat{\mathcal{S}}}_{\dff \Delta}$\sss
by an\sss infinite sequence of\dss push-outs\qss
(see\sss Section\qss \ref{simplicial-spaces})\sss
shows\sss that\sss $\bm{\pi}$\sss is\dss a\sss locally\sss trivial\sss bundle
with\sss the fiber\sss $\gr$\nnsp.\oss
As we will\sss see,\oss the\sss total\sss space $\mathbf{G}$ of\dss this
bundle\dss is\dss contractible
and\sss $\gr$\sss is\dss
homotopy equivalent\sss to\sss the\sss loop space of\dss
$\hat{\mathcal{F}}$\nnsp.\oss

\myuppar{The\dss Grassmannian\dss quasi-fibration.}
Now we are going\sss to replace\sss the spaces $\gr\trf(\trf P\trf)$
by\sss the spaces $\num{\mathcal{G}\dff (\trf P\trf)}$\sss
in\sss the above construction.\oss
One cannot\sss get\sss a\sss locally\sss trivial\sss bundle in\sss
this way\sss because\dss
$\num{\mathcal{G}\dff (\trf P\fff'\trf)}
\off\dff \neq\off\dff
\num{\mathcal{G}\dff (\trf P\trf)}$\trs
if\qss $P\fff'\off <\off P$\dnsp.\oss
Nevertheless,\oss one can adapt\sss this construction\sss to
get\sss a quasi-fibration\sss with similar\sss homotopy\sss properties.\oss
The following adaptation\dss is\dss based on\sss the ideas of\trs
Quillen\qss \cite{q}.\oss

Suppose\sss that\sss $P\fff'\off \leq\off P$\sss
and\sss let\sss
$u\dff \colon\dff
P\fff'\qff \ttoo\qff P$\sss
be\sss the morphism corresponding\sss to\sss this inequality.\oss
Recall\sss that\sss an object\sss of\dss
$\mathcal{G}\dff (\trf P\trf)$\sss
is\dss a diagram in\sss 
$\mathcal{P}{\nsp}\hat{\mathcal{S}}$\sss
of\dss the form\sss
$P\qff \ttoo\qff M\off \longleftarrow\off N$\nnsp,\oss
where\sss $N$\sss is\dss an object\sss of\sss
$\mathcal{P}$\dnsp.\oss
By\sss taking\sss the composition of\dss the morphism\sss
$P\qff \ttoo\qff M$\sss with $u$
we get\sss the diagram\sss
$P\fff'\qff \ttoo\qff M\off \longleftarrow\off N$\nnsp,\oss
which\dss is\dss an object\sss of\dss
$\mathcal{G}\dff (\trf P\fff'\trf)$\nnsp.\oss
This construction extends in an obvious way\sss to morphisms
and defines a continuous functor
\[
\quad
u^{\dff *}\dff \colon\qff
\mathcal{G}\dff (\trf P\trf)
\qff \ttoo\pff
\mathcal{G}\dff (\trf P\fff'\trf)
\pff.
\]

\vspace{-12pt}
In\sss terms of\dss adapted subspaces\sss $u^{\dff *}$ can\sss be described
as follows.\oss
A subspace adapted\sss to\sss $P$\sss
is\dss obviously adapted\sss to\sss $P\fff'$\dnsp.\oss
Given a subspace\sss $W$\sss adapted\sss to\sss $P$\sss
together with a splitting\sss
$W\off =\off W_{\dff -}\dff \oplus\trf W_{\dff +}$\nsp,\oss
the functor\sss $u^{\dff *}$\sss takes\sss them\sss to\sss
the same subspace\sss together with\sss the same splitting,\oss
but\sss now considered as a subspace adapted\sss to\sss $P\fff'$
with a splitting.\oss
From\sss the point\sss of\dss view of\dss adapted subspaces\sss $u^{\dff *}$\sss
is\dss simply\sss the inclusion of\dss the subcategory\sss
$\mathcal{G}\dff (\trf P\trf)$\sss
of\dss the category\sss
$\mathcal{G}\dff (\trf P\fff'\trf)$\sss
into\sss $\mathcal{G}\dff (\trf P\fff'\trf)$\nnsp.\oss
By\sss passing\sss to geometric realizations
we get\sss a map
\[
\quad
\num{u^{\dff *}}\dff \colon\qff
\num{\mathcal{G}\dff (\trf P\trf)}
\qff \ttoo\pff
\num{\mathcal{G}\dff (\trf P\fff'\trf)}
\pff.
\]

\vspace{-12pt}
The categories\sss $\mathcal{G}\dff (\trf P\trf)$
and\sss $\mathcal{G}\dff (\trf P\fff'\trf)$
are associated\sss with\sss partial\sss orders,\oss
and\dss the functor\sss $u^{\dff *}$\sss
is\dss strictly order-preserving as map of\dss ordered sets,\oss
i.e.\dss 
$u^{\dff *}\dff(\dff a\trf)
\qff <\qff
u^{\dff *}\dff(\dff b\trf)$\sss
is\dss equivalent\sss to\sss $a\qff <\qff b$\nnsp.\oss
In\sss particular,\pss $u^{\dff *}$\sss
takes non-degenerate simplices\sss to non-degenerate simplices.\oss

The construction of\dss the bundle\sss
$\bm{\pi}\dff \colon\dff
\mathbf{G}
\qff \ttoo\qff 
\num{\mathcal{P}{\nsp}\hat{\mathcal{S}}}$\sss
implicitly used\sss the fact\sss that\dss all\dss Grassmannians\sss
are contained\sss in a single space,\oss
the unrestricted\dss Grassmannian\sss $\mathbf{G{\fff}r}$\nnsp.\oss
See\dss Section\qss \ref{spaces-operators}\qss for\sss the\sss latter.\oss
Now we need a category containing\sss all\sss categories\sss
$\mathcal{G}\dff (\trf P\trf)$ as subcategories.\oss
It\dss is\dss clear\sss that\sss the space of\dss objects of\dss
this category should\sss be\sss the space of\dss diagrams in\sss 
$\mathcal{P}{\nsp}\hat{\mathcal{S}}$\sss
of\dss the form\sss
$P\qff \ttoo\qff M\off \longleftarrow\off N$\nnsp,\oss
where\sss $P$\sss is\dss an arbitrary object\sss of\sss
$\mathcal{P}{\nsp}\hat{\mathcal{S}}$\sss
and\sss $N$\sss is\dss an object\sss of\sss
$\mathcal{P}$\dnsp.\oss
Similarly,\oss morphisms of\dss all\sss categories\sss
$\mathcal{G}\dff (\trf P\trf)$ should\sss be morphisms
of\dss this category.\oss
Following\dss Quillen\qss \cite{q},\oss
we will\sss also incorporate\sss morphisms\sss
$P\fff'\qff \ttoo\qff P$\sss into\sss its structure.\oss
There\dss is\dss a\sss topological\sss category\sss
having diagrams as above as objects,\oss and commutative diagrams\vspace{0.75pt}
\begin{equation}
\label{sub-morphisms}
\quad
\begin{tikzcd}[column sep=boom, row sep=normal]
P
\arrow[r, "\dis v"]
&
M
\arrow[dd, "\dis w"]
&
\\
&&
N
\arrow[lu]
\arrow[ld]
\\
P\fff'
\arrow[r, "\dis v'"]
\arrow[uu, "\dis u"]
&
M\fff'
&
\end{tikzcd}
\end{equation}

\vspace{-12pt}\vspace{-0pt}
as morphisms from\sss
$P\qff \ttoo\qff M\off \longleftarrow\off N$\sss
to\sss
$P\fff'\qff \ttoo\qff M\fff'\off \longleftarrow\off N$\nnsp.\oss
In\sss particular,\oss the objects $N$ are always\sss
the same for\sss the source and\sss the\sss target\sss of\dss a morphism.\oss
The commutativity of\dss the square,\oss of\dss course,\oss
means\sss that\sss
$v'\off =\off w\qff \circ\qff v\qff \circ\qff u$\nnsp.\oss
We will\sss denote\sss this category\sss by\sss $S\dff(\trf \phi\trf)$\sss
because it\dss is\dss a special\sss case of\trs Quillen's\dss categories\sss
$S\dff(\trf f\trf)$\sss applied\sss to\sss the forgetting\sss functor\sss
$\phi\dff \colon\dff
s\dff \hat{\mathcal{S}}
\qff \ttoo\qff
\mathcal{P}{\nsp}\hat{\mathcal{S}}$\sss
from\trs Section\qss \ref{polarizations-splittings}\qss
in\sss the role of\sss $f$\dnsp.\oss
There\dss is\dss an obvious\qss \emph{contravariant}\pss
forgetting\sss functor\sss
$p\dff \colon\dff
S\dff(\trf \phi\trf)
\qff \ttoo\qff
\mathcal{P}{\nsp}\hat{\mathcal{S}}$\nnsp,\oss
which can\sss be considered as a covariant\sss functor\sss
from\sss $S\dff(\trf \phi\trf)$\sss to\sss the category
opposite\sss to\sss
$\mathcal{P}{\nsp}\hat{\mathcal{S}}$\nnsp.\oss
Since\sss the geometric realizations do not\sss change
when a category\dss is\dss replaced\sss by\sss its opposite,\pss
$p$\sss induces a continuous map\vspace{-0.125pt}
\[
\quad
\num{p}\dff \colon\dff
\num{S\dff(\trf \phi\trf)}
\qff \ttoo\qff
\num{\mathcal{P}{\nsp}\hat{\mathcal{S}}}
\pff.
\]

\vspace{-12pt}
This\dss is\dss our categorical\sss analogue of\dss the bundle\sss
$\bm{\pi}\dff \colon\dff
\mathbf{G}
\qff \ttoo\qff 
\bbnum{\mathcal{P}{\nsp}\hat{\mathcal{S}}}
\off =\off
\num{\mathcal{P}{\nsp}\hat{\mathcal{S}}}$\nnsp.\oss

\myuppar{The partial\sss order defining\sss $S\dff(\trf \phi\trf)$\nnsp.}
Let\sss us describe\sss the category\sss $S\dff(\trf \phi\trf)$\sss
in\sss terms of\dss the partial\sss order defining\sss
$\mathcal{P}{\nsp}\hat{\mathcal{S}}$\sss
and adapted subspaces.\oss
An object\sss of\sss $S\dff(\trf \phi\trf)$\sss
can\sss be identified\sss with an object\sss $P$\sss of\sss
$\mathcal{P}{\nsp}\hat{\mathcal{S}}$\sss
together with a subspace\sss $W$\sss adapted\sss to $P$\sss
and a splitting\sss
$W\off =\off W_{\dff -}\dff \oplus\trf W_{\dff +}$\nsp.\oss
Let\sss us\sss denote such an object\sss by\sss
$(\qff P\fff,\off
W\off =\off W_{\dff -}\dff \oplus\trf W_{\dff +}\qff)$\nnsp.\oss
A morphism\vspace{0.5pt}
\[
\quad
(\qff P\fff,\off
W\off =\off W_{\dff -}\dff \oplus\trf W_{\dff +}\qff)
\off \ttoo\off
(\qff P\fff'\fff,\off
W\fff'\off =\off W\fff'_{\fff -}\dff \oplus\trf W\fff'_{\fff +}\qff)
\]

\vspace{-12pt}\vspace{0.5pt}
exists\dss if\trs and\dss only\trs if\trs
$P\fff'\off \leq\off P$\sss
and\dss $W\off \leq\off W\fff'$\dnsp.\oss
Moreover,\oss if\dss such a morphism exists,\oss it\dss is\dss unique.\oss
In\sss fact,\oss
the category\sss $S\dff(\trf \phi\trf)$\sss
is\dss associated\sss with\sss the partial\sss order\sss $\leq$\nnsp,\oss
where\vspace{0.5pt}
\[
\quad
(\qff P\fff,\off
W\off =\off W_{\dff -}\dff \oplus\trf W_{\dff +}\qff)
\off \leq\off
(\qff P\fff'\fff,\off
W\fff'\off =\off W\fff'_{\fff -}\dff \oplus\trf W\fff'_{\fff +}\qff)
\]

\vspace{-12pt}\vspace{0.5pt}
if\trs
$P\fff'\off \leq\off P$\dnsp,\oss
$W\off \leq\off W\fff'$\sss
and\sss the splitting\sss
\[
\quad
W\fff'
\off =\off 
W\fff'_{\fff -}\dff \oplus\trf W\fff'_{\fff +}
\]

\vspace{-12pt}
is\dss determined\sss by\sss $W\fff'$\sss and\dss the splitting\sss
\[
\quad
W
\off =\off 
W_{\dff -}\dff \oplus\trf W_{\dff +}
\]

\vspace{-12pt}
as\sss in\qss (\ref{induced-splitting}).\oss
Equivalently,\oss 
$P\fff'\off \leq\off P$\dnsp,\oss
$W\off \leq\off W\fff'$\sss and\sss
the splitting\sss
\[
\quad
W
\off =\off 
W_{\dff -}\dff \oplus\trf W_{\dff +}
\]

\vspace{-12pt}
is\dss determined\sss by\sss $W$\sss and\dss the splitting\sss
$W\fff'
\off =\off
W\fff'_{\fff -}\dff \oplus\trf W\fff'_{\fff +}$\sss
as\sss in\qss (\ref{induced-splitting}).\oss

This description of\sss $S\dff(\trf \phi\trf)$\sss
shows,\oss in\sss particular,\oss that\sss the category\sss
$\mathcal{G}\dff (\trf P\trf)$\sss
is\dss the\sss full\sss subcategory\sss of\sss $S\dff(\trf \phi\trf)$\sss
defined\sss by\sss objects of\dss the form\sss
$(\qff P\fff,\off
W\off =\off W_{\dff -}\dff \oplus\trf W_{\dff +}\qff)$\nnsp.\oss

The\sss partially\sss ordered\sss space\sss $S\dff(\trf \phi\trf)$\sss
is\dss free\sss by\sss the same reasons as\sss 
$\mathcal{P}{\nsp}\hat{\mathcal{S}}$\nnsp.\oss
See\sss the proof\dss of\trs 
Lemma\qss \ref{free-degeneracies}.\oss
Therefore\trs Lemma\qss \ref{free-orders}\qss
implies\sss that\sss $S\dff(\trf \phi\trf)$
has free degeneracies,\oss
and\trs Corollary\qss \ref{free-full-realizations}\qss
implies\sss that\sss the geometric realizations of\sss
$S\dff(\trf \phi\trf)$
as a\sss topological\sss category and 
as a\sss topological\sss simplicial\sss complex
are\sss the same.\oss

\mypar{Lemma.}{fibers-are-he}
\emph{For every\sss morphism\dss
$u\dff \colon\dff
P\fff'\qff \ttoo\qff P$\sss
of\qss $\mathcal{P}{\nsp}\hat{\mathcal{S}}$\sss
the induced\sss map}\vspace{-0.5pt}
\[
\quad
\num{u^{\dff *}}\dff \colon\qff
\num{\mathcal{G}\dff (\trf P\trf)}
\qff \ttoo\pff
\num{\mathcal{G}\dff (\trf P\fff'\trf)}
\]

\vspace{-12pt}\vspace{-0.5pt}
\emph{is\dss a\sss homotopy\sss equivalence.\oss}

\proof
The morphism $u$\sss has\sss the form\vspace{-0.5pt}
\[
\quad
(\trf W\fff,\qff K_{\dff -}\fff,\qff K_{\dff +}\dff)
\qff \ttoo\qff
(\trf V\fff,\qff H_{\dff -}\fff,\qff H_{\dff +}\dff)
\qff,
\]

\vspace{-12pt}
where\sss
$V
\off =\off
U_{\dff -}\qff \oplus\qff
W\qff \oplus\qff U_{\dff +}$\sss
for some subspaces\sss 
$U_{\dff -}\qff \subset\pff K_{\dff -}$\sss
and\dss
$U_{\dff +}\qff \subset\pff K_{\dff +}$\nsp.\oss
We may assume\sss that\sss
the\sss decomposition\qss (\ref{spectral-decomposition})\qss
from\dss Section\qss \ref{restricted-grassmannians}\qss is\dss chosen\sss
in such a way\sss that
\[
\quad
H_{\dff -\dff 1}\off =\off U_{\dff -}
\qff,\quad
H_{\trf 0}\off =\off W
\qff,\quad
\mbox{and}\quad
H_{\dff 1}\off =\off U_{\dff +}
\pff.
\]

\vspace{-12pt}
Then\sss the functor\sss $u^{\dff *}$\sss induces a functor\dss 
$u^{\dff *}_{\qff \nnn}\qff \colon\qff
\mathcal{G}_{\qff \nnn}\trf (\trf P\trf)
\qff \ttoo\qff
\mathcal{G}_{\qff \nnn}\trf (\trf P\fff'\trf)$\sss
such\sss that\sss the diagram\vspace{-3pt}
\[
\quad
\begin{tikzcd}[column sep=boom, row sep=normal]
\mathcal{G}_{\qff \nnn}\trf (\trf P\trf)
\arrow[dd, "\dis u^{\dff *}_{\qff \nnn}\dff"']
\arrow[rd]
&
\\
&
\mathbb{G}\trf(\dff \infty\dff)
\\
\mathcal{G}_{\qff \nnn}\trf (\trf P\fff'\trf)
\arrow[ru]
\end{tikzcd}
\]

\vspace{-12pt}\vspace{-1.5pt}
is\dss commutative.\oss
By\sss passing\sss to\sss the geometric realizations and\sss using\trs
Theorem\qss \ref{grassmannian}\qss we conclude\sss that\sss
$\num{\mathcal{G}_{\qff \nnn}\trf (\trf P\trf)}
\qff \ttoo\qff 
\num{\mathcal{G}_{\qff \nnn}\trf (\trf P\fff'\trf)}$\sss
is\dss a homotopy equivalence.\oss
Clearly,\oss the square\vspace{-3pt}
\[
\quad
\begin{tikzcd}[column sep=boom, row sep=boom]
\mathcal{G}_{\qff \nnn}\trf (\trf P\trf)
\arrow[d, "\dis u^{\dff *}_{\qff \nnn}\dff"']
\arrow[r]
&
\mathcal{G}\dff (\trf P\trf)
\arrow[d, "\dis u^{\dff *}\dff"']
\\
\mathcal{G}_{\qff \nnn}\trf (\trf P\fff'\trf)
\arrow[r]
&
\mathcal{G}\dff (\trf P\fff'\trf)\dff,
\end{tikzcd}
\]

\vspace{-12pt}\vspace{-1.5pt}
where\sss the horizontal\sss arrows are\sss the inclusions,\oss
is\dss commutative.\oss
It\dss remains\sss to pass\sss to\sss the geometric realizations
and apply\trs Theorem\qss \ref{framing}.\oss  \eproof\vspace{-0.5pt}

\mypar{Theorem.}{p-is-quasi-fibration}
\emph{The map\dss
$\num{p}\dff \colon\dff
\num{S\dff(\trf \phi\trf)}
\qff \ttoo\qff
\num{\mathcal{P}{\nsp}\hat{\mathcal{S}}}$\dss
is\dss a\sss quasi-fibration.\oss}

\proof
The proof\dss is\dss fairly\sss similar\sss to\sss the proof\dss of\trs 
Theorem\qss \ref{to-models}.\oss
Let\sss $\pi\off =\off \num{p}$\nnsp.\oss
For every\sss $n\qff \in\qff \nnn$\sss and every\sss object\sss $P$\sss of\sss
$\mathcal{P}{\nsp}\hat{\mathcal{S}}$\sss
let\sss us consider\sss the product\sss of\dss categories\sss
$[\halfff n\dff]\dff \times\dff \mathcal{G}\dff (\trf P\trf)$\nnsp.\oss
The category\sss
$[\halfff n\dff]\dff \times\dff \mathcal{G}\dff (\trf P\trf)$\sss
is\dss also associated\sss with a partial\sss order,\oss
namely,\oss with\sss the order\sss
$(\dff i\fff,\qff P\trf)
\qff \leq\qff
(\dff i'\fff,\qff P\fff'\trf)$\sss
if\dss $i\qff \leq\qff i'$\sss and\dss $P\qff \leq\qff P\fff'$\dnsp.\oss
As a\sss partially ordered space\sss
$[\halfff n\dff]\dff \times\dff \mathcal{G}\dff (\trf P\trf)$\sss
is\dss free\sss together with\sss $\mathcal{G}\dff (\trf P\trf)$\sss
and\sss hence can\sss be considered as a\sss topological\sss
simplicial\sss complex.\oss

We will\sss also consider\sss the categories $S\dff(\trf \phi\trf)$ and\sss
$\mathcal{P}{\nsp}\hat{\mathcal{S}}$\sss as\sss topological\sss
simplicial\sss complexes associated\sss with\sss 
the partial\sss orders defining\sss them.\oss
Let\sss us\sss fix an {\nsp}$n$\dnsp-simplex $\sigma$
of\sss $\mathcal{P}{\nsp}\hat{\mathcal{S}}$\dnsp,\oss
and\sss let\sss
$P_{\dff n}\off <\off
P_{\dff n\dff -\dff 1}\off <\off
\ldots\off <\off
P_{\dff 0}$\sss
be\sss the corresponding sequence of\dss objects of\sss 
$\mathcal{P}{\nsp}\hat{\mathcal{S}}$\dnsp.\oss
The numbering of\dss the objects\sss $P_{\dff i}$\sss
reflects\sss the reversal\sss of\dss their order\sss in\sss $S\dff(\trf \phi\trf)$\nnsp.\oss

Let\sss $\tau$\sss be an {\nsp}$l$\dnsp-simplex\sss of\dss 
$S\dff(\trf \phi\trf)$\nnsp.\oss
It\dss is\dss determined\sss by\sss a\sss sequence\dss
$Q_{\dff l}\off \leq\off
Q_{\dff l\trf -\dff 1}\off \leq\off
\ldots\off \leq\off
Q_{\dff 0}$\sss
of\dss objects of\sss 
$\mathcal{P}{\nsp}\hat{\mathcal{S}}$\dnsp,\oss
a\sss sequence\sss
$W_{\dff 0}\off \leq\off
W_{\dff 1}\off \leq\off
\ldots\off \leq\off
W_{\dff l}$\sss
of\dss subspaces adapted\sss to\sss $P_{\dff 0}$\nsp,\oss
and a splitting of\dss $W_{\dff 0}$\nsp,\oss
such\sss that\sss for every\sss
$i\off =\off 1\fff,\qff 2\fff,\qff \ldots\fff,\qff l$\sss
either\sss
$Q_{\dff i\dff -\dff 1}\off >\off Q_{\dff i}$\dss
or\trs
$W_{\dff i\dff -\dff 1}\off <\off W_{\dff i}$\nsp.\oss
The\sss last\sss condition ensures\sss that\sss the corresponding sequence of\dss
objects of\sss $S\dff(\trf \phi\trf)$\sss is\dss strictly\sss increasing.\oss
The functor $p$\nnsp,\oss considered as a simplicial\sss map,\oss
takes\sss the {\nsp}$l$\dnsp-simplex\sss $\tau$\sss to\sss the simplex of\dss
$\mathcal{P}{\nsp}\hat{\mathcal{S}}$
obtained\sss from\sss the non-decreasing\sss sequence\dss
$Q_{\dff l}\off \leq\off
Q_{\dff l\trf -\dff 1}\off \leq\off
\ldots\off \leq\off
Q_{\dff 0}$\dss
by\sss removing\sss repetitions.\oss
Therefore\sss $p\dff(\trf \tau\trf)\off =\off \sigma$\sss
if\trs and\dss only\trs if\dss
removing\sss the repetitions from\sss the sequence\dss
$Q_{\dff l}\off \leq\off
Q_{\dff l\trf -\dff 1}\off \leq\off
\ldots\off \leq\off
Q_{\dff 0}$\dss
results\sss in\sss the sequence\sss
$P_{\dff n}\off <\off
P_{\dff n\dff -\dff 1}\off <\off
\ldots\off <\off
P_{\dff 0}$\nsp.\oss

If\dss $i\qff \in\qff [\halfff n\dff]$\sss and\sss 
$(\qff P_{\dff 0}\fff,\off
W\off =\off W_{\dff -}\dff \oplus\trf W_{\dff +}\qff)$\sss
is\dss an object\sss of\dss
$\mathcal{G}\dff (\trf P_{\dff 0}\trf)$\nnsp,\oss
then\sss $W$\sss is\dss adapted\sss to\sss $P_{\dff i}$\sss
and\sss
$(\qff P_{\dff i}\fff,\off
W\off =\off W_{\dff -}\dff \oplus\trf W_{\dff +}\qff)$\sss
is\dss an object\sss of\trs
$\mathcal{G}\dff (\trf P_{\dff i}\trf)$
and\sss hence an object\sss of\dss
$S\dff(\trf \phi\trf)$\nnsp.\oss
The rule\vspace{1.5pt}
\[
\quad
\bigl(\qff
i\dff,\pff
(\qff P_{\dff 0}\fff,\off
W\off =\off W_{\dff -}\dff \oplus\trf W_{\dff +}\qff)
\qff\bigr)
\off \longmapsto\off
(\qff P_{\dff i}\fff,\off
W\off =\off W_{\dff -}\dff \oplus\trf W_{\dff +}\qff)
\]

\vspace{-12pt}\vspace{1.5pt}
extends\sss to morphisms of\sss
$[\halfff n\dff]\dff \times\dff \mathcal{G}\dff (\trf P_{\dff 0}\trf)$\sss 
and defines a simplicial\sss map\vspace{1.5pt}
\[
\quad
t_{\trf \sigma}\dff \colon\dff
[\halfff n\dff]\dff \times\dff \mathcal{G}\dff (\trf P_{\dff 0}\trf)
\qff \ttoo\qff
S\dff(\trf \phi\trf)
\]

\vspace{-12pt}\vspace{1.5pt}
The reversed\sss partial\sss order of\sss $\mathcal{P}{\nsp}\hat{\mathcal{S}}$\sss
leads\sss to an\sss isomorphism\sss 
$[\halfff n\dff]\qff \ttoo\qff \sigma$
such\sss that\sss the square\vspace{0pt}
\[
\quad
\begin{tikzcd}[column sep=boom, row sep=boomm]
[\halfff n\dff]\dff \times\dff \mathcal{G}\dff (\trf P_{\dff 0}\trf)
\arrow[r, "\dis t_{\trf \sigma}"]
\arrow[d, "\dis \pr\dff"']
&
S\dff(\trf \phi\trf)
\arrow[d, "\dis \phi"]
\\
\protect{[\halfff n\dff]}
\arrow[r]
&
\sigma\dff.
\end{tikzcd}
\]

\vspace{-10.5pt}
is\dss commutative.\oss
Clearly,\pss 
$p\trf(\trf \tau\trf)\off =\off \sigma$\sss
if\dss and\dss only\dss if\dss
$\tau\off =\off t_{\trf \sigma}\trf(\trf \alpha\trf)$\sss
for a simplex $\alpha$ of\dss the product\sss
$[\halfff n\dff]\dff \times\dff \mathcal{G}\dff (\trf P_{\dff 0}\trf)$
such\sss that\sss
$\pr\trf(\trf \alpha\trf)\off =\off [\halfff n\dff]$\nnsp,\oss
and\dss if\trs such $\alpha$ exists,\oss
it\dss is\dss unique.\oss
By\sss passing\sss to\sss the geometric realizations we see\sss
that\sss $\num{t_{\trf \sigma}}$\sss induces a homeomorphism\vspace{1.5pt}
\[
\quad
h_{\dff \sigma}\qff \colon\qff
\pi^{\dff -\dff 1}\qff \bigl(\qff \inte \num{\sigma}\qff\bigr)
\off \ttoo\off
\bigl(\qff \inte \Delta^n\qff\bigr)
\qff \times\qff
\num{\mathcal{G}\dff (\trf P_{\dff 0}\trf)}
\pff,
\]

\vspace{-10.5pt}
and\sss that\sss there\dss is\dss a neighborhood\sss $U_{\dff \sigma}$\sss
of\trs $\partial\trf \num{\sigma}$
in $\num{\sigma}$ 
such\sss that\sss 
$\pi^{\dff -\dff 1}\trf(\trf \partial\trf \num{\sigma}\trf)$\sss
is\dss a deformation\sss retract\sss of\dss
$\pi^{\dff -\dff 1}\trf(\trf U_{\dff \sigma}\trf)$\nnsp.\oss 
In\sss particular,\pss $\pi$\sss is\dss a\sss trivial\sss bundle
over\sss $\inte \Delta^n$\dnsp.\oss

Let\sss us\sss allow $\sigma$ vary among $n$\dnsp-simplices of\dss 
the simplicial\sss complex\sss $\mathcal{P}{\nsp}\hat{\mathcal{S}}$\dnsp.\oss
Clearly,\oss the union of\dss the interiors $\inte \sigma$\sss is\dss
equal\sss to\sss the difference\sss
$\bbnum{\ssk_{\dff n}\dff \mathcal{P}{\nsp}\hat{\mathcal{S}}}
\off \smallsetminus\off
\bbnum{\ssk_{\dff n\dff -\dff 1}\dff \mathcal{P}{\nsp}\hat{\mathcal{S}}}$\nnsp.\oss
Also,\oss the polarized subspace models\dss $P_{\dff 0}$
continuously depends on\sss $\sigma$\nnsp,\oss
as are\sss the maps $\num{t_{\dff \sigma}}$
and\sss homeomorphisms $h_{\dff \sigma}$\nsp.\oss
It\sss follows\sss that\sss the map\vspace{3pt}
\[
\quad
\pi^{\dff -\dff 1}\qff \bigl(\qff 
\bbnum{\ssk_{\dff n}\qff \mathcal{P}{\nsp}\hat{\mathcal{S}}}
\off \smallsetminus\off
\bbnum{\ssk_{\dff n\dff -\dff 1}\qff \mathcal{P}{\nsp}\hat{\mathcal{S}}}
\qff\bigr)
\off \ttoo\off
\bbnum{\ssk_{\dff n}\qff \mathcal{P}{\nsp}\hat{\mathcal{S}}}
\off \smallsetminus\off
\bbnum{\ssk_{\dff n\dff -\dff 1}\qff \mathcal{P}{\nsp}\hat{\mathcal{S}}}
\]

\vspace{-12pt}\vspace{3pt}
induced\sss by\sss
$\pi$\sss is\dss a\sss locally\sss trivial\sss bundle.\oss
Clearly,\oss the neighborhoods\sss $U_{\dff \sigma}$\sss can\sss be chosen\sss
to continuously depend\sss on $\sigma$\nnsp,\oss
as also\sss the deformations of\dss $U_{\dff \sigma}$\sss 
into\sss $\partial\trf \num{\sigma}$\nnsp.\oss
Then\sss the union\sss\vspace{1.75pt} 
\[
\quad
\mathcal{U}_{\dff n\dff -\dff 1}
\off =\off
\bigcup\nolimits_{\dff \sigma}\qff U_{\dff \sigma}
\]

\vspace{-12pt}\vspace{1.5pt}
is\dss an open\sss neighborhood\sss of\dss
$\bbnum{\ssk_{\dff n\dff -\dff 1}\qff \mathcal{P}{\nsp}\hat{\mathcal{S}}}$\dss
in\dss
$\bbnum{\ssk_{\dff n}\qff \mathcal{P}{\nsp}\hat{\mathcal{S}}}$\dss
and\dss the deformations of\dss $U_{\dff \sigma}$\sss 
into\sss $\partial\trf \num{\sigma}$ define\sss a\sss deformation\sss retraction of\trs
$\mathcal{U}_{\dff n\dff -\dff 1}$\sss to\dss
$\bbnum{\ssk_{\dff n\dff -\dff 1}\qff \mathcal{P}{\nsp}\hat{\mathcal{S}}}$\nnsp.\oss
Moreover,\oss this deformation\dss is\dss covered\sss by\sss
a deformation of\dss the preimage\sss
$\pi^{\dff -\dff 1}\qff (\trf \mathcal{U}_{\dff n\dff -\dff 1}\trf)$\sss
into\sss the preimage\vspace{3pt}
\[
\quad
\pi^{\dff -\dff 1}\qff 
\bigl(\qff 
\bbnum{\ssk_{\dff n\dff -\dff 1}\qff \mathcal{P}{\nsp}\hat{\mathcal{S}}}
\qff\bigr)
\pff.
\]

\vspace{-12pt}\vspace{3pt}
Suppose\sss that\sss a point\sss $x\qff \in\qff U_{\dff \sigma}$\sss
is\dss deformed\sss to\sss a point\sss
$y\qff \in\qff \partial\trf \num{\sigma}$\nnsp.\oss
Then\sss $y\qff \in\qff \inte \num{\sigma'}$\sss
for a face $\sigma'$\sss of\sss $\sigma$\nnsp,\oss
i.e.\qss for a simplex $\sigma'$\sss corresponding\sss to a subsequence\vspace{1.75pt}
\[
\quad
P_{\dff i_{\trf m}}\off <\off
P_{\dff i_{\trf m\dff -\dff 1}}\off <\off
\ldots\off <\off
P_{\dff i_{\trf 0}}
\]

\vspace{-12pt}\vspace{1.75pt}
of\dss the sequence\dss
$P_{\dff n}\off <\off
P_{\dff n\dff -\dff 1}\off <\off
\ldots\off <\off
P_{\dff 0}$\nsp.\oss
It\dss is\dss easy\sss to see\sss that\vspace{3pt}
\[
\quad
\pi^{\dff -\dff 1}\dff(\trf x\trf)
\off =\off
\num{\mathcal{G}\dff (\trf P_{\dff 0}\trf)}
\quad
\mbox{and}\dff\quad
\pi^{\dff -\dff 1}\dff(\trf y\trf)
\off =\off
\num{\mathcal{G}\dff (\trf P_{\dff k}\trf)}
\pff,
\]

\vspace{-12pt}\vspace{3pt}
where\sss $k\off =\off i_{\trf 0}$\nsp,\oss
and\sss the map\dss
$\pi^{\dff -\dff 1}\dff(\trf x\trf)
\qff \ttoo\qff
\pi^{\dff -\dff 1}\dff(\trf y\trf)$\dss
induced\sss by\sss the deformation\dss is\dss the map\vspace{1.75pt}
\[
\quad
\num{\mathcal{G}\dff (\trf P_{\dff 0}\trf)}
\qff \ttoo\pff
\num{\mathcal{G}\dff (\trf P_{\dff k}\trf)}
\]

\vspace{-12pt}\vspace{1.75pt}
induced\sss by\sss the unique morphism\sss
$P_{\dff k}\qff \ttoo\qff P_{\dff 0}$\nsp,\oss
which exists because\sss
$P_{\dff k}\qff <\qff P_{\dff 0}$\nsp.\oss
Lemma\qss \ref{fibers-are-he}\qss
implies\sss that\sss this map\dss is\dss a\sss homotopy equivalence.\oss\vspace{1.5pt}

Now\sss we are ready\sss to use\sss the standard\sss way\sss of\dss proving\sss
that\sss a map\dss is\dss a\sss quasi-fibration.\oss
Arguing\sss by\sss induction\sss we can assume\sss that\sss $\pi$\sss
is\dss a quasi-fibration over\sss
$\bbnum{\ssk_{\dff n\dff -\dff 1}\qff \mathcal{P}{\nsp}\hat{\mathcal{S}}}$\nnsp.\oss
The deformation constructed above shows\sss that\sss $\pi$\sss
is\dss a quasi-fibration over\sss $\mathcal{U}_{\dff n\dff -\dff 1}$\nsp.\oss
Over\sss the difference\vspace{3pt}
\[
\quad
\bbnum{\ssk_{\dff n}\qff \mathcal{P}{\nsp}\hat{\mathcal{S}}}
\off \smallsetminus\off
\bbnum{\ssk_{\dff n\dff -\dff 1}\qff \mathcal{P}{\nsp}\hat{\mathcal{S}}}
\]

\vspace{-12pt}\vspace{3pt}
the map\sss $\pi$\sss is\dss a\sss locally\sss trivial\sss bundle and\sss
hence\dss is\dss a quasi-fibration.\oss
Together\sss with\sss the fact\sss that\sss the maps\sss
$\pi^{\dff -\dff 1}\dff(\trf x\trf)
\qff \ttoo\qff
\pi^{\dff -\dff 1}\dff(\trf y\trf)$\sss
are homotopy equivalences\sss this implies\sss that\sss
$\pi$\sss is\dss a quasi-fibration over\sss
$\bbnum{\ssk_{\dff n}\qff \mathcal{P}{\nsp}\hat{\mathcal{S}}}$\nnsp.\oss
It\sss follows\sss that\sss $\pi$\sss is\dss a quasi-fibration over\sss
$\bbnum{\ssk_{\dff n}\qff \mathcal{P}{\nsp}\hat{\mathcal{S}}}$\sss
for every $n$\nnsp.\oss
In\sss turn,\oss this implies\sss that\sss $\pi$\sss
is\dss a quasi-fibration.\oss
See\dss Section\qss \ref{introduction}\qss and\sss the
references\sss in\dss that\dss section.\oss  \eproof

\mypar{Corollary.}{homotopy-fiber-found}
\emph{The homotopy\sss fiber of\trs the map\dss
$\pi
\off =\off
\num{p}\qff \colon\dff
\num{S\dff(\trf \phi\trf)}
\qff \ttoo\qff
\num{\mathcal{P}{\nsp}\hat{\mathcal{S}}}$\dss
is\dss homotopy equivalent\dss to\trs
$\gr$\sss and\dss
$\gr\trf(\dff \infty\dff)$\nnsp.\oss}

\proof
The spaces\sss $\gr$\sss and\sss $\gr\trf(\dff \infty\dff)$
are homotopy equivalent\sss by\trs Theorem\qss \ref{two-grassmannians}.\oss
Let\sss us\sss consider an object\sss $P$ of\dss
$\mathcal{P}{\nsp}\hat{\mathcal{S}}$\sss
as a point\sss in\sss $\num{\mathcal{P}{\nsp}\hat{\mathcal{S}}}$\nnsp.\oss
Clearly,\oss the fiber\sss $\pi^{\dff -\dff 1}\dff(\trf P\trf)$\sss
is\dss nothing else but\sss the geometric realization\sss
$\num{\mathcal{G}\dff (\trf P\trf)}$\sss
of\dss $\mathcal{G}\dff (\trf P\trf)$\nnsp.\oss
But\dss by\trs Theorem\qss \ref{s-sigma-homotopy-type}\qss
$\num{\mathcal{G}\dff (\trf P\trf)}$\sss
is\dss homotopy equivalent\sss to\sss
$\gr\trf(\dff \infty\dff)$\nnsp.\oss
Since\sss $\pi$\sss is\dss a quasi-fibration,\oss
it\sss follows\sss that\sss the homotopy fiber\dss 
is\dss weakly homotopy equivalent\sss to\sss
$\gr\trf(\dff \infty\dff)$\nnsp.\oss
Since\sss the homotopy fiber and\sss
$\gr\trf(\dff \infty\dff)$ are homotopy equivalent\sss
to\dss CW-complexes,\oss
this implies\sss the\sss theorem.\oss  \eproof

\mypar{Theorem.}{category-contractibility}
\emph{The classifying\sss space\sss
$\num{S\dff(\trf \phi\trf)}$\sss
is\dss contractible.\oss}

\proof
There\dss is\dss a forgetting\sss functor\sss
$\varphi\dff \colon\dff
S\dff(\trf \phi\trf)
\qff \ttoo\qff
s\dff \hat{\mathcal{S}}$\sss
taking an object\sss
$P\qff \ttoo\qff M\off \longleftarrow\off N$
of\dss $S\dff(\trf \phi\trf)$\sss
to\sss the object\sss
$N\qff \ttoo\qff M$
of\sss $s\dff \hat{\mathcal{S}}$\nnsp,\oss
and\sss taking a morphism\qss (\ref{sub-morphisms})\qss
to\sss the morphism\vspace{-6.5pt}\vspace{-1.25pt}\vspace{-0.375pt}
\[
\quad
\begin{tikzcd}[column sep=boom, row sep=shqboom]
&
M
\arrow[dd, "\dis w"]
\\
N
\arrow[ru]
\arrow[rd]
&
\\
&
M\fff'\qff.
\end{tikzcd}
\]

\vspace{-12pt}\vspace{-4.5pt}\vspace{-1.25pt}\vspace{-0.375pt}
For an\sss object\sss
$s\dff \colon\dff N\qff \ttoo\qff M$\sss
of\dss $s\dff \hat{\mathcal{S}}$\dss
let\sss $\mathcal{F}\trf(\dff s\trf)$\sss
be\sss the full\sss subcategory of\dss $S\dff(\trf \phi\trf)$
having as objects\sss diagrams of\dss the form\dss
$P\qff \ttoo\qff M\off \longleftarrow\off N$\nnsp,\oss
where\sss the right\sss arrow\dss is\sss $s$\nnsp.\oss
The category\sss $\mathcal{F}\trf(\dff s\trf)$\sss
has a\sss terminal\sss object,\oss
namely,\oss the object\sss
$M\qff \ttoo\qff M\off \longleftarrow\off N$\nnsp,\oss
where\sss the\sss left\sss arrow\dss is\dss the identity
and\sss right\sss arrow\dss is\sss $s$\nnsp.\oss
In\sss particular,\pss
$\num{\mathcal{F}\trf(\dff s\trf)}$\sss
is\dss contractible.\oss
Arguing as in\sss the proof\dss of\trs 
Theorem\qss \ref{to-models}\qss one can\sss prove\sss that\sss
$\num{\varphi}\dff \colon\dff
\num{S\dff(\trf \phi\trf)}
\qff \ttoo\qff
\num{s\dff \hat{\mathcal{S}}}$\sss
is\dss a homotopy equivalence.\oss
Since\sss $\num{s\dff \hat{\mathcal{S}}}$\sss is\dss contractible\sss
by\trs Theorem\qss \ref{split-contractible},\oss 
this proves\sss the\sss theorem.\oss  \eproof

\myuppar{Remark.}
In\sss the above proof\dss one can also argue as in\sss the proof\dss of\trs
Theorem\qss \ref{p-is-quasi-fibration}\qss and\sss prove\sss that\sss
$\num{\varphi}$\sss
is\dss a quasi-fibration.\oss
Since\sss the spaces\sss $\num{\mathcal{F}\trf(\dff s\trf)}$\sss
are fibers of\sss $\num{\varphi}$\sss and are contractible,\oss
it\sss follows\sss that\sss $\num{\varphi}$\sss is\dss
a weak\sss homotopy equivalence.\oss
Similarly\sss to our other classifying spaces,\oss
the spaces\sss $\num{S\dff(\trf \phi\trf)}$\sss
and\sss $\num{s\dff \hat{\mathcal{S}}}$\sss are 
homotopy equivalent\sss to\sss CW-comp\-lex\-es.\oss
This again\sss implies\sss that\sss
$\num{\varphi}$\sss is\dss a homotopy equivalence
and\sss hence proves\sss the\sss theorem.\oss

\myuppar{Comparing\sss the bundle $\bm{\pi}$ and\sss 
the quasi-fibration\sss $\num{p}$\nnsp.}
Let\sss $P$\sss be an object\sss of\sss $\mathcal{P}{\nsp}\hat{\mathcal{S}}$\nnsp.\oss
There\dss is\dss a canonical\sss map\sss
$\num{\mathcal{G}\dff(\trf P\trf)}
\qff \ttoo\qff
\gr\trf(\trf P\trf)$\nnsp,\oss
which can\sss be defined as follows.\oss
Consider\sss $\gr\trf(\trf P\trf)$ as a\sss topological\sss category\sss
having\sss $\gr\trf(\trf P\trf)$ as\sss the space of\dss objects
and only\sss identity morphisms.\oss
An object\sss of\dss $\mathcal{G}\dff(\trf P\trf)$
defines a polarized subspace model\sss
$(\trf W\fff,\qff K_{\dff -}\fff,\qff K_{\dff +}\dff)$\sss
together with a splitting\sss
$W\off =\off W_{\dff -}\dff \oplus\trf W_{\dff +}$\nsp.\oss
Let\sss us assign\sss to such an object\sss the subspace\sss
$K_{\dff -}\dff \oplus\dff W_{\dff -}
\qff \in\qff 
\gr\trf(\trf P\trf)$\nnsp.\oss
A\sss trivial\sss verification shows\sss that\sss
this rule assigns\sss the same subspace\sss to\sss two objects
of\dss $\mathcal{G}\dff(\trf P\trf)$ related\sss by a morphism.\oss
Hence we can extend\sss this rule\sss to morphisms 
of\dss $\mathcal{G}\dff(\trf P\trf)$\sss
by assigning\sss to every morphism an\sss identity morphism 
of\dss $\gr\trf(\trf P\trf)$ and\sss
get\sss a functor\sss
$\mathcal{G}\dff(\trf P\trf)
\qff \ttoo\qff
\gr\trf(\trf P\trf)$\nnsp.\oss
By\sss passing\sss to\sss the geometric realizations we get\sss a map\sss
$h\trf(\trf P\trf)\dff \colon\dff
\num{\mathcal{G}\dff(\trf P\trf)}
\qff \ttoo\qff
\num{\gr\trf(\trf P\trf)}
\off =\off
\gr\trf(\trf P\trf)$\nnsp.\oss

The proof\dss of\trs Theorem\qss \ref{p-is-quasi-fibration}\qss shows\sss that\sss
the fiber of\sss $\num{p}$ over a point\sss 
$x\qff \in\qff \num{\mathcal{P}{\nsp}\hat{\mathcal{S}}}$ contained\sss in\sss
the interior of\dss a simplex\sss represented\sss by\sss a sequence\sss
$P_{\dff n}\off <\off
P_{\dff n\dff -\dff 1}\off <\off
\ldots\off <\off
P_{\dff 0}$\sss
is\dss equal\sss to\sss $\num{\mathcal{G}\dff (\trf P_{\dff 0}\trf)}$\nnsp.\oss
At\sss the same\sss time\sss the fiber of\sss $\bm{\pi}$ over $x$\sss
is\dss equal\sss to\sss $\gr\trf(\trf P_{\dff 0}\trf)$\sss by\sss the definition.\oss
Hence $h\trf(\trf P_{\dff 0}\trf)$\sss is\dss a map\sss
$\pi^{\dff -\dff 1}\dff(\trf x\trf)
\qff \ttoo\qff
\num{p}^{\dff -\dff 1}\dff(\trf x\trf)$\nnsp.\oss
Together\sss these maps define a map\sss
$h\dff \colon\dff
\num{S\dff(\trf \phi\trf)}
\qff \ttoo\qff
\mathbf{G}$\sss
such\sss that\sss
$\bm{\pi}\dff \circ\dff h\off =\off \num{p}$\nnsp,\oss
i.e.\qss a map from\sss the quasi-fibration $\num{p}$\sss to\sss 
the\sss bundle $\bm{\pi}$\nnsp.\oss
By\trs Theorem\qss \ref{s-sigma-homotopy-type}\qss
the map $h$\sss induces homotopy equivalences of\dss fibers.\oss
Comparing\sss the homotopy sequences of\dss 
$\num{p}$\sss and\sss $\bm{\pi}$\sss 
shows\sss that\sss $h$\sss is\dss a weak\sss homotopy equivalence.\vspace{0.75pt}

\mypar{Lemma.}{g-is-cw}
\emph{The space\dss $\mathbf{G}$ is\dss homotopy equivalent\sss to a\dss CW-complex.\oss}\vspace{0.75pt}

\proof
The space\sss
$\num{\mathcal{P}{\nsp}\hat{\mathcal{S}}}$\sss
is\dss homotopy equivalent\sss to a\dss CW-complex\sss
by\sss the same reasons as\sss the spaces from\trs Theorem\qss \ref{framing}.\oss
It\dss is\dss easy\sss to see\sss that\sss the skeletons\sss
$\bbnum{\ssk_{\dff n}\qff \mathcal{P}{\nsp}\hat{\mathcal{S}}}$\sss
are metrizable spaces.\oss
Since\sss  
$\num{\mathcal{P}{\nsp}\hat{\mathcal{S}}}
\off =\off
\bbnum{\ssk_{\dff n}\qff \mathcal{P}{\nsp}\hat{\mathcal{S}}}$\sss 
is\dss the direct\sss limit\sss of\dss skeletons\sss
$\bbnum{\ssk_{\dff n}\qff \mathcal{P}{\nsp}\hat{\mathcal{S}}}$\nnsp,\oss
it\dss is\dss paracompact\dss by a\sss theorem of\qss Morita\qss \cite{mo}.\oss
This implies\sss that\sss $\bm{\pi}$\sss is\dss a\dss Hurewicz\dss fibration\qss
(see\qss \cite{td3},\oss Theorem\qss 14.3.5)\qss
and,\oss in\sss particular,\oss is\dss an $h$\dnsp-fibration.\oss
By\trs Lemma\qss \ref{gr-is-cw}\qss the fibers\sss $\gr$\sss of\sss $\bm{\pi}$ are homotopy
equivalent\sss to\dss CW-complexex.\oss
By a\sss theorem of\trs tom\dss Dieck\qss \cite{td1}\qss
this implies\sss that\sss the\sss total\sss space\sss $\mathbf{G}$\sss
has\sss the homotopy\sss type of\dss a\dss CW-complex.\oss
See also\qss \cite{td3},\oss Section\qss 13.4,\oss
Problem\qss 2.\oss  \eproof\vspace{0.75pt}

\mypar{Theorem.}{g-is-contractible}
\emph{The space\dss $\mathbf{G}$ is\dss contractible.\oss}\vspace{0.75pt}

\proof
Since $\mathbf{G}$\sss is\dss weakly\sss homotopy equivalent\sss 
to\sss $\num{S\dff(\trf \phi\trf)}$\nnsp,\oss
Theorem\qss \ref{category-contractibility}\qss implies\sss that\sss $\mathbf{G}$\sss
is\dss weakly contractible.\oss
In view of\trs Lemma\qss \ref{g-is-cw}\qss this implies\sss that\sss $\mathbf{G}$\sss
is\dss contractible.\oss  \eproof\vspace{0.75pt}

\mypar{Theorem\qss ({\fff}Bott{\fff}).}{loops}
\emph{The\sss loop space of\qss $U\trf(\dff \infty\dff)$
is\dss homotopy equivalent\sss to\dss 
$\gr\trf(\dff \infty\dff)$\nnsp.\oss}\vspace{0.75pt}

\proof
Theorems\qss \ref{forgetting-equivalence}\qss and\qss \ref{harris-h}\qss
imply\sss that\sss
$\num{\mathcal{P}{\nsp}\hat{\mathcal{S}}}$\sss
is\dss homotopy\sss equivalent\sss to\sss 
$U\ffin$\dnsp.\oss
Together\sss with\dss Corollary\qss \ref{classical-subspace-unitary}\qss
this implies\sss that\sss
$\num{\mathcal{P}{\nsp}\hat{\mathcal{S}}}$\sss
is\dss homotopy equivalent\sss to\sss $U\trf(\dff \infty\dff)$\nnsp.\oss
Theorem\qss \ref{g-is-contractible}\qss implies\sss that\sss
the\sss loop space of\sss
$\num{\mathcal{P}{\nsp}\hat{\mathcal{S}}}$\sss 
is\dss homotopy equivalent\sss 
to\sss the fiber $\gr$ of\sss $\bm{\pi}$\nnsp.\oss
Finally,\oss by\trs Theorem\qss \ref{two-grassmannians}\qss
the spaces\sss $\gr$ and\sss $\gr\trf(\dff \infty\dff)$
are homotopy equivalent.\oss  \eproof\vspace{0.75pt}

\myuppar{Bott\dss periodicity.}
Theorem\qss \ref{g-is-contractible}\qss is\dss one of\dss
the most\sss classical\sss forms of\dss the\trs
Bott\dss periodicity.\oss
The\dss Grassmannian\sss $\gr\trf(\dff \infty\dff)$\sss
is\dss one of\dss the standard\sss forms of\dss the classifying space\sss
$\zzz\dff \times\dff B\fff U$\nnsp.\oss
At\sss the same\sss time $U$ in\sss the\sss last\sss formula\dss
is\dss nothing else but\sss our\sss $U\trf(\dff \infty\trf)$\nnsp.\oss 
This gives us another classical\sss
form of\trs Bott\dss periodicity\fff:\oss
the\sss loop space\sss of\sss $U$\sss is\dss homotopy
equivalent\sss to\sss $\zzz\dff \times\dff B\fff U$\nnsp.\oss
This proof\dss of\qss Bott\trs periodicity\dss is\dss partially\dss inspired\dss by\trs
Atiyah--Singer\trs proof\pss \cite{as},\oss
but,\oss in contract\sss with\sss the\sss latter,\oss
is\dss independent\sss of\dss any\sss results about\sss spaces of\trs
Fredholm\dss operators.\vspace{0.75pt}

The key\sss steps in\dss this proof\dss of\qss Bott\trs periodicity\dss are\dss
Theorem\qss \ref{p-is-quasi-fibration}\qss and\trs
Corollary\qss \ref{homotopy-fiber-found}.\oss 
They are\sss inspired\sss by\trs Theorem\qss B\qss of\qss Quillen\qss \cite{q},\oss
and\sss their\sss proofs\sss follow\dss Quillen's\dss ideas.\oss
Similarly,\oss the proof\dss of\trs Theorem\qss \ref{category-contractibility}\qss
is\dss inspired\sss by\sss the proof\dss of\trs 
Theorem\qss A\qss of\qss Quillen\qss \cite{q}.\oss

\myuppar{The\dss Grassmannian\sss bundle over $\hat{\mathcal{F}}$\dnsp.}
Suppose\sss that\sss $A\qff \in\qff \hat{\mathcal{F}}$\dnsp.\oss
Then\sss there exists $\varepsilon\qff >\qff 0$\sss such\sss that\sss
$(\trf A\fff,\qff \varepsilon\trf)$\sss is\dss an enhanced\sss operator.\oss 
In\sss particular,\pss $\varepsilon\qff \not\in\qff \sigma\trf(\trf A\trf)$\sss and\vspace{1.5pt}
\[
\quad
H
\off =\off
\image P_{\dff \leq\dff \varepsilon}\dff(\trf A\trf)
\qff \oplus\qff
\image P_{\dff \geq\dff \varepsilon}\dff(\trf A\trf)
\]

\vspace{-12pt}\vspace{1.5pt}
is\dss a polarization.\oss
Clearly,\oss the restricted\dss Grassmannian\sss corresponding\sss to\sss this
polarizations does not\sss depend on\sss the choice of\sss $\varepsilon$\nnsp.\oss
We will\sss denote\sss it\sss by\sss $\gr\trf(\trf A\trf)$\nnsp.\oss
Let\sss $\mathbf{F}$\sss be\sss the space of\dss pairs\sss
$(\trf A\fff,\qff K\trf)$\sss 
such\sss that\sss
$A\qff \in\qff \hat{\mathcal{F}}$\sss
and\dss $K\qff \in\qff \gr\trf(\trf A\trf)$\nnsp,\oss
and\sss let\sss\vspace{1.5pt} 
\[
\quad
\mathbf{p}\dff \colon\dff
\mathbf{F}
\qff \ttoo\qff 
\hat{\mathcal{F}}
\]

\vspace{-12pt}\vspace{1.5pt}
be\sss the projection\sss
$(\trf A\fff,\qff K\trf)
\off \longmapsto\off
A$\nnsp.\oss
The stability of\dss half-line projections\qss (see\dss Section\qss \ref{spaces-operators})\qss
implies\sss that\sss
$\mathbf{p}\dff \colon\dff
\mathbf{F}
\qff \ttoo\qff 
\hat{\mathcal{F}}$\sss
is\dss a\sss locally\sss trivial\sss bundle with\sss the fiber $\gr$\nnsp.\oss
As we will\sss see now,\oss this bundle\dss is\dss essentially\sss
the bundle\sss 
$\bm{\pi}\dff \colon\dff
\mathbf{G}
\qff \ttoo\qff 
\num{\mathcal{P}{\nsp}\hat{\mathcal{S}}}$\sss
transferred\sss to $\hat{\mathcal{F}}$\nnsp.\oss

By\trs Theorem\qss \ref{forgetting-enhancement}\qss
the canonical\sss map\sss
$\num{\hat{\mathcal{E}}}
\qff \ttoo\qff
\hat{\mathcal{F}}$\sss 
is\dss a homotopy equivalence.\oss
Let\sss us\sss consider\sss the composition\sss
$\hat{\mathcal{E}}
\qff \ttoo\qff
\mathcal{P}{\nsp}\hat{\mathcal{S}}$\sss
of\dss forgetting\sss functors\dss
$\mathcal{P}\fff\psi$\nnsp,\qss
$\mathcal{P}m$\nnsp,\oss
and\sss
$\mathcal{P}\varphi$\dss
from\dss Section\qss \ref{polarizations-splittings}.\oss
By\trs Theorem\qss \ref{forgetting-equivalence}\qss the induced\sss map\sss
$\num{\hat{\mathcal{E}}}
\qff \ttoo\qff
\num{\mathcal{P}{\nsp}\hat{\mathcal{S}}}$\sss
is\dss a homotopy equivalence.\oss

\mypar{Theorem.}{bundles-f-ps}
\emph{The\sss two bundles over\dss $\num{\hat{\mathcal{E}}}$\sss
induced\sss from\sss the bundles}\vspace{3pt}
\[
\quad
\mathbf{p}\dff \colon\dff
\mathbf{F}
\qff \ttoo\qff 
\hat{\mathcal{F}}
\quad
\mbox{\emph{and}}\quad\dff
\bm{\pi}\dff \colon\dff
\mathbf{G}
\qff \ttoo\qff 
\num{\mathcal{P}{\nsp}\hat{\mathcal{S}}}
\]

\vspace{-12pt}\vspace{3pt}
\emph{by\sss the maps\sss
$\num{\hat{\mathcal{E}}}
\qff \ttoo\qff
\hat{\mathcal{F}}$\sss
and\qss
$\num{\hat{\mathcal{E}}}
\qff \ttoo\qff
\num{\mathcal{P}{\nsp}\hat{\mathcal{S}}}$\sss
respectively,\oss
are\sss the same.\oss}

\proof
Since we are using\sss the discrete\sss topology\sss for\sss the controlling
parameters $\varepsilon$ of\dss enhanced operators,\oss
the partially ordered space $\hat{\mathcal{E}}$\sss has free equalities.\oss
By\dss Corollary\qss \ref{free-full-realizations}\qss the geometric realizations\sss
$\num{\hat{\mathcal{E}}}$\sss and\sss $\bbnum{\hat{\mathcal{E}}}$\sss
are canonically\sss homeomorphic.\oss

Let\sss $x\qff \in\qff \bbnum{\hat{\mathcal{E}}}$\nnsp.\oss
Then $x$ can\sss be represented\sss 
by\sss a weighted sums\vspace{1.5pt}
\[
\quad
t_{\trf 0}\trf E_{\dff 0}\pff +\pff
t_{\dff 1}\trf E_{\dff 1}\pff +\pff
\ldots\pff +\pff
t_{\dff n}\trf E_{\dff n}
\pff,
\]

\vspace{-12pt}\vspace{1.5pt}
where\sss 
$(\trf t_{\dff 0}\dff,\qff t_{\dff 1}\dff,\qff \ldots\dff,\qff t_{\dff n}\trf)
\qff \in\qff
\Delta^n$\dss
and\dss 
$E_{\dff 0}\fff,\qff
E_{\dff 1}\fff,\qff
\ldots\fff,\qff
E_{\dff n}$\sss 
are objects of\dss 
$\hat{\mathcal{E}}$
such\sss that\vspace{1.5pt}
\[
\quad
E_{\dff 0}\off <\off
E_{\dff 1}\off <\off
\ldots\off <\off
E_{\dff n}
\pff.
\]

\vspace{-12pt}\vspace{1.5pt}
These inequalities imply\sss that\sss there exists $A\qff \in\qff \hat{\mathcal{F}}$\sss
and\sss positive numbers\sss
$\varepsilon_{\dff 0}\qff <\qff
\varepsilon_{\dff 1}\qff <\qff
\ldots\qff <\qff
\varepsilon_{\dff n}$
such\sss that\dss 
$E_{\dff i}
\off =\off 
(\trf A\fff,\qff \varepsilon_{\dff i}\trf)$\sss
for every $i$\nnsp.\oss
A direct\sss verification shows\sss that\sss the fibers over\sss the point\sss $x$ of\dss
the bundles over\sss $\num{\hat{\mathcal{E}}}$\sss 
induced\dss from\sss the bundles $\mathbf{p}$ and\sss $\bm{\pi}$\sss
are both equal\dss to $\gr\trf(\trf A\trf)$\nnsp.\oss
It\sss follows\sss that\sss the induced\sss 
bundles are\sss the same.\oss  \eproof

\mypar{Corollary.}{f-tot-contractible.}
\emph{The space\sss $\mathbf{F}$\sss is\dss contractible.\oss}  \eproof

\mysection{Classifying\qss spaces\qss for\qss 
Fredholm\qss operators}{classifying-spaces-odd-saf}

\myuppar{Fredholm\dss and\sss self-adjoint\dss Fredholm\dss operators.}
It\dss is\dss technically\sss convenient\sss to\sss
interpret\dss Fredholm\sss operators in\sss $H$\sss
as\qss \emph{odd}\pss self-adjoint\dss Fredholm\sss operators in\sss
$H\dff \oplus\dff H$\sss and\sss then
adapt\dss the\sss technique of\trs Section\qss \ref{classifying-spaces-saf}\qss
to\sss such operators.\oss
Let\sss us\sss begin\sss with\sss the main definitions.\oss

Let\dss 
$\gamma\dff \colon\dff 
H\dff \oplus\dff H
\qff \ttoo\qff
H\dff \oplus\dff H$\sss
be\sss the operator acting as\sss $\id$\sss on\sss $H\dff \oplus\dff 0$
and as\sss $-\qff \id$\sss on\sss $0\dff \oplus\dff H$\nnsp.\oss
Let\sss $V$\sss be a $\gamma$\dnsp-invariant\sss subspace of\sss $H$\nnsp.\oss
Then\sss $V\pff =\pff V_{\dff -}\dff \oplus\dff V_{\dff +}$\sss
for some\sss $V_{\dff -}\dff,\off V_{\dff +}\qff \subset\qff H$\nnsp.\oss
An\qss \emph{odd\sss operator}\qss
$B\dff \colon\dff V\qff \ttoo\qff V$\sss
is\dss defined as an operator 
anti-commuting\sss with\sss $\gamma$\nnsp,\oss
i.e.\qss is\dss such\sss that\sss
$B\dff \circ\dff \gamma\off =\off -\qff \gamma\dff \circ\dff B$\dnsp.\oss
An odd\sss operator\qss 
$B\dff \colon\dff 
H\dff \oplus\dff H
\qff \ttoo\qff
H\dff \oplus\dff H$
can\sss be represented\sss by a matrix
\[
\begin{pmatrix}
0 & \qff A\qff
\\
\pff A'\dff & 0
\end{pmatrix}
\qff \colon\qff
H\dff \oplus\dff H
\qff \ttoo\qff
H\dff \oplus\dff H
\]

\vspace{-12pt}
for some operators\sss
$A\fff,\qff A'\dff \colon\dff H\qff \ttoo\qff H$\nnsp.\oss
Such operator\sss is\dss Fredholm\dss if\trs and\dss only\trs if\trs
both\sss $A\fff,\qff A'$\sss are\dss Fredholm,\oss
and\dss is\dss self-adjoint\dss if\trs and\dss only\trs if\dss
$A'\off =\off A^*$\dnsp,\oss
i.e.\qss $A'$\sss is\dss the adjoint\sss operator of\dss $A$\nnsp.\oss
Let\sss us\sss assign\sss to an operator\sss
$A\dff \colon\dff H\qff \ttoo\qff H$\sss the operator
\[
\quad
A^\sa
\off =\off
\begin{pmatrix}
0 & \qff A\qff
\\
\pff A^* & 0
\end{pmatrix}
\qff \colon\qff
H\dff \oplus\dff H
\qff \ttoo\qff
H\dff \oplus\dff H
\pff.
\]

\vspace{-12pt}
Then\sss the map\sss $A\qff \longmapsto\qff A^\sa$\sss
establishes a\sss homeomorphism\sss between\sss the space $\mathcal{F}$
of\trs bounded\dss Fredholm\sss operators\sss $H\qff \ttoo\qff H$\sss
and\sss the space $\hat{\mathcal{F}}^{\dff \odd}$ of\dss 
odd\sss self-adjoint\dss bounded\trs Fredholm\dss operators\sss
$H\dff \oplus\dff H
\qff \ttoo\qff
H\dff \oplus\dff H$\nnsp.\oss
We will\sss use\sss this map\sss
to identify\sss these\sss two spaces.\oss

\myuppar{The eigenspaces of\sss $A^\sa$\dnsp.}
Suppose\sss that\sss $x\qff \in\qff H\dff \oplus\dff H$\sss
and\sss 
$A^\sa\dff (\dff x\trf)
\off =\off
\lambda\dff x$\sss
for some $\lambda\qff \in\qff \rrr$\nnsp.\oss
Then\vspace{-0.55pt}
\[
\quad
A^\sa\trf \bigl(\qff \gamma\dff (\dff x\trf)\dff \bigr)
\off =\off
-\qff \gamma\qff \bigl(\qff A^\sa\dff (\dff x\trf)\dff \bigr)
\off =\off
-\qff \gamma\qff (\trf \lambda\dff x\trf) 
\off =\off
-\qff \lambda\qff \gamma\dff(\dff x\trf)
\pff.
\]

\vspace{-12pt}\vspace{-0.55pt}
It\dss follows\sss that\sss if\dss $x$\sss is\dss an eigenvector
corresponding\sss to\sss the eigenvalue $\lambda$\nnsp,\oss
then\dss $\gamma\dff (\dff x\trf)$\sss is\dss an eigenvector
corresponding\sss to\sss the eigenvalue $-\qff \lambda$\nnsp.\oss
In\sss particular,\oss the kernel\sss $\kernel A^\sa$\sss
is\sss $\gamma$\dnsp-invariant.\oss
In\sss fact,\oss obviously,\qss
$\kernel A^\sa
\off =\off
\kernel A\dff \oplus\dff \kernel A^*$\dnsp.\oss

Suppose now\sss that\sss
$x
\off =\off
(\dff u\fff,\qff v\trf)
\qff \in\qff 
H\dff \oplus\dff H$\nnsp,\oss
and\sss 
$A^\sa\dff (\dff x\trf)
\off =\off
\lambda\dff x$\sss
with\sss $\lambda\off \neq\off 0$\nnsp.\oss
Then\sss\vspace{-0.55pt}
\[
\quad 
A^\sa\dff (\dff u\fff,\qff v\trf)
\off =\off
\left(\trf A\trf(\dff v\trf)\fff,\qff  A^*\dff(\dff u\trf)\trf\right)
\]

\vspace{-12pt}\vspace{-0.55pt}
and\sss hence\sss
$A\trf(\dff v\trf)\off =\off \lambda\dff u$\nnsp,\pss
$A^*\dff(\dff u\trf)\off =\off \lambda\dff v$\nnsp.\oss
It\sss follows\sss that\sss
\[
\quad
A\trf A^*\dff(\dff u\trf)
\off =\off 
\lambda^{\fff 2}\dff u
\quad
\mbox{and}\quad
A^* A\trf(\dff v\trf)\off =\off \lambda^{\fff 2}\dff v
\pff.
\]

\vspace{-12pt}
Therefore\sss $v$\sss is\dss an eigenvector\sss of\dss the self-adjoint\sss operator\sss
$A^* A$\sss with\sss the eigenvalue\sss $\lambda^{\fff 2}$\sss and\vspace{-0.55pt}
\[
\quad
(\dff u\fff,\qff v\trf)
\off =\off 
\left(\trf
\lambda^{\fff -\dff 1}\dff A\trf(\dff v\trf)\dff,\qff v
\trf\right)
\pff.
\]

\vspace{-12pt}
Conversely,\dff\oss if\qss
$A^* A\trf(\dff v\trf)\off =\off \lambda^{\fff 2}\dff v$\nnsp,\oss
then\dss\vspace{0.5pt}
\[
\quad
\left(\trf
\lambda^{\fff -\dff 1}\dff A\trf(\dff v\trf)\dff,\qff v
\trf\right)
\quad
\mbox{and}\quad
\left(\trf
\lambda^{\fff -\dff 1}\dff A\trf(\dff v\trf)\dff,\qff -\qff v
\trf\right)
\]

\vspace{-12pt}\vspace{0.5pt}
are eigenvectors of\sss $A^\sa$\sss with\sss the eigenvalues\sss
$\lambda$\sss and\sss $-\qff \lambda$\sss respectively,\oss
and\sss $\gamma$\sss interchanges\sss them.\oss
These observations\sss take a nice form\sss in\sss terms of\dss
polar\sss decompositions.\oss
Recall\dss that\sss the\qss 
\emph{polar\sss decomposition}\pss of\sss $A$\sss is\dss the
unique presentation\sss of\dss the form\sss
$A\off =\off U\trf \num{A}$\nnsp,\oss
where\vspace{0.5pt}
\[
\quad
\num{A}
\off =\off
\sqrt{\dff A^* A\dff}
\]

\vspace{-12pt}\vspace{0.5pt}
and\sss $U$\sss is\dss a partial\sss isometry\sss
with\sss $\kernel U\off =\off \kernel A$\sss
and\dss
$\image U\off =\off \image A$\dss
(since $A$\sss is\trs Fredholm,\oss
the image\sss $\image A$\sss is\dss closed\fff).\oss
Then\sss $\lambda$\sss is\dss an eigenvalue of\sss $A^\sa$\sss if\trs and\dss only\trs if\dss
$\num{\lambda}$\sss is\dss an eigenvalue of\sss $\num{A}$\nnsp.\oss
If\dss $\lambda\qff \geq\qff 0$\sss and\sss
$\num{A}\qff(\dff v\trf)\off =\off \lambda\dff v$\nnsp,\oss
then\vspace{0.5pt}
\[
\quad
\bigl(\trf
U\trf(\dff v\trf)\dff,\qff v
\trf\bigr)
\quad
\mbox{and}\quad
\bigl(\trf
U\trf(\dff v\trf)\dff,\qff -\qff v
\trf\bigr)
\]

\vspace{-12pt}\vspace{0.5pt}
are eigenvectors of\sss $A^\sa$ with\sss the eigenvalues\sss
$\lambda$ and $-\qff \lambda$ respectively.\oss
Moreover,\oss every eigenvector of\sss $A^\sa$\sss
with non-zero eigenvalue has\sss this form.\oss

\myuppar{Categories related\dss to odd self-adjoint\trs Fredholm\dss operators.}
The categories\sss $\hat{\mathcal{E}}$\dnsp,\qss 
$\mathcal{E}\hat{\mathcal{O}}$\dnsp,\qss
$\hat{\mathcal{O}}$\dnsp,\oss
and\sss $\hat{\mathcal{S}}$\sss have natural\sss analogues\sss in\sss the context\sss
of\dss odd self-adjoint\trs Fredholm\dss operators.\oss
Let\sss us\sss
define\qss \emph{enhanced\sss odd}\qss (self-adjoint\trs Fredholm)\qss operator\dss as\dss
a\sss pair $(\trf A\dff,\qff \varepsilon\trf)$\nnsp,\oss
where $A\qff \in\qff \hat{\mathcal{F}}^{\dff \odd}$
is\dss an odd\sss self-adjoint\dss bounded\trs Fredholm\dss operator,\oss and 
$\varepsilon\qff \in\qff \rrr$\sss
is\dss such\sss that\sss $0\qff <\qff \varepsilon$\nnsp,\oss the\sss interval\sss
$[\dff -\qff \varepsilon\fff,\qff 
\varepsilon\trf]$\sss 
is\dss disjoint\sss from\sss the essential\sss
spectrum of\sss $A$\nnsp,\oss
and\sss 
$-\qff \varepsilon\fff,\qff 
\varepsilon\qff \not\in\qff \sigma\dff(\trf A\trf)$\nnsp.\oss
The space $\hat{\mathcal{E}}^{\dff \odd}$ of\dss enhanced odd operators\dss
is\dss ordered\sss by\sss the relation\sss
$(\dff A\dff,\qff \varepsilon\trf)
\qff \leq\qff
(\dff A'\dff,\qff \varepsilon'\trf)$\sss
if\trs $A\off =\off A'$\sss and\sss 
$\varepsilon
\qff \leq\qff 
\varepsilon'$\nnsp.\oss
As usual,\oss we consider\sss this order as\sss a structure of\dss a\sss topological\sss
category on $\hat{\mathcal{E}}^{\dff \odd}$\dnsp.\oss
There\dss is\dss an obvious forgetting\sss functor\sss
$\hat{\varphi}^{\dff \odd}\dff \colon\dff 
\hat{\mathcal{E}}^{\dff \odd}
\qff \ttoo\qff 
\hat{\mathcal{F}}^{\dff \odd}$\dnsp,\oss
where $\hat{\mathcal{F}}^{\dff \odd}$\sss is\dss considered as a\sss topological\sss
category\sss having only\sss identity\sss morphisms.\oss

An\qss \emph{odd\dss enhanced\dss operator\dss model}\pss is\dss a\sss triple\sss
$(\trf V,\qff F\dff,\qff \varepsilon\trf)$\nnsp,\oss
where\sss $V$\sss is\dss a $\gamma$\dnsp-invariant\sss finitely 
dimensional\sss subspace of\sss $H\dff \oplus\dff H$\nnsp,\oss
$F\dff \colon\dff V\qff \ttoo\qff V$\sss is\dss an odd\sss self-adjoint\sss operator,\oss
and $\varepsilon$\sss is\dss a real\sss number
such\sss that\sss $0\qff <\qff \varepsilon$ and\sss 
$\sigma\dff(\trf F\trf)
\qff \subset\qff
(\dff -\qff \varepsilon\fff,\qff \varepsilon\trf)$\nnsp.\oss
The space\sss $\mathcal{E}\hat{\mathcal{O}}^{\dff \odd}$\sss
of\dss odd\sss enhanced operator\sss models\dss
is\dss ordered\sss by\sss the relation\sss 
$\leq$\sss 
defined exactly\sss as\sss the order of\dss enhanced operators models
was defined\sss in\dss Section\qss \ref{classifying-spaces-saf}.\oss
As usual,\oss we consider\sss this order as\sss a structure of\dss a\sss topological\sss
category on $\mathcal{E}\hat{\mathcal{O}}^{\dff \odd}$\nnsp.\oss
There\dss is\dss an obvious functor\sss
$\hat{\psi}^{\dff \odd}\dff \colon\dff 
\hat{\mathcal{E}}^{\dff \odd}
\qff \ttoo\qff 
\mathcal{E}\hat{\mathcal{O}}^{\dff \odd}$\sss
taking an enhanced\sss odd\sss operator\sss
$(\dff A\dff,\qff \varepsilon\trf)$\sss
to\qss \emph{its\sss odd\sss enhanced\sss operator\sss model}\qss
$(\trf V,\qff F\dff,\qff \varepsilon\trf)$\nnsp,\oss
where\vspace{1.5pt}
\[
\quad
V
\off =\off 
\image\dff P_{\dff [\dff -\qff \varepsilon\fff,\qff \varepsilon\trf]}\trf(\trf A \trf)
\off =\off 
\image\dff P_{\dff (\dff -\qff \varepsilon\fff,\qff \varepsilon\trf)}\trf(\trf A \trf)
\]

\vspace{-12pt}\vspace{1.5pt}
and\sss the operator\sss 
$F\dff \colon\dff V\qff \ttoo\qff V$\sss
is\dss induced\sss by $A$\nnsp.\oss

An\qss \emph{odd\dss operator\dss model}\oss is\dss a\sss pair\sss
$(\trf V\fff,\pff F\trf)$\nnsp,\oss
where\sss $V$\sss is\dss a $\gamma$\dnsp-invariant\sss finitely 
dimensional\sss subspace of\sss $H\dff \oplus\dff H$\sss and\sss
$F\dff \colon\dff V\qff \ttoo\qff V$\sss 
is\dss an odd\sss self-adjoint\sss operator.\oss
Let\sss $\hat{\mathcal{O}}^{\dff \odd}$\sss 
be\sss the set\sss of\dss odd\sss operator\sss models\dss
with\sss the obvious\sss topology  
and ordered\sss by\sss the relation\dss $\leq$\sss 
defined exactly\sss as\sss the order of\dss operators models
was defined\sss in\dss Section\qss \ref{classifying-spaces-saf}.\oss
There\dss is\dss an obvious forgetting functor\sss
$\hat{o}^{\dff \odd}\dff \colon\dff 
\mathcal{E}\hat{\mathcal{O}}^{\dff \odd}
\qff \ttoo\qff 
\hat{\mathcal{O}}^{\dff \odd}$\dss
taking\sss enhanced operator model\sss
$(\trf V\fff,\pff F\fff,\pff \varepsilon\trf)$\sss
to\sss the operator\sss model\dss
$(\trf V\fff,\pff F\trf)$\nnsp.\oss

A\qss \emph{odd\sss subspace\sss model}\pss is\dss defined
as a {\nsp}$\gamma$\dnsp-invariant\sss finitely dimensional\sss subspace of\sss 
$H\dff \oplus\dff H$\nnsp.\oss
A\qss \emph{morphism\sss of\dss odd\sss subspace models}\dss $V\qff \ttoo\qff V\fff'$\sss 
is\dss defined as a subspace\sss $U\qff \subset\qff H\dff \oplus\dff H$\sss
such\sss that\vspace{1.5pt}
\[
\quad
V\fff'
\off =\off
\gamma\trf(\trf U \trf)\qff \oplus\qff
V\qff \oplus\qff U
\qff.
\]

\vspace{-12pt}\vspace{1.5pt}
The composition of\dss morphisms\dss is\dss defined\sss by\sss taking\sss the sum of\dss
the corresponding subspaces $U$\dnsp.\oss 
The\sss topology on\sss the set\sss of\dss
morphisms\dss is\dss defined\sss in\sss the obvious manner.\oss
This defines a\sss topological\sss category\sss $\hat{\mathcal{S}}^{\dff \odd}$\sss
having odd\sss subspace models as\sss objects.\oss
Let\sss us\sss assign\sss to an odd\sss operator\sss model\sss 
$(\trf V,\qff F\dff,\qff \varepsilon\trf)$\sss 
the odd\sss subspace model\sss
$V$\dnsp,\oss
and\sss to a morphism\vspace{1.5pt}
\[
\quad
(\trf V,\qff F\dff,\qff \varepsilon\trf)
\qff \ttoo\qff
\left(\trf V\fff',\qff F\fff'\dff,\qff \varepsilon'\trf\right)
\qff
\]

\vspace{-12pt}\vspace{1.5pt}
of\sss odd\sss operator models\sss
the morphism of\dss odd\sss subspace models
defined\sss by\sss the subspace\vspace{0.625pt}
\[
\quad
U
\off =\off
\image\dff 
P_{\dff [\dff \varepsilon\fff,\qff \varepsilon'\trf]}\dff\left(\trf F\fff'\trf\right)
\qff.
\]

\vspace{-12pt}\vspace{0.625pt}
The symmetry of\dss eigenvectors\sss implies\sss that\sss
\[
\quad
\image\dff 
P_{\dff [\dff -\qff \varepsilon'\fff,\qff -\qff \varepsilon\trf]}\dff\left(\trf F\fff'\trf\right)
\off =\off
\gamma\trf\left(\trf
\image\dff 
P_{\dff [\dff \varepsilon\fff,\qff \varepsilon'\trf]}\dff\left(\trf F\fff'\trf\right)
\trf\right)
\]

\vspace{-12pt}
and\sss hence\sss
$V\fff'
\off =\off
\gamma\trf(\trf U \trf)\qff \oplus\qff
V\qff \oplus\qff U$\dnsp,\oss
i.e.\qss the subspace $U$\sss indeed defines a morphism.\oss
Clearly,\oss these rules define\sss a\sss forgetting\sss functor\sss
$\hat{\omega}^{\dff \odd}\dff \colon\dff
\hat{\mathcal{O}}^{\dff \odd}
\qff \ttoo\qff
\hat{\mathcal{S}}^{\dff \odd}$\dnsp.\oss

\mypar{Lemma.}{non-sa-contractible}
\emph{Suppose\sss that\sss $a\qff >\qff 0$\nnsp.\oss
The space of\trs bounded\sss operators\sss
$A\dff \colon\dff H\qff \ttoo\qff H$\sss
such\sss that\qss
$\norm{A\trf(\dff v\trf)}
\qff >\qff
a\trf \norm{v}$\qss
for every\sss $v\off \neq\off 0$\sss
is\dss contractible.\oss}

\proof
Every such operator $A$\sss is\dss invertible,\oss
and\sss the map\sss
$A\off \longmapsto\off A^{\dff -\dff 1}$\sss
is\dss a homeomorphism\sss between our space and\sss the space of\dss
invertible bounded operators\sss
$B\dff \colon\dff H\qff \ttoo\qff H$\sss
such\sss that\qss
$\norm{B\trf(\dff v\trf)}
\qff <\qff
a^{\dff -\dff 1}\trf \norm{v}$\qss
for every $v\qff \in\qff H$\nnsp.\oss
Let\sss us choose some\sss $b\qff >\qff 0$\sss such\sss that\sss $b\qff <\qff a^{\dff -\dff 1}$\dnsp.\oss
Using\sss the polar\sss decomposition one can deform\sss the space of\dss
such operators $B$\sss into\sss the subspace of\dss operators $B$\sss such\sss that\sss 
$\norm{B\trf(\dff v\trf)}
\off =\off
b\trf \norm{v}$\qss
for every $v\qff \in\qff H$\nnsp.\oss
The\sss latter\dss is\dss homeomorphic\sss to\sss the
unitary\sss group of\sss $H$\sss and\dss is\dss contractible
by\trs Kuiper's\dss theorem.\oss 
The lemma\sss follows.\oss  \eproof

\mypar{Lemma.}{odd-contractibility-prelim}
\emph{Suppose\sss that\trs $\varepsilon\qff >\qff 0$\nnsp.\oss
The space of\trs odd\dss bounded\dss self-adjoint\trs Fredholm\dss operators\qss
$B\dff \colon\dff
H\dff \oplus\dff H
\qff \ttoo\qff
H\dff \oplus\dff H$\dss
such\sss that\trs
$\sigma\dff(\trf B\trf)
\dff \cap\trf 
[\dff -\qff \varepsilon\fff,\qff \varepsilon\trf]
\off =\off
\varnothing$\dss
is\dss contractible.\oss}

\proof
Every\sss such operator $B$ has\sss the form\sss
$B\off =\off A^\sa$\sss for a unique bounded\dss Fredholm\sss operator\sss
$A\dff \colon\dff H\qff \ttoo\qff H$\nnsp.\oss
The square\sss $B^{\dff 2}$\sss is\dss a positive operator and\dss 
is\dss represented\sss by\sss the diagonal\sss
matrix of\dss operators with\sss the diagonal\sss entries
$A\dff A^*$\sss and $A^* A$\nnsp.\oss
The condition\sss imposed on\sss the spectrum\sss 
$\sigma\dff(\trf B\trf)$\sss
implies\sss that\sss
$\sigma\dff(\trf B^{\dff 2}\trf)
\dff \cap\trf 
[\dff 0\fff,\qff \varepsilon^{\dff 2}\trf]
\off =\off
\varnothing$\sss
and\sss hence\sss\vspace{1.5pt}
\[
\quad
\sco{B^{\dff 2}\dff(\dff x\trf)\fff,\qff x}
\off >\off
\varepsilon^{\dff 2}\dff
\sco{x\fff,\qff x}
\]

\vspace{-12pt}\vspace{1.5pt}
for every\sss 
$x\qff \in\qff H\dff \oplus\dff H$\nnsp,\qss
$x\off \neq\off 0$\nnsp.\oss
It\sss follows\sss that\sss\vspace{1.5pt}
\[
\quad\sco{A\trf(\dff v\trf)\fff,\qff A\trf(\dff v\trf)}
\off =\off
\sco{A^* A\trf(\dff v\trf)\fff,\qff v}
\off =\off
\sco{B^{\dff 2}\dff(\dff v\trf)\fff,\qff v}
\off >\off
\varepsilon^{\dff 2}\dff
\sco{v\fff,\qff v}
\]

\vspace{-10.5pt}
for every\sss $v\qff \in\qff H$\nnsp,\qss $v\off \neq\off 0$\nnsp,\oss
where we identify\sss $H$\sss with\sss $H\dff \oplus\dff 0$\nnsp.\oss
Hence\sss\vspace{1.5pt}
\[
\quad
\norm{A\trf(\dff v\trf)}
\qff >\qff
\varepsilon\qff \norm{v}
\]

\vspace{-10.5pt}
for every\sss $v\off \neq\off 0$\nnsp.\oss
In\sss particular,\pss $A$\sss is\dss an\sss invertible operator.\oss
Conversely,\oss if\sss $A$\sss is\dss a bounded\sss invertible 
operator such\sss that\sss
$\norm{A\trf(\dff v\trf)}
\qff >\qff
\varepsilon\qff \norm{v}$\sss
for every\sss $v\off \neq\off 0$\nnsp,\oss
then\sss
$B\off =\off A^\sa$\sss
is\dss an odd\sss self-adjoint\trs Fredholm\dss operator such\sss that\sss
$\sco{B^{\dff 2}\dff(\dff x\trf)\fff,\qff x}
\qff >\qff
\varepsilon^{\dff 2}\dff
\sco{x\fff,\qff x}$\sss
for every\sss $x\off \neq\off 0$\nnsp.\oss
The last\sss condition\sss implies\sss that\sss
$\sigma\dff(\trf B\trf)
\dff \cap\trf 
[\dff -\qff \varepsilon\fff,\qff \varepsilon\trf]
\off =\off
\varnothing$\nnsp.\oss

It\sss follows\sss that\sss the space of\dss operators from\sss the\sss lemma\dss
is\dss homeomorphic\sss to\sss the space of\dss bounded operators $A$ such\sss that\sss
$\norm{A\trf(\dff v\trf)}
\qff >\qff
\varepsilon\qff \norm{v}$\sss
for every\sss $v\off \neq\off 0$\nnsp.\oss
But\sss the\sss latter space\dss is\dss contractible by\trs
Lemma\qss \ref{non-sa-contractible}.\oss  \eproof

\mypar{Lemma.}{odd-fredholm-contractibility}
\emph{Let\sss
$(\trf V,\qff F\dff,\qff \varepsilon\trf)$\sss
be an odd\sss operator\sss model.\oss
The space of\dss bounded enhanced\sss odd\sss operators\sss $(\dff A\dff,\qff \varepsilon\trf)$\sss
such\sss that\dss
$\hat{\psi}^{\dff \odd}\trf(\dff A\dff,\qff \varepsilon\trf)
\off =\off
(\trf V,\qff F\dff,\qff \varepsilon\trf)$\sss
is\dss contractible.\oss}

\proof
This space can\sss be identified\sss with\sss the space of\dss
odd\sss self-adjoint\trs Fredholm\dss operators\vspace{1.5pt}
\[
\quad
B\dff \colon\dff
(\trf H\dff \oplus\dff H\trf)
\dff \ominus\dff
V
\qff \ttoo\qff
(\trf H\dff \oplus\dff H\trf)
\dff \ominus\dff
V
\]

\vspace{-10.5pt}
such\sss that\sss
$\sigma\dff(\trf B\trf)
\dff \cap\trf 
[\dff -\qff \varepsilon\fff,\qff \varepsilon\trf]
\off =\off
\varnothing$\nnsp.\oss
Therefore\sss the\sss lemma\sss follows\sss from\trs
Lemma\qss \ref{odd-contractibility-prelim}\qss
applied\sss to\sss
 $(\trf H\dff \oplus\dff H\trf)
\dff \ominus\dff
V$\sss
in\sss the role of\sss $H\dff \oplus\dff H$\nnsp.\oss  \eproof

\mypar{Theorem.}{forgetting-odd}
\emph{For\sss bounded operators\sss the\sss maps\sss of\dss
classifying\sss spaces}\dss\vspace{1.25pt}
\[
\quad
\num{\hat{\mathcal{E}}^{\dff \odd}}
\qff \ttoo\qff 
\num{\hat{\mathcal{F}}^{\dff \odd}}\quad
\mbox{and}\quad
\]

\vspace{-38pt}
\[
\quad
\num{\hat{\mathcal{E}}^{\dff \odd}}
\qff \ttoo\qff 
\num{\mathcal{E}\hat{\mathcal{O}}^{\dff \odd}}
\qff \ttoo\qff 
\num{\hat{\mathcal{O}}^{\dff \odd}}
\qff \ttoo\qff
\num{\hat{\mathcal{S}}^{\dff \odd}}
\]

\vspace{-10.75pt}
\emph{induced\dss by\dss the\sss functors defined\sss above,\oss
are\sss homotopy equivalences.\oss}

\proof
The proof\dss that
$\num{\hat{\mathcal{E}}^{\dff \odd}}
\qff \ttoo\qff 
\num{\hat{\mathcal{F}}^{\dff \odd}}$
is\dss a\sss homotopy\sss equivalence\dss is\dss
completely similar\sss to\sss the proof\dss of\qss
Theorem\qss \ref{forgetting-enhancement}.\oss
The proof\dss that\sss
$\num{\hat{\mathcal{E}}^{\dff \odd}}
\qff \ttoo\qff 
\num{\mathcal{E}\hat{\mathcal{O}}^{\dff \odd}}$\sss
is\dss a\sss homotopy\sss equivalence\dss is\dss
completely similar\sss to\sss the proof\dss of\qss
Theorem\qss \ref{to-enhanced-models}.\oss
One needs only\sss to use\dss Lemma\qss \ref{odd-fredholm-contractibility}\qss
instead of\trs Proposition\qss \ref{invertible-contractible}.\oss
The proof\dss that\sss
$\num{\mathcal{E}\hat{\mathcal{O}}^{\dff \odd}}
\qff \ttoo\qff 
\num{\hat{\mathcal{O}}^{\dff \odd}}$\sss
is\dss a\sss homotopy\sss equivalence\dss is\dss
similar\sss to\sss the proof\dss of\qss
Theorem\qss \ref{to-models}.\oss
Finally,\oss the proof\dss that
$\num{\hat{\mathcal{O}}^{\dff \odd}}
\qff \ttoo\qff 
\num{\hat{\mathcal{S}}^{\dff \odd}}$
is\dss a\sss homotopy\sss equivalence\dss is\dss
similar\sss to\sss the proof\dss of\qss
Theorem\qss \ref{forgetting-operators}.\oss
One only\sss needs\sss to keep in\sss mind\dss that\sss in\sss
the odd case\sss the positive part\sss of\dss a morphism
determines\sss its negative part.\oss  \eproof

\myuppar{Fredholm\dss subspace models.}
A\qss \emph{Fredholm\dss subspace model}\oss is\dss a pair\sss
$(\trf E_{\dff 1}\fff,\qff E_{\dff 2}\trf)$ of\dss finitely\sss
dimensional\sss subspaces of\sss $H$\nnsp.\oss
Let\sss $\mathcal{S}$\sss be\sss the category\sss having\trs such\sss pairs\sss
as\sss objects,\oss
with\sss morphisms\sss
$(\trf E_{\dff 1}\fff,\qff E_{\dff 2}\trf)
\qff \ttoo\qff 
(\trf E\fff'_{\dff 1}\fff,\qff E\fff'_{\dff 2}\trf)$\sss
being\sss pairs of\dss subspaces\sss
$(\trf F_{\dff 1}\fff,\qff F_{\dff 2}\trf)$
such\sss that\vspace{1.5pt}
\[
\quad
E_{\dff 1}\dff \oplus\dff F_{\dff 1}
\off =\off
E\fff'_{\dff 1}
\quad
\mbox{and}\quad
E_{\dff 2}\dff \oplus\dff F_{\dff 2}
\off =\off
E\fff'_{\dff 2}
\]

\vspace{-12pt}\vspace{1.5pt}
together with an\sss isometry\sss
$f\dff \colon\dff
F_{\dff 1}\qff \ttoo\qff F_{\dff 2}$\nsp.\oss
The composition\dss is\dss defined\sss by\sss taking\sss the direct\sss sums
of\dss the corresponding subspaces\sss $F_{\dff 1}\fff,\qff F_{\dff 2}$
and of\dss the isometries.\oss
The category $\mathcal{S}$\sss is\dss a\sss topological\sss category\sss 
in\sss an obvious way.\oss
The category\sss $\mathcal{S}$\sss appears\sss to be a more intuitive
analogue of\trs the category\sss $\hat{\mathcal{S}}$\sss
than\sss the category\sss $\hat{\mathcal{S}}^{\dff \odd}$\dnsp.\oss 
Our next\sss goal\dss is\dss to relate\sss the category\sss
$\hat{\mathcal{S}}^{\dff \odd}$ with\sss the category\sss $\mathcal{S}$\sss
of\qss Fredholm\dss subspace models
by constructing\sss functors\sss\vspace{1.5pt}
\[
\quad
p\dff \colon\dff
\hat{\mathcal{S}}^{\dff \odd}
\qff \ttoo\qff 
\mathcal{S}
\quad
\mbox{and}\quad\qff
q\dff \colon\dff
\mathcal{S}
\qff \ttoo\qff 
\hat{\mathcal{S}}^{\dff \odd}
\pff.
\]

\vspace{-12pt}\vspace{1.5pt}
\myuppar{The functor\sss 
$p\dff \colon\dff
\hat{\mathcal{S}}^{\dff \odd}\qff \ttoo\qff \mathcal{S}$\nsp\dnsp.}
Let\dss us\sss begin\sss with defining\sss the action of\sss $p$ on objects.\oss
Let\sss $V$\sss be an object\sss of\sss $\hat{\mathcal{S}}^{\dff \odd}$\dnsp,\oss
i.e.\qss a {\nsp}$\gamma$\dnsp-invariant\sss 
finitely dimensional\sss subspace of\sss $H\dff \oplus\dff H$\nnsp.\oss
The operator\sss $\gamma$\sss induces a self-adjoint\sss unitary operator\sss
$V\qff \ttoo\qff V$\sss
with eigenvalues equal\sss to $-\qff 1$ or $1$
and\sss hence induces a decomposition\sss
$V\off =\off V_{\dff -\dff 1}\dff \oplus\trf V_{\dff 1}$ of\dss $V$\sss
into\sss the corresponding eigenspaces.\oss
Clearly,\qss 
$V_{\dff -\dff 1}
\off =\off
V\dff \cap\dff (\trf 0\dff \oplus\dff H\trf)$
and\sss
$V_{\dff 1}
\off =\off
V\dff \cap\dff (\trf H\dff \oplus\dff 0\trf)$\nnsp.\oss
By\sss identifying\sss $0\dff \oplus\dff H$
and\sss $H\dff \oplus\dff 0$\sss with\sss $H$\sss we can consider\sss
$V_{\dff -\dff 1}$ and\sss $V_{\dff 1}$ as subspaces of\dss $H$\nnsp.\oss
The functor $p$\sss assigns\sss to\sss $V$\sss 
the object\sss 
$(\trf V_{\dff -\dff 1}\fff,\qff V_{\dff 1}\trf)$ 
of\dss $\mathcal{S}$\nsp\dnsp.\oss

In order\sss to define\sss the action of\sss $p$ on\sss morphisms,\oss
let\sss us\sss consider a morphism\sss
$V\qff \ttoo\qff V\fff'$\sss 
of\sss $\hat{\mathcal{S}}^{\dff \odd}$\nsp\dnsp.\dff\oss 
It\dss corresponds\sss to\sss an orthogonal\sss decomposition\sss
$V\fff'
\off =\off
\gamma\trf(\trf U \trf)\qff \oplus\qff
V\qff \oplus\qff U$\nnsp.\oss
The orthogonal\sss complement\sss 
$W\off =\off V\fff'\dff \ominus\dff V$\sss
is\dss $\gamma$\dnsp-invariant\sss and\dss hence\dss is\dss a direct\sss sum\sss
$W\off =\off W_{\dff -\dff 1}\dff \oplus\trf W_{\dff 1}$ of\dss eigenspaces 
of\dss the operator\sss $W\qff \ttoo\qff W$\sss induced\sss by\sss
$\gamma$\nnsp.\oss
Clearly,\oss\vspace{2.5pt}
\[
\quad
V\fff'_{\fff -\dff 1}
\off =\off
V_{\dff -\dff 1}\dff \oplus\dff W_{\dff -\dff 1}
\quad
\mbox{and}\dff\quad
V\fff'_{\fff 1}
\off =\off
V_{\dff 1}\dff \oplus\dff W_{\dff 1}
\qff.
\]

\vspace{-12pt}\vspace{2.5pt}
At\sss the same\sss time\sss
$W\off =\off \gamma\trf(\trf U \trf)\dff \oplus\dff U$\dnsp.\oss
By\sss the definition of\dss morphisms,\pss
$\gamma\trf(\trf U \trf)$\sss is\dss orthogonal\sss to $U$\dnsp.\oss
Therefore\sss
$(\dff u\fff,\qff v\trf)\fff,\off
(\dff u'\fff,\qff v'\trf)
\qff \in\qff U$\sss
implies\sss that\vspace{1.5pt}
\[
\quad
\sco{u\fff,\qff u'\dff}
\qff -\qff
\sco{v\fff,\qff v'\dff}
\off =\off
\sco{\fff(\dff u\fff,\qff v\trf)\fff,\pff
(\dff u'\fff,\qff -\qff v'\trf)\fff}
\]

\vspace{-37.5pt}\vspace{-0.25pt}
\[
\quad
\phantom{\sco{u\fff,\qff u'\dff}
\qff -\qff
\sco{v\fff,\qff v'\dff}
\off }
=\off
\sco{\fff(\dff u\fff,\qff v\trf)\fff,\pff
\gamma\trf(\dff u'\fff,\qff v'\trf)\fff}
\off =\off
0
\]

\vspace{-12pt}\vspace{1.5pt}
and\sss hence\sss
$\sco{u\fff,\qff u'\dff}
\off =\off
\sco{v\fff,\qff v'\dff}$\dnsp.\oss
In\sss particular,\oss
if\dss $(\dff u\fff,\qff v\trf)\qff \in\qff U$\dnsp,\oss
then\sss
$\sco{u\fff,\qff u\dff}
\off =\off
\sco{v\fff,\qff v\dff}$\dnsp.\oss
It\sss follows\sss that\sss $U$\sss is\dss the graph
of\dss an\sss isometric\sss isomorphism\sss
between\sss two subspaces of\dss $H$\nnsp,\oss
namely,\oss between\sss the images of\dss the projections of\dss $U$\sss
onto\sss the summands\sss $W_{\dff -\dff 1}$ and\sss $W_{\dff 1}$ of\dss
$W\off =\off W_{\dff -\dff 1}\dff \oplus\trf W_{\dff 1}$\nnsp.\oss
Clearly,\oss the images of\dss the projections of\dss $U$\sss and of\dss $\gamma\dff(\trf U\trf)$
are equal\sss and\sss hence are equal\sss to\sss the images of\dss the projections of\dss
$W\off =\off W_{\dff -\dff 1}\dff \oplus\trf W_{\dff 1}$\nsp.\oss
It\dss follows\sss that\sss $U$\sss is\dss the graph of\dss an\sss
isometric\sss isomorphism\sss
$W_{\dff -\dff 1}\qff \ttoo\qff W_{\dff 1}$\nnsp.\oss

Now\sss we are ready\sss to define\sss the action of\dss $p$ on\sss morphisms.\oss
Namely,\oss to\sss the morphism\sss
$V\qff \ttoo\qff V\fff'$\sss
as\sss in\sss the previous paragraph\sss the functor\sss $p$\sss
assigns\sss the morphism\vspace{0pt}
\[
\quad
\bigl(\trf V_{\dff -\dff 1}\fff,\qff V_{\dff 1}\trf\bigr)
\qff \ttoo\qff
\bigl(\trf V\fff'_{\fff -\dff 1}\fff,\qff V\fff'_{\fff 1}\trf\bigr)
\]

\vspace{-12pt}
defined\sss by\sss the pair\sss
$(\trf W_{\dff -\dff 1}\fff,\qff W_{\dff 1}\trf)$
and\sss the isometric\sss isomorphism\sss
$W_{\dff -\dff 1}\qff \ttoo\qff W_{\dff 1}$
constructed above.\oss
A direct\sss verification shows\sss that\sss $p$\sss is\dss 
indeed a\sss functor.\oss

\myuppar{The functor\sss 
$q\dff \colon\dff
\mathcal{S}\qff \ttoo\qff \hat{\mathcal{S}}^{\dff \odd}$\nsp\dnsp.}
The functor\sss 
$q$\sss
assigns\sss to an object\sss
$(\trf E_{\dff 1}\fff,\qff E_{\dff 2}\trf)$ of\dss $\mathcal{S}$\sss
the $\gamma$\dnsp-invariant\sss subspace\sss
$E_{\dff 1}\dff \oplus\dff E_{\dff 2}
\pff \subset\pff
H\dff \oplus\dff H$\nnsp.\oss
In order\sss to define\sss the action of\sss $q$\sss on\sss morphisms,\oss
let\dss us\sss consider an\sss isometry\sss
$f\dff \colon\dff
F_{\dff 1}\qff \ttoo\qff F_{\dff 2}$\sss
between\sss two finitely dimensional\sss subspaces of\dss $H$\nnsp.\oss
Let\sss $U\qff \subset\qff F_{\dff 1}\dff \oplus\dff F_{\dff 2}$\dss
be\sss the graph of\sss $f$\dnsp.\oss
Since\sss $f$\sss is\dss injective,\qss
$\gamma\dff(\trf U\trf)\dff \cap\trf U
\off =\off
0$\nnsp,\oss
and\sss hence\sss
$F_{\dff 1}\qff \oplus\qff F_{\dff 2}
\off =\off
\gamma\dff(\trf U\trf)\qff +\qff U$\dnsp.\oss
Moreover,\oss
$\gamma\dff(\trf U\trf)$ and\sss $U$ are orthogonal.\oss
Indeed,\oss if\dss $u\fff,\qff u'\qff \in\qff U$
and\sss
$v\off =\off f\dff(\dff u\trf)$\nnsp,\qss
$v'\off =\off f\dff(\dff u'\trf)$\nnsp,\oss
then\sss
$\sco{u\fff,\qff u'\dff}
\off =\off
\sco{v\fff,\qff v'\dff}$\sss
and\sss hence\vspace{0pt}
\[
\quad
\sco{\fff(\dff u\fff,\qff v\trf)\fff,\pff
\gamma\trf(\dff u'\fff,\qff v'\trf)\fff}
\off =\off
\sco{\fff(\dff u\fff,\qff v\trf)\fff,\pff
(\dff u'\fff,\qff -\qff v'\trf)\fff}
\]

\vspace{-37.5pt}
\[
\quad
\phantom{\sco{\fff(\dff u\fff,\qff v\trf)\fff,\pff
\gamma\trf(\dff u'\fff,\qff v'\trf)\fff}
\off }
=\off
\sco{u\fff,\qff u'\dff}
\qff -\qff
\sco{v\fff,\qff v'\dff}
\off =\off
0
\qff.
\]

\vspace{-12pt}
It\sss follows\sss that\sss
$\gamma\dff(\trf U\trf)$\sss is\dss orthogonal\sss to $U$
and\sss hence\sss
$F_{\dff 1}\qff \oplus\qff F_{\dff 2}
\off =\off
\gamma\dff(\trf U\trf)\dff \oplus\dff U$\dnsp.\oss

Now\sss we are ready\sss to define\sss the action of\dss $q$\sss on\sss morphisms.\oss
Suppose\sss that\sss\vspace{1.5pt}
\[
\quad
\mu\dff \colon\dff
(\trf E_{\dff 1}\fff,\qff E_{\dff 2}\trf)
\qff \ttoo\qff 
(\trf E\fff'_{\dff 1}\fff,\qff E\fff'_{\dff 2}\trf)
\]

\vspace{-10.5pt}
is\dss a\sss morphism
of\dss $\mathcal{S}$\dnsp.\oss
By\sss the definition,\qss
$\mu$\sss corresponds\sss to a pair of\dss subspaces\sss
$(\trf F_{\dff 1}\fff,\qff F_{\dff 2}\trf)$
and an\sss isometry\sss
$f\dff \colon\dff
F_{\dff 1}\qff \ttoo\qff F_{\dff 2}$\sss
such\sss that\vspace{3pt}
\[
\quad
E\fff'_{\dff 1}\dff \oplus\dff E\fff'_{\dff 2}
\off =\off
\bigl(\trf
E_{\dff 1}\dff \oplus\dff F_{\dff 1}
\trf\bigr)
\qff \oplus\qff
\bigl(\trf
E_{\dff 2}\dff \oplus\dff F_{\dff 2}
\trf\bigr)
\off =\off
\bigl(\trf
E_{\dff 1}\dff \oplus\dff E_{\dff 2}
\trf\bigr)
\qff \oplus\qff
\bigl(\trf
F_{\dff 1}\dff \oplus\dff F_{\dff 2}
\trf\bigr)
\pff.
\]

\vspace{-9pt}
By\sss the previous paragraph\sss
$F_{\dff 1}\qff \oplus\qff F_{\dff 2}
\off =\off
\gamma\dff(\trf U\trf)\dff \oplus\dff U$\dnsp,\oss
where\sss $U$\sss is\dss the graph of\sss $f$\dnsp.\oss
Hence\vspace{3pt}
\[
\quad
E\fff'_{\dff 1}\dff \oplus\dff E\fff'_{\dff 2}
\off =\off
\gamma\trf(\trf U\trf)
\dff \oplus\dff
\bigl(\trf
E_{\dff 1}\dff \oplus\dff E_{\dff 2}
\trf\bigr)
\dff \oplus\dff
U
\]

\vspace{-9pt}
and\sss this decomposition defines a morphism\sss\vspace{1.5pt}\vspace{-0.125pt}
\[
\quad
E_{\dff 1}\dff \oplus\dff E_{\dff 2}
\qff \ttoo\qff
E\fff'_{\dff 1}\dff \oplus\dff E\fff'_{\dff 2}
\]

\vspace{-10.5pt}
of\dss $\hat{\mathcal{S}}^{\dff \odd}$\dnsp.\oss
We\sss define\sss $q\dff(\trf \mu\trf)$\sss
as\sss this morphism.\oss
A routine verification shows\sss that\sss these rules indeed define a functor.\oss
We\sss leave\sss the routine proof\dss of\dss the following\sss theorem\sss to\sss the reader.\oss

\mypar{Theorem.}{two-fredholm}
\emph{The\dss functor\dss
$p$\sss and\sss $q$\sss
are mutually\sss inverse\sss isomorphisms\sss of\dss categories.\oss}  \eproof\vspace{1.25pt}

\mypar{Theorem.}{s-fredholm}
\emph{The classifying\sss space\sss
$\num{\mathcal{S}}$\sss
is\dss canonically\sss homotopy\sss equivalent\dss to\sss
the space\sss $\mathcal{F}$\sss of\pss bounded\trs Fredholm\dss
operators\sss
$H\qff \ttoo\qff H$\nnsp.\oss}\vspace{1.25pt}

\proof
Since $\hat{\mathcal{F}}^{\dff \odd}$ is\dss homeomorphic\sss to
$\mathcal{F}$\dnsp\dnsp,\oss
this follows from\dss Theorems\dss \ref{forgetting-odd}\dss
and\qss \ref{two-fredholm}.\oss  \eproof\vspace{1.25pt}

\myuppar{Fredholm\dss operators\sss and\sss the category\sss $\mathcal{S}$\nsp\dnsp.}
Let\sss us\sss describe more directly\sss the relation\sss between\sss
$\mathcal{F}$\sss and\sss $\mathcal{S}$\sss arising\sss from\sss
the homeomorphism\sss
$\mathcal{F}\qff \ttoo\qff \hat{\mathcal{F}}^{\dff \odd}$\sss
and\sss the isomorphism\sss
$\hat{\mathcal{S}}^{\dff \odd}\qff \ttoo\qff \mathcal{S}$\nsp\dnsp.\oss
To\sss this end\sss we need an analogue\sss $\mathcal{E}$\sss of\dss
the category\sss $\hat{\mathcal{E}}$\sss of\dss enhanced\sss self-adjoint\sss operators.\oss
The definition of\dss $\hat{\mathcal{S}}^{\dff \odd}$\sss and\dss the above discussion 
of\dss eigenspaces and eigenvalues of\sss $A^\sa$\sss
suggest\sss the following definitions.\oss

An\qss \emph{enhanced\trs Fredholm\dss operator}\pss is\dss
a pair\sss $(\trf A\dff,\qff \varepsilon\trf)$\nnsp,\oss
where $A\qff \in\qff \mathcal{F}$ and $\varepsilon\qff \in\qff \rrr$\sss
are such\sss that\sss $\varepsilon\qff >\qff 0$\nnsp,\oss the\sss interval\sss
$[\trf 0\fff,\qff \varepsilon\trf]$\sss is\dss disjoint\sss from\sss the essential\sss
spectrum of\sss $\num{A}$\nnsp,\oss
and\sss $\varepsilon\qff \not\in\qff \sigma\dff(\trf \num{A}\trf)$\nnsp.\oss
Let\sss $\mathcal{E}$ be\sss the space\sss of\dss enhanced\trs Fredholm\dss operators.\oss
The\sss topology\dss is\dss defined\sss by\sss
the\sss topology of\sss $\mathcal{F}$\sss
and\sss the discrete\sss topology on\sss $\rrr$\nnsp.\oss
The space\sss $\mathcal{E}$\sss 
is\dss ordered\sss by\sss the relation\sss\vspace{1.5pt}
\[
\quad
(\trf A\dff,\qff \varepsilon\trf)
\off \leq\off
(\trf A'\dff,\qff \varepsilon'\trf)
\quad
\mbox{if}\quad 
A\off =\off A'
\quad
\mbox{and}\quad 
\varepsilon\qff \leq\qff \varepsilon'
\pff.
\]

\vspace{-12pt}\vspace{1.5pt}
As in\sss the case of\sss $\hat{\mathcal{E}}$\dnsp,\oss
this order allows\sss to consider\sss $\mathcal{E}$ as a\sss
topological\sss category.\oss
It\dss is\dss easy\sss to see\sss that\sss $(\trf A\dff,\qff \varepsilon\trf)$\sss
is\dss an enhanced\dss Fredholm\dss operator\dss if\trs and\dss only\trs if\dss
$(\trf A^\sa\dff,\qff \varepsilon\trf)$\sss is\dss an enhanced odd\sss operator.\oss
Clearly,\oss the rule\sss
$(\trf A\fff,\qff \varepsilon\trf)
\off \longmapsto\off
(\trf A^\sa\dff,\qff \varepsilon\trf)$\sss
defines an\sss isomorphism of\dss categories\sss
$\sa\dff \colon\dff
\mathcal{E}
\qff \ttoo\qff 
\hat{\mathcal{E}}^{\dff \odd}$\dnsp.\oss
Let\sss us\sss assign\sss to an enhanced\dss Fredholm\dss operator\sss
$(\trf A\dff,\qff \varepsilon\trf)$\sss the object\vspace{3pt}
\[
\quad
\iota\qff(\trf A\dff,\qff \varepsilon\trf)
\off =\off
\left(\qff
\image\dff P_{\qff [\dff 0\fff,\qff \varepsilon\dff]}\trf(\trf \num{A}\trf)\fff,\off
\image\dff P_{\qff [\dff 0\fff,\qff \varepsilon\dff]}\trf(\trf \num{A^*}\trf)
\qff\right)
\] 

\vspace{-12pt}\vspace{3pt}
of\sss $\mathcal{S}$\nsp\dnsp,\oss
and assign\sss to a morphism\sss
$(\trf A\dff,\qff \varepsilon\trf)
\qff \ttoo\qff
(\trf A\dff,\qff \varepsilon'\trf)$\nnsp,\oss
where\sss $\varepsilon\qff \leq\qff \varepsilon'$\nnsp,\oss
the morphism\vspace{3pt}
\[
\quad
\iota\qff(\trf A\dff,\qff \varepsilon\trf)
\off \ttoo\off
\iota\qff(\trf A\dff,\qff \varepsilon'\trf)
\] 

\vspace{-12pt}\vspace{3pt}
of\sss $\mathcal{S}$\sss
defined\sss by\sss the pair\dss 
$\left(\qff
\image\dff P_{\qff [\dff \varepsilon\fff,\qff \varepsilon'\dff]}\trf(\trf \num{A}\trf)\fff,\off
\image\dff P_{\qff [\dff \varepsilon\fff,\qff \varepsilon'\dff]}\trf(\trf \num{A^*}\trf)
\qff\right)$\dss
together\sss with\sss the isometry\sss\vspace{3pt}
\[
\quad
\image\dff P_{\qff [\dff \varepsilon\fff,\qff \varepsilon'\dff]}\trf(\trf \num{A}\trf)
\off \ttoo\off
\image\dff P_{\qff [\dff \varepsilon\fff,\qff \varepsilon'\dff]}\trf(\trf \num{A^*}\trf)
\] 

\vspace{-12pt}\vspace{3pt}
induced\sss by\sss $U$\nnsp,\oss where\sss
$A\off =\off U\trf \num{A}$\sss is\dss the polar decomposition of\sss $A$\nnsp.\oss
A routine verification shows\sss that\sss these rules define a functor\sss
$\iota\trf \colon\dff
\mathcal{E}\qff \ttoo\qff \mathcal{S}$\nsp\dnsp.\oss
It\dss is\dss easy\sss to see\sss that\vspace{3pt}
\[
\quad
P_{\qff [\dff 0\fff,\qff \varepsilon\dff]}\trf(\trf \num{A}\trf)
\off \subset\off
\bigl\{\pff 
v\qff \in\qff H
\off \bigl|\off 
\norm{A\trf(\trf v\trf)}\qff \leq\qff \varepsilon\dff \norm{v}
\pff\bigr\}
\pff.
\] 

\vspace{-12pt}\vspace{3pt}
The set\sss at\sss the right\sss hand side of\dss this 
inclusion\dss is\dss almost\sss never a vector space.\oss
It\sss may be suggestive\sss to\sss think about\sss
$P_{\qff [\dff 0\fff,\qff \varepsilon\dff]}\trf(\trf \num{A}\trf)$\sss
as\sss the canonical\sss vector space replacement\sss of\dss this set.

\newpage
Let\sss
$\hat{\iota}^{\dff \odd}\dff \colon\dff 
\hat{\mathcal{E}}^{\dff \odd}
\qff \ttoo\qff 
\hat{\mathcal{S}}^{\dff \odd}$\dss
be\sss the composition\sss\vspace{3.625pt}
\[
\quad
\hat{\varphi}^{\dff \odd}\qff \circ\qff
\hat{o}^{\dff \odd}\qff \circ\qff
\hat{\psi}^{\dff \odd}
\qff \colon\qff
\hat{\mathcal{E}}^{\dff \odd}
\qff \ttoo\qff
\mathcal{E}\hat{\mathcal{O}}^{\dff \odd}
\qff \ttoo\qff 
\hat{\mathcal{O}}^{\dff \odd}
\qff \ttoo\qff 
\hat{\mathcal{S}}^{\dff \odd}
\pff.
\]

\vspace{-12pt}\vspace{3.625pt}
The functors $\iota$ and\sss $\hat{\iota}^{\dff \odd}$\sss appear\sss in\sss 
the following diagram.\vspace{0pt}
\[
\quad
\begin{tikzcd}[column sep=boom, row sep=boomm]
\mathcal{E}
\arrow[r, "\dis \sa"]
\arrow[d, "\dis \iota\qff"']
&
\hat{\mathcal{E}}^{\dff \odd}
\arrow[d, "\dis \pff \hat{\iota}^{\dff \odd}"]
\\
\mathcal{S}
\arrow[r, "\dis q"]
&
\hat{\mathcal{S}}^{\dff \odd}
\end{tikzcd}
\]

\vspace{-10.5pt}
The above description of\dss eigenvectors of\dss operators $A^\sa$\sss implies\sss
that\sss this diagram\dss is\dss commutative.\oss
We\sss leave\sss the details\sss to\sss the reader.\oss

\myuppar{Fredholm\dss vector\sss space models.}
A\qss \emph{Fredholm\dss vector\sss space\sss model}\oss
is\dss a\sss pair 
$(\trf V_{\fff 1}\dff,\pff V_{\fff 2}\dff)$ 
of\dss finite\-ly\sss dimensional\sss vector spaces with a scalar product.\oss
Let\sss $S$\sss be\sss the category\sss having\trs such\sss pairs\sss
as\sss objects,\oss
with\sss morphisms\sss
$(\trf V_{\dff 1}\dff,\pff V_{\dff 2}\trf)
\qff \ttoo\qff 
(\trf W_{\dff 1}\dff,\pff W_{\dff 2}\trf)$\sss
being\sss triples\sss
$(\trf i_{\trf 1}\dff,\pff i_{\trf 2}\dff,\pff g\trf)$\nnsp,\oss
where\vspace{3pt}
\[
\quad
i_{\trf 1}\dff \colon\dff V_{\fff 1}\qff \ttoo\qff W_{\fff 1}
\quad
\mbox{and}\quad
i_{\trf 2}\dff \colon\dff V_{\fff 2}\qff \ttoo\qff W_{\fff 2}
\]

\vspace{-12pt}\vspace{3pt}
are isometric embeddings and\vspace{3pt}
\[
\quad
g\dff \colon\dff
W_{\fff 1}\qff \ominus\qff i_{\trf 1}\dff(\trf V_{\fff 1}\trf)
\qff \ttoo\qff
W_{\fff 2}\qff \ominus\qff i_{\trf 2}\dff(\trf V_{\fff 2}\trf)
\]

\vspace{-12pt}\vspace{3pt}
is\dss an\sss isometry.\oss
The composition\dss is\dss defined\sss in\sss an obvious way
and amounts\sss to\sss taking\sss the direct\sss sum of\dss
the isometric\sss isomorphisms $g$\nnsp.\oss
This definition\dss is\dss slightly\sss different\sss from,\oss
but\dss is\dss trivially\sss equivalent\sss to\trs Segal's\qss \cite{s4}\qss definition of\dss
his category $\hat{C}$\nnsp.\oss

Similarly\sss to\dss Section\qss \ref{classifying-spaces-saf},\oss
there\dss is\dss a category\sss 
$S/\fff H$\sss serving as an\sss intermediary\sss
between $S$ and\sss $\mathcal{S}$\dnsp.\oss
The objects of\dss $S/\fff H$\sss are quadruples\sss
$(\trf V_{\fff 1}\dff,\pff V_{\fff 2}\trf;\pff h_{\dff 1}\dff,\pff h_{\trf 2}\trf)$\sss
such\sss that\sss
$(\trf V_{\fff 1}\dff,\pff V_{\fff 2}\trf)$\sss
is\dss an object\sss of\dss $S$\sss 
and\sss 
$h_{\dff 1}\dff \colon\dff V_{\fff 1}\qff \ttoo\qff H$\nnsp,\qss 
$h_{\dff 2}\dff \colon\dff V_{\fff 2}\qff \ttoo\qff H$\sss
are isometric\sss embeddings.\oss
A morphism\vspace{3pt}
\begin{equation}
\label{morphism-pairs}
\quad
\bigl(\trf V_{\fff 1}\dff,\pff V_{\fff 2}\trf;\pff h_{\dff 1}\pff,\pff h_{\trf 2}
\trf\bigr)
\qff \ttoo\qff
\left(\trf W_{\fff 1}\dff,\pff W_{\fff 2}\trf;\pff k_{\dff 1}\dff,\pff k_{\trf 2}
\trf\right)
\end{equation}

\vspace{-12pt}\vspace{3pt}
is\dss a\sss triple\sss
$(\trf i_{\trf 1}\dff,\pff i_{\trf 2}\dff,\pff g\trf)$\sss
defining a morphism\sss
$(\trf V_{\fff 1}\dff,\pff V_{\fff 2}\trf)
\qff \ttoo\qff
(\trf W_{\fff 1}\dff,\pff W_{\fff 2}\trf)$
of\trs $S$\sss such\sss that\vspace{3pt}
\[
\quad
k_{\dff 1}\qff \circ\qff i_{\trf 1}
\off =\off
h_{\dff 1}
\qff,\qquad
k_{\trf 2}\qff \circ\qff i_{\trf 2}
\off =\off
h_{\trf 2}
\pff.
\] 

\vspace{-12pt}\vspace{3pt}
Ignoring\sss embeddings $h_{\dff 1}\fff,\qff h_{\trf 2}$ 
defines a forgetting\sss functor\sss
$S/\fff H\qff \ttoo\qff S$\nnsp.\oss
Let\sss $S/\fff H\qff \ttoo\qff \mathcal{S}$\sss
be\sss the functor\sss
acting\sss on\sss the objects by assigning\sss to\sss a quadruple\sss 
$(\trf V_{\fff 1}\dff,\pff V_{\fff 2}\trf;\pff h_{\dff 1}\dff,\pff h_{\trf 2}\trf)$\sss
the pair\sss
$(\trf h_{\dff 1}\dff(\trf V_{\fff 1}\trf)\fff,\pff 
h_{\trf 2}\dff(\trf V_{\fff 2}\trf)\trf)$\nnsp,\oss
and acting\sss on\sss morphisms by\sss assigning\sss to 
a morphism\qss (\ref{morphism-pairs})\qss
defined\sss by\sss the\sss triple
$(\trf i_{\trf 1}\fff,\qff i_{\trf 2}\dff,\qff g\trf)$\sss
the morphism\vspace{3pt}
\[
\quad
\bigl(\trf h_{\dff 1}\dff(\trf V_{\fff 1}\trf)\fff,\qff 
h_{\trf 2}\dff(\trf V_{\fff 2}\trf)\trf\bigr)
\qff \ttoo\qff
\bigl(\trf k_{\dff 1}\dff(\trf W_{\fff 1}\trf)\fff,\qff 
k_{\trf 2}\dff(\trf W_{\fff 2}\trf)\trf\bigr)
\]

\vspace{-12pt}\vspace{3pt}
defined\sss by\sss the pair
$(\trf F_{\dff 1}\fff,\qff F_{\dff 2}\trf)$\sss
together with\sss the isometry\sss 
$f\dff \colon\dff
F_{\dff 1}\qff \ttoo\qff F_{\dff 2}$\nsp,\dff\oss
where\vspace{3pt}
\[
\quad
F_{\dff 1}
\off =\off
k_{\dff 1}\trf\bigl(\trf
W_{\fff 1}\qff \ominus\qff i_{\trf 1}\dff(\trf V_{\fff 1}\trf)
\trf\bigr)
\qff,\quad
\]

\vspace{-39pt}
\[
\quad
F_{\dff 2}
\off =\off
k_{\trf 2}\qff\bigl(\trf
W_{\fff 2}\qff \ominus\qff i_{\trf 2}\dff(\trf V_{\fff 2}\trf)
\trf\bigr)
\qff,\quad
\]

\vspace{-16.5pt}
and\vspace{-9.25pt}
\[
\quad
f
\off =\off 
k_{\trf 2}\dff \circ\dff g\dff \circ\dff (\dff k_{\dff 1}\dff)^{\dff -\dff 1}
\pff.
\]

\vspace{-12pt}\vspace{2.75pt}
Clearly,\oss these rules indeed define a functor\sss
$i\dff \colon\dff
S/\fff H\qff \ttoo\qff \mathcal{S}$\dnsp.\oss

\mypar{Theorem.}{intermediary-f}
\emph{The maps\sss
$\num{S}
\off \longleftarrow\off 
\num{S/\fff H}
\qff \ttoo\qff
\num{\mathcal{S}}$\sss
induced\sss by\sss these functors are homotopy equivalences.\oss}

\proof
The proof\dss is\dss similar\sss to\sss the proof\dss of\qss
Theorem\qss \ref{intermediary}.\oss 

Let\sss us consider\sss
$\num{S/\fff H}
\qff \ttoo\qff
\num{S}$\nnsp.\oss
The space of\dss objects of\sss $S$\sss is\dss discrete,\oss
and\sss the space of\dss objects of\sss $S/\fff H$\sss is\dss
the disjoint\sss union over\sss 
$(\trf V_{\fff 1}\dff,\pff V_{\fff 2}\dff)\qff \in\qff \ob\trf S$ of\dss the
spaces of\dss pairs of\dss isometric embeddings 
$V_{\fff 1}\dff,\pff V_{\fff 2}\qff \ttoo\qff H$\nnsp.\oss
Such spaces of\dss pairs are products of\dss two contractible spaces
and\sss hence are contractible.\oss 
It\dss follows\sss that\sss $S/\fff H\qff \ttoo\qff S$\sss 
induces a homotopy equivalence of\dss spaces of\dss objects.\oss
The space of\dss morphisms of\sss $S/\fff H$\sss 
is\dss a\sss locally\sss trivial\dss bundle over\sss the space of\dss
morphisms of\dss $S$\nnsp,\oss
with\sss the fiber over a morphism\sss
$(\trf V_{\dff 1}\dff,\pff V_{\dff 2}\trf)
\qff \ttoo\qff 
(\trf W_{\dff 1}\dff,\pff W_{\dff 2}\trf)$\sss 
of\dss $S$ being\sss the space 
of\dss pairs of\dss isometric em\-bed\-dings
$W_{\fff 1}\dff,\pff
W_{\fff 2}\qff \ttoo\qff H$\nnsp.\oss
Since such spaces of\dss pairs are contractible,\oss
the functor\sss $S/\fff H\qff \ttoo\qff S$\sss 
induces a homotopy equivalence of spaces of\dss morphisms.\oss
A similar argument\sss applies\sss to\sss the spaces
of\sss $n$\dnsp-simplices.\oss
Obviously,\oss our categories have free units,\oss
and\sss hence\sss
$\num{S/\fff H\halfff}
\qff \ttoo\qff
\num{S}$\sss
is\dss a\sss homotopy equivalence.\oss  

The proof\dss for\sss
$\num{S/\fff H}
\qff \ttoo\qff
\num{\mathcal{S}}$\sss
also closely\sss follows\sss the corresponding\sss part\sss of\dss the proof\dss 
of\qss Theorem\qss \ref{intermediary}.\oss 
The role of\dss the categories $\vect\downarrow W$\sss
is\dss played\sss by\sss the categories\sss
$\vect\downarrow (\trf E_{\dff 1}\fff,\qff E_{\dff 2}\trf)$\sss
defined as\sss the products\sss
$\vect\downarrow E_{\dff 1}
\dff \times\dff
\vect\downarrow E_{\dff 2}$.\oss
The key fact\dss is\dss that\sss classifying space\vspace{1.5pt}\vspace{0.625pt}
\[
\quad
\num{\vect\downarrow (\trf E_{\dff 1}\fff,\qff E_{\dff 2}\trf)}
\off\qff =\off\qff
\num{\vect\downarrow E_{\dff 1}}
\qff \times\qff
\num{\vect\downarrow E_{\dff 2}}
\]

\vspace{-12pt}\vspace{1.5pt}\vspace{0.625pt}
is\dss contractible.\oss 
This ensures\sss that\sss the same arguments as\sss in\sss the proof\dss of\trs
Theorem\qss \ref{intermediary}\qss work.\oss
We\sss leave\sss the details\sss to\sss the reader.\oss  \eproof

\mypar{Theorem.}{qs-fredholm}
\emph{The classifying\sss space\sss
$\num{S}$\sss
is\dss canonically\sss homotopy\sss equivalent\dss to\sss
the space\sss $\mathcal{F}$\sss of\pss bounded\trs Fredholm\dss
operators\sss
$H\qff \ttoo\qff H$\nnsp.\oss}

\proof
This immediately\sss follows\sss from\trs Theorems\qss \ref{s-fredholm}\qss
and\qss \ref{intermediary-f}.\oss  \eproof

\myuppar{Remark.}
Theorem\qss \ref{qs-fredholm}\qss 
is\dss due\sss to\dss Segal\qss \cite{s4}.\oss
The above proof\dss does not\sss follows\dss Segal's\dss outline.\oss
Segal\dss works directly\sss with\dss Fredholm operators and\sss uses
as an approximation\sss to\sss $\mathcal{F}$\sss the category
defined\sss by\sss the ordered set\sss of\dss pairs\sss $(\trf f\fff,\qff V\trf)$\nnsp,\oss
where\sss $f\dff \colon\dff H\qff \ttoo\qff H$\sss is\dss
a\dss Fredholm\sss operator and\sss $V\qff \subset\qff H$\sss is\dss
a finitely dimensional\sss subspace such\sss that\sss
$V\qff +\qff \image f\off =\off H$\nnsp.\oss
By\sss the definition,\qss
$(\trf f\fff,\qff V\trf)\off \leq\off (\trf f\fff'\fff,\qff V\fff'\trf)$\sss
if\trs and\dss only\trs if\dss
$f\off =\off f\fff'$ and\sss $V\qff \subset\qff V\fff'$\nnsp.\oss
This approach encounters some\sss technical\sss difficulties related\sss to\sss
the fact\sss that\sss the map\sss
$(\trf V\fff,\qff W\trf)\off \longrightarrow\off V\qff +\qff W$\nnsp,\oss
where\sss $V\fff,\qff W\qff \subset\qff H$\sss are finitely dimensional\sss subspaces,\oss
is\dss not\sss continuous,\oss and\sss by\sss this reason\trs
Proposition\qss 2.7\qss from\qss \cite{s2}\qss
does not\sss apply,\oss at\dss least\sss not\sss directly.\oss

\myuppar{A subcategory of\sss $\mathcal{S}$\nsp\dnsp.}
Let\sss $\mathcal{S}_{\qff \id}$\sss be\sss the\sss topological\sss category\sss
with\sss the same objects as $\mathcal{S}$ and\sss with\sss morphisms defined as follows.\oss
For an\sss object\sss $(\trf E_{\dff 1}\fff,\pff E_{\dff 2}\trf)$\sss of\sss 
$\mathcal{S}_{\qff \id}$\sss and\sss a finitely dimensional\sss subspace\sss 
$F\qff \subset\qff H$\sss
orthogonal\dss to both\sss $E_{\dff 1}$\sss and\sss $E_{\dff 2}$\nsp,\oss
there\dss is\dss a corresponding\sss morphism\vspace{3pt}
\[
\quad
(\trf E_{\dff 1}\fff,\pff E_{\dff 2}\trf)
\off \ttoo\off
(\trf E_{\dff 1}\dff \oplus\dff F\fff,\pff E_{\dff 2}\dff \oplus\dff F\trf)
\pff.
\]

\vspace{-12pt}\vspace{3pt}
There are no other\sss morphisms.\oss
The composition\dss is\dss defined\sss by\sss taking\sss the sums of\dss
the corresponding subspaces {\nsp}$F$\nnsp.\oss
The category\sss $\mathcal{S}_{\qff \id}$\dss is,\oss 
in\sss a natural\sss way,\oss
a subcategory of\sss $\mathcal{S}$\nsp\dnsp.\oss
Indeed,\oss these categories have\sss the same objects,\oss
and a morphism\sss
$(\trf E_{\dff 1}\fff,\qff E_{\dff 2}\trf)
\qff \ttoo\qff 
(\trf E\fff'_{\dff 1}\fff,\qff E\fff'_{\dff 2}\trf)$\sss
of\sss the category\sss 
$\mathcal{S}$\sss
defined\sss by\sss a\sss pair of\dss subspaces\sss
$(\trf F_{\dff 1}\fff,\qff F_{\dff 2}\trf)$\sss
and\sss an\sss isometry\sss
$f\dff \colon\dff
F_{\dff 1}\qff \ttoo\qff F_{\dff 2}$\sss
is\dss a morphism of\sss the category\sss 
$\mathcal{S}_{\qff \id}$\dss
if\trs and\dss only\trs if\dss
$F_{\dff 1}\off =\off F_{\dff 2}$\sss
and\sss $f$\sss is\dss the identity.\oss
Let\sss
$i\dff \colon\dff
\mathcal{S}_{\qff \id}\qff \ttoo\qff \mathcal{S}$\sss
be\sss the inclusion.\oss

The category\sss $\mathcal{S}_{\qff \id}$\dss
can\sss be also defined as\sss the category resulting\sss from\sss
a partial\sss order on\sss its space of\dss objects.\oss
Namely,\oss the space of\dss objects\dss is\dss ordered\sss by\sss
the relation\sss $\leq$\nnsp,\oss where\vspace{3pt}
\[
\quad
(\trf E_{\dff 1}\fff,\pff E_{\dff 2}\trf)
\off \leq\off
(\trf E\fff'_{\dff 1}\fff,\pff E\fff'_{\dff 2}\trf)
\pff
\]

\vspace{-34pt}
\[
\quad
\mbox{if}\dff\quad
E_{\dff 1}\off \subset\off E\fff'_{\dff 1}\qff,\quad
E_{\dff 2}\off \subset\off E\fff'_{\dff 2}\qff,\quad
\mbox{and}\quad\dff
E\fff'_{\dff 1}\dff \ominus\dff E_{\dff 1}
\off =\off
E\fff'_{\dff 2}\dff \ominus\dff E_{\dff 2}
\pff.
\]

\vspace{-12pt}\vspace{3pt}
Clearly,\oss the category defined\sss by\sss this order\dss is\dss
nothing else but\sss $\mathcal{S}_{\qff \id}$\nsp.\oss

The objects\sss $(\trf E_{\dff 1}\fff,\qff E_{\dff 2}\trf)$ of\dss
$\mathcal{S}_{\qff \id}$ can\sss be\sss thought\sss as formal\sss differences
$E_{\dff 1}\qff -\qff E_{\dff 2}$ of\dss finitely dimensional\sss
vector subspaces of\sss $H$\nnsp,\oss
with morphisms imposing,\oss in a moral\sss sense,\oss the relations\sss
$(\trf E_{\dff 1}\dff \oplus\dff F\trf)\qff -\qff (\trf E_{\dff 2}\dff \oplus\dff F\trf)
\off =\off
E_{\dff 1}\qff -\qff E_{\dff 2}$
when\sss $F$\sss is\dss orthogonal\sss to both\sss $E_{\dff 1}$ and\sss $E_{\dff 2}$\nsp.\oss
Of\dss course,\pss $\mathcal{S}$\sss admits a similar,\oss
but\dss less straightforward,\oss interpretation.\oss

\mypar{Theorem.}{categories-of-pairs}
\emph{The map\dss
$\num{i\fff}\dff \colon\dff
\num{\mathcal{S}_{\qff \id}}\qff \ttoo\qff \num{\mathcal{S}}$\dss
is\dss a\sss homotopy\sss equivalence.\oss}

\proof
A morphisms of\dss 
$\mathcal{S}$\sss is\dss determined\sss
by\sss a quadruple\sss
$(\trf E_{\dff 1}\fff,\pff E_{\trf 2}\fff,\pff F\fff,\pff f\trf)$\sss
such\sss that\sss 
$E_{\dff 1}\fff,\pff E_{\trf 2}$\nsp,\pss 
$F
\off \subset\off
H$\sss
are finitely dimensional\sss subspaces,\pss
$F$\dss is\dss orthogonal\sss to $E_{\dff 1}$\nsp,\oss and\sss 
$f\dff \colon\dff
F
\qff \ttoo\qff
H\dff \ominus\dff E_{\trf 2}$\sss
%
is\dss an\sss isometric embedding.\oss
A morphism corresponding\sss to\sss
$(\trf E_{\dff 1}\fff,\qff E_{\dff 2}\fff,\qff F\fff,\qff f\trf)$\sss
is\dss a morphism of\sss $\mathcal{S}_{\qff \id}$\sss 
if\trs and\dss only\trs if\qss 
$f$\sss is\dss the inclusion\sss map and,\oss
in\sss particular,\pss
$F\qff \subset\qff  H\dff \ominus\dff E_{\dff 2}$\nsp.\oss
Since\sss the inclusion\sss of\dss the space of\dss subspaces\sss
$F\qff \subset\qff H\dff \ominus\dff (\trf E_{\dff 1}\qff +\qff E_{\trf 2}\trf)$\sss
into\sss the space of\dss subspaces\sss
$F\qff \subset\qff H\dff \ominus\dff E_{\dff 1}$\sss
is\dss a homotopy equivalence,\oss and\dss
the space of\dss isometric embeddings of\dss 
$F
\qff \ttoo\qff
H\dff \ominus\dff E_{\dff 2}$\sss
is\dss contractible,\oss
the inclusion\dss
$\mor\dff \mathcal{S}_{\qff \id}
\qff \ttoo\qff
\mor\dff \mathcal{S}$\sss is\dss a\sss homotopy equivalence.\oss
More generally,\oss an $n$\dnsp-simplex of\dss 
$\mathcal{S}$\sss is\dss determined\sss
by an object\sss $(\trf E_{\dff 1}\fff,\qff E_{\dff 2}\trf)$\nnsp,\oss
a sequence\sss
$F_{\dff 1}\dff,\off
F_{\dff 2}\dff,\off
\ldots\dff,\off
F_{\dff n}$\sss
of\dss subspaces orthogonal\dss to\sss $E_{\dff 1}$ and\dss 
pair-wise orthogonal,\oss
and\dss isometric embeddings\vspace{3pt}
\[
\quad
f_{\dff i}\qff \colon\qff
F_{\dff i}
\off \ttoo\off
H\qff \ominus\qff
\bigl(\qff 
E_{\dff 2}\dff \oplus\dff
f_{\dff 1}\dff(\trf F_{\dff 1}\trf)\dff \oplus\dff
\ldots\dff \oplus\dff
f_{\dff i\dff -\dff 1}\dff(\trf F_{\dff i\dff -\dff 1}\trf)
\qff\bigr)
\pff,
\]

\vspace{-12pt}\vspace{3pt}
where\sss $i\off =\off 1\fff,\qff 2\fff,\qff \ldots\fff,\qff n$\nnsp.\oss
The corresponding $n$\dnsp-simplex of\sss $\mathcal{S}$\sss is\dss
an $n$\dnsp-simplex of\sss $\mathcal{S}_{\qff \id}$\sss
if\trs and\dss only\trs if\dss the embedding\sss 
$f_{\dff i}$\sss is\dss the inclusion\sss map for every $i$\nnsp.\oss
Since\sss the space of\dss isometric embeddings such as\sss
$f_{\dff i}$\sss is\dss contractible for every $i$\nnsp,\oss
the inclusion of\dss the space of\sss $n$\dnsp-simplices of\sss $\mathcal{S}_{\qff \id}$\sss
into\sss the space of\sss $n$\dnsp-simplices of\sss $\mathcal{S}$\sss
is\dss a homotopy equivalence.\oss
It\dss remains\sss to observe\sss that\sss 
$\mathcal{S}_{\qff \id}$\sss and\sss $\mathcal{S}$\sss
are categories with\sss free units.\oss
By\trs Lemma\qss \ref{free-d-categories}\qss
this implies\sss that\sss their\sss nerves have free degeneracies,\oss
and\sss hence\sss the\sss theorem\sss follows from\trs
Proposition\qss \ref{level-heq}.\oss  \eproof

\myuppar{Remark.}
Given\sss two objects
$(\trf E_{\dff 1}\fff,\pff E_{\dff 2}\trf)$
and\sss
$(\trf E\fff'_{\dff 1}\fff,\pff E\fff'_{\dff 2}\trf)$
of\dss $\mathcal{S}_{\qff \id}$,\oss
the spaces of\dss morphisms\sss
$(\trf E_{\dff 1}\fff,\pff E_{\dff 2}\trf)
\qff \ttoo\qff 
(\trf E\fff'_{\dff 1}\fff,\pff E\fff'_{\dff 2}\trf)$\sss
in $\mathcal{S}_{\qff \id}$\dss and\dss in $\mathcal{S}$\sss
do not\dss need\sss to be homotopy equivalent.\oss
Also,\oss there\dss is\dss no functor
$\mathcal{S}\qff \ttoo\qff \mathcal{S}_{\qff \id}$\sss
inducing\sss a homotopy inverse\sss to $\num{i\fff}$\nnsp.\oss

\myuppar{Small\sss versions of\dss the category $S$\nnsp.}
As in\sss the case of\dss the category $Q$\nnsp,\oss
in\sss the definition of\sss $S$\sss one can use only\sss the pairs of\dss
vector spaces with a scalar\sss product\sss from any set\sss
containing\sss a representative from each\sss isomorphism class.\oss
The resulting version of\sss $S$\sss is\dss 
equivalent\sss to\sss the original\sss one,\oss
and\sss hence\sss its\sss classifying space\dss is\dss
homotopy equivalent\sss to $\num{S}$\nnsp.\oss
The smallest\sss version\sss
$S^{\fff \mathrm{st}}$\sss
has as objects\sss the pairs\sss
$(\trf \ccc^{\dff n}\fff,\qff \ccc^{\dff m}\trf)$\nnsp.\oss
This\dss is\dss exactly\sss the category\sss $PC$\sss defined\dss by\dss Harris\qss \cite{h}.\oss
Cf.\qss \cite{h},\oss Section\qss 4.\oss
As\dss Harris\dss points out,\oss the category\sss
$S^{\fff \mathrm{st}}\off =\off PC$\sss is\dss a special\sss case of\dss
a construction used\sss by\trs Quillen\sss for\sss general\sss rings.\oss
See\qss \cite{q},\oss \cite{gr}.\oss

\myuppar{Harris's\dss category\sss $PG$\nnsp.}
The category\sss $\mathcal{S}_{\qff \id}$\sss is\dss a\dss Hilbert\dss space 
version of\sss a category $PG$
defined\sss by\dss Harris\qss \cite{h}.\oss
See\qss \cite{h},\oss Section\qss 4.\oss
Let\sss $\ccc^{\dff \infty}$\sss be a vector space of\dss countable infinite dimension
over $\ccc$ equipped\sss with a\dss Hermitian\dss scalar\sss product.\oss
We may assume\sss that\sss
$\ccc^{\dff \infty}$\sss is\dss a subspace of\dss $H$\nnsp.\oss
Then\sss the category\sss $PG$\sss considered\sss by\dss Harris\dss
is\dss equal\dss to\sss the full\sss subcategory of\dss
$\mathcal{S}_{\qff \id}$\sss having as objects\sss pairs\sss
$(\trf E_{\dff 1}\fff,\pff E_{\dff 2}\trf)$\sss of\dss subspaces\sss
$E_{\dff 1}\fff,\pff E_{\dff 2}$\sss of\dss $\ccc^{\dff \infty}$\dnsp.\oss

Harris\qss \cite{h}\qss indicated an approach\sss to determining\sss the homotopy\sss type
of\dss the classifying space of\dss $PG$ and\sss then wrote\fff:\oss
\emph{``Instead we will\dss just remark\sss that\dss the classifying\sss space of\qss $PG$ 
at\dss least\dss intuitively classifies pairs of\dss vector bundles modulo 
addition of\dss a common bundle\fff:\oss i.e.,\oss we claim\dss 
$BPG\off =\off \zzz\dff \times\dff B\halfff U$\nnsp.''}\oss
His space\sss $\zzz\dff \times\dff B\fff U$\sss is\dss the classifying\sss space for $K$\dnsp-theory.\oss
It\dss is\dss not\sss hard\sss to see\sss that\sss the inclusion 
$\num{PG}\qff \ttoo\qff \num{\mathcal{S}_{\qff \id}}$
is\dss a homotopy equivalence.\oss
Theorem\qss \ref{fu-quasi-fibration}\qss below\sss implies\sss that\sss
$\num{\mathcal{S}}$\sss is\dss also a classifying space for $K$\dnsp-theory.\oss
Together\sss with\dss Theorem\qss \ref{categories-of-pairs}\qss 
this provides a\sss justification of\trs Harris's\dss claim,\oss
certainly,\oss not\sss the first\sss one.\oss

\newpage
\mysection{Finite-unitary\qss bundle\qss and\qss finite-unitary\qss quasi-fibration}{unitary-quasi-fibrations}

\myuppar{Polarized\trs odd\dss subspace models.}
The goal\sss of\dss this section\dss is\dss to prove analogues of\dss
results of\trs Section\qss \ref{grassmannian-fibrations}\qss 
for categories related\sss to\sss general\trs Fredholm\dss operators\qss
(as opposed\sss to self-adjoint\sss ones).\oss
To\sss this end\sss we need a\sss category\sss
to play a role similar\sss to\sss that\sss of\sss 
$\mathcal{P}{\nsp}\hat{\mathcal{S}}$\dnsp,\oss
i.e.\qss a notion of\dss polarization for\sss objects of\sss $\mathcal{S}$\dnsp.\oss
From\sss the viewpoint\sss of\dss odd\sss operators\sss there\dss is\dss a
natural\sss candidate.\oss
Namely,\oss a\qss \emph{polarized\dss object}\pss of\sss
$\hat{\mathcal{S}}^{\dff \odd}$\sss should\sss be defined\sss as an object\sss
of\sss $\hat{\mathcal{S}}^{\dff \odd}$\dnsp,\oss
i.e.\qss a finitely dimensional\sss $\gamma$\dnsp-invariant\sss subspace\sss
$V\qff \subset\qff H\dff \oplus\dff H$\sss together\sss with a polarization\sss
$(\trf H\dff \oplus\dff H\trf)\dff \ominus\dff V
\off =\off
K_{\dff -}\dff \oplus\dff K_{\dff +}$\sss
which\dss is\qss \emph{$\gamma$\dnsp-invariant}\pss in\sss the sense\sss that\sss
$\gamma\trf(\trf K_{\dff +}\trf)\off =\off K_{\dff -}$\sss
and\sss hence\sss
$\gamma\trf(\trf K_{\dff -}\trf)\off =\off K_{\dff +}$\nsp.\oss
A\qss \emph{morphism}\pss of\dss polarized objects\vspace{1.25pt}
\[
\quad
(\trf V\fff,\off K_{\dff -}\dff,\off K_{\dff +}\trf)
\qff \ttoo\qff
(\trf V\fff',\off K\fff'_{\dff -}\dff,\off K\fff'_{\dff +}\trf)
\]

\vspace{-12pt}\vspace{1.25pt}
is\dss defined as a morphism\sss $V\qff \ttoo\qff V\fff'$\sss such\sss that\trs
if\dss $V\fff'\off =\off \gamma\trf(\trf U\trf)\dff \oplus\dff V\dff \oplus\dff U$\sss
is\dss the corresponding decomposition,\oss
then\sss $U\qff \subset\pff K_{\dff +}$\sss 
and\sss hence\sss
$\gamma\trf(\trf U\trf)\qff \subset\pff K_{\dff -}$.\oss
In\sss this way we get\sss a category\sss 
$\mathcal{P}{\nsp}\hat{\mathcal{S}}^{\dff \odd}$
of\dss polarized objects of\sss $\hat{\mathcal{S}}^{\dff \odd}$\dnsp.\oss
There\dss is\dss a\sss forgetting\sss functor\sss
$\pi\dff \colon\dff
\mathcal{P}{\nsp}\hat{\mathcal{S}}^{\dff \odd}
\qff \ttoo\qff 
\hat{\mathcal{S}}^{\dff \odd}$\nnsp.

Like\sss $\mathcal{P}{\nsp}\hat{\mathcal{S}}$\dnsp,\oss
the category\sss $\mathcal{P}{\nsp}\hat{\mathcal{S}}^{\dff \odd}$\sss 
can\sss be defined\sss in\sss terms of\dss an order.\oss
The space of\dss objects of\sss $\mathcal{P}{\nsp}\hat{\mathcal{S}}^{\dff \odd}$\sss
is\dss partially ordered\sss by\sss 
the relation\sss $\leq$\nnsp,\oss where\vspace{0pt}
\[
\quad
(\trf V\fff,\off K_{\dff -}\dff,\off K_{\dff +}\trf)
\off \leq\off
(\trf V\fff'\fff,\off K\fff'_{\dff -}\dff,\off K\fff'_{\dff +}\trf)
\quad
\]

\vspace{-38pt}
\[
\quad
\mbox{if}\dff\quad
V\off \subset\off V\fff'
\quad
\mbox{and}\quad
K\fff'_{\dff +}\qff \subset\pff K_{\dff +}
\pff.
\]

\vspace{-12pt}\vspace{0pt}
The $\gamma$\dnsp-invariance implies\sss that\sss then also\sss
$K\fff'_{\dff -}\qff \subset\pff K_{\dff -}$.\oss
By assigning\sss to\sss such an\sss inequality\dss the morphism defined\sss by\sss the subspace\sss
$U\off =\off K_{\dff +}\dff \ominus\dff K\fff'_{\dff +}$\sss
we get\sss an\sss isomorphism\sss between\sss
the category associated\sss with\sss this order and\sss
$\mathcal{P}{\nsp}\hat{\mathcal{S}}^{\dff \odd}$\dnsp.\oss
The morphisms of\dss $\mathcal{P}{\nsp}\hat{\mathcal{S}}^{\dff \odd}$\sss 
have\sss the form\vspace{1.25pt}
\[
\quad
(\trf V\fff,\off K_{\dff -}\dff,\off K_{\dff +} \trf)
\qff \ttoo\qff
(\trf\gamma\trf(\trf U\trf)\qff \oplus\qff
V\qff \oplus\qff U\dff,\off 
K_{\dff -}\dff \ominus\dff \gamma\trf(\trf U\trf)\dff,\off 
K_{\dff +}\dff \ominus\dff U
\trf)
\qff,
\]

\vspace{-12pt}\vspace{1.25pt}
where\sss $U$\sss
is\dss a\sss finitely dimensional\sss subspace of\dss  $K_{\dff +}$.\oss

\myuppar{Polarized\trs Fredholm\dss subspace models.}
Section\qss \ref{classifying-spaces-odd-saf}\qss suggests\sss the corresponding\sss
notion for\trs Fredholm\dss subspace models.\oss
A\qss \emph{polarized\trs Fredholm\dss subspace model}\oss is\dss a\sss triple
$(\trf E_{\dff 1}\dff,\qff E_{\trf 2}\dff,\qff i\trf)$ 
such\sss that\sss $E_{\dff 1}\dff,\qff E_{\trf 2}$ are\sss finitely\sss
dimensional\sss subspaces of\sss $H$\sss and\vspace{0pt}
\[
\quad
i\dff \colon\dff
H\dff \ominus\dff E_{\dff 1}
\qff \ttoo\qff
H\dff \ominus\dff E_{\trf 2}
\]

\vspace{-12pt}\vspace{0pt}
is\dss an\sss isometry.\oss
Let\sss $\mathcal{P}{\nsp}\mathcal{S}$\sss be\sss the category\sss having\trs 
such\sss triples\sss as\sss objects,\oss
with\sss morphisms\sss\vspace{0pt}
\[
\quad
\bigl(\trf E_{\dff 1}\dff,\off E_{\dff 2}\dff,\off i\trf\bigr)
\qff \ttoo\qff 
\bigl(\trf E\fff'_{\dff 1}\dff,\off E\fff'_{\dff 2}\dff,\off i\fff'\trf\bigr)
\]

\vspace{-12pt}\vspace{0pt}
being\dss triples\sss
$(\trf F_{\dff 1}\dff,\qff F_{\dff 2}\dff,\qff f\qff)$\sss
defining a morphism\sss
$(\trf E_{\dff 1}\dff,\qff E_{\trf 2}\trf)
\qff \ttoo\qff 
(\trf E\fff'_{\dff 1}\dff,\qff E\fff'_{\trf 2}\trf)$
of\sss $\mathcal{S}$\sss
and such\sss that\sss
$i\trf\bigl(\trf F_{\dff 1}\trf\bigr)
\off =\off
F_{\dff 2}$\dss
and\dss the maps\sss 
$f\dff \colon\dff F_{\dff 1}\qff \ttoo\qff F_{\dff 2}$\sss
and\vspace{1.5pt}
\[
\quad
i'\dff \colon\dff
H\trf \ominus\dff \bigl(\trf E_{\dff 1}\dff \oplus\dff F_{\dff 1}\trf\bigr)
\qff \ttoo\qff
H\trf \ominus\dff \bigl(\trf E_{\dff 2}\dff \oplus\dff F_{\dff 2}\trf\bigr)
\off =\off
i\trf \left(\qff H\trf \ominus\dff \bigl(\trf E_{\dff 1}\dff \oplus\dff F_{\dff 1}\trf\bigr)\qff \right)
\]

\vspace{-12pt}\vspace{1.5pt}
are induced\sss by\sss $i$\nnsp.\oss
Clearly,\oss a morphism\dss is\dss determined\sss by\sss 
the subspace\sss $F_{\dff 1}$\sss and\sss has\sss the form\vspace{1.5pt}
\[
\quad
\bigl(\trf E_{\dff 1}\dff,\off E_{\dff 2}\dff,\off i\trf\bigr)
\qff \ttoo\qff 
\bigl(\trf E_{\dff 1}\dff \oplus\dff F_{\dff 1}\dff,\off 
E_{\trf 2}\dff \oplus\dff i\trf(\trf F_{\dff 1}\trf)\dff,\off i\fff'\trf\bigr)
\pff,
\]

\vspace{-12pt}\vspace{1.5pt}
where\sss $i'$\sss is\dss induced\sss by\sss $i$\nnsp.\oss
The composition\dss is\dss defined\sss by\sss taking\sss the sums of\dss
the subspaces\sss $F_{\dff 1}$\nsp.\oss
Clearly,\oss the category\sss $\mathcal{P}{\nsp}\mathcal{S}$\sss
is\dss associated\sss with\sss the order\sss $\leq$\nnsp,\oss
where\vspace{1.5pt}\vspace{0.75pt}
\[
\quad
\bigl(\trf E_{\dff 1}\dff,\off E_{\dff 2}\dff,\off i\trf\bigr)
\off \leq\off 
\bigl(\trf E\fff'_{\fff 1}\dff,\off E\fff'_{\dff 2}\dff,\off i\fff'\trf\bigr)
\]

\vspace{-36pt}
\[
\quad
\mbox{if}\dff\quad
E_{\dff 1}\off \subset\off E\fff'_{\fff 1}
\off,
\qff\quad
E\fff'_{\dff 2}
\off =\off
E_{\trf 2}\dff \oplus\dff i\trf(\trf E\fff'_{\dff 1}\dff \ominus\dff E_{\dff 1}\trf) 
\pff,
\]

\vspace{-12pt}\vspace{1.5pt}\vspace{0.75pt}
and\sss $i'$\sss is\dss induced\sss by $i$\nnsp.\oss
The\sss topology on\sss the space of\dss objects\dss is\dss induced\sss
by\sss the obvious\sss topology on\sss the space of\dss pairs
$(\trf E_{\dff 1}\dff,\qff E_{\dff 2}\trf)$
and\sss the norm\sss topology on\sss the space of\dss isometries $i$\nnsp.\oss
Clearly,\oss the partial\sss order\sss $\leq$\sss on\sss the space
of\dss objects of\sss $\mathcal{P}{\nsp}\mathcal{S}$\sss 
has free equalities.\oss
Hence\sss $\mathcal{P}{\nsp}\mathcal{S}$ can\sss be 
considered as a\sss topological\sss simplicial\sss complex.\oss

One can\sss lift\sss the functor\sss $\iota\dff \colon\dff \mathcal{E}\qff \ttoo\qff \mathcal{S}$\sss
from\dss Section\qss \ref{classifying-spaces-odd-saf}\pss
to a functor\sss
$\mathcal{P}\dff\iota\dff \colon\dff 
\mathcal{E}
\qff \ttoo\qff 
\mathcal{P}{\nsp}\mathcal{S}$\dnsp.\oss
Let\sss $(\trf A\fff,\qff \varepsilon\trf)$\sss be an enhanced\dss Fredholm\sss operator,\oss
and\dss let\sss
$(\trf E_{\dff 1}\dff,\pff E_{\trf 2}\trf)
\off =\off
\iota\qff(\trf A\fff,\qff \varepsilon\trf)$\nnsp.\oss
Let\sss
$A\off =\off U\trf \num{A}$\sss be\sss the polar decomposition of\sss $A$\nnsp.\oss
Then\sss $U$\sss induces an\sss isometry\sss
$i\dff \colon\dff 
H\dff \ominus\dff E_{\dff 1}
\qff \ttoo\qff
H\dff \ominus\dff E_{\trf 2}$\nsp.\oss
Let\sss us\sss set\sss
$\mathcal{P}\dff\iota\trf (\trf A\fff,\qff \varepsilon\trf)
\off =\off
(\trf E_{\dff 1}\dff,\pff E_{\trf 2}\dff,\pff i\trf)$\nnsp.\oss
This defines\sss the action of\sss $\mathcal{P}\dff\iota$\sss on objects.\oss
We\sss leave\sss to\sss the reader\sss to define\sss the action on\sss morphisms
and\sss to check\sss that\sss $\mathcal{P}\dff\iota$\dss lifts $\iota$\nnsp.\oss

There\dss is\dss a\sss forgetting\sss functor\sss
$\pi\dff \colon\dff
\mathcal{P}{\nsp}\mathcal{S}
\qff \ttoo\qff 
\mathcal{S}$\nnsp.\oss
The isomorphism\sss between\sss
$\mathcal{S}$\sss and\sss $\hat{\mathcal{S}}^{\dff \odd}$\sss
from\dss Section\qss \ref{classifying-spaces-odd-saf}\qss
can\sss be\sss lifted\sss to\sss an\sss isomorphism\sss
between\sss
$\mathcal{P}{\nsp}\mathcal{S}$\sss 
and\sss $\mathcal{P}{\nsp}\hat{\mathcal{S}}^{\dff \odd}$\dnsp.\oss
This\sss isomorphism\sss takes\sss an object\sss
$(\trf E_{\dff 1}\dff,\qff E_{\trf 2}\dff,\qff i\trf)$\sss
of\sss $\mathcal{P}{\nsp}\mathcal{S}$\sss
to\sss the object\sss
$(\trf V\fff,\off K_{\dff -}\dff,\off K_{\dff +}\trf)$\sss
of\sss $\mathcal{P}{\nsp}\hat{\mathcal{S}}^{\dff \odd}$\dnsp,\oss
where\vspace{1.5pt}
\[
\quad
V
\off =\off
E_{\dff 1}\dff \oplus\dff E_{\trf 2}
\off \subset\off
H\dff \oplus\dff H
\]

\vspace{-12pt}\vspace{1.5pt}
$K_{\dff +}
\qff \subset\qff
H\dff \oplus\dff H$\sss
is\dss the graph of\dss the isometry\sss $i$\nnsp,\oss
and\sss
$K_{\dff -}\off =\off \gamma\trf(\trf K_{\dff +}\trf)$\nnsp.\oss

\mypar{Theorem.}{fredholm-polarization-forget}
\emph{The map\sss
$\num{\pi}\dff \colon\dff
\num{\mathcal{P}{\nsp}\mathcal{S}}
\qff \ttoo\qff 
\num{\mathcal{S}}$\sss
is\dss a\sss homotopy\sss equivalence.\oss}

\proof
The proof\trs is\dss similar\sss to\sss the proofs of\trs
Theorems\qss \ref{forgetting-polarization}\qss and\qss 
\ref{forgetting-equivalence}.\oss  \eproof

\myuppar{The unitary\sss bundle.}
For an object
$P\off =\off (\trf E_{\dff 1}\dff,\qff E_{\dff 2}\dff,\qff i\trf)$
of\trs $\mathcal{P}{\nsp}\mathcal{S}$
let\sss $U\ffin\dff(\trf P\trf)$\sss
be\sss the space of\dss isometries\sss
$H\qff \ttoo\qff H$\sss
equal\sss to $i$\sss on a subspace of\dss finite codimension\sss 
in\sss $H\dff \ominus\dff E_{\dff 1}$\nsp.\oss
We consider\sss $U\ffin\dff(\trf P\trf)$\sss
with\sss the norm\sss topology.\oss
Clearly,\pss $U\ffin\dff(\trf P\trf)$\sss is\dss non-empty\dss if\trs and\dss only\trs if\dss
$\dim\trf E_{\dff 1}\off =\off \dim\trf E_{\dff 2}$\nsp.\oss
If\trs $P\off =\off (\trf 0\dff,\qff 0\dff,\qff \id_{\trf H}\trf)$\nnsp,\oss
then\sss $U\ffin\dff(\trf P\trf)\off =\off U\ffin$\nsp.\oss
In\sss general,\oss if\dss $U\ffin\dff(\trf P\trf)$\sss is\dss nonempty,\oss then\sss
$U\ffin$ acts simply\sss transitively on\sss
$U\ffin\dff(\trf P\trf)$\sss from\sss the right.\oss

Let\sss $\mathcal{S}_{\qff 0}$\sss and\sss $\mathcal{P}{\nsp}\mathcal{S}_{\qff 0}$\sss
be\sss the full\sss categories of\dss
$\mathcal{S}$\sss and\sss $\mathcal{P}{\nsp}\mathcal{S}$\sss
having as objects,\oss respectively,\oss
pairs\sss $(\trf E_{\dff 1}\dff,\qff E_{\dff 2}\trf)$\sss
and\sss triples\sss
$(\trf E_{\dff 1}\dff,\qff E_{\dff 2}\dff,\qff i\trf)$\sss
such\sss that\sss
$\dim\trf E_{\dff 1}\off =\off \dim\trf E_{\dff 2}$\nsp.\oss
Let\sss us\sss consider\sss the space of\dss pairs\sss
$(\trf P\fff,\qff u\trf)$\sss
such\sss that\sss $P$\sss is\dss an object\sss of\dss
$\mathcal{P}{\nsp}\mathcal{S}_{\qff 0}$\sss
and\sss $u\qff \in\qff U\ffin\trf(\trf P\trf)$\nnsp.\oss
Clearly,\oss the projection\sss
$(\trf P\fff,\qff u\trf)
\off \longmapsto\off
P$\sss
from\sss this space of\dss pairs\sss to\sss the space of\dss
objects of\sss $\mathcal{P}{\nsp}\mathcal{S}_{\qff 0}$\sss
is\dss a\sss locally\sss trivial\sss bundle
with\sss the fiber\sss $U\ffin$\nsp.\oss
Similarly\sss to\dss Section\qss \ref{grassmannian-fibrations},\oss
this bundle can\sss be extended\sss to a bundle over\sss
$\num{\mathcal{P}{\nsp}\mathcal{S}_{\qff 0}}
\off =\off
\bbnum{\mathcal{P}{\nsp}\mathcal{S}_{\qff 0}}$\nnsp.\oss
Indeed,\oss if\dss there exists a morphism\sss
$P\qff \ttoo\qff P\fff'$\nnsp,\oss then\sss
$U\ffin\dff(\trf P\trf)\off =\off U\ffin\dff(\trf P\fff'\trf)$\nnsp,\oss
and\sss the construction of\trs Section\qss \ref{grassmannian-fibrations}\qss
is\dss based only\sss on\sss the similar\sss property of\dss $\gr\trf(\trf P\trf)$\sss 
and\sss the interpretation of\dss  
$\num{\mathcal{P}{\nsp}\hat{\mathcal{S}}}$\sss
as\sss 
$\bbnum{\mathcal{P}{\nsp}\hat{\mathcal{S}}}$\nnsp.\oss
Let\vspace{1.5pt}\vspace{0.575pt}
\[
\quad
\bm{\pi}\dff \colon\dff
\mathbf{U}
\qff \ttoo\qff 
\bbnum{\mathcal{P}{\nsp}\mathcal{S}_{\qff 0}}
\off =\off
\num{\mathcal{P}{\nsp}\mathcal{S}_{\qff 0}}
\]

\vspace{-12pt}\vspace{1.5pt}\vspace{0.575pt}
be\sss the resulting\sss bundle.\oss
Our\sss next\sss goal\dss is\dss to construct\sss a categorical\sss
analogue of\dss this bundle.\vspace{0.575pt}

\myuppar{Splittings.}
Let\sss $\mathcal{U}$\sss be\sss
be\sss the full\sss subcategory of\sss
$\mathcal{P}{\nsp}\mathcal{S}$\sss
having as objects\sss the objects of\sss
$\mathcal{P}{\nsp}\mathcal{S}$\sss
of\dss the form\sss
$(\trf 0\fff,\qff 0\fff,\qff i\trf)$\nnsp.\oss
This subcategory\sss has only\sss identity\sss morphisms,\oss
and\dss its space of\dss objects\dss is\dss nothing else but\sss the
space of\dss isometries\sss $H\qff \ttoo\qff H$\sss with\sss
the norm\sss topology and\sss hence\dss is\dss contractible by\trs Kuiper's\dss theorem.\oss
It\dss follows\sss that\sss 
$\num{\mathcal{U}}$\sss
is\dss contractible.\oss
A\qss \emph{split\dss object}\pss of\sss $\mathcal{P}{\nsp}\mathcal{S}$\sss 
is\dss defined as a morphism\sss
$N\qff \ttoo\qff M$\sss 
of\trs $\mathcal{P}{\nsp}\mathcal{S}$
such\sss that\sss $N$\sss is\dss an object\sss of\trs
$\mathcal{U}$\dnsp.\oss
A split\sss object\sss can\sss be identified\sss with an object\sss
$M
\off =\off
(\trf E_{\dff 1}\dff,\qff E_{\trf 2}\dff,\qff i\trf)$ 
of\sss
$\mathcal{P}{\nsp}\mathcal{S}$\sss
together\sss with an\sss isometry\sss
$e\dff \colon\dff
E_{\dff 1}\qff \ttoo\qff E_{\trf 2}$\nsp.\oss
Under\sss this identification\sss $N$\sss corresponds\sss to\sss
the isometry\sss equal\sss to $e$ on\sss $E_{\dff 1}$ and\sss to\sss $i$\sss on\sss
$H\dff \ominus\dff E_{\dff 1}$\nsp,\oss
and\sss
$N\qff \ttoo\qff M$\sss
is\dss the unique morphism\sss from\sss $N$\sss to\sss $M$\nnsp.\oss 
\emph{Morphisms of\dss split\sss objects}\pss are\sss the same commutative\sss
triangles as used\sss in\dss Section\qss \ref{polarizations-splittings}\qss
to define split\sss models.\oss 
Hence\sss the split\sss objects of\sss $\mathcal{P}{\nsp}\mathcal{S}$\sss
are\sss objects of\dss a\sss topological\sss category,\oss
which we denote by
$s\dff \mathcal{S}$\nnsp.\oss
There\dss is\dss a canonical\dss forgetting functor\sss
$\phi\dff \colon\dff
s\dff \mathcal{S}
\qff \ttoo\qff
\mathcal{P}{\nsp}\mathcal{S}$\sss
taking\sss
$N\qff \ttoo\qff M$\sss
to\sss $M$\nnsp.\oss\vspace{0.575pt}

\mypar{Theorem.}{split-contractible-u}
\emph{The classifying\sss space\sss 
$\num{s\dff \mathcal{S}}$\sss
is\dss contractible.\oss}\vspace{0.575pt}

\proof
The proof\dss is\dss completely\sss similar\sss
to\sss the proof\dss of\trs Theorem\qss \ref{split-contractible}.\oss  \eproof\vspace{0.575pt}

\myuppar{Categories related\sss to\sss finite-unitary\sss groups.}
Let\sss
$P\off =\off (\trf E_{\dff 1}\dff,\qff E_{\dff 2}\dff,\qff i\trf)$\sss
be an object\sss of\sss $\mathcal{P}{\nsp}\mathcal{S}_{\qff 0}$\nsp.\oss
Let\sss $\mathcal{U}\dff (\trf P\trf)$
be\sss the category\sss
having as objects diagrams in\sss 
$\mathcal{P}{\nsp}\mathcal{S}_{\qff 0}$\sss
of\dss the form\sss
$P\trf \ttoo\trf M\pff \longleftarrow\off N$\nnsp,\oss
where\sss $N$\sss is\dss an object\sss of\sss
$\mathcal{U}$\dnsp,\oss
and\sss with\sss morphisms\sss being\sss the same commutative diagrams
as in\sss the definition of\dss morphisms of\dss $\mathcal{G}\dff (\trf P\trf)$\sss
is\dss Section\qss \ref{categories-grassmannians}.\oss
The objects of\sss $\mathcal{U}\dff (\trf P\trf)$\sss
can\sss be identified\sss with\sss pairs\sss $(\trf F\fff,\qff f\trf)$\nnsp,\oss 
where\sss $F\qff \subset\qff H\dff \ominus\dff E_{\dff 1}$\sss 
is\dss a\sss finitely\sss dimensional\sss subspace\sss and\sss\vspace{1.5pt}\vspace{0.575pt} 
\[
\quad
f\dff \colon\dff
E_{\dff 1}\dff \oplus\dff F\qff \ttoo\qff E_{\trf 2}\dff \oplus\dff i\trf(\trf F\trf)
\pff
\]

\vspace{-12pt}\vspace{1.5pt}\vspace{0.575pt}
is\dss an\sss isometry.\oss
In\sss these\sss terms\sss morphisms have\sss the form\sss\vspace{1.5pt}\vspace{0.575pt}
\[
\quad
\bigl(\qff 
F\fff,\pff f
\qff\bigr)
\off \ttoo\off
\left(\qff 
F\dff \oplus\dff F\fff'\fff,\pff f\dff \oplus\trf 
\left(\trf i\qff \bigl|_{\qff F\fff'} \trf\right)
\qff\right)
\pff.
\]

\vspace{-12pt}\vspace{1.5pt}\vspace{0.575pt}
Like\sss the categories\sss $\mathcal{P}{\nsp}\mathcal{S}$\sss
and\sss $\mathcal{G}\dff (\trf P\trf)$\nnsp,\oss
the category\sss $\mathcal{U}\dff (\trf P\trf)$\sss
is\dss associated\sss with a partial\sss order.\oss
As in\dss Section\qss \ref{grassmannian-fibrations},\oss
every morphism\dss 
$u\dff \colon\dff
P\fff'\qff \ttoo\qff P$\dss 
of\qss $\mathcal{P}{\nsp}\mathcal{S}_{\qff 0}$\sss
defines a continuous functor\vspace{1.5pt}
\[
\quad
u^{\dff *}\dff \colon\qff
\mathcal{U}\dff (\trf P\trf)
\qff \ttoo\pff
\mathcal{U}\dff (\trf P\fff'\trf)
\pff,
\]

\vspace{-12pt}\vspace{1.5pt}
and\sss $u^{\dff *}$\sss is\dss strictly order-preserving\sss 
as a map of\dss partially ordered sets.\oss
More explicitly,\oss one can describe\sss $u^{\dff *}$\sss as follows.\oss
If\qss
$P\fff'
\off =\off 
(\trf E\fff'_{\dff 1}\dff,\qff E\fff'_{\dff 2}\dff,\qff i'\trf)$\nnsp,\oss
then\sss 
$E_{\dff 1}\off =\off E\fff'_{\dff 1}\dff \oplus\dff F\fff'_{\dff 1}$\sss
for some\sss $F\fff'_{\dff 1}$\nsp,\oss
and\vspace{1.5pt}
\[
\quad
u^{\dff *}\trf
\bigl(\qff 
F\fff,\pff f
\qff\bigr)
\off =\off
\left(\qff 
F\fff'_{\dff 1}\dff \oplus\dff F\fff,\pff f 
\qff\right)
\pff.
\]

\vspace{-12pt}\vspace{1.5pt}
As an example,\oss suppose\sss that\sss
$E_{\dff 1}\off =\off E_{\trf 2}\off =\off E$\sss
and\sss $i$\dss is\dss the identity map\sss
$H\dff \ominus\dff E
\qff \ttoo\qff
H\dff \ominus\dff E$\nnsp.\oss
Such an object\sss $P$\sss is\dss determined\sss by\sss the subspace\sss $E$\nnsp,\oss
and\sss we will\sss denote\sss the category\sss $\mathcal{U}\dff (\trf P\trf)$\sss
also by\sss $\mathcal{U}\dff (\trf E\trf)$\nsp.\oss
The objects of\sss $\mathcal{U}\dff (\trf E\trf)$\sss
are\sss the pairs\sss $(\trf F\fff,\qff f\trf)$\nnsp,\oss 
where\sss $F\qff \subset\qff H\dff \ominus\dff E$\sss
and\sss $f$\sss is\dss an\sss isometry\sss of\dss
$E\dff \oplus\dff F$\nnsp,\oss
and\sss morphisms\sss have\sss the form\sss
$(\qff F\fff,\qff f\qff)
\qff \ttoo\qff
(\qff F\dff \oplus\dff F\fff'\fff,\qff f\dff \oplus\dff \id_{\trf F\fff'}\qff)$\nnsp.\oss
In\sss fact,\oss this\dss is\dss a\sss typical\sss example\fff:\oss
clearly,\oss every category\sss $\mathcal{U}\dff (\trf P\trf)$\sss
is\dss isomorphic\sss to some\sss $\mathcal{U}\dff (\trf E\trf)$\nnsp.\oss
Let\sss us\sss fix a subspace $E$ and consider\sss the category\sss
$\mathcal{U}\dff (\trf E\trf)$\nnsp.\oss
We need also analogues of\dss categories\sss
$\mathcal{G}_{\qff \infty}\trf (\trf P\trf)$\nnsp,\pss
$\mathcal{G}_{\qff \nnn}\trf (\trf P\trf)$\nnsp,\oss
and\sss
$\mathcal{G}_{\qff \nnn,\qff n}\trf (\trf P\trf)$\sss
from\dss Section\qss \ref{categories-grassmannians}.\oss
It\dss is\dss more natural\sss
to define\sss these analogues\sss in\sss terms of\dss 
a fixed decomposition\vspace{1.5pt}
\[
\quad
H
\off =\off
\bigoplus\nolimits_{\dff n\qff \geq\qff 0}\qff H_{\dff n}
\] 

\vspace{-12pt}\vspace{1.5pt}
rather\sss than\sss a decomposition of\dss the form\qss
(\ref{spectral-decomposition})\qss into subspaces indexed\dss by all\dss integers.\oss
As in\dss Section\qss \ref{restricted-grassmannians},\oss the subspaces\sss $H_{\dff n}$\sss
are assumed\dss to be finitely dimensional.\oss
In order\sss to deal\sss with\sss $\mathcal{U}\dff (\trf E\trf)$\sss
we will\sss assume\sss that\sss $H_{\dff 0}\off =\off E$\nnsp.\oss
Let\sss 
$H_{\trf \leq\trf n}
\off =\off
H_{\dff 0}\dff \oplus\dff
H_{\dff 1}\dff \oplus\dff
\ldots\dff \oplus\dff
H_{\dff n}$\nsp.\oss
Let\sss
$\mathcal{U}_{\qff \infty}\trf (\trf E\trf)$\sss
be\sss the full\sss subcategory
of\dss $\mathcal{U}\dff (\trf E\trf)$\sss having as objects\sss
the pairs $(\trf F\fff,\qff f\trf)$ such\sss that\sss
$E\dff \oplus\dff F
\qff \subset\pff
H_{\trf \leq\trf n}$\sss
for some $n$\nnsp,\oss
and\dss let\sss
$\mathcal{U}_{\qff \nnn}\trf (\trf E\trf)$\sss
be\sss the full\sss subcategory
of\dss $\mathcal{U}\dff (\trf E\trf)$\sss having as objects\sss
the pairs\sss $(\trf F\fff,\qff f\trf)$\sss such\sss that\sss
$E\dff \oplus\dff F
\off =\off
H_{\trf \leq\trf n}$\sss
for some $n$\nnsp.\oss
Also,\oss for each\sss $n\qff \in\qff \nnn$\dss we will\sss need\sss
the full subcategories\sss
$\mathcal{U}_{\dff n}\trf (\trf E\trf)$\sss
and\sss
$\mathcal{U}_{\qff \nnn,\qff n}\trf (\trf E\trf)$\sss
of\dss $\mathcal{U}\dff (\trf E\trf)$\sss having as objects\sss
the pairs\sss $(\trf F\fff,\qff f\trf)$\sss such\sss that\sss
$E\dff \oplus\dff F
\qff \subset\pff
H_{\trf \leq\trf n}$\sss
and\dss
$E\dff \oplus\dff F
\off =\off
H_{\trf \leq\trf n}$\sss
respectively.\oss
Of\dss course,\oss the category\sss 
$\mathcal{U}_{\qff \nnn,\qff n}\trf (\trf E\trf)$\sss
has only\sss identity morphisms
and\sss its space of\dss objects\dss is\dss the usual\sss unitary\sss
group of\sss $H_{\trf \leq\trf n}$\nsp.\oss

With\sss these notions at\sss hand one can\sss prove analogues of\dss
all\sss results of\trs Section\qss \ref{categories-grassmannians}.\oss
In\sss the statements one needs simply\sss to\sss replace\sss
$\mathcal{G}_{\qff *}\trf (\trf P\trf)$\sss
by\sss $\mathcal{U}_{\qff *}\trf (\trf E\trf)$\nnsp,\oss
where $*$ stands for one of\dss the possible subscripts\qss
(or\sss the absence of\dss a subscript).\oss
In\sss the proof\dss one needs\sss to replace\sss terms such as\qss
``splitting of\dss $W$\nnsp''\pss by\dss terms\sss like\pss 
``isometry of\dss $E\dff \oplus\dff F$\nnsp''.\oss
Otherwise\sss the proofs are either completely similar\sss to\sss the proofs in\dss
Section\qss \ref{categories-grassmannians},\oss or even simpler,\oss
as\dss is\dss the case for\sss the analogue of\dss the key\trs
Lemma\qss \ref{using-compactness}.\oss
The main\sss results are summarized\sss in\sss the following\dss theorem.\oss

\mypar{Theorem.}{unitary-categories}
\emph{The canonical\sss maps\dss 
$\num{\mathcal{U}_{\qff \nnn}\trf (\trf E\trf)}
\qff \ttoo\qff
\num{\mathcal{U}_{\qff \infty}\trf (\trf E\trf)}
\qff \ttoo\qff
\num{\mathcal{U}\dff (\trf E\trf)}$\dss
are\sss homotopy\sss equivalences.\oss
There\sss are canonical\sss homotopy\sss
equivalences\sss 
$\num{\mathcal{U}_{\qff \nnn}\trf (\trf E\trf)}
\qff \ttoo\qff 
U\trf(\dff \infty\dff)$\sss
and\dss
$\num{\mathcal{U}_{\qff \infty}\trf (\trf E\trf)}
\qff \ttoo\qff
\num{\mathcal{U}\dff (\trf E\trf)}
\qff \ttoo\qff 
U\trf(\dff \infty\dff)$\hnsp\dnsp.\oss}  \eproof

\myuppar{The\sss finite-unitary\sss quasi-fibration.}
Similarly\sss to\dss Section\qss \ref{grassmannian-fibrations},\oss
let\sss us consider\sss the category\sss having as objects 
diagrams in\sss $\mathcal{P}{\nsp}\mathcal{S}_{\qff 0}$ of\dss the form\sss
$P\trf \ttoo\trf M\pff \longleftarrow\off N$\nnsp,\oss
where\sss $N$\sss is\dss an object\sss of\sss
$\mathcal{U}$\dnsp,\oss
and\sss with\sss morphisms\sss being\sss diagrams of\dss the form\qss
(\ref{sub-morphisms}).\oss
By\sss the same reasons as in\dss Section\qss \ref{grassmannian-fibrations}\qss
we will\sss denote\sss this category\sss by\sss $S\dff(\trf \phi\trf)$\nnsp,\oss
where now\sss $\phi$\sss is\dss the forgetting\sss functor\sss
$s\dff \mathcal{S}
\qff \ttoo\qff
\mathcal{P}{\nsp}\mathcal{S}_{\qff 0}$\nsp.\oss
There\dss is\dss an obvious\qss \emph{contravariant}\pss
forgetting\sss functor\sss
$p\dff \colon\dff
S\dff(\trf \phi\trf)
\qff \ttoo\qff
\mathcal{P}{\nsp}\mathcal{S}_{\qff 0}$\nsp,\oss
which can\sss be considered as a covariant\sss functor\sss
from\sss $S\dff(\trf \phi\trf)$\sss to\sss the category
opposite\sss to\sss
$\mathcal{P}{\nsp}\mathcal{S}_{\qff 0}$\nsp.\oss
Since\sss the geometric realizations do not\sss change
when a category\dss is\dss replaced\sss by\sss its opposite,\pss
$p$\sss induces a continuous map\sss
$\num{p}\dff \colon\dff
\num{S\dff(\trf \phi\trf)}
\qff \ttoo\qff
\num{\mathcal{P}{\nsp}\mathcal{S}_{\qff 0}}$\dnsp.\oss
This\dss is\dss a categorical\sss analogue of\dss the bundle\sss
$\bm{\pi}\dff \colon\dff
\mathbf{U}
\qff \ttoo\qff 
\bbnum{\mathcal{P}{\nsp}\mathcal{S}_{\qff 0}}
\off =\off
\num{\mathcal{P}{\nsp}\mathcal{S}_{\qff 0}}$\nnsp.\oss
The results of\trs Section\qss \ref{grassmannian-fibrations}\qss concerned\sss with\sss
$\num{p}\dff \colon\dff
\num{S\dff(\trf \phi\trf)}
\qff \ttoo\qff
\num{\mathcal{P}{\nsp}\hat{\mathcal{S}}}$\nnsp,\oss
with\sss the obvious modifications,\oss remain valid\sss in\sss the present\sss
context,\oss with essentially\sss the same proofs.\oss
We summarize\sss them\sss in\sss the following\sss theorem.\oss

\mypar{Theorem.}{fu-quasi-fibration}
\emph{The map\dss
$\num{p}\dff \colon\dff
\num{S\dff(\trf \phi\trf)}
\qff \ttoo\qff
\num{\mathcal{P}{\nsp}\mathcal{S}_{\qff 0}}$\dss
is\dss a\sss quasi-fibration.\oss
Its\sss homotopy\sss fiber\dss is\dss
homotopy equivalent\dss to\sss
$U\ffin$\sss and\dss
$U\trf(\dff \infty\dff)$\nnsp,\oss
and\sss the\sss total\sss space\sss
$\num{S\dff(\trf \phi\trf)}$\sss
is\dss contractible.\oss}  \eproof

\myuppar{Comparing\sss the bundle $\bm{\pi}$ and\sss 
the quasi-fibration\sss $\num{p}$\nnsp.}
Let\sss 
$P\off =\off (\trf E_{\dff 1}\dff,\qff E_{\dff 2}\dff,\qff i\trf)$\sss 
be an object\sss of\sss $\mathcal{P}{\nsp}\mathcal{S}_{\qff 0}$\nsp.\oss
There\dss is\dss a canonical\sss map\sss
$\num{\mathcal{U}\dff(\trf P\trf)}
\qff \ttoo\qff
U\ffin\trf(\trf P\trf)$\nnsp,\oss
which can\sss be defined as follows.\oss
Let\sss us\sss consider\sss $U\ffin\trf(\trf P\trf)$ as a\sss topological\sss category\sss
having\sss $U\ffin\trf(\trf P\trf)$ as\sss the space of\dss objects
and only\sss identity morphisms.\oss
An object\sss $(\trf F\fff,\qff f\trf)$\sss of\dss $\mathcal{U}\dff(\trf P\trf)$
defines a map\sss $H\qff \ttoo\qff H$\sss equal\dss to\sss $f$\sss on\sss
$E_{\dff 1}\dff \oplus\dff F$\sss and\sss to\sss $i$\sss on\sss 
$H\dff \ominus\dff (\trf E_{\dff 1}\dff \oplus\dff F\trf)$\nnsp.\oss
Clearly,\oss this map\dss is\dss an\sss isometry of\sss $H$\sss and\dss is\dss
equal\dss to\sss $i$\sss on a subspace of\dss finite codimension\sss 
in\sss $H\dff \ominus\dff E_{\dff 1}$\nsp,\oss
i.e.\qss is\dss an element\sss of\sss $U\ffin\trf(\trf P\trf)$\nnsp.\oss
Let\sss us assign\sss to\sss $(\trf F\fff,\qff f\trf)$\sss this isometry.\oss
A\sss trivial\sss verification shows\sss that\sss
this rule assigns\sss the same\sss isometry\sss to\sss two objects
of\dss $\mathcal{U}\dff(\trf P\trf)$ related\sss by a morphism.\oss
Hence we can extend\sss this rule\sss to morphisms 
of\dss $\mathcal{U}\dff(\trf P\trf)$\sss
by assigning\sss to every morphism an\sss identity morphism 
of\dss $U\ffin\trf(\trf P\trf)$ and\sss
get\sss a functor\sss
$\mathcal{U}\dff(\trf P\trf)
\qff \ttoo\qff
U\ffin\trf(\trf P\trf)$\nnsp.\oss
By\sss passing\sss to\sss geometric realizations 
we get\sss a map\sss\vspace{1.5pt}\vspace{-0.8pt}
\[
\quad
h\trf(\trf P\trf)\dff \colon\dff
\num{\mathcal{U}\dff(\trf P\trf)}
\qff \ttoo\qff
\num{U\ffin\trf(\trf P\trf)}
\off =\off
U\ffin\trf(\trf P\trf)
\pff.
\]

\vspace{-12pt}\vspace{1.5pt}\vspace{-0.8pt}
As in\dss Section\qss \ref{grassmannian-fibrations},\oss
the maps\sss $h\trf(\trf P\trf)$\sss define a map\sss
$h\dff \colon\dff
\num{S\dff(\trf \phi\trf)}
\qff \ttoo\qff
\mathbf{U}$\sss
such\sss that\sss
$\bm{\pi}\dff \circ\dff h\off =\off \num{p}$\nnsp,\oss
i.e.\qss a map from\sss the quasi-fibration $\num{p}$\sss to\sss 
the\sss bundle $\bm{\pi}$\nnsp.\oss
Theorem\qss \ref{unitary-categories}\qss implies\sss that\sss $h$\sss
induces homotopy equivalences on\sss the fibers.\oss
By comparing\sss the homotopy sequences of\dss 
$\num{p}$\sss and\sss $\bm{\pi}$\sss we see\sss that\sss
$\mathbf{U}$\sss is\dss weakly\sss homotopy equivalent\sss 
to\sss $\num{S\dff(\trf \phi\trf)}$\nnsp.\oss
Now\sss the arguments of\dss the proof\dss of\trs
Theorem\qss \ref{g-is-contractible}\qss apply\sss
and\dss lead\sss to\sss the following\sss theorem.\oss

\mypar{Theorem.}{u-is-contractible}
\emph{The space\dss $\mathbf{U}$ is\dss contractible.\oss}  \eproof

\mypar{Theorem\qss ({\fff}Atiyah--Singer{\fff}).}{f-loops}
\emph{The\sss loop space of\qss $\mathcal{F}$\dss
is\dss homotopy equivalent\dss to 
$\hat{\mathcal{F}}$\nnsp.}

\proof
The\sss loop spaces of\dss
$\num{\mathcal{P}{\nsp}\mathcal{S}}$\sss
and\sss
$\num{\mathcal{P}{\nsp}\mathcal{S}_{\qff 0}}$\sss
are\sss the same,\oss
and\trs Theorem\qss \ref{u-is-contractible}\qss
implies\sss that\sss the\sss loop space of\dss
$\num{\mathcal{P}{\nsp}\mathcal{S}_{\qff 0}}$\sss
is\dss homotopy equivalent\sss to\sss the fiber\sss
$U\ffin$ of\sss $\bm{\pi}$\nnsp.\oss
By\trs Theorem\qss \ref{s-fredholm}\qss the space\sss $\mathcal{F}$\sss
is\dss homotopy equivalent\sss to $\num{\mathcal{S}}$\nnsp,\oss
and\sss by\trs Theorem\qss\ref{fredholm-polarization-forget}\qss 
the space\sss $\num{\mathcal{S}}$\dss
is\dss homotopy equivalent\sss to $\num{\mathcal{P}{\nsp}\mathcal{S}}$\nnsp.\oss
Hence\sss the\sss loop space of\dss $\mathcal{F}$\sss is\dss
homotopy equivalent\sss to\sss $U\ffin$\dnsp.\oss
On\sss the other\sss hand,\oss 
by\trs Theorem\qss \ref{operators-categories}\qss 
the space\sss $\hat{\mathcal{F}}$\sss is\dss homotopy equivalent\sss
to\sss $\num{\hat{\mathcal{S}}}$\nnsp,\oss
and\sss by\trs Theorem\qss \ref{harris-h}\qss
the space\sss $\num{\hat{\mathcal{S}}}$\sss is\dss homeomorphic\sss
to\sss $U\ffin$\dnsp.\oss
The\sss theorem\sss follows.\oss  \eproof

\mysection{Atiyah--Singer\qss map}{as-map}

\myuppar{Atiyah--Singer\dss map\sss
$\hat{\mathcal{F}}
\fff \ttoo\dff
\Omega\qff \mathcal{F}$\nsp\dnsp.}
By\trs Theorem\qss \ref{operators-categories}\qss
and\dss Corollary\qss \ref{f-loops}\qss
the\sss loop space\sss $\Omega\qff \num{\mathcal{F}}$\sss
of\trs $\num{\mathcal{F}}$\sss is\dss 
homotopy equivalent\sss to\sss $\num{\hat{\mathcal{F}}}$\nnsp.\oss
In\sss fact,\oss the proofs provide a canonical\dss homotopy equivalence.\oss
This result\dss is\dss due\sss to\dss
Atiyah\dss and\trs Singer\qss \cite{as},\oss
who also constructed an explicit\dss homotopy equivalence\sss
$\alpha\dff \colon\dff
\hat{\mathcal{F}}
\qff \ttoo\qff
\Omega\qff \mathcal{F}$\dnsp.\oss
We are going\sss to prove\sss that\sss these\sss two homotopy equivalences
are\sss the same up\sss to homotopy.\oss
To\sss this end we need\sss to discuss\sss the\dss Atiyah--Singer\dss proof\dss
in some details and\sss to make\sss the proof\dss of\dss the present\sss paper
more explicit.\oss 

The\dss Atiyah--Singer\qss \cite{as}\qss map\sss
$\alpha\dff \colon\dff
\hat{\mathcal{F}}
\qff \ttoo\qff
\Omega\qff \mathcal{F}$\sss
is\dss defined as follows.\oss
First\sss of\dss all,\oss Atiyah\dss and\dss Singer\dss work\sss with bounded operators.\oss
By\sss this reason\sss the space\sss $\hat{\mathcal{F}}$\sss will\sss be understood
as\sss the space of\dss bounded self-adjoint\trs Fredholm\sss operators with\sss
the norm\sss topology,\oss subject\sss to\sss the usual\sss condition of\dss not\sss
being essentially\sss positive or\sss negative.\oss
In\qss \cite{as}\qss this space\dss is\dss denoted\sss by\sss $\hat{\mathcal{F}}_{\dff *}$\nsp.\oss
Second,\oss
Atiyah\dss and\dss Singer\dss understand\sss
the\sss loop space\sss $\Omega\qff \mathcal{F}$\sss
as\sss the space of\dss paths connecting\sss 
$\id$\sss with\sss $-\qff \id$\nnsp.\oss
To simplify\sss the notations,\oss we will\sss assume\sss that\sss
loops\qss (or,\oss rather,\oss paths)\qss 
are defined on\sss $[\trf 0\fff,\qff \pi \trf]$\nnsp.\oss
Given an operator\sss $A\qff \in\qff \hat{\mathcal{F}}$\dnsp,\oss
the\dss loop\sss $\alpha\trf(\trf A\trf)$\sss is\dss the map\vspace{1.5pt}
\[
\quad
t
\off \longmapsto\off
\id_{\trf H}\dff \cos t\qff +\qff i\trf A\trf \sin t\qff,
\quad
0\qff \leq\qff t\qff \leq\qff \pi
\pff.
\]

\vspace{-12pt}\vspace{1.5pt}
The proof\dss that\sss $\alpha$\sss is\dss a homotopy equivalence\dss
is\dss not\sss direct,\oss and\sss in order\sss to compare $\alpha$\sss
with\sss the homotopy equivalence implicit\sss in\dss 
Section\qss \ref{unitary-quasi-fibrations},\oss
we need some details of\dss this proof.\oss
Let\sss us begin with a general\sss construction\sss behind\sss both proofs.\oss

\myuppar{Loops and\sss bundles.}
Let\sss $p\dff \colon\dff E\qff \ttoo\qff B$\sss 
be a\sss locally\sss trivial\sss bundle.\oss 
Suppose\sss that\qss \emph{two base points}\dss $s\fff,\qff t\qff \in\qff B$\sss are fixed,\oss
and\sss let\sss $\widetilde{s}\qff \in\qff E$\sss be such\sss that\sss
$p\trf(\trf \widetilde{s}\trf)\off =\off s$\nnsp.\oss
Let\sss $F\off =\off p^{\dff -\dff 1}\dff(\trf t\trf)$\nnsp.\oss
We will\sss denote by\sss $\Omega\trf(\trf E\fff,\qff F\trf)$\sss
the space of\dss all\sss paths in $E$ starting at $\widetilde{s}$\sss
and ending\sss in $F$\dnsp.\oss

Let\sss $\Pi\dff B$\sss be\sss the space of\dss paths\sss
$[\trf 0\fff,\qff \pi\trf]\qff \ttoo\qff B$\sss starting at\sss $s$\nnsp,\oss
and\sss let\sss $\ev\dff \colon\dff \Pi\dff B\qff \ttoo\qff B$\sss
be\sss the evaluation\sss map assigning\sss to a path\sss its endpoint.\oss
Let\sss us\sss interpret\sss the loops space $\Omega\dff B$\sss
as\sss the space of\dss paths 
starting at\sss $s$\sss and\sss ending at\sss $t$\nnsp.\oss
Then\sss $\Omega\dff B$\sss is\dss the fiber of\dss
$\ev\dff \colon\dff \Pi\dff B\qff \ttoo\qff B$\sss over\sss $t$\nnsp.\oss
The maps\sss $p$\sss and\sss $\ev$\sss
are related\sss by\sss 
the commutative diagram of\dss the form\vspace{-3pt}
\[
\quad
\begin{tikzcd}[column sep=sboom, row sep=sboom]
\Pi\dff B
\arrow[d, "\dis \ev\dff"']
\arrow[r, "\dis e"]
&
E
\arrow[d, "\dis \dff p"]
\\
B
\arrow[r, "\dis \id"]
&
B\dff.
\end{tikzcd}
\]

\vspace{-12pt}
The map\sss
$e\dff \colon\dff
\Pi\dff B\qff \ttoo\qff E$\sss is\dss constructed as follows.\oss
Let\sss  
$\mathbf{s}\dff \colon\dff \Pi\dff B\qff \ttoo\qff B$\sss 
be\sss the map\sss taking\sss all\sss paths\sss
to\sss the point $s$ and\trs 
$\widetilde{\mathbf{s}}\dff \colon\dff \Pi\dff B\qff \ttoo\qff E$\sss
be\sss the map\sss
taking\sss all\sss paths\sss to\sss the point\sss $\widetilde{s}\qff \in\qff E$\nnsp.\oss
Then\sss $p\dff \circ\qff \widetilde{\mathbf{s}}\off =\off \mathbf{s}$\nnsp.\oss
For\sss $u\qff \in\qff [\trf 0\fff,\qff \pi \trf]$\dss  let\dss 
$\ev_{\fff u}\dff \colon\dff \Pi\dff B\qff \ttoo\qff B$\sss
be\sss the evaluation\sss map\sss
$k\off \longmapsto\off k\trf(\dff u\trf)$\nnsp.\oss
The maps\sss $\ev_{\fff u}$\sss form\sss a homotopy\sss
between\sss 
$\ev_{\dff 0}
\off =\off
\mathbf{s}$\sss 
and\sss
$\ev_{\dff \pi}
\off =\off
\ev\dff \colon\dff
\Pi\dff B\qff \ttoo\qff B$\nnsp.\oss
By\sss the covering\sss homotopy property of\sss
$p\dff \colon\dff E\qff \ttoo\qff B$\sss
there exists a homotopy
$h_{\fff u}
\dff \colon\dff 
\Pi\dff B\trf \ttoo\qff E\fff,\qff u\qff \in\qff [\trf 0\fff,\qff \pi\trf]$\sss
covering\sss the homotopy\sss
$\ev_{\fff u}\dff \colon\dff \Pi\dff B\qff \ttoo\qff B$\nnsp,\oss
starting at\sss $\widetilde{\mathbf{s}}\off =\off h_{\dff 0} $\sss
and ending at\sss some map\sss 
$e\off =\off h_{\dff \pi}\dff \colon\dff \Pi\dff B\qff \ttoo\qff E$\nnsp.\oss
Then\sss 
$p\dff \circ\dff e
\off =\off 
p\dff \circ\dff h_{\dff \pi}
\off =\off 
\ev_{\dff \pi}
\off =\off
\ev$\sss
and\sss hence\sss the above diagram\dss is\dss commutative.\oss
The map $e$\sss induces a map of\dss fibers\sss
$f\dff \colon\dff
\Omega\dff B\qff \ttoo\qff F$\nnsp.\oss
For every\sss 
$\lambda\qff \in\qff \Omega\trf B$\sss
the path\sss 
$g\trf(\trf \lambda\trf)\dff \colon\dff
u\off \longmapsto\off h_{\fff u}\trf(\trf \lambda\trf)$\sss
covers\sss $\lambda$\sss and ends in\sss $F$\dnsp.\oss
This defines a map\vspace{0.3pt}
\[
\quad
g\dff \colon\dff
\Omega\trf B
\qff \ttoo\qff
\Omega\trf(\trf E\fff,\qff F\trf)
\qff.
\]

\vspace{-12pt}\vspace{0.3pt}
Clearly,\pss
$\ev\dff \circ\dff g\off =\off f$\nnsp.\oss
If\dss $E$\sss is\dss contractible,\oss
then\sss $f$\sss and\sss $g$\sss are weak\sss homotopy equivalences
by\sss the elementary\sss homotopy\sss theory.\oss
If,\oss in addition,\pss $B$\sss and\sss $F$\sss are homotopy equivalent\sss
to\sss CW-complexes,\oss then\sss $f$\sss and\sss $g$\sss are homotopy equivalence.\oss
When\sss we will\sss need\sss to stress\sss the dependence on\sss
the bundle\sss $p$\nnsp,\oss we will\sss use\sss the notations\sss
$e^{\dff p}\nsp,\off f^{\dff p}$\dnsp,\oss etc.\oss

We will\sss need also a generalization of\dss this construction.\oss
Suppose\sss that\sss 
$\alpha\dff \colon\dff X\qff \ttoo\qff \Pi\dff B$\sss
is\dss a continuous map.\oss
We may consider $\alpha$ as a homotopy\sss 
$\ev_{\dff u}\dff \circ\qff \alpha\fff,\qff u\qff \in\qff [\trf 0\fff,\qff \pi\trf]$\sss 
between\sss the constant\sss map\sss 
$\ev_{\dff 0}\dff \circ\qff \alpha
\off =\off
\mathbf{s}\dff \circ\trf \alpha$\sss
and\sss
$\ev_{\dff \pi}\dff \circ\qff \alpha
\off =\off
\ev\dff \circ\qff \alpha$\nnsp.\oss
By\sss the covering\sss homotopy property\sss
there exists a homotopy
$h_{\fff u}\trf(\trf \alpha\trf)
\dff \colon\dff 
X\trf \ttoo\qff E\fff,\qff u\qff \in\qff [\trf 0\fff,\qff \pi\trf]$\sss
covering\sss the homotopy\sss
$\ev_{\fff u}\dff \circ\qff \alpha
\dff \colon\dff 
X\qff \ttoo\qff B$\nnsp,\oss
starting at\dss $\widetilde{\mathbf{s}}\dff \circ\qff \alpha$\sss
and ending at\sss some map\sss 
$e\trf(\trf \alpha\trf)
\off =\off 
h_{\dff \pi}\trf(\trf \alpha\trf)\dff \colon\dff 
X\qff \ttoo\qff E$\nnsp.\oss
If\sss $\alpha$\sss is\dss the identity\sss map\sss
$\Pi\dff B\qff \ttoo\qff \Pi\dff B$\nnsp,\oss
then\sss $e\trf(\trf \alpha\trf)\off =\off e$\nnsp.\oss
In our\sss applications\sss $\alpha$\sss will\sss be a map\sss from\sss $X$\sss to\sss
$\Omega\dff B\qff \subset\qff \Pi\dff B$\nnsp.\oss
In\sss this case\sss the image of\sss $e\trf(\trf \alpha\trf)$\sss
is\dss contained\sss in\sss $F$\sss and\sss $e\trf(\trf \alpha\trf)$\sss
induces a map\sss
$f\trf(\trf \alpha\trf)\dff \colon\dff
X\qff \ttoo\qff F$\nnsp.\oss
Also,\oss in\sss the case\sss the homotopy\sss
$h_{\fff u}\trf(\trf \alpha\trf)
\dff \colon\dff 
X\trf \ttoo\qff E\fff,\qff u\qff \in\qff [\trf 0\fff,\qff \pi\trf]$\sss
defines a map\vspace{0.3pt}
\[
\quad
g\trf(\trf \alpha\trf)\dff \colon\dff
X
\qff \ttoo\qff
\Omega\trf(\trf E\fff,\qff F\trf)
\qff.
\]

\vspace{-12pt}\vspace{0.3pt}
Of\dss course,\oss one can always\sss take\sss
$h_{\fff u}\trf(\trf \alpha\trf)
\off =\off
h_{\fff u}\trf \circ\qff \alpha$\nnsp.\oss 
The point\sss of\dss this generalization\dss is\dss the following\sss uniqueness property,\oss
which\dss is\dss stronger\sss in\sss this more general\sss situation.\oss

\mypar{Lemma.}{loops-bundles-lifts}
\emph{Suppose\sss that\sss $\alpha$ is\dss a map\qss
$X\qff \ttoo\qff \Omega\trf B$\nnsp.\oss
Then up\dss to\sss homotopy\sss the maps\sss
$f\trf(\trf \alpha\trf)\dff \colon\dff
X\qff \ttoo\qff F$\sss
and\trs
$g\trf(\trf \alpha\trf)\dff \colon\dff
X
\qff \ttoo\qff
\Omega\trf(\trf E\fff,\qff F\trf)$\sss
do not\sss depend on\sss the choice of\dss a covering\sss homotopy\dss
$h_{\fff u}\trf(\trf \alpha\trf)
\dff \colon\dff 
X\trf \ttoo\qff E\fff,\qff u\qff \in\qff [\trf 0\fff,\qff \pi\trf]$\nnsp.\oss}

\proof
Two such choices\sss $h_{\dff 0}\dff,\qff h_{\dff 1}$\sss 
can\sss be considered as a\sss homotopy covering\sss the\sss
composition of\dss the projection\sss
$X\dff \times\dff \{\qff 0\fff,\qff 1\qff\}
\qff \ttoo\qff
X$\sss
with\sss the homotopy\sss
$\ev_{\dff u}\dff \circ\qff \alpha\fff,\qff u\qff \in\qff [\trf 0\fff,\qff \pi\trf]$\nnsp.\oss
This\sss covering\sss homotopy\sss can\sss be extended\sss to a\sss 
homotopy\sss covering\sss the composition of\dss the projection\sss
$X\dff \times\dff [\trf 0\fff,\qff 1\trf]
\qff \ttoo\qss 
X$\sss
with\sss $\ev_{\dff u}\dff \circ\qff \alpha$\nnsp.\oss
By considering\sss the factor\sss $[\trf 0\fff,\qff 1\trf]$\sss in\sss
$X\dff \times\dff [\trf 0\fff,\qff 1\trf]$\sss as\sss the parameter space
of\dss a homotopy\sss
we can consider\sss the extended\dss covering\sss homotopy 
as a homotopy\sss between\sss $h_{\dff 0}$\sss and\sss $h_{\dff 1}$\sss 
in\sss the class of\dss homotopies covering\sss $\ev_{\dff u}\dff \circ\qff \alpha$\nnsp.\oss
This homotopy\sss induces homotopies between\sss the maps\sss
$X\qff \ttoo\qff F$\sss
resulting from\sss $h_{\dff 0}$\sss and\sss $h_{\dff 1}$\nsp,\oss
as also between\sss the maps\sss
and\sss
$X\qff \ttoo\qff \Omega\trf(\trf E\fff,\qff F\trf)$\nnsp.\oss  \eproof

\myuppar{Functorial\dss properties.}
The above constructions have good functorial\sss properties with respect\sss to bundle maps.\oss
Let\sss $p'\dff \colon\dff E\fff'\qff \ttoo\qff B\fff'$\dss be another\sss locally\sss
trivial\sss bundle,\oss
and\dss let\sss $s'\fff,\qff t'$\sss be\sss the base points of\dss $B\fff'$\dnsp.\oss
Let\sss $\widetilde{s}\qff'\qff \in\qff E\fff'$\sss be a point\sss such\sss that\sss
$p'\trf(\trf \widetilde{s}\qff'\trf)\off =\off s'$\nnsp.\oss
Suppose\sss that\vspace{-3pt}\vspace{-1.5pt}
\[
\quad
\begin{tikzcd}[column sep=sboom, row sep=sboom]
E\fff'
\arrow[d, "\dis p'\dff"']
\arrow[r]
&
E
\arrow[d, "\dis \dff p"]
\\
B\fff'
\arrow[r]
&
B
\end{tikzcd}
\]

\vspace{-12pt}\vspace{3pt}\vspace{-1.5pt}
is\dss a bundle map and\sss the map\sss $B\fff'\qff \ttoo\qff B$\sss
takes\sss 
$s'\fff,\qff t'$\sss 
to\sss 
$s\fff,\qff t$\sss  
respectively.\oss
Suppose also\sss that\sss the map\sss $E\fff'\qff \ttoo\qff E$\sss
takes\sss $\widetilde{s}\qff'$\sss to\sss $\widetilde{s}$\nnsp.\oss
Let\sss $F\fff'\off =\off (\trf p'\trf)^{\dff -\dff 1}\dff(\trf t'\trf)$\nnsp.\oss

We will\sss need\sss two special\sss cases of\dss this situation.\oss
In\sss the first\sss case\sss $B\fff'\off =\off B$\sss
and\sss $B\off =\off B\fff'\qff \ttoo\qff B$\dss is\dss the identity\sss map.\oss
Then we can choose as\sss the covering\sss homotopy\sss for\sss
$p\dff \colon\dff E\qff \ttoo\qff B$\dss
the composition of\dss the covering\sss homotopy\sss for\sss
$p'\dff \colon\dff E\fff'\qff \ttoo\qff B\fff'$\sss
with\sss the map\sss $E\fff'\qff \ttoo\qff E$\nnsp.\oss
In\sss this case we have\sss the following\sss commutative diagram\oss\vspace{-3pt}
\[
\quad
\begin{tikzcd}[column sep=boom, row sep=boom]
\Omega\dff B
\arrow[d, "\dis ="']
\arrow[r, "\dis g^{\dff p'}"]
&
\Omega\trf(\trf E\fff'\fff,\qff F\fff'\trf)
\arrow[d]
\arrow[r, "\dis \ev"]
&
F\fff'
\arrow[d]
\\
\Omega\dff B
\arrow[r, "\dis g^{\dff p}"]
&
\Omega\trf(\trf E\fff,\qff F\trf)
\arrow[r, "\dis \ev"]
&
F
\end{tikzcd}
\]

\vspace{-12pt}\vspace{3pt}
In\sss the second case\sss the bundle\sss
$p'\dff \colon\dff E\fff'\qff \ttoo\qff B\fff'$\sss
is\dss assumed\sss to be\sss the bundle induced\sss from\sss
$p\dff \colon\dff E\qff \ttoo\qff B$\sss
by\sss the map\sss $B\fff'\qff \ttoo\qff B$\nnsp.\oss
In other\sss words,\oss the canonical\sss map\sss\vspace{3pt}
\[
\quad
E\fff'
\qff \ttoo\qff B\fff'\qff \times_{\dff B}\qff E
\]

\vspace{-12pt}\vspace{3pt}
is\dss a homeomorphism.\oss
In\sss this case we can\sss identify\sss the fibers\sss $F\fff'$\sss and\sss $F$\dnsp.\oss
In\sss this situation\sss we can\sss take as\sss the covering\sss homotopy\sss for\sss
$p'\dff \colon\dff E\fff'\qff \ttoo\qff B\fff'$\sss
the pull-back\sss of\sss a covering homotopy\sss for\sss
$p\dff \colon\dff E\qff \ttoo\qff B$\nnsp.\oss
Then\sss the diagram\oss\vspace{-3pt}
\[
\quad
\begin{tikzcd}[column sep=boom, row sep=boom]
\Omega\dff B\fff'
\arrow[d]
\arrow[r, "\dis g^{\dff p'}"]
&
\Omega\trf(\trf E\fff'\fff,\qff F\trf)
\arrow[d]
\arrow[r, "\dis \ev"]
&
F
\arrow[d, "\dis ="]
\\
\Omega\dff B
\arrow[r, "\dis g^{\dff p}"]
&
\Omega\trf(\trf E\fff,\qff F\trf)
\arrow[r, "\dis \ev"]
&
F\dff.
\end{tikzcd}
\]

\vspace{-12pt}\vspace{3pt}
is\dss commutative.\oss
If\dss we choose another covering\sss homotopy\sss for\sss $p'$\nnsp,\oss
then\sss this diagram\dss is\dss not\sss necessarily\sss commutative,\oss
but\dss is\dss still\sss commutative up\sss to homotopy.\oss

\myuppar{The\dss Atiyah--Singer\dss proof.}
Let\sss $\mathcal{A}$\sss be\sss the $C^{\dff *}$\nsp\dnsp-algebra of\dss
bounded\sss operators in\sss $H$\nnsp,\oss
and\sss let\sss $\mathcal{K}$\sss be\sss the closed\sss two-sided\sss ideal\sss
of compact\sss operators.\oss
The quotient\sss $\mathcal{A}/\mathcal{K}$\sss is\dss the\qss \emph{Calkin\sss algebra}.\oss
It\dss is\dss also a $C^{\dff *}$\nsp\dnsp-algebra.\oss
Let\sss $p\dff \colon\dff \mathcal{A}\qff \ttoo\qff \mathcal{A}/\mathcal{K}$\sss
be\sss the quotient\sss map,\oss
and\sss let\sss $\mathcal{G}$\sss be\sss the group of\dss invertible elements of\sss
$\mathcal{A}/\mathcal{K}$\dnsp.\oss
Then\sss $\mathcal{F}\off =\off p^{\dff -\dff 1}\dff(\trf \mathcal{G}\trf)$\nnsp.\oss
Let\sss $G$\sss be\sss the subgroup\sss of\dss unitary elements of\sss $\mathcal{G}$\dnsp,\oss
i.e.\qss elements\sss 
$x\qff \in\qff \mathcal{G}$\sss such\sss that\sss
$x\dff x^{\dff *}\off =\off 1$\nnsp.\oss
By a slight\sss abuse of\dss notations we will\sss denote
the induced\sss maps\sss
$\mathcal{F}\qff \ttoo\qff \mathcal{G}$\sss
and\sss
$p^{\dff -\dff 1}\dff(\trf G\trf)\qff \ttoo\qff G$\sss
by\sss $p$\nnsp.\oss
Both of\dss them are homotopy equivalences.\oss
See\qss \cite{as},\oss the discussion\sss preceding\trs Lemma\qss (2.3).\oss

There\dss is\dss a deformation\sss retraction\sss
$r\dff \colon\dff
\mathcal{G}\qff \ttoo\qff G$\nnsp,\oss
which can\sss be described as follows.\oss
If\sss $A\qff \in\qff \mathcal{F}$\sss
and\sss $A\off =\off U\trf \num{A}$\sss is\dss the polar decomposition,\oss
then\sss 
$r\trf(\trf p\trf(\trf A\trf)\trf)\off =\off p\trf(\trf U\trf)$\nnsp.\oss
The map\sss $r$\sss is\dss connected\sss to\sss 
the identity\sss by\sss the\sss linear\sss homotopy,\oss which\dss is\dss fixed on\sss $G$\nnsp.\oss
See\qss \cite{as},\oss the formula\qss (2.1).\oss
It\sss follows\sss that\sss
$\Omega\dff r\dff \circ\trf \Omega\dff p
\dff \colon\dff
\Omega\trf \mathcal{F}
\qff \ttoo\qff 
\Omega\trf G$\sss
is\dss a homotopy equivalence.\oss

Let\sss us\sss consider\sss $\mathcal{U}$\sss 
simply\sss as\sss the space of\dss isometries\sss 
$u\dff \colon\dff
H\qff \ttoo\qff H$\sss with\sss 
the norm\sss topology.\oss
Let\sss $-\qff U\comp$\sss be\sss the subspace of\dss
isometries $u$ differing\sss from\sss $-\qff \id$\sss by a compact\sss operator,\oss
i.e.\qss such\sss that\sss $u\qff +\qff \id$\sss is\dss compact.\oss
The image $p\trf(\trf \mathcal{U}\trf)$\sss is\dss the 
component\sss of\sss the identity\sss $G_{\dff 0}$\sss of\sss $G$\nnsp,\oss
and\sss $p$\sss induces a map\sss 
$\mathbf{p}\dff \colon\dff
\mathcal{U}\qff \ttoo\qff G_{\dff 0}$\nsp,\oss
which\sss turns out\sss to be\sss a\sss locally\sss trivial\sss bundle.\oss
See\qss \cite{as},\oss Proposition\qss 3.2.\oss
Trivially,\pss
$\Omega\trf G\off =\off \Omega\trf G_{\dff 0}$\nsp.\oss
Let\sss us choose\sss  
$1\off =\off p\trf(\trf \id\trf)$\sss 
and\sss 
$-\qff 1\off =\off p\trf(\trf -\qff \id \trf)$\sss 
as\sss the base points of\sss $G$\nnsp.\oss
Clearly,\pss 
$\mathbf{p}^{\dff -\dff 1}\dff(\trf -\qff 1\trf)
\off =\off 
-\qff U\comp$\nnsp.\oss
Since $\mathcal{U}$\sss is\dss contractible,\oss the bundle\sss $\mathbf{p}$\sss
leads\sss to  homotopy equivalences\vspace{1.5pt}
\[
\quad
g^{\trf \mathbf{p}}\dff \colon\trf 
\Omega\trf G
\qff \ttoo\qff 
\Omega\trf(\trf \mathcal{U}\fff,\qff -\qff U\comp\trf)
\quad
\mbox{and}\quad
f^{\trf \mathbf{p}}\dff \colon\trf 
\Omega\trf G
\qff \ttoo\qff
-\qff U\comp
\pff.
\]

\vspace{-12pt}\vspace{1.5pt}
In order\sss to prove\sss that\sss $\alpha$\sss is\dss
a homotopy equivalence,\oss
it\dss is\dss sufficient\sss to prove\sss that\sss\vspace{1.5pt}
\[
\quad
f^{\trf \mathbf{p}}\dff \circ\qff 
\Omega\dff r\dff \circ\qff 
\Omega\dff p\dff \circ\trf \alpha
\]

\vspace{-12pt}\vspace{1.5pt}
is\dss
a homotopy equivalence.\oss
The proof\dss of\dss this\sss fact\sss by\trs
Atiyah\dss and\dss Singer\qss \cite{as}\qss 
involves\sss two preliminary steps.\oss

First,\oss the space\sss $\hat{F}\comp$\sss from\dss
Section\qss \ref{polarizations-splittings}\qss is\dss 
a deformation\sss retract\sss of\sss $\hat{\mathcal{F}}$\dnsp.\oss
Therefore\sss it\dss is\dss sufficient\sss to prove\sss that\sss the 
restriction of\dss
$f^{\trf \mathbf{p}}\dff \circ\qff 
\Omega\dff r\dff \circ\qff 
\Omega\dff p\dff \circ\trf \alpha$\sss
to\sss $\hat{F}\comp$\dss is\dss a homotopy equivalence.\oss
It\dss turns out\sss that\sss
$\alpha$\sss 
maps\sss $\hat{F}\comp$\sss
into\sss $\Omega\trf p^{\dff -\dff 1}\dff(\trf G\trf)$\sss
and\sss hence\sss the restrictions of\dss
$\Omega\dff r\dff \circ\qff 
\Omega\dff p\dff \circ\trf \alpha$\sss
and\dss
$\Omega\dff p\dff \circ\trf \alpha$\sss to\sss $\hat{F}\comp$\sss
are equal.\oss
See\qss \cite{as},\oss the discussion around\trs Lemma\qss (2.4).\oss
Therefore\sss it\dss is\dss sufficient\sss to prove\sss that\sss
$f^{\trf \mathbf{p}}\dff \circ\pff
\mathbf{a}$\sss
is\dss a homotopy equivalence,\oss
where\vspace{1.5pt}
\[
\quad
\mathbf{a}
\qff \colon\qff
\hat{F}\comp
\qff \ttoo\qff
\Omega\trf G
\pff
\]

\vspace{-12pt}\vspace{1.5pt}
is\dss the map induced\sss by\sss
$\Omega\dff p\dff \circ\dff \alpha$\nnsp.\oss
On\sss the other\sss hand,\oss
Lemma\qss \ref{loops-bundles-lifts}\qss
and\sss the functoriality\sss of\dss maps\sss $f^{\dff p}$\sss
imply\sss that\sss up\sss to homotopy\sss the maps\vspace{1.5pt}
\[
\quad
f^{\trf \mathbf{p}}\dff \circ\pff
\mathbf{a}
\quad\dff
\mbox{and}\quad
f^{\trf \mathbf{p}}\dff(\trf \mathbf{a}\trf)
\qff \colon\qff
\hat{F}\comp
\qff \ttoo\qff
-\qff U\comp
\pff
\]

\vspace{-12pt}\vspace{1.5pt}
are equal.\oss
The second\sss preliminary step\dss is\dss the observation\sss
that\dss there is\dss a natural\sss choice of\dss a covering\sss homotopy\sss
$h_{\fff u}\trf(\trf \mathbf{a}\trf)
\dff \colon\dff 
\hat{F}\comp\trf \ttoo\qff \mathcal{U}\fff,\qff 
u\qff \in\qff [\trf 0\fff,\qff \pi\trf]$\nnsp.\oss
Namely,\oss the homotopy\vspace{1.5pt}
\[
\quad
h_{\fff u}\trf(\trf \mathbf{a}\trf)
\dff \colon\dff
A
\off \longmapsto\off 
\exp\trf(\trf u\fff i\dff A\trf)
\fff,\quad 
u\qff \in\qff [\trf 0\fff,\qff \pi\trf]
\]

\vspace{-12pt}\vspace{1.5pt}
is\dss a covering\sss homotopy.\oss
See\qss \cite{as},\oss the discussion around\trs Lemma\qss (2.4).\oss
The map\vspace{1.5pt}
\[
\quad
f^{\trf \mathbf{p}}\trf(\trf \mathbf{a}\trf)
\dff \colon\dff 
\hat{F}\comp
\qff \ttoo\qff 
-\qff U\comp
\]

\vspace{-12pt}\vspace{1.5pt}
corresponding\sss to\sss this covering\sss homotopy\dss 
is\dss the map\vspace{1.5pt}
\[
\quad
\exp\dff \pi\fff i
\qff \colon\qff
A
\off \longmapsto\off 
\exp\trf(\trf \pi\fff i\dff A\trf)
\pff.
\]

\vspace{-12pt}\vspace{1.5pt}
It\sss follows\sss that\sss
$f^{\trf \mathbf{p}}\dff \circ\pff
\mathbf{a}$\sss
is\dss homotopic\sss to $\exp\dff \pi\fff i$\nnsp.\oss
It\dss remains\sss to prove\sss that\sss
$\exp\dff \pi\fff i$\sss is\dss a homotopy equivalence.\oss
This\sss proof\trs is\dss the heart\sss of\dss the paper\qss \cite{as}.\oss
See\qss \cite{as},\oss Proposition\qss (3.3).\oss
But\sss we are more concerned\sss with\sss these preliminary steps.\oss

\myuppar{A characterization of\dss the homotopy class of\sss $\alpha$\nnsp.}
The preliminary steps show\sss that\sss up\sss to homotopy\sss
the map\sss
$\alpha$\sss 
can\sss be characterized as a map\sss 
$\hat{\mathcal{F}}
\qff \ttoo\qff
\Omega\qff \mathcal{F}$\sss
such\sss that\sss the composition\sss 
$f^{\trf \mathbf{p}}\dff \circ\qff 
\Omega\dff p\dff \circ\trf \alpha$\sss
induces a map\dss 
$\hat{F}\comp
\qff \ttoo\qff 
-\qff U\comp$\sss
homotopic\sss to $\exp\dff \pi\fff i$\nnsp.\oss
By\trs Corollary\qss \ref{subspace-polarized}\qss 
the inclusion\sss map\sss
$\hat{F}\ffin\qff \ttoo\qff \hat{\mathcal{F}}$\sss
is\dss a homotopy equivalence.\oss
Hence\sss 
$\alpha$\sss 
can\sss be also characterized\sss as a map\sss
$\hat{\mathcal{F}}
\qff \ttoo\qff
\Omega\qff \mathcal{F}$\sss 
such\sss that\sss the restriction\sss of\sss
$f^{\trf \mathbf{p}}\dff \circ\qff 
\Omega\dff p\dff \circ\trf \alpha$\sss
to\sss $\hat{F}\ffin$
is\dss a map\dss 
$\hat{F}\ffin
\qff \ttoo\qff 
-\qff U\comp$\sss
homotopic\sss to $\exp\dff \pi\fff i$\nnsp.\oss

\myuppar{Another\sss map\sss
$\hat{\mathcal{F}}
\qff \ttoo\qff
\Omega\qff \mathcal{F}$\dnsp.}
Now\sss we will\sss discuss\sss the homotopy\sss equivalence\sss
$\hat{\mathcal{F}}
\qff \ttoo\qff
\Omega\qff \mathcal{F}$\sss 
implicit\sss in\sss the proof\dss of\trs Theorem\qss \ref{f-loops}.\oss
The starting\sss point\dss is\dss the bundle\sss
$\bm{\pi}\dff \colon\dff
\mathbf{U}
\qff \ttoo\qff 
\num{\mathcal{P}{\nsp}\mathcal{S}_{\qff 0}}$\nnsp,\oss
leading\sss to\sss the homotopy equivalence\sss\vspace{1.5pt}
\[
\quad
f^{\trf \bm{\pi}}
\dff \colon\dff
\Omega\trf \num{\mathcal{P}{\nsp}\mathcal{S}}
\off =\off
\Omega\trf \num{\mathcal{P}{\nsp}\mathcal{S}_{\qff 0}}
\qff \ttoo\qff 
-\qff U\ffin
\pff.
\]

\vspace{-12pt}\vspace{1.5pt}
On\sss the other\sss hand,\oss the map\sss
$\num{\hat{\pi}}\dff \colon\dff
\num{\mathcal{P}{\nsp}\hat{\mathcal{S}}}
\qff \ttoo\qff
\num{\hat{\mathcal{S}}}$\sss
is\dss a homotopy equivalence
and\sss the map\sss 
$h\dff \colon\dff
\num{\hat{\mathcal{S}}}
\qff \ttoo\qff
-\qff U\ffin$\sss
is\dss a homeomorphism.\oss
Hence\sss for some\sss $\bm{\gamma}$\sss the diagram\vspace{-1.5pt}
\[
\quad
\begin{tikzcd}[column sep=sboom, row sep=sboom]
\protect{\num{\mathcal{P}{\nsp}\hat{\mathcal{S}}}}
\arrow[r, "\dis \bm{\gamma}\vphantom{f}\dff"]
\arrow[d, "\dis \protect{\num{\hat{\pi}}}\qff"']
&
\Omega\trf \protect{\num{\mathcal{P}{\nsp}\mathcal{S}}}
\arrow[d, "\dis \qff f^{\trf \bm{\pi}}\dff"]
\\
\protect{\num{\hat{\mathcal{S}}}}
\arrow[r, "\dis h"]
&
-\qff U\ffin
\end{tikzcd}
\]

\vspace{-9pt}
is\dss commutative up\sss to homotopy.\oss
Clearly,\oss such a map\sss 
$\bm{\gamma}$\sss is\dss unique up\sss to homotopy 
and\dss is\dss a homotopy equivalence.\oss
Since\sss there are canonical\sss homotopy equivalences\sss
$\hat{\mathcal{F}}
\qff \ttoo\qff 
\num{\mathcal{P}{\nsp}\hat{\mathcal{S}}}$\sss
and\sss
$\mathcal{F}
\qff \ttoo\qff
\num{\mathcal{P}{\nsp}\mathcal{S}}$\nnsp,\oss
the map\sss $\bm{\gamma}$\sss defines a homotopy equivalence\sss
$\bm{\alpha}\dff \colon\dff
\hat{\mathcal{F}}
\qff \ttoo\qff
\Omega\trf \mathcal{F}$\nnsp.\oss

\myuppar{A canonical\sss map\sss 
$\num{\mathcal{P}{\nsp}\mathcal{S}}\qff \ttoo\qff G$\nnsp.}
Let\sss us consider\sss $G$\sss as a\sss topological\sss category\sss
having only\sss identity\sss morphisms.\oss
Somewhat\sss surprisingly,\oss
there\dss is\dss a\sss continuous\sss functor\vspace{1.5pt}
\[
\quad
\kappa\dff \colon\dff
\mathcal{P}{\nsp}\mathcal{S}
\qff \ttoo\qff
G
\pff.
\]

\vspace{-12pt}\vspace{1.5pt}
For an object\sss\sss
$P
\off =\off 
(\trf E_{\dff 1}\dff,\pff E_{\dff 2}\dff,\pff i\trf)$\sss
of\trs $\mathcal{P}{\nsp}\mathcal{S}$\dss
let\sss us\sss define $\kappa\trf(\trf P\trf)$\sss
as\sss the image in\sss
$\mathcal{A}/\mathcal{K}$\sss 
of\dss an arbitrary\sss linear operator\sss
$k\dff \colon\dff H\qff \ttoo\qff H$\sss 
equal\dss to $i$ on\sss 
$H\dff \ominus\dff E_{\dff 1}$\nsp.\oss
Any\sss two such operators differ by an operator of\dss finite rank
and\sss hence\sss their\sss images in\sss
$\mathcal{A}/\mathcal{K}$\sss
are equal.\oss
Therefore\sss $\kappa\trf(\trf P\trf)$\sss is\dss well\sss defined.\oss
If\dss we\sss take as $k$\sss the operator equal\sss to $0$ on\sss $E_{\dff 1}$\nsp,\oss
then\sss $k\dff k^{\dff *}$\sss and\sss $k^{\dff *}\fff k$\sss
are projections on subspaces of\dss finite codimension and\sss hence\sss
their\sss images in\sss
$\mathcal{A}/\mathcal{K}$\sss 
are equal\sss to $1$\nnsp.\oss
It\sss follows\sss that\sss
$\kappa\trf(\trf P\trf)\qff \in\qff G$\sss 
for every object\sss $P$\sss of\sss $\mathcal{P}{\nsp}\mathcal{S}$\dnsp.\oss

Therefore\sss the map\sss
$\kappa\dff \colon\dff
P\pff \longmapsto\pff \kappa\trf(\trf P\trf)$\sss
defines a map,\oss still\sss denoted\sss by\sss $\kappa$\nnsp,\oss 
from\sss the space of\dss objects of\dss
$\mathcal{P}{\nsp}\mathcal{S}$\sss to\sss $G$\nnsp.\oss
Clearly,\oss $\kappa$\sss is\dss continuous,\oss
and\dss if\dss there\dss is\dss a morphism\sss
$P\qff \ttoo\qff P\fff'$\nnsp,\oss then\sss
$\kappa\trf(\trf P\trf)\off =\off \kappa\trf(\trf P\fff'\trf)$\nnsp.\oss
Hence\sss $\kappa$\sss can\sss be considered as a\sss functor\sss
$\mathcal{P}{\nsp}\mathcal{S}
\qff \ttoo\qff
G$\nnsp.\oss
The geometric realization of\dss the functor\sss $\kappa$\sss 
is\dss a canonical\sss map\sss
\vspace{1.5pt}
\[
\quad
\num{\kappa}\dff \colon\dff
\num{\mathcal{P}{\nsp}\mathcal{S}}
\qff \ttoo\qff 
\num{G}
\off =\off 
G
\pff.
\]

\vspace{-12pt}\vspace{1.5pt}
Clearly,\pss $\num{\kappa}$\sss induces a map\sss
$\num{\mathcal{P}{\nsp}\mathcal{S}_{\qff 0}}\qff \ttoo\qff G_{\dff 0}$\nsp,\oss
which\sss we will\sss also denote by\sss $\num{\kappa}$\nnsp.\oss

\mypar{Lemma.}{two-maps-to-g}
\emph{The diagram}\vspace{-3pt}
\[
\quad
\begin{tikzcd}[column sep=sboom, row sep=sboom]
\mathcal{E}
\arrow[r, "\dis \mathcal{P}\dff\iota\vphantom{a}\dff"]
\arrow[d, "\dis \varphi\qff"']
&
\mathcal{P}{\nsp}\mathcal{S}
\arrow[d, "\dis \trf \kappa\dff"]
\\
\mathcal{F}
\arrow[r, "\dis r\dff \circ\dff p\phantom{\mathcal{P}}"]
&
G\dff,
\end{tikzcd}
\]

\vspace{-10.5pt}
\emph{where\sss $\mathcal{F}$\sss and\sss $G$\sss are considered as categories
with only identity morphisms,\oss
is\dss commutative.\oss}

\proof
Let\sss
$(\trf A\fff,\qff \varepsilon\trf)$\sss
be an object\sss of\sss $\mathcal{E}$\dnsp,\oss
and\sss let\sss $A\off =\off U\trf \num{A}$\sss 
be\sss the polar decomposition of\sss $A$\nnsp.\oss
Then\sss
$r\dff \circ\dff p\trf \circ\trf \varphi$\sss
takes\sss $(\trf A\fff,\qff \varepsilon\trf)$\sss
to\sss $p\trf(\trf U\trf)$\nnsp.\oss
On\sss the other\sss hand,\pss
$\mathcal{P}\dff\iota$\sss
takes\sss $(\trf A\fff,\qff \varepsilon\trf)$\sss
to a\sss triple\sss
$(\trf E_{\dff 1}\dff,\pff E_{\dff 2}\dff,\pff i\trf)$\sss
such\sss that\sss the isometry\sss $i$\sss is\dss induced\sss by\sss $U$\nnsp.\oss
Therefore\sss $\kappa\dff \circ\trf \mathcal{P}\dff\iota$\sss
also\sss takes\sss $(\trf A\fff,\qff \varepsilon\trf)$\sss
to\sss $p\trf(\trf U\trf)$\nnsp.\oss
This proves\sss that\sss the diagram\dss is\dss commutative on objects.\oss
Since\sss $G$\sss has only\sss identity\sss morphisms,\oss
it\dss is\dss commutative on\sss morphisms also.\oss  \eproof

\mypar{Corollary.}{kappa-equivalence}
\emph{The map\sss
$\num{\kappa}\dff \colon\dff
\num{\mathcal{P}{\nsp}\mathcal{S}}
\qff \ttoo\qff
G$\sss
is\dss a homotopy equivalence.\oss}

\proof
Since\sss $\num{\varphi}$\nnsp,\qss $\num{\mathcal{P}\dff\iota}$\nnsp,\oss
and\sss $r\dff \circ\dff p$\sss are homotopy\sss equivalences,\oss
the\sss corollary\sss follows\sss from\trs Lemma\qss \ref{two-maps-to-g}.\oss  \eproof

\myuppar{The induced\sss bundle.}
Let\sss us\sss consider\sss the bundle
$\mathbf{q}\dff \colon\dff
\mathcal{V}
\qff \ttoo\qff
\num{\mathcal{P}{\nsp}\mathcal{S}_{\qff 0}}$\sss
induced\sss from\sss the bundle\sss
$\mathbf{p}\dff \colon\dff
\mathcal{U}
\qff \ttoo\qff
G_{\dff 0}$\sss
by\sss the map\sss
$\num{\kappa}\dff \colon\dff
\num{\mathcal{P}{\nsp}\mathcal{S}_{\qff 0}}
\qff \ttoo\qff 
G_{\dff 0}$\nsp.\oss
Then we have a pull-back diagram\vspace{-1.5pt}
\[
\quad
\begin{tikzcd}[column sep=sboom, row sep=sboom]
\mathcal{V}
\arrow[d, "\dis \mathbf{q}\qff"']
\arrow[r]
&
\mathcal{U}
\arrow[d, "\dis \fff\qff\mathbf{p}"]
\\
\protect{\num{\mathcal{P}{\nsp}\mathcal{S}_{\qff 0}}}
\arrow[r, "\dis \protect{\num{\kappa}}"]
&
G_{\dff 0}\qff.
\end{tikzcd}
\]

\vspace{-9pt}
\mypar{Lemma.}{v-contractible}
\emph{The space\dss $\mathcal{V}$\sss
is\dss contractible.\oss}

\proof
Corollary\qss \ref{kappa-equivalence}\qss implies\sss that\sss
$\num{\kappa}\dff \colon\dff
\num{\mathcal{P}\fff\mathcal{O}_{\dff 0}}
\qff \ttoo\trf
G_{\dff 0}$\sss
is\dss a homotopy equivalences.\oss
Since\sss $\mathcal{U}$\sss is\dss contractible,\oss 
comparing\sss the homotopy sequences of\dss
the bundles\sss $\mathbf{q}$\sss and\sss $\mathbf{p}$\sss
and of\dss the bundles\sss $\bm{\rho}$\sss and\sss $\bm{\pi}$\nnsp,\oss
shows\sss that\sss $\mathcal{V}$\sss is\dss weakly\sss contractible.\oss
Therefore\sss it\dss is\dss sufficient\sss to prove\sss that\sss
$\mathcal{V}$\sss is\dss homotopy equivalent\sss to a\dss CW-complex.\oss
Since\sss the base\sss 
$\num{\mathcal{P}{\nsp}\mathcal{S}_{\qff 0}}$\sss 
of\dss $\mathbf{q}$\sss has\sss
the homotopy\sss type of\dss a\dss CW-complex,\oss
a\sss theorem of\trs tom\dss Dieck\qss \cite{td1}\qss
implies\sss that\dss
it\dss is\dss sufficient\sss to prove\sss that\sss the fibers of\dss
the bundle $\mathbf{q}$\sss are homotopy equivalent\sss to\dss CW-complexes.\oss
These fibers are homeomorphic\sss to\sss $U\comp$\dnsp.\oss
The space\sss $U\comp$\sss is\dss a deformation\sss retract\sss
of\dss an open subset\sss of\dss a\dss Banach\dss space and\sss
hence\dss is\dss homotopy equivalent\sss to a\dss CW-complex.\oss
See\qss \cite{as},\oss the beginning of\trs Section\qss 3.\oss
The\sss open subset\sss in question\dss is\dss 
the space of\dss invertible operators\sss $T$\sss 
such\sss that\sss $\id\qff +\qff T$\sss is\dss compact,\oss
and\sss the deformation\sss retraction\dss is\dss 
defined\sss by\sss the standard\sss homotopy\qss (2.10)\qss from\qss \cite{as}.\oss
The\sss lemma\sss follows.\oss  \eproof

\myuppar{The induced\dss bundle and\dss loop spaces.}
By\sss the functoriality of\dss maps\sss $f^{\dff p}$\sss  
the diagram\vspace{-1.5pt}
\begin{equation}
\label{r-p}
\quad
\begin{tikzcd}[column sep=sboom, row sep=sboom]
\Omega\trf \protect{\num{\mathcal{P}{\nsp}\mathcal{S}_{\qff 0}}}
\arrow[d, "\dis \Omega\trf \protect{\num{\kappa}}\qff"']
\arrow[r, "\dis f^{\dff \mathbf{q}}"]
&
-\qff U\comp
\arrow[d]
\\
\Omega\trf G_{\dff 0}
\arrow[r, "\dis f^{\dff \mathbf{p}}"]
&
-\qff U\comp
\end{tikzcd}
\end{equation}

\vspace{-9pt}
is\dss commutative.\oss

On\sss the other\sss hand,\oss the bundle\sss
$\bm{\pi}\dff \colon\dff
\mathbf{U}
\qff \ttoo\qff 
\num{\mathcal{P}{\nsp}\mathcal{S}_{\qff 0}}$\sss
is\dss almost\sss a subbundle\sss of\dss the bundle
$\mathbf{q}
\dff \colon\dff
\mathcal{V}
\qff \ttoo\qff 
\num{\mathcal{P}{\nsp}\mathcal{S}_{\qff 0}}$\nnsp.\oss
Name\-ly,\oss the fiber\sss of\sss $\bm{\pi}$\sss
over\sss
$P
\off =\off 
(\trf E_{\dff 1}\dff,\pff E_{\dff 2}\dff,\pff i\trf)$\sss
is\dss the space of\dss isometries\sss
$H\qff \ttoo\qff H$\sss equal\sss to\sss $i$\sss on a subspace of\dss
finite codimension\sss in\sss $E_{\dff 1}$\nsp.\oss
If\sss $u$\sss is\dss a\sss such\sss isometry,\oss
then\sss the fiber\sss of\sss $\mathbf{q}$\sss over\sss $P$\sss
is\dss the space of\dss isometries\sss
$u'$\sss such\sss that\sss $u'\qff -\qff u$\sss is\dss
a compact\sss operator.\oss
It\sss follows\sss that\sss fibers of\sss $\bm{\pi}$\sss
are contained\sss in\sss fibers of\sss $\mathbf{q}$\nnsp.\oss
While\sss the\sss topology of\sss $\mathbf{U}$\sss is\dss
different\sss from\sss the\sss topology\sss induced\sss from\sss $\mathcal{V}$\nnsp,\oss
the inclusion\sss 
$\mathbf{U}\qff \ttoo\qff \mathcal{V}$\sss 
is\dss continuous,\oss
and\sss hence\dss is\dss a map\sss of\dss bundles from\sss the bundle\sss
$\bm{\pi}\dff \colon\dff
\mathbf{U}
\qff \ttoo\qff 
\num{\mathcal{P}{\nsp}\mathcal{S}_{\qff 0}}$\sss 
to\sss the bundle\sss
$\mathbf{q}\dff \colon\dff
\mathcal{V}
\qff \ttoo\qff 
\num{\mathcal{P}{\nsp}\mathcal{S}_{\qff 0}}$\nnsp.\oss
Arguing\sss as in\sss the case of\dss actual\sss subbundles,\oss we see\sss that\sss
the diagram\vspace{-3pt}
\begin{equation}
\label{rho-r}
\quad
\begin{tikzcd}[column sep=sboom, row sep=sboom]
\Omega\trf \protect{\num{\mathcal{P}{\nsp}\mathcal{S}_{\qff 0}}}
\arrow[d, "\dis \dff ="']
\arrow[r, "\dis f^{\dff \bm{\pi}}"]
&
-\qff U\ffin
\arrow[d]
\\
\Omega\trf \protect{\num{\mathcal{P}{\nsp}\mathcal{S}_{\qff 0}}}
\arrow[r, "\dis f^{\dff \mathbf{q}}"]
&
-\qff U\comp\dff,
\end{tikzcd}
\end{equation}

\vspace{-10.5pt}
where\sss the right\sss vertical\sss arrow\dss is\dss the inclusion,\oss
is\dss commutative.\oss

\myuppar{A convention.}
By an abuse of\dss notations,\oss 
if\dss $f\dff \colon\dff X\qff \ttoo\qff Y$\sss is\dss a map and\sss 
$f\trf(\trf A\trf)\qff \subset\qff B$\sss for some\sss
$A\qff \subset\qff X$\sss and\dss $B\qff \subset\qff Y$\dnsp,\oss
we will\sss denote the map\sss $A\qff \ttoo\qff B$\sss
induced\sss by\sss $f$\sss simply\sss by\sss  
$f\dff \colon\dff A\qff \ttoo\qff B$\nnsp.\oss

\myuppar{The map\sss $\alpha$\sss and classifying\sss spaces.}
Since\sss the map\sss $\num{\hat{\varphi}\ffin}$\sss
from\trs Theorem\qss \ref{finite-to-enhanced-models}\qss
is\dss a homotopy equivalence,\oss
we can\qss ``lift''\qss the above characterization of\dss the map\sss $\alpha$\sss
up\sss to a homotopy\sss to\sss $\num{\hat{\mathcal{E}}\ffin}$\nnsp.\oss
In more details,\oss 
let\sss us\sss consider\sss 
the following diagram.\oss\vspace{-1.5pt}\vspace{-0.25pt}
\[
\quad
\begin{tikzcd}[column sep=sboom, row sep=sboom]
\protect{\num{\hat{\mathcal{E}}\ffin}}
\arrow[r, "\dis \beta"]
\arrow[d, "\dis \protect{\num{\hat{\varphi}\ffin}}\qff"']
&
\Omega\trf \protect{\num{\mathcal{E}}}
\arrow[r, "\dis \Omega\trf \protect{\num{\mathcal{P}\dff\iota}}\vphantom{a}"]
\arrow[d]
&
\Omega\trf \protect{\num{\mathcal{P}{\nsp}\mathcal{S}_{\qff 0}}}
\arrow[d, "\dis \qff\Omega\trf \protect{\num{\kappa}}"]
&
\\
\hat{F}\ffin
\arrow[r, "\dis \alpha\vphantom{f}"]
&
\Omega\trf p^{\dff -\dff 1}\dff(\trf G\trf)
\arrow[r, "\dis \Omega\trf p"]
&
\Omega\trf G
\arrow[r, "\dis f^{\dff \mathbf{p}}"]
&
-\qff U\comp
\end{tikzcd}
\]

\vspace{-9pt}
Here\sss the middle vertical\sss arrow\dss is\dss the map\dss
$\Omega\dff r\dff \circ\qff 
\Omega\dff p\dff \circ\qff
\Omega\trf \num{\varphi}$\nnsp,\oss
and as\sss the arrow\sss $\beta$\sss one can\sss take any map such\sss that\sss the\sss
left\sss square\dss is\dss commutative up\sss to homotopy.\oss
Such a map\sss $\beta$\sss exists and\dss is\dss unique up\sss to homotopy\sss
because\sss
$\num{\hat{\varphi}\ffin}\qff \colon\qff
\num{\hat{\mathcal{E}}\ffin}
\qff \ttoo\qff 
\hat{F}\ffin$\dss 
and\vspace{1.5pt}
\[
\quad
\Omega\dff r\dff \circ\qff 
\Omega\dff p\dff \circ\qff
\Omega\trf \num{\varphi}
\off =\off
\Omega\trf (\trf r\dff \circ\dff p\dff \circ\dff \num{\varphi}\trf)
\]

\vspace{-12pt}\vspace{1.5pt}
are homotopy equivalences.\oss
Clearly,\oss the map\sss $\beta$\sss is\dss a homotopy equivalence.\oss
Lemma\qss \ref{two-maps-to-g}\qss implies\sss that\sss the right\sss square\dss
is\dss commutative.\oss
Up\sss to homotopy\sss the map\sss $\alpha$\sss is\dss determined\sss by\sss requiring\sss that\sss
the map\sss $\hat{F}\ffin\qff \ttoo\qff -\qff U\comp$\dnsp,\oss
obtained\sss by following\sss the arrows of\dss the diagram,\oss
is\dss homotopic\sss to $\exp\dff \pi\fff i$\nnsp.\oss
Since\sss the vertical\sss arrows in\sss the\sss left\sss square are homotopy equivalences,\oss
up\sss to homotopy\sss $\beta$\sss is\dss determined\sss by\sss requiring\sss that\sss
the map\sss \vspace{1.5pt}
\[
\quad
\num{\hat{\mathcal{E}}\ffin}\qff \ttoo\qff -\qff U\comp 
\]

\vspace{-12pt}\vspace{1.5pt}
obtained\sss by following\sss the arrows of\dss 
the diagram,\oss
is\dss homotopic\sss to 
$\exp\dff \pi\fff i
\qff \circ\qff
\num{\hat{\varphi}\ffin}$\nnsp.\oss

Let\sss us consider\sss now the diagram\vspace{-4.5pt}
\[
\quad
\begin{tikzcd}[column sep=sboom, row sep=sboom]
\protect{\num{\hat{\mathcal{E}}\ffin}}
\arrow[r, "\dis \beta"]
\arrow[d]
&
\Omega\trf \protect{\num{\mathcal{E}}}
\arrow[d, "\dis \dff \Omega\trf \protect{\num{\mathcal{P}\dff\iota}}\vphantom{a}"]
&
&
\\
\protect{\num{\mathcal{P}{\nsp}\hat{\mathcal{S}}}}
\arrow[r, "\dis \gamma\vphantom{f}"]
\arrow[d, "\dis \dff ="]
&
\Omega\trf \protect{\num{\mathcal{P}{\nsp}\mathcal{S}_{\qff 0}}}
\arrow[r, "\dis \Omega\trf \protect{\num{\kappa}}"]
\arrow[d, "\dis \dff ="]
&
\Omega\trf G
\arrow[d, "\dis \qff f^{\dff \mathbf{p}}"]
\\
\protect{\num{\mathcal{P}{\nsp}\hat{\mathcal{S}}}}
\arrow[r, "\dis \gamma\vphantom{f}"]
&
\Omega\trf \protect{\num{\mathcal{P}{\nsp}\mathcal{S}_{\qff 0}}}
\arrow[r, "\dis f^{\dff \mathbf{q}}"]
&
-\qff U\comp\dff.
\end{tikzcd}
\]

\vspace{-12pt}
Here one can\sss take as\sss the arrow\sss $\gamma$\sss any\sss map such\sss that\sss the
upper\sss left\sss square\dss is\dss commutative up\sss to homotopy.\oss
Such a map\sss $\gamma$\sss exists and\dss is\dss unique up\sss to homotopy\sss
because\sss the vertical\sss arrows in\sss this square are homotopy equaivalences.\oss
Indeed,\oss Theorems\qss \ref{forgetting-odd}\qss and\qss \ref{two-fredholm}\qss
together\sss with\trs Theorem\qss \ref{fredholm-polarization-forget}\qss
imply\sss that\sss
$\num{\mathcal{P}\dff\iota}\dff \colon\dff
\num{\mathcal{E}}
\qff \ttoo\qff 
\num{\mathcal{P}{\nsp}\mathcal{S}}$\sss
is\dss a homotopy equivalence,\oss 
and\trs Theorems\qss \ref{operators-categories}\qss and\qss \ref{finite-to-enhanced-models}\qss 
imply\sss that\sss
$\num{\hat{\mathcal{E}}\ffin}
\qff \ttoo\qff
\num{\mathcal{P}{\nsp}\hat{\mathcal{S}}}$\sss
is\dss a homotopy equivalence.\oss
Clearly,\oss $\gamma$\sss is\dss a homotopy equivalence.\oss
The bottom\sss right\sss square\dss is\dss simply\sss
the diagram\qss (\ref{r-p})\qss drawn differently,\oss
The bottom\sss left\sss square\dss is\dss trivially commutative.\oss
The characterization of\sss $\beta$\sss implies\sss that\sss
up\sss to homotopy $\gamma$\sss is\dss determined\sss by\sss requiring\sss that\sss
the map\sss $\num{\hat{\mathcal{E}}\ffin}\qff \ttoo\qff -\qff U\comp$\dnsp,\oss
obtained\sss by following\sss the arrows of\dss 
this diagram,\oss
is\dss homotopic\sss to 
$\exp\dff \pi\fff i
\qff \circ\qff
\num{\hat{\varphi}\ffin}$\nnsp.\oss

\mypar{Theorem.}{gamma-gamma}
\emph{The maps\sss $\gamma$\sss and\pss
$\bm{\gamma}\trf \colon\dff
\num{\mathcal{P}{\nsp}\hat{\mathcal{S}}}
\qff \ttoo\qff
\Omega\trf \num{\mathcal{P}{\nsp}\mathcal{S}_{\qff 0}}$\trs
are homotopic.\oss}

\proof
In view of\trs Theorem\qss \ref{exp}\qss the map\sss $\bm{\gamma}$\sss
can\sss be also defined as a map such\sss that\sss\vspace{-3pt}
\[
\quad
\begin{tikzcd}[column sep=sboom, row sep=sboom]
\protect{\num{\mathcal{P}{\nsp}\hat{\mathcal{S}}}}
\arrow[r, "\dis \bm{\gamma}\vphantom{f}\dff"]
\arrow[d, "\dis \mathcal{P}h\dff"']
&
\Omega\trf \protect{\num{\mathcal{P}{\nsp}\mathcal{S}_{\qff 0}}}
\arrow[d, "\dis \qff f^{\trf \bm{\pi}}\dff"]
\\
\hat{F}\ffin
\arrow[r, "\dis \exp\ffin"]
&
-\qff U\ffin
\end{tikzcd}
\]

\vspace{-10.5pt}
is\dss a commutative diagram,\oss
where\sss $\mathcal{P}h$\sss is\dss the homeomorphism\sss from\trs
Theorem\qss \ref{polarized-homeo}.\oss
By\trs Theorem\qss \ref{enhanced-and-finite}\qss the composition\sss
\[
\quad
\num{\hat{\mathcal{E}}\ffin}
\qff \ttoo\qff
\num{\mathcal{P}{\nsp}\hat{\mathcal{S}}}
\qff \ttoo\qff
\hat{F}\ffin 
\pff,
\]

\vspace{-12pt}
where\sss the second\sss map\dss is\sss $\mathcal{P}h$\nnsp,\oss
is\dss homotopic\sss to\sss $\num{\hat{\varphi}\ffin}$\nnsp.\oss

Therefore\sss $\bm{\gamma}$\sss 
can\sss be defined as\sss the map such\sss that\sss 
the composition\vspace{-3pt}
\[
\quad
\begin{tikzcd}[column sep=large, row sep=large]
\protect{\num{\hat{\mathcal{E}}\ffin}}
\arrow[r]
&
\protect{\num{\mathcal{P}{\nsp}\hat{\mathcal{S}}}}
\arrow[r, "\dis \bm{\gamma}"]
&
\Omega\trf \protect{\num{\mathcal{P}{\nsp}\mathcal{S}_{\qff 0}}}
\arrow[r, "\dis f^{\dff \bm{\pi}}"]
&
-\qff U\ffin
\end{tikzcd}
\]

\vspace{-10.5pt}
is\dss homotopic\sss to\sss
$\exp\ffin\dff \circ\pff \num{\hat{\varphi}\ffin}$\nnsp.\oss
But\sss $\exp\ffin$\sss differs\sss from\sss from\sss $\exp\dff \pi\fff i$\sss
only\sss by\sss the domain of\dss definition,\oss which\dss is\dss smaller\sss
for\sss $\exp\ffin$\dnsp.\oss
Hence\dss\vspace{3pt}
\[
\quad
\exp\ffin\dff \circ\pff \num{\hat{\varphi}\ffin}
\off =\off
\exp\dff \pi\fff i\qff \circ\pff \num{\hat{\varphi}\ffin}
\]

\vspace{-12pt}\vspace{3pt}
and\sss the map\sss $\bm{\gamma}$\sss can\sss be defined as\sss the map such\sss that\sss 
the above composition\dss is\dss homotopic\sss to\sss
$\exp\dff \pi\fff i\qff \circ\pff \num{\hat{\varphi}\ffin}$\nnsp.\oss
Now\sss the commutativity of\dss the square\qss (\ref{rho-r})\qss implies\sss that\sss
the composition\vspace{-3pt}
\[
\quad
\begin{tikzcd}[column sep=large, row sep=large]
\protect{\num{\hat{\mathcal{E}}\ffin}}
\arrow[r]
&
\protect{\num{\mathcal{P}{\nsp}\hat{\mathcal{S}}}}
\arrow[r, "\dis \bm{\gamma}"]
&
\Omega\trf \protect{\num{\mathcal{P}{\nsp}\mathcal{S}_{\qff 0}}}
\arrow[r, "\dis f^{\dff \mathbf{q}}"]
&
-\qff U\comp
\end{tikzcd}
\]

\vspace{-9pt}
is\dss homotopic\sss to\sss
$\exp\dff \pi\fff i\qff \circ\pff \num{\hat{\varphi}\ffin}$\nnsp.\oss
As we saw above,\oss this property characterizes\sss the map\sss $\gamma$\sss
up\sss to homotopy.\oss
It\sss follows\sss that\sss $\bm{\gamma}$\sss 
is\dss homotopic\sss to\sss $\gamma$\nnsp.\oss  \eproof

\mypar{Theorem.}{alpha-alpha-beta-beta}
\emph{The maps\dss $\alpha$\sss and\pss
$\bm{\alpha}\dff \colon\dff
\hat{\mathcal{F}}
\qff \ttoo\qff
\Omega\qff \mathcal{F}$ 
are homotopic.\oss}

\proof
The square\vspace{-1.5pt}
\[
\quad
\begin{tikzcd}[column sep=sboom, row sep=sboom]
\hat{\mathcal{F}}
\arrow[r, "\dis \alpha\dff"]
\arrow[d]
&
\Omega\qff \mathcal{F}
\arrow[d]
\\
\protect{\num{\mathcal{P}{\nsp}\hat{\mathcal{S}}}}
\arrow[r, "\dis \gamma"]
&
\Omega\trf \protect{\num{\mathcal{P}{\nsp}\mathcal{S}_{\qff 0}}}\dff,
\end{tikzcd}
\]

\vspace{-9pt}
where\sss the vertical\sss arrows are\sss the canonical\sss homotopy equivalences,\oss
is\dss commutative\sss by\sss the definition of\sss $\gamma$\nnsp.\oss
The square\vspace{0pt}
\[
\quad
\begin{tikzcd}[column sep=sboom, row sep=sboom]
\hat{\mathcal{F}}
\arrow[r, "\dis \bm{\alpha}\dff"]
\arrow[d]
&
\Omega\qff \mathcal{F}
\arrow[d]
\\
\protect{\num{\mathcal{P}{\nsp}\hat{\mathcal{S}}}}
\arrow[r, "\dis \bm{\gamma}"]
&
\Omega\trf \protect{\num{\mathcal{P}{\nsp}\mathcal{S}_{\qff 0}}}\dff,
\end{tikzcd}
\]

\vspace{-9pt}
where\sss the vertical\sss arrows are\sss the same,\oss
is\dss commutative\sss by\sss the definition of\dss $\bm{\alpha}$\nnsp.\oss
The\sss theorem\sss follows from\sss the commutativity of\dss
these squares\sss together\sss with\trs
Theorem\qss \ref{gamma-gamma}.\oss  \eproof

\newpage
\mysection{Hilbert\qss bundles}{hilbert-bundles}

\myuppar{Categories related\sss to\sss Hilbert\sss bundles.}
The following modification of\dss the constructions of\trs 
Section\qss \ref{classifying-spaces-saf}\qss is\dss intended\sss for
applications\sss to families of\dss self-adjoint\dss Fredholm\sss operators.\oss
Let\sss $X$\sss be a\qss (compactly\sss generated)\qss space,\oss
and\sss let\sss $\mathbb{H}$\sss be a\sss 
locally\sss trivial\dss Hilbert\dss bundle with separable fibers over $X$\nnsp.\oss
It\dss is\dss convenient\sss to\sss think\sss that\sss $\mathbb{H}$\sss is\dss a family\sss
$H_{\dff x}\dff,\qff x\qff \in\qff X$ of\dss Hilbert\sss spaces parametrized\sss by $X$\nnsp.\oss
We will\dss replace\sss finitely dimensional\sss
subspaces of\dss the fixed\dss Hilbert\sss space $H$ 
by\sss finitely dimensional\sss subspaces of\dss all\dss fibers 
$H_{\dff x}$\nsp,\dss $x\qff \in\qff X$\nnsp.\oss
Let\sss us define a\sss topology on\sss the set\sss of\dss such subspaces.\oss
A\sss local\sss trivialization of\sss $\mathbb{H}$\sss over an open subset\sss
$U\qff \subset\qff X$\sss allows\sss to identify\sss fibers $H_{\dff x}$\sss
with\sss $x\qff \in\qff U$\sss with\sss $H$\nnsp,\oss
and\sss hence\sss to identify subspaces contained\sss in\sss these fibers\sss
with\sss pairs\sss $(\trf x\fff,\qff V\trf)$\nnsp,\oss where\sss
$x\qff \in\qff U$\sss and\sss $V$\sss is\dss a finitely dimensional\sss
subspace of\sss $H$\nnsp.\oss
This\sss leads\sss to a\sss topology\sss on\sss the set\sss of\dss subspaces
of\dss fibers\sss $H_{\dff x}$\sss with\sss $x\qff \in\qff U$\nnsp.\oss
Since we consider only finitely dimensional\sss subspaces,\oss
this\sss topology\dss is\dss independent\sss on\sss the choice of\dss
the\sss trivialization.\oss
Therefore\sss these\sss locally defined\sss topologies\sss lead\sss to a\sss
topology on\sss the set\sss of\dss finitely dimensional\sss subspaces of\dss
all\sss fibers.\oss\vspace{-0.25pt}

Let\sss us\sss define a\sss topological\sss category\sss 
$\hat{\mathcal{S}}\dff(\trf \mathbb{H}\trf)$\sss
as follows.\oss
Its objects are\sss finitely dimensional\sss subspaces of\dss the fibers 
$H_{\dff x}\dff,\qff x\qff \in\qff X$\nnsp.\oss
Morphisms\sss $V\qff \ttoo\qff V\fff'$ exists only\dss if\sss
$V\fff,\qff V\fff'$\sss are contained\sss in\sss the same fiber\sss $H_{\dff x}$\nsp,\oss
and\sss in\sss this case morphisms are defined exactly as morphisms of\sss
$\hat{\mathcal{S}}$\sss with\sss $H$\sss replaced\sss by\sss $H_{\dff x}$\nsp.\oss
If\dss we consider $X$ as\sss topological\sss category\sss with $X$\sss being\sss
the space of\dss objects and only\sss the identity\sss morphisms,\oss
then\sss there\dss is\dss a\sss functor\sss
$\eta\dff \colon\dff 
\hat{\mathcal{S}}\dff(\trf \mathbb{H}\trf)
\qff \ttoo\qff
X$\sss
assigning\sss to an object\sss $V$\sss the point\sss $x\qff \in\qff X$\sss
such\sss that\sss $V\qff \subset\qff H_{\dff x}$\nsp.\oss\vspace{-0.05pt}

Let\sss us\sss define a\sss topological\sss category\sss 
$Q\dff(\trf \mathbb{H}\trf)$\sss 
as follows.\oss
Its objects are pairs\sss $(\trf x\fff,\qff V\trf)$\sss
such\sss that\sss $x\qff \in\qff X$\sss and\sss $V$\sss is\dss a finitely\sss dimensional\dss
Hilbert\dss space.\oss
The\sss topology on\sss the set\sss of\dss objects\dss is\dss the product\sss
of\dss the\sss topology of\sss $X$\sss and\sss the discrete\sss topology on\sss
the set\sss of\dss spaces.\oss
Morphisms\sss
$(\trf x\fff,\qff V\trf)
\qff \ttoo\qff 
(\trf x'\fff,\qff V\fff'\trf)$\sss
exist\sss only\sss if\sss $x\off =\off x'$\nnsp,\oss
and\sss in\sss this case\sss they are\sss the same as morphisms\sss
$V\qff \ttoo\qff V\fff'$\sss in\sss the category $Q$\nnsp.\oss
The\sss topology on\sss the set\sss of\dss morphisms\dss is\dss defined\sss by\sss
the\sss topology of\sss $X$\sss and\sss the\sss topology on\sss the set\sss
of\dss morphisms of\sss $Q$\nnsp.\oss
Clearly,\oss the topological\sss category\sss $Q\dff(\trf \mathbb{H}\trf)$\sss  
depends only\sss on $X$ and\sss not\sss on\sss the bundle $\mathbb{H}$\nnsp.\oss
There\dss is\dss a functor\sss
$\pi\dff \colon\dff 
Q\dff(\trf \mathbb{H}\trf)
\qff \ttoo\qff
X$\sss
assigning\sss to\sss $(\trf x\fff,\qff V\trf)$\sss the point $x$\nnsp.\oss
Clearly,\oss there\dss is\dss a canonical\sss bijection\sss
$\num{Q\dff(\trf \mathbb{H}\trf)}
\qff \ttoo\qff
X\dff \times\dff \num{Q}$\nnsp.\oss
Since $X$\sss is\dss compactly\sss generated,\oss this bijection\dss
is\dss actually a homeomorphism.\oss\vspace{-0.05pt}

Finally,\oss there\dss is\dss also an analogue\sss $Q/\fff\mathbb{H}$
of\dss the category\sss $Q/\fff H$\nnsp.\oss
The objects of\sss $Q/\fff\mathbb{H}$ are pairs\sss 
$(\trf V\fff,\qff h\trf)$\nnsp,\oss
where $V$\sss is\dss an object\sss of\sss $Q$\sss and\sss
$h$\sss is\dss an\sss isometric embedding\sss $V\qff \ttoo\qff H_{\dff x}$\sss
for some $x\qff \in\qff X$\nnsp.\oss
Morphisms\sss
$(\trf V\fff,\qff h\trf)
\qff \ttoo\qff 
(\trf V\fff'\fff,\qff h'\trf)$\sss
exist\sss only\sss if\dss the maps $h$ and\sss $h'$\sss are maps\sss
to\sss the same fiber\sss $H_{\dff x}$\nsp,\oss
and\sss in\sss this case morphisms are defined as morphisms of\sss
$Q/\fff H_{\dff x}$\nsp.\oss
Ignoring\sss isometric embeddings defines a functor\sss
$Q/\fff \mathbb{H}\qff \ttoo\qff Q\dff(\trf \mathbb{H}\trf)$\nnsp.\oss
Also,\oss by assigning\sss to an object\sss
$(\trf V\fff,\qff h\trf)$\sss of\sss $Q/\fff\mathbb{H}$\sss
the image\sss $h\dff(\trf V\trf)$\sss we get\sss a functor\sss
$i\dff \colon\dff
Q/\fff\mathbb{H}
\qff \ttoo\qff
\hat{\mathcal{S}}\dff(\trf \mathbb{H}\trf)$\nnsp.\oss\vspace{-0.25pt}

\mypar{Theorem.}{intermediary-bundle}
\emph{If\pss $X$ is\dss paracompact,\oss
then\sss the maps\sss
$\num{Q\dff(\trf \mathbb{H}\trf)}
\off \longleftarrow\off 
\num{Q/\fff \mathbb{H}}
\qff \ttoo\qff
\num{\hat{\mathcal{S}}\dff(\trf \mathbb{H}\trf)}$\sss
induced\sss by\sss these functors are homotopy equivalences.\oss}

\proof
The proof\dss is\dss similar\sss to\sss the proof\dss of\trs Theorem\qss \ref{intermediary}.\oss
Let\sss us consider\sss first\sss the map\sss
$\num{Q/\fff \mathbb{H}}
\qff \ttoo\qff
\num{Q\dff(\trf \mathbb{H}\trf)}$\nnsp.\oss
Clearly,\oss the composition\sss
$\ob\trf Q/\fff \mathbb{H}
\qff \ttoo\qff
\ob\trf Q\dff(\trf \mathbb{H}\trf)
\qff \ttoo\qff
X$\nnsp,\oss
where\sss the second\sss map\dss is\sss $\pi$\nnsp,\oss
is\dss equal\sss to $\eta$\nnsp.\oss
For every $x\qff \in\qff X$\sss the preimage\sss $\pi^{\trf -\dff 1}\dff(\dff x\trf)$\sss
is\dss the discrete space\sss $\mathbf{V}$ with\sss points corresponding\sss to 
finitely dimensional\dss Hilbert\dss spaces.\oss
In\sss particular,\pss $\pi$\sss is\dss a\sss locally\sss trivial\sss bundle.\oss
The preimage\sss $\eta^{\dff -\dff 1}\dff(\dff x\trf)$\sss
is\dss the disjoint\sss union over\sss $V\qff \in\qff \mathbf{V}$\sss
of\dss spaces of\dss isometric embeddings\sss $V\qff \ttoo\qff H_{\dff x}$\nsp.\oss
It\sss follows\sss that\sss $\ob\trf Q/\fff \mathbb{H}$\sss is\dss is\dss
the disjoint\sss union over\sss $V\qff \in\qff \mathbf{V}$\sss
of\dss spaces of\dss pairs\sss $(\dff x\fff,\qff h\trf)$\nnsp,\oss
where\sss $x\qff \in\qff X$\sss and\sss 
$h\dff \colon\dff V\qff \ttoo\qff H_{\dff x}$\sss
is\dss an\sss isometric embedding.\oss
Every such space\dss is\dss a\sss locally\sss trivial\sss bundle over $X$
with contractible fibers.\oss
Since $X$\sss is\dss paracompact,\oss its projection\sss to $X$\sss
is\dss a homotopy equivalence.\oss
It\sss follows\sss that\sss the map\sss
$\ob\trf Q/\fff \mathbb{H}
\qff \ttoo\qff
\ob\trf Q\dff(\trf \mathbb{H}\trf)$\sss
is\dss a homotopy equivalence,\oss
i.e.\qss the functor\sss
$Q/\fff \mathbb{H}
\qff \ttoo\qff
Q\dff(\trf \mathbb{H}\trf)$\sss
induces a homotopy equivalence of\dss the spaces of\dss objects.\oss\vspace{-0.16pt}

The space of\dss morphisms of\sss $Q/\fff \mathbb{H}$\sss is\dss 
a\sss locally\sss trivial\dss bundle over\sss the space of\dss
morphisms of\dss $Q\dff(\trf \mathbb{H}\trf)$\nnsp,\oss
with\sss the fiber over a morphism\sss
$(\trf x\fff,\qff V\trf)
\qff \ttoo\qff 
(\trf x\fff,\qff  V\fff'\trf)$\sss 
of\dss $Q\dff(\trf \mathbb{H}\trf)$\sss being\sss the space 
of\dss isometric embeddings
$V\fff'\qff \ttoo\qff H_{\dff x}$\nsp.\oss
Since\sss the spaces of\dss isometric embeddings
$V\fff'\qff \ttoo\qff H_{\dff x}$\sss are contractible,\oss
$Q/\fff \mathbb{H}\qff \ttoo\qff Q\dff(\trf \mathbb{H}\trf)$\sss 
induces a homotopy equivalence of\dss the spaces of\dss morphisms also.\oss
A similar argument\sss applies\sss to $n$\dnsp-simplices.\oss
Since all\dss involved categories have free units,\oss
it\sss follows\sss that\sss
$\num{Q/\fff \mathbb{H}}
\qff \ttoo\qff
\num{Q\dff(\trf \mathbb{H}\trf)}$\sss
is\dss a\sss homotopy equivalence.\oss\vspace{-0.16pt}  

Let\sss us consider\sss the map
$\iota
\off =\off
\num{i}\dff \colon\dff
\num{Q/\fff \mathbb{H}}
\qff \ttoo\qff
\num{\hat{\mathcal{S}}\dff(\trf \mathbb{H}\trf)}$\nnsp.\oss
Since every simplex\sss involves only one fiber,\oss
the corresponding\sss part\sss of\dss the proof\dss of\trs Theorem\qss \ref{intermediary}\qss
applies without\sss any changes.\oss
The paracompactness assumption allows\sss to apply\sss the\sss theorem of\trs
tom\dss Dieck\qss \cite{td1}\qss as in\sss the proof\dss of\trs Theorem\qss \ref{to-models}\qss
({\fff}this step was omitted\sss in\dss 
the proof\dss of\trs Theorem\qss \ref{intermediary}).\oss  \eproof\vspace{-0.16pt}

\mypar{Theorem.}{harris-h-bundles}
\emph{There\dss is\dss a canonical\dss homeomorphism\dss
$\num{\hat{\mathcal{S}}\dff(\trf \mathbb{H}\trf)}
\qff \ttoo\qff 
-\qff U^{\dff \mathrm{fin}}\dff(\trf \mathbb{H}\trf)${\nsp},\oss
where\sss $-\qff U^{\dff \mathrm{fin}}\dff(\trf \mathbb{H}\trf)$\sss
is\dss the\sss total\sss space of\trs the bundle with\sss
the fiber\sss $-\qff U^{\dff \mathrm{fin}}$\sss associated\sss with\sss $\mathbb{H}$\nnsp.\oss}\vspace{-0.16pt}

\proof
The proof\trs follows\sss the proof\dss of\trs 
Theorem\qss \ref{harris-h}\qss with\sss parameters $x\qff \in\qff X$\sss added.\oss  \eproof\vspace{-0.16pt}

\myuppar{Remarks.}
One can also define categories\sss 
$\mathcal{E}\hat{\mathcal{O}}\dff(\trf \mathbb{H}\trf)$\sss
and\dss
$\hat{\mathcal{O}}\dff(\trf \mathbb{H}\trf)$\sss
related\sss to\sss
$\mathcal{E}\hat{\mathcal{O}}$\sss
and\dss
$\hat{\mathcal{O}}$\sss
in\sss the same way as\sss
$\hat{\mathcal{S}}\dff(\trf \mathbb{H}\trf)$\sss
is\dss related\sss to\sss
$\hat{\mathcal{S}}$\dnsp.\oss
Moreover,\oss one can\sss prove for\sss these categories\sss
analogues of\trs Theorems\qss \ref{to-models}\qss and\qss \ref{forgetting-operators}.\oss
We\sss leave\sss these\sss tasks\sss to\sss the interested\sss readers.\oss\vspace{-0.16pt}

On\sss the other\sss hand\sss there\dss is\dss no bundle analogue of\dss
$\hat{\mathcal{F}}$\dnsp,\oss
neither as a space nor as a\sss topological\sss category,\oss
and no bundle analogue of\dss $\hat{\mathcal{E}}$\dnsp.\oss
The reason\dss is\dss that\sss the\sss topologies of\dss
$\hat{\mathcal{F}}$\sss and\sss $\hat{\mathcal{E}}$\sss are\sss too
strong\sss to construct\sss bundles associated\sss with\sss $\mathbb{H}$\sss
with\sss fibers\sss
$\hat{\mathcal{F}}$\sss and\sss $\hat{\mathcal{E}}$\dnsp.\oss
One can construct\sss associated\sss bundles if\sss 
$\hat{\mathcal{F}}$\sss is\dss equipped\sss with\sss the compact-open\sss
topology on\sss the space of\dss operators\sss $H\qff \ttoo\qff H$\nnsp,\oss
but\sss for\sss this\sss topology\sss the stability\sss property\sss
from\dss Section\qss \ref{spaces-operators}\qss fails.\oss\vspace{-0.16pt}

In a more plain\sss language,\oss there are no good\sss notions
of\dss continuous families\sss of\qss Fredholm\sss and of\dss self-adjoint\trs Fredholm\sss
operators acting on\sss fibers of\dss a\sss general\dss Hilbert\dss bundle.\oss
At\sss the same\sss time\sss families of\dss operators arising\sss in\sss the index\sss theory\sss
are good enough\sss to apply\sss the\sss theory of\dss this section\sss without\sss
defining\sss analogues of\dss
$\hat{\mathcal{F}}$\sss and\sss $\hat{\mathcal{E}}$\dnsp.\oss

\vspace{\bigskipamount}

\begin{flushright}

November\qss 28\fff,\oss 2021
 
https\halfff:/\!/\!nikolaivivanov.com

E-mail\halfff:\oss nikolai.v.ivanov{\fff}@{\dff}icloud.com\vspace{12pt}

Department\sss of\trs Mathematics,\oss Michigan\sss State\sss University

\end{flushright}

\end{document}